TECHNISCHE UNIVERSITÄT DRESDEN

FAKULTÄT MATHEMATIK / INSTITUT FÜR ALGEBRA

Masterarbeit

*zur Erlangung des Hochschulgrades*

Master of Science (M.Sc.)

# On Abstract and Concrete Minions

Vorgelegt von **Lukas Juhrich**

*geboren am 15. März 1998 in Berlin*

| Eingereicht am | Betreut von | Zweitgutachter |
|---|---|---|
| 3. Juni 2024 | Prof. Dr. Manuel Bodirsky | Prof. Dr. Ulrich Krähmer |

*Für Jasmin*

## Selbstständigkeitserklärung

Hiermit erkläre ich, Lukas Juhrich, dass ich diese Arbeit unter dem Titel "On abstract and concrete minions" unter Betreuung von Prof. Dr. Manuel Bodirsky selbstständig erarbeitet, verfasst, und Zitate kenntlich gemacht habe. Andere als die angegebenen Hilfsmittel wurden nicht benutzt.

Ort, Datum

Matrikelnummer

Unterschrift

# Contents

# ABSTRACT

This thesis is expository in nature. We analyze the connection between abstract minions, which can be described as functors from the category of finite ordinals to sets, and concrete minions, which are sets $\mathrm{Pol}(\underline{A}, \underline{B})$ of polymorphisms $\underline{A}^k \to \underline{B}$ between relational structures $\underline{A}$ and $\underline{B}$. The functorial structure arises because a function $\alpha\colon n \to k$ transforms a $n$-ary polymorphism $f\colon \underline{A}^n \to \underline{B}$ to the polymorphism $(x_0, ..., x_{k-1}) \mapsto f\left(x_{\alpha(0)}, ..., x_{\alpha(n-1)}\right)$. This data is relevant to the constraint satisfaction problems over finite domain, because minion homomorphisms $\mathrm{Pol}(\underline{A}, \underline{B}) \to \mathrm{Pol}(\underline{A'}, \underline{B'})$ give rise to log-space reductions from the promise constraint satisfaction problem $\mathrm{PCSP}(\underline{A'}, \underline{B'})$ to $\mathrm{PCSP}(\underline{A}, \underline{B})$. Thus, they are valuable to understand the homomorphism order of polymorphism minions over finite-domain structures, especially since it is unknown whether the P/NP complexity dichotomy of constraint satisfaction problems extends to the more generalized promise setting.

Crucially, although implicitly used in some papers to solve nontrivial problems, the concept of an abstract minion has not been exposed on so far, even though much is known about functor categories. The aim of this thesis is to start closing this gap, and to apply well-known constructions from category theory to minions. Furthermore, we identify a condition under which an abstract minion can arise as a concrete one over a finite domain, and deduce that many constructions remain stable under some of these concreteness assumptions. We will conclude that some of the homomorphism orders in question are uncountable distributive bounded lattices, and moreover (bi-) Heyting algebras. Along the way, we collect open questions related to abstract and concrete minions.

## Foreword to the revised edition

This document presents a significant revision of the author's master's thesis, originally submitted on June 3rd, 2024.

The author acknowledges the support of the European Research Council grant 101071674, POCOCOP, during the revision process. Views and opinions expressed are however those of the author only and do not necessarily reflect those of the European Union or the European Research Council Executive Agency. Neither the European Union nor the granting authority can be held responsible for them.

## Foreword

It is an old mathematical wisdom that symmetries can tell us a lot about an object we are studying. No less is true for the study of constraint satisfaction problems ("CSPs"), such as the problem of determining whether a graph can be 3-colored: Whether the problem is "hard" (NP-hard) or "easy" (in P) is completely determined by the generalized, higher-dimensional symmetries of the so-called CSP template. More formally, given a relational structure $\underline{H}$, the CSP is concerned with deciding on a given input structure $\underline{G}$ whether there exists a homomorphism $\underline{G} \to \underline{H}$. The symmetries we alluded to are the *polymorphisms* of the template $\underline{H}$, i.e., homomorphisms from an arbitrary power $\underline{H}^k$ to $\underline{H}$. Let us denote the set of polymorphisms on $\underline{H}$ by $\mathrm{Pol}(\underline{H})$.

In 2017, after decades of research on the topic, A. Bulatov and D. Zhuk independently managed to establish a *complexity dichotomy* for CSPs (Theorem 1.2.i), which states that every CSP is either NP-complete or in P.

They proved that a long known condition sufficient for hardness [1, Cor. 7.3] is also *necessary*: By providing polynomial-time algorithms which work when this condition is violated, they demonstrated that all other templates must be tractable.

A key tool to the analysis of the CSP is the so-called "algebraic approach". Its core insight is that maps $h\colon \mathrm{Pol}(\underline{H}) \to \mathrm{Pol}(\underline{H'})$ which preserve arities and identities give rise to reductions $\mathrm{CSP}(\underline{H'}) \to \mathrm{CSP}(\underline{H})$.

This follows from the notion of primitive positive gadget reductions [2, Thm. 5.15] and Birkhoff's theorem of universal algebra [2, Prop. 8.37].

However, that insight has been sharpened in the meantime, in that it suffices to consider maps which preserve arities and *height-one identities*, i.e., identities of the form $f\left(x_{\alpha(0)}, \dots, x_{\alpha(n-1)}\right) = g\left(x_{\beta(0)}, \dots, x_{\beta(k-1)}\right)$. Taking a function $\alpha\colon \{0, \dots, n-1\} \to \{0, \dots, l-1\}$ and a $n$-ary polymorphism $f$ to produce the associated function $f\alpha(x_0, \dots, x_{l-1}) \coloneqq f\left(x_{\alpha(0)}, \dots, x_{\alpha(n-1)}\right)$ is called *taking a minor*, where we call $\alpha$ a *minor operation*. We can then see that such a map $h$ preserves height-one identities if and only if it is compatible with minor operations. Hence, it is of interest to analyze the structure with base set $\mathrm{Pol}(\underline{H})$ sorted by arities and equipped with the action of minor operations. This is, roughly speaking, an *abstract minion*.

The thesis is structured as follows:

1. We shall introduce the necessary background in the "Preliminaries" chapter.
2. Then, we shall continue to define abstract and concrete minions, and give the most basic constructions, in the chapter "Minions".
3. Later, we will prove a "Cayley theorem" for minions, by which we mean that every abstract minion is isomorphic to a polymorphism minion (albeit not necessarily over finite domain).
4. We shall analyze the question of when the domains can be taken to be finite more in-depth in Chapter 4.

5. Following up, we are ready to make some basic statements about the homomorphism order of certain classes of minions in Chapter 5.
6. In a brief digression, we will then analyze the "growth" of a minion, and use this as a criterion to show a certain "smallness" condition.
7. In the last (non-summary) chapter, we will construct the exponential minion, and show how this allows us to make even more statements about the homomorphism orders.

It is at this moment that we should warn the reader about some notational conventions, each of which has a specific rationale.

- Function composition is most often denoted in diagrammatic order, i.e., $fg = f;g = g \circ f$.
- Numbers $n$ are considered to be ordinals, i.e., $0 = \emptyset$ and $n = \{0, ...n-1\}$.
- Powers $X^Y$ of sets are to be understood as functions $Y \to X$, even when the exponent is a finite ordinal like $X^n$.
- We use placeholders to signify "curried" functions, e.g. $H(_,_) = (x,y) \mapsto H(x,y)$. If the placeholder is in a sub- or superscript, we shall use a bullet • instead as in $r_{\alpha(\bullet)}$.
- The set of homomorphisms from $\underline{A}$ to $\underline{B}$ will not be denoted $\mathrm{Hom}(\underline{A}, \underline{B})$, but rather $\underline{A} \to \underline{B}$.
- The $\alpha$-minor of a function $f$ is denoted by $f\alpha$ instead of $f^\alpha$ or $f_\alpha$, and
- The arity-n component of a clone or minion is denoted by $\mathcal{M}_n$ instead of $\mathcal{M}^{(n)}$ to emphasize the covariance of the index.

The author is grateful to his advisor, Manuel Bodirsky, for his continued support, expertise, and patience, and also for including me in the meetings of the research project POCOCOP. The author is also grateful to Sebastian Meyer, without whom many arguments would be much more complicated, and whose proofreading provided valuable insights; to Andrew Moorhead for some enlightening and fun talks, and to Jakub Opršal, who replied to mails even though they only consisted of trivial insights about abstract nonsense. Additionaly, the author is indebted to his mother for her continued support and patience.

Most importantly however, the author would like to express his deepest gratitude to Jasmin, without whom the world would be a much dimmer place.

# 1 Preliminaries

## 1.1 Categories and Homomorphisms

Since abstract minions are most simply described as certain functors, it is necessary to introduce some basic category-theoretical concepts. Moreover, since the concept of a category is an abstract way to talk about any context in which there is a notion of "homomorphism" and of "composition", this is a perfectly natural context in which to define *homomorphism orders*, some minion-related instances of which will be of central importance of this thesis.

For a more general and comprehensive introduction, the reader is referred to [3] for a classic, [4], [5] for a modern introduction, or [6] as a special recommendation for the reader fluent in reading German.

### 1.1.1 Categories, Functors, and Natural Transformations

The notion of "homomorphism" is pervasive in modern mathematics: There are homomorphisms of algebraic or relational structures, like groups, graphs, or sets; homomorphisms of geometric structures such as metric spaces, topological spaces or banach spaces; Or even homomorphisms of mixed structures, such as topological groups or $C^*$-algebras.

All of these instances have one key property in common: homomorphisms can be *composed*.[1] A category is defined in such a way that it captures this compositionality. Surprisingly, a lot of properties of the structures we consider can be recovered by just looking at homomorphisms and the way in which homomorphisms compose. For instance:

- elements of a group $G$ are in bijection with group homomorphisms $\mathbb{Z} \to G$,
- edges in a graph $H$ are in bijection with graph homomorphisms $\vec{P}_1 \to H$, where $\vec{P}_1$ is the directed graph with one edge, i.e., $\bullet \to \bullet$,
- homomorphisms $h$ can often be characterized as *injective* if they satisfy a certain cancellation property with respect to all other compatible homomorphisms,
- the usual product of groups or graphs is uniquely characterized by a universal property, which we will introduce in Section 1.1.3.

The following is a formal way to define a category.

**Definition 1.1.1**: A *category* $\mathcal{C}$ consists of

- a class $\mathcal{C}_0$ whose elements are called *objects*
- a class $\mathcal{C}_1$ whose elements are called *morphisms* or *arrows*
- maps $s, t\colon \mathcal{C}_1 \to \mathcal{C}_0$ called *source* and *target*
- for every $X, Y, Z \in \mathcal{C}_0$ a *composition map* $c_{XYZ}\colon \mathcal{C}(X,Y) \times \mathcal{C}(Y,Z) \to \mathcal{C}(X,Z)$ where $\mathcal{C}(X,Y) := \{f \in \mathcal{C}_1 \mid s(f) = X \wedge t(f) = Y\}$

such that the following properties are upheld:

1. Identities exist: for every $X \in \mathcal{C}_0$ there is a morphism $\mathrm{id}_X \in \mathcal{C}(X,X)$ called the *identity morphism*, such that all of them together satisfy

[1] Just like we can define multiple monoid structures on the same set, we can define multiple notions of composition on the same set of homomorphisms. Therefore, it is more accurate to think of "composition" as extra structure than to think of "composability" as an intrinsic property.

$$c_{XXY}(\mathrm{id}_X, f) = f = c_{XYY}(f, \mathrm{id}_Y)$$

for all morphisms $f\colon \mathcal{C}(X, Y)$.

2. Composition is associative: we have

$$c_{XZW}(c_{XYZ}(f, g), h)) = c_{XYW}(f, c_{YZW}(g, h))$$

whenever $f \in \mathcal{C}(X, Y)$, $g \in \mathcal{C}(Y, Z)$ and $h \in \mathcal{C}(Z, W)$.

*Remark*: That a morphism $f \in \mathcal{C}_1$ has source $x$ and target $y$ is commonly denoted in one of the following ways:

$$f\colon X \to_{\mathcal{C}} Y \qquad X \overset{f}{\to}_{\mathcal{C}} Y$$

$$f \in \mathcal{C}(X, Y) \qquad f \in \mathrm{Mor}_{\mathcal{C}}(X, Y) \qquad f \in \mathrm{Hom}_{\mathcal{C}}(X, Y)$$

If $\mathcal{C}$ is clear from the context, we will omit the subscript, and write $f\colon X \to Y$ or $X \overset{f}{\to} Y$.

Furthermore, we will call $\mathcal{C}(X, Y) = X \to_{\mathcal{C}} Y$ the *hom-set* from $X$ to $Y$, even though in general it might not be a set but a proper class.[2]

When defining a category, we will often just state what the objects and morphisms are; the composition is usually implicit. For instance, when talking about the "category of graphs and edge-preserving maps", it is clear that we consider functions between the vertex sets and compose them accordingly.

Sometimes it even suffices to say what the objects are because it is clear what notion of homomorphism we care about: If the objects $\mathcal{C}_0$ are groups, it is usually clear that $\mathcal{C}_1$ will consist of group homomorphisms. As a rule of thumb, this usually works well in algebraic categories[3], but not so much in geometrical or relational ones. For instance, a homomorphism of metric spaces might be either a continuous map or a non-expansive one, and it is similarly unclear whether a homomorphism of graphs should preserve the presence of edges or the absence of edges as well[4].

In practice, instead of referring to the composition $c_{XYZ}$, we will use an infix operator which in addition omits the objects $X, Y$ and $Z$, because every morphism has a well-defined source and target. For such an infix operator there are two options (see Figure 1): The *diagrammatic* composition $f; g$ (read: "$f$ then $g$") and the opposite one coming from function evaluation, where $f(g(x)) =: (f \circ g)(x)$.

While the latter one is more common, in the scope of this thesis it often makes more sense to use the former one, mainly because the application of minor operations $\alpha$ and $\beta$ on an element $f$ of a minion (all to be defined in due time) satisfies

$$(f\alpha)\beta = f(\alpha; \beta),$$

so choosing $\alpha\beta$ to mean $\alpha; \beta$ allows us to write the above expression as

[2] A category is called *small* if $\mathcal{C}_1$ is a set, and *locally small* if $X \to_{\mathcal{C}} Y$ is a set for all $X, Y \in \mathcal{C}_0$.

[3] Here, the predicate *algebraic* is to be understood in a colloquial and not a formal way.

[4] The latter is often called a *strong graph homomorphism*.

$$f\alpha\beta,$$

which prevents mistakes and confusion. Therefore:

**In this thesis, we will mostly use $\alpha;\beta$, and abbreviate it as $\alpha\beta$.**

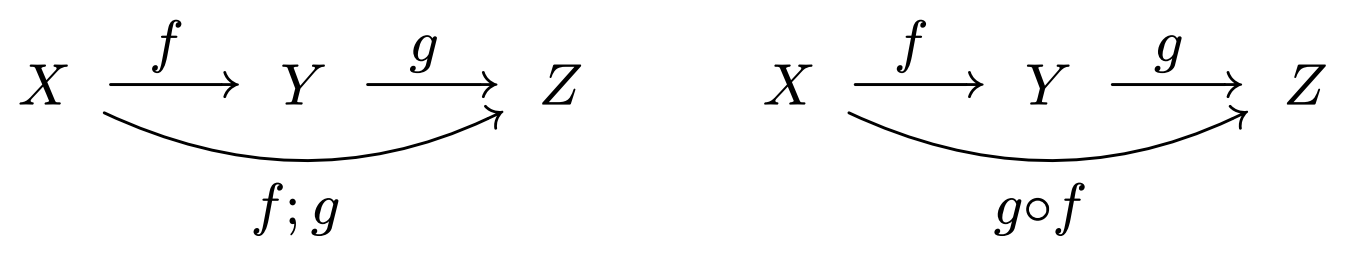

Figure 1: The two ways to denote arrow composition in a category

The definition of a category allows us to formalize when two objects are structurally equivalent, i.e., *isomorphic*. In fact, it is central to category theory that we rarely refer to a pair of objects as being equal on-the-nose, but only to them being isomorphic. For more discussion on this topic, the reader is referred to [7].

**Definition 1.1.2**: A morphism $f\colon X \to_{\mathcal{C}} Y$ is called an *isomorphism* if there exists a $g\colon Y \to_{\mathcal{C}} X$ such that $f;g = \mathrm{id}_X$ and $g;f = \mathrm{id}_Y$. The morphism $g$ is called a (two-sided, bi-) *inverse* to $f$. Objects $X$ and $Y$ are called *isomorphic*, written $X \cong Y$ , if there exists an isomorphism between them.

**Observation 1.1.1**: Being isomorphic is an equivalence relation on $\mathcal{C}_0$: Indeed,

- identities are isomorphisms,
- the inverse of an isomorphism has the original morphism as an inverse, so it is an isomorphism as well, and
- the composition of two isomorphisms has an inverse given by composing the individual inverses, so they remain to be isomorphisms.

This proves reflexivity, symmetry, and transitivity, respectively.

**Definition 1.1.3**: For a category $\mathcal{C}$ and we define $\mathrm{End}_{\mathcal{C}}(X)$ to be the monoid of *endomorphisms* consisting of all morphisms $X \to X$, equipped with composition, and $\mathrm{Aut}_{\mathcal{C}}(X)$ to be the group of *automorphisms* consisting of all isomorphisms $X \to X$.

**Definition 1.1.4** (Opposite category): Let $\mathcal{C}$ be a category with composition map $c$. The *opposite category* $\mathcal{C}^{\mathsf{op}}$ of $\mathcal{C}$ is defined as the category having the same objects, but with morphisms $\mathcal{C}^{\mathsf{op}}(X,Y) := \mathcal{C}(Y,X)$, and composition $c'_{XYZ}(f,g) := c_{ZYX}(g,f)$ whenever $g \in \mathcal{C}^{\mathsf{op}}(X,Y) = \mathcal{C}(Y,X)$ and $f \in \mathcal{C}^{\mathsf{op}}(Y,Z) = \mathcal{C}(Z,Y)$.

**Definition 1.1.5** (Functor): Let $\mathcal{C}$ and $\mathcal{D}$ be categories. A (covariant) *functor $F$* from $\mathcal{C}$ to $\mathcal{D}$, denoted by $F\colon \mathcal{C} \to \mathcal{D}$ or $F\colon [\mathcal{C}, \mathcal{D}]$ is given by

- a map $F_0\colon \mathcal{C}_0 \to \mathcal{D}_0$ on objects (denoted by $X \mapsto FX$)
- a map $F_1\colon \mathcal{C}_1 \to \mathcal{D}_1$ on morphisms (denoted by $f \mapsto Ff$)

and shall have the property that it commutes with source and target maps – i.e., $F_1$ restricts to a family of maps $\mathcal{C}(X, Y) \to \mathcal{D}(FX, FY)$ – and preserves composition (see Figure 2) and identities.

A functor $K\colon [\mathcal{C}^{\mathsf{op}}, \mathcal{D}]$ is also called a *contravariant functor* from $\mathcal{C}$ to $\mathcal{D}$.

$$X \xrightarrow{f} Y \xrightarrow{g} Z \quad (f;g\colon X \to Z) \qquad\qquad FX \xrightarrow{Ff} FY \xrightarrow{Fg} FZ \quad (F(f;g)\colon FX \to FZ)$$

Figure 2: A functor must preserve composition

**Definition 1.1.6**: A *natural transformation* between $F, G\colon [\mathcal{C}, \mathcal{D}]$ is a sequence[5] $\eta = (\eta_X)_X$ of morphisms $\eta_X\colon FX \to_{\mathcal{D}} GX$ for every $X \in \mathcal{C}_0$ such that for every $f\colon X \to_{\mathcal{C}} Y$ we have $\eta_X; Gf = Ff; \eta_Y$ – which means that the following diagram commutes:

$$\begin{array}{ccc} FX & \xrightarrow{\eta_X} & GX \\ \downarrow{\scriptstyle Ff} & & \downarrow{\scriptstyle Gf} \\ FY & \xrightarrow{\eta_Y} & GY. \end{array}$$

*Remark*: If $\eta$ is a natural transformation, this is often denoted by $\eta\colon F \Rightarrow G$, or visually:

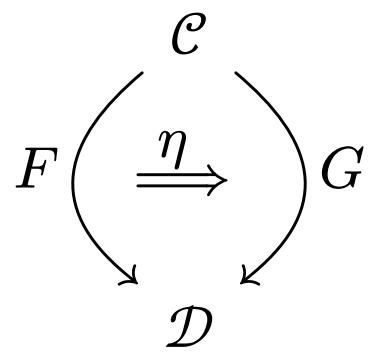

However, in this thesis, we will just use a normal arrow, i.e., $F \to G$. This notation is not ambiguous, because the only kind of "morphism between functors" one could refer to is a natural transformation.

[5] "Dependently typed sequence" might be a more accurate description, because $\eta$ does not live in some power, but rather in the product $\prod_{X \in \mathcal{C}_0} \mathcal{D}(FX, GX)$ which has different factors in every index.

Functors between fixed categories $\mathcal{C}$ and $\mathcal{D}$, together with natural transformations, form a *category*. Many categories take this form, and as we will see in Definition 2.1.1, the category of minions will be no exception.

**Definition 1.1.7** (functor category): Let $\mathcal{C}$ and $\mathcal{D}$ be categories. The *functor category* $[\mathcal{C}, \mathcal{D}]$ is defined to be the category whose objects are functors and whose morphisms are natural transformations $F \Rightarrow G$, where the composition is given object-wise as $\eta;\kappa := (\eta_X;\kappa_X)_X$.

*Proof that the functor category is a category*: Due to the pointwise definition, the associativity is inherited by the associativity in $\mathcal{D}$, and the identity morphism is given by the identity transformation $(\mathrm{id}_X)_X$. □

$$\begin{array}{ccccc} FX & \xrightarrow{\eta_X} & GX & \xrightarrow{\kappa_X} & HX \\ \big\downarrow{\scriptstyle Ff} & & \big\downarrow{\scriptstyle Gf} & & \big\downarrow{\scriptstyle Hf} \\ FY & \xrightarrow{\eta_Y} & GY & \xrightarrow{\kappa_Y} & HY \end{array}$$

Figure 3: Composition of $\eta\colon F \Rightarrow G$ and $\kappa\colon G \Rightarrow H$ is given componentwise

The prototypical exmple of a functor is the (co- or contravariant) *hom functor*. The name is historical, and arises from the fact that the morphism sets $X \to_{\mathcal{C}} Y$ were usually known as “hom-sets”.

**Definition 1.1.8** (Co- and contravariant hom-functors): Let $\mathcal{C}$ be a category. For every object $X$, we define a functor $X \to_{\mathcal{C}} _\colon \mathcal{C} \to \underline{\mathrm{Set}}$, given by

$$X \to_{\mathcal{C}} _\colon \quad \begin{array}{c} Y \\ \big\downarrow{\scriptstyle f} \\ Y' \end{array} \quad \mapsto \quad \begin{array}{c} X \to_{\mathcal{C}} Y \\ \big\downarrow{\scriptstyle f_* := _;f} \\ X \to_{\mathcal{C}} Y' \end{array}$$

called *covariant hom-functor* or *pushforward functor*. The associated map $f_*$, which sends a morphism $g\colon X \to Y$ to the composite $g;f\colon X \to Y'$ is called the *pushforward* of $f$.

Analogously, for every object $Y$, we define a functor $_ \to_{\mathcal{C}} Y\colon \mathcal{C}^{\mathrm{op}} \to \underline{\mathrm{Set}}$, given by

$$_ \to_{\mathcal{C}} Y\colon \quad \begin{array}{c} X \\ \big\downarrow{\scriptstyle f} \\ X' \end{array} \quad \mapsto \quad \begin{array}{c} X \to_{\mathcal{C}} Y \\ \big\uparrow{\scriptstyle f^* := f;_} \\ X' \to_{\mathcal{C}} Y \end{array}$$

called *contravariant hom-functor* or *pullback functor*. The associated map $f^*$, which sends a morphism $g\colon X' \to Y$ to the composite $f;g\colon X \to Y$ is called the *pullback* of $f$.

We can combine these notions by considering a *product* of categories.

**Definition 1.1.9** (Product category): Let $\mathcal{C}$, $\mathcal{D}$ be categories. We define a *product category* $\mathcal{C} \times \mathcal{D}$ by taking the product of the underlying multigraphs $(\mathcal{C}_0, \mathcal{C}_1)$ resp. $(\mathcal{D}_0, \mathcal{D}_1)$, and defining identities and composition componentwise.

Visually speaking, an arrow in the product category is of the form $(f, g)\colon (X, Y) \to (X', Y')$, where $f\colon X \to X'$ is an arrow in $\mathcal{C}$ and $g\colon Y \to Y'$ is an arrow in $\mathcal{D}$.

It is routine to verify that this definition satisfies the axioms of a category.

**Definition 1.1.10** (Hom-functor): Let $\mathcal{C}$ be a category. We define a unified *hom-functor*

$$_ \to_{\mathcal{C}} _\colon \mathcal{C}^{\mathsf{op}} \times \mathcal{C} \to \underline{\mathsf{Set}}$$

given by

$$_ \to_{\mathcal{C}} _\colon \quad \begin{array}{c} (X', Y) \\ \Big\downarrow (f^{\mathsf{op}}, g) \\ (X, Y') \end{array} \mapsto \begin{array}{c} X' \to_{\mathcal{C}} Y \\ \Big\downarrow f;_;g \\ X \to_{\mathcal{C}} Y' \end{array}$$

where a morphism $h\colon X' \to_{\mathcal{C}} Y$ gets "reflected" to a morphism

$$\begin{array}{ccc} X' & \xrightarrow{h} & Y \\ f\uparrow & & \downarrow g \\ X & & Y' \end{array}$$

### 1.1.2 (Homomorphism) Orders and Lattices

In this section we will introduce a concept most central to this thesis: associated to a category is its *homomorphism preorder*, which places $A \leq B$ if and only if there is a morphism $A \to B$. Let us develop this concept in some more detail, since it will be important later on how the relationship between different categories is reflected in their respective homomorphism preorders.

**Definition 1.1.11**: A *preorder*[6] $(P, \leq)$ is a set $P$ together with a binary relation $\leq$ which is reflexive and transitive. If the relation is also antisymmetric, i.e., if $x \leq y \leq x$ implies $x = y$, then it is called a *partial order*.

Note that a preorder can always be viewed as a category, where we take the set of arrows $p \to q$ to be populated with a distinguished element whenever $p \leq q$, and let it be empty otherwise:

[6] Often also called *quasi-order*

as there is at most one arrow between objects, there is only one way to define compostion, and it is trivial to check that the identities – the unique arrows $p \to p$ – satisfy the identity axioms.

**Definition 1.1.12**: Let $\mathcal{C}$ be a category. The *homomorphism preorder* of $\mathcal{C}$, denoted by $\mathcal{C}/\rightrightarrows$, is defined to have the objects $\mathcal{C}_0$ of $\mathcal{C}$ as elements, and to have $X \leq Y$ whenever there is a homomorphism $X \to Y$.

The *homomorphism order* of $\mathcal{C}$, denoted by $\mathcal{C}/\leftrightarrow$, is the poset quotient of the homomorphism preorder, which considers $X$ and $Y$ to be equivalent if both $X \leq Y$ and $Y \leq X$, and orders equivalence classes $[X] \leq [Y]$ if and only if $X \leq Y$.

Furthermore, we will often write $X \to Y$ for the statement that there exists a homomorphism from $X$ to $Y$.

The rationale for the notation $\mathcal{C}/\rightrightarrows$ is that moving to the homomorphism preorder collapses parallel arrows ($X \rightrightarrows Y$) to one arrow ($X \to Y$), whereas the notation $\mathcal{C}/\leftrightarrow$ alludes to collapsing objects $X \leftrightarrow Y$ to a single equivalence class $[X]$.

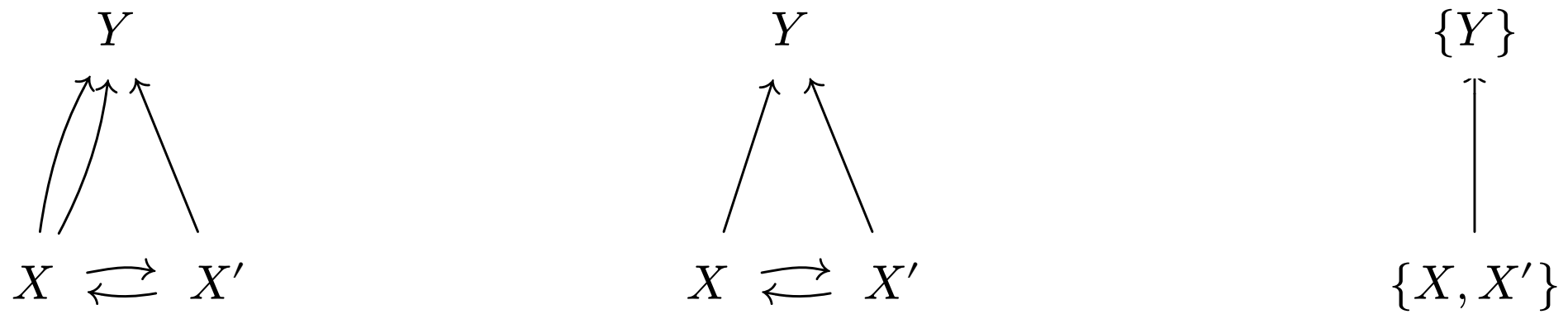

Figure 4: Collapsing a category $\mathcal{C}$ to its homomorphism preorder $\mathcal{C}/\rightrightarrows$ and then its homomorphism order $\mathcal{C}/\leftrightarrow$. Identity arrows are omitted.

**Lemma 1.1.1**: Let $\mathcal{C}$ and $\mathcal{D}$ be categories.

1. If $F\colon \mathcal{C} \to \mathcal{D}$ is a functor, then $X \mapsto FX$ induces a preorder homomorphism from $\mathcal{C}/\rightrightarrows$ to $\mathcal{D}/\rightrightarrows$ and thus $[X] \mapsto [FX]$ is a well-defined poset homomorphism $\mathcal{C}/\leftrightarrow$ to $\mathcal{D}/\leftrightarrow$.
2. If $F\colon \mathcal{C} \to \mathcal{D}$ is full – i.e., the induced map $(c \to_{\mathcal{C}} c') \to (Fc \to_{\mathcal{D}} Fc')$ is surjective for all $c, c' \in \mathcal{C}_0$ – then the aforementioned maps are *strong homomorphisms*, i.e., $X \mapsto FX$ is a strong preorder homomorphism $\mathcal{C}/\rightrightarrows \to \mathcal{D}/\rightrightarrows$ and $[X] \mapsto [FX]$ is a strong poset homomorphism $\mathcal{C}/\leftrightarrow \to \mathcal{C}/\leftrightarrow$.

*Proof*: For point one, note that if $f\colon X \to Y$, then $Ff\colon FX \to FY$, so the relation $\to$ which makes up the homomorphism preorder is preserved. Similarly, if $[X] \leq [Y]$, that means that there is an $f\colon X \to Y$, hence $FX \to FY$, hence $[FX] \leq [FY]$.

For the second point, assume now that $FX \leq FY$ in $\mathcal{D}/\rightrightarrows$, i.e., $FX \to FY$, as witnessed by a homomorphism $g$. Since $F$ is full, this map $g$ must be of the form $Ff$ for some $f\colon X \to Y$, implying that $X \leq Y$ in $\mathcal{C}/\rightrightarrows$. The argument for the second map works analogously. □

**Definition 1.1.13** (lattices): A poset $(P, \leq)$ is called a *lattice order* if any two elements have a least upper bound ("supremum") and a greatest lower bound ("infimum"). It is called *bounded* if it has a greatest element $1$ and a least element $0$.

An algebra $(L, \wedge, \vee)$ is called a *lattice* if $\wedge$ (called "meet") and $\vee$ (called "join") are both associative, commutative, and idempotent binary operations satisfying the absorptive laws $l \wedge (l \vee m) = l \vee (l \wedge m) = l$.

Such a lattice is called *distributive* if for every $l \in L$, joining with $l$ ($_ \vee l$) and meeting with $l$ ($_ \wedge l$) are both lattice homomorphisms[7].

An algebra $(L, \wedge, \vee, 0, 1)$ is called *bounded lattice* if $(L, \wedge, \vee)$ is a lattice and if all $l \in L$ satisfy $0 \wedge l = 0$, $0 \vee l = l$, $1 \wedge l = l$ and $1 \vee l = 1$.

For every lattice order we can obtain a lattice structure by setting $\wedge$ to the infimum and $\vee$ to the supremum, and from every lattice we can recover a lattice order by letting $x \leq y$ if and only if $x \wedge y = x$ (or, equivalently, $x \vee y = y$). It is not hard to check that this correspondence is well-defined and one-to-one. We will use both concepts interchangeably when we talk about certain homomorphism orders. For a reference, the reader can consult [8] or [9]

### 1.1.3 Products and Sums

To motivate the following definitions, the reader is urged to think about the following question: Assuming we work with Zermelo-Fraenkel set theory as our foundation, how is the third cartesian power $X^3$ of a set $X$ actually defined?

One answer would be: As the set of three-tuples $(x, y, z)$ where $x, y, z \in X$, of course! But we will have to continue to define what $(x, y, z)$ actually means, and this can be answered by saying $(x, y, z) := ((x, y), z)$, where $(x, y) := \{\{x\}, \{x, y\}\}$. This actually gives two answers, because we also could have chosen $(x, y, z) := (x, (y, z))$

An alternative construction would be to say that $X^3$ contains *functions* from the ordinal $3 := \{0, 1, 2\}$ to $X$, for which we also need to define what an ordered pair is.

While $X^3$ could be defined as $(X \times X) \times X$, $X \times (X \times X)$, or $3 \to X$, intuitively, all of these constructions provide the "same" kind of product. But what does that mean formally? Certainly, all of these sets are in bijection. But that is not sufficient, since for instance $A \times \mathbb{N}$ and $\mathbb{N} \times \mathbb{N}$ are in bijection for every countable set $A$, but clearly they are meant to describe different products when $A$ is finite.

The resolution here is that for each notion of product we have associated *projection maps* $\pi_1, \pi_2, \pi_3$, and the bijective correspondences between them, like $((x, y), z) \mapsto (x, (y, z))$, respect the projection maps in a manner which we will make precise shortly. So in some ways, there is not "the" product, but rather various things behaving like a product, but as we will see later, there is some universal property that pins them down up to the above established notion of "same product".

The dual notion is that of a *coproduct*, more intuitively called a *sum*.

[7] In the universal algebra sense, i.e., they have to commute with all operations, in this case $\wedge$ and $\vee$.

**Definition 1.1.14**: Let $(X_i)_{i \in I}$ be an indexed collection of objects in a category $\mathcal{C}$. A *product* of $(X_i)_i$ is a tuple $\left(P, (p_i)_i\right)$ with $P \in \mathcal{C}_0$ and $p_i \colon P \to C_i$ (called "cone over $(X_i)_i$"), such that the following universal property holds: For every other such tuple $\left(Q, (q_i)_i\right)$, there exists a unique morphism $h \colon Q \to P$ such that $h;p_i = q_i$ for every $i$. This unique morphism is called the *tupling* of the $q_i$ and denoted by $\langle q_i \rangle_i$ .

**Lemma 1.1.2**: Any two products $\left(P, (p_i)_i\right)$ and $\left(Q, (q_i)_i\right)$ of $(X_i)_i$ are isomorphic as products, i.e., there is an isomorphism $h \colon P \to Q$ such that $p_i = h;q_i$.

*Proof*: Since $\left(Q, (q_i)_i\right)$ is a product of $(X_i)_i$, there must be a homomorphism $h \colon P \to Q$ such that $h;q_i = p_i$ for every $i$. Analogously, there exists a $h' \colon Q \to P$ such that $h';p_i = q_i$ for every $i$. Thus, their composites $h;h' \colon P \to P$ and $h';h \colon Q \to Q$ are maps between products which commute with the projections, but by the uniqueness assumption they must be identities. □

**Definition 1.1.15**: We will thus denote – if it exists – "the product" of $(X_i)_i$ as

$$\left(\Pi_i X_i, (\pi_i)_i\right),$$

knowing that it is only defined up to canonical isomorphism. Giving factors explicitly, we may also denote the object as $X_1 \times X_2 \times \ldots \times X_n$. By slight abuse of notation, we will just refer to $\Pi_i X_i$ as "the product" and leave the projections implicit.

Similarly, for every $k > 0$, we will denote "the $k$th power" of $X$ by $\left(X^k, (\pi_i)_{i=0}^{k-1}\right)$. If we consider multiple powers $X^k$, $X^l$ etc., we will add a superscript $\left(\pi_i^k\right)$ to remove ambiguity.

**Definition 1.1.16**: We say a category has *finite products* if for every finite collection $X_1, \ldots, X_n$ of objects, the product $\Pi_i X_i$ exists. We say a category has *finite powers* if for every object $X$ and $n \in \mathbb{N}$, the power $X^n := \Pi_{i=1}^n X$ exists.

**Definition 1.1.17** (comparison morphism): Let $F \colon \mathcal{C} \to \mathcal{D}$ be a functor between two categories with finite powers. We define

$$\eta_{A,n} := \langle F\pi_0, \ldots, F\pi_{n-1} \rangle \colon F(A^n) \to (FA)^n.$$

**Definition 1.1.18** (argument operation): Let $\mathcal{C}$ be a category with finite powers and $A \in \mathcal{C}_0$. Given a map $\alpha \colon n \to k$, we define the pullback $\alpha^* \colon A^k \to A^n$ as the tupling $\alpha^* := \left\langle \pi_{\alpha(0)}, \ldots, \pi_{\alpha(n-1)} \right\rangle$, where $\pi_i \colon A^k \to A$ is the $i$th projection associated to the power $A^k$.

A technical observation we will need later is that argument operations "commute" with the comparison morphisms.

**Lemma 1.1.3** (cf. [3, Exercise V.2.4]): Let $F\colon \mathcal{C} \to \mathcal{D}$ be a functor between categories with finite powers, $\alpha\colon n \to n'$, and $A \in \mathcal{C}_0$. Then the following diagram commutes:

$$\begin{array}{ccc} F\left(A^{n'}\right) & \xrightarrow{F\alpha^*} & F(A^n) \\ \Big\downarrow \eta_{A,n'} & & \Big\downarrow \eta_{A,n} \\ (FA)^{n'} & \xrightarrow[\alpha^*]{} & (FA)^n \end{array}$$

*Proof*: By the universal property of the product, it suffices to show that commutativity of the square in each component of $(FA)^n$. Indeed, we have

$$\begin{aligned} F\alpha^*;\eta_{A,n};\pi_i &= F\alpha^*;\langle F\pi_0, ..., F\pi_{n-1}\rangle;\pi_i && \text{Definition of } \eta \\ &= F\alpha^*;F\pi_i && \text{Definition of tupling} \\ &= F(\alpha^*;\pi_i) && F \text{ is functor} \\ &= F\left(\left\langle \pi_{\alpha(0)}, ..., \pi_{\alpha(n-1)}\right\rangle;\pi_i\right) && \text{Definition of } \alpha^* \\ &= F\pi_{\alpha(i)} && \text{Definition of tupling} \\ &= \langle F\pi_0, ..., F\pi_{n'-1}\rangle;\pi_{\alpha(i)} && \text{Definition of tupling} \\ &= \eta_{A,n'};\pi_{\alpha(i)} && \text{Definition of } \eta \\ &= \eta_{A,n'};\left\langle \pi_{\alpha(0)}, ..., \pi_{\alpha(n-1)}\right\rangle;\pi_i && \text{Definition of tupling} \\ &= \eta_{A,n'};\alpha^*;\pi_i. && \text{Definition of } \alpha^* \end{aligned}$$

□

Analogously, we can define the *coproduct*, more intuitively called a *sum*:

**Definition 1.1.19**: Let $(X_i)_{i\in I}$ be an indexed collection of objects in a category $\mathcal{C}$. A *coproduct* or *sum* of $(X_i)_i$ is a tuple $\left(S, (s_i)_i\right)$ with $S \in \mathcal{C}_0$ and $s_i\colon S \to C_i$ (called "cone under $(X_i)_i$"), such that the following universal property holds: For every other such tuple $\left(T, (t_i)_i\right)$, there exists a unique morphism $h\colon S \to T$ such that $s_i;h = t_i$. This unique morphism is called the *cotupling* of the $s_i$ and denoted by $[s_i]_i$ .

Completely analogously[8] to products, we have

[8]Indeed, we could have defined sums as products in the opposite category – omitting the need to restate – at the expense of readability.

**Lemma 1.1.4**: Any two coproducts $\left(S, (s_i)_i\right)$ and $\left(T, (t_i)_i\right)$ of $(X_i)_i$ are isomorphic as coproducts, i.e., there is an isomorphism $h\colon S \to T$ such that $s_i;h = t_i$ for every $i$.

*Proof*: *Omitted.* □

**Definition 1.1.20**: We will thus denote – if it exists – "the sum" of $(X_i)_i$ as

$$\left(\Sigma_i X_i, (\iota_i)_i\right),$$

knowing that it is only defined up to canonical isomorphism. Giving summands explicitly, we may also denote the object as $X_1 + X_2 + \ldots + X_n$. By slight abuse of notation, we will just refer to $\Sigma_i X_i$ as "the sum" and leave the injections $\iota_i$ implicit.

**Definition 1.1.21**: We say a category has *finite sums* if for every finite collection $X_1, \ldots, X_n$ of objects, the sum $\Sigma_i X_i$ exists.

### 1.1.4 Adjunctions: Free Objects, Galois Connections, and More

The concept of an adjunction formalizes a phenomenon that is pervasive in mathematics, namely whenever some kind of "free structure" – for instance a free group – is involved. It should be noted however that the concept is vastly more general and encompasses cases which do not quite fit the "free" naming: For instance, a *Galois connection* between posets (see Section A.1) is a special case where the involved categories are posets. For another example, consider the tensor product of modules over a commutative ring. More examples shall follow after the formal definition.

**Definition 1.1.22** (Adjunction): Let $\mathcal{C}, \mathcal{D}$ be categories and $F\colon \mathcal{D} \to \mathcal{C}$ and $G\colon \mathcal{C} \to \mathcal{D}$ be functors. We say that $(F, G)$ is an adjoint pair, or an *adjunction* – in symbols $F \dashv G$ – if we have a natural isomorphism

$$\Phi\colon (F_ \to_{\mathcal{C}} _) \to (_ \to_{\mathcal{D}} G_)$$

of functors $\mathcal{D}^{\mathsf{op}} \times \mathcal{C} \to \underline{\mathsf{Set}}$. If $F \dashv G$, we also say that *$F$ is left adjoint to $G$* and *$G$ is right adjoint to $F$*. We call $F$ a *left adjoint* if it admits a right adjoint $G$ and vice versa.

We may also denote functors visually in a diagram like this:

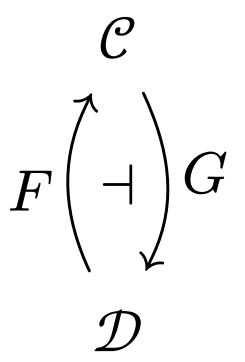

The compact definition above deserves to be spelt out concretely. To give a natural isomorphism is to give a family $\left(\Phi_{X,A}\right)_{X\in\mathcal{D}_0, A\in\mathcal{C}_0}$ of bijections

$$\Phi\colon (FX \to_{\mathcal{C}} A) \to (X \to_{\mathcal{D}} GA)$$

such that for all $l\colon X \to_{\mathcal{D}} X'$ and all $r\colon A \to_{\mathcal{C}} A'$, the following diagrams commute:

$$\begin{array}{ccc} (FX \to_{\mathcal{C}} A) & \xrightarrow{\Phi_{X,A}} & (X \to_{\mathcal{D}} GA) \\ \Big\downarrow r_* & & \Big\downarrow (Gr)_* \\ (FX \to_{\mathcal{C}} A') & \xrightarrow{\Phi_{X,A'}} & (X \to_{\mathcal{D}} GA') \end{array} \qquad \begin{array}{ccc} (FX' \to_{\mathcal{C}} A) & \xrightarrow{\Phi_{X',A}} & (X' \to_{\mathcal{D}} GA) \\ \Big\downarrow (Fl)^* & & \Big\downarrow l^* \\ (FX \to_{\mathcal{C}} A) & \xrightarrow{\Phi_{X,A}} & (X \to_{\mathcal{D}} GA) \end{array}$$

Naturality of $\Phi$ is enough, because if all components of a natural transformation are isomorphisms, then the collection of all the inverses will be natural as well. The argument for that is analogous to the well-known argument that bijective homomorphisms between algebraic structures are also isomorphisms.

The first diagram says that given $h\colon FX \to_{\mathcal{C}} A$, first postcomposing with $r\colon A \to A'$ and then applying the translation $\Phi$ to get a homomorphism $X \to_{\mathcal{D}} GA'$ should yield the same result as applying the correspondence first and then postcomposing with $Gr\colon GA \to GA'$. The second diagram gives a dual condition with regards to precomposing $l\colon X \to X'$. Intuitively speaking, the naturality condition on $\Phi$ ensures that the construction does not depend on specific information about the objects $X$ or $A$ involved, but is completely canonical.

The importance of these naturality conditions can be witnessed in many ways, most prominently by the uniqueness theorem:

**Proposition 1.1.5**:
- If $F \dashv G$ and $F \dashv G'$, then $G$ and $G'$ are naturally isomorphic.
- If $F \dashv G$ and $F' \dashv G$, then $F$ and $F'$ are naturally isomorphic.

A proof will be provided in Section A.3.

Another powerful statement is that *left adjoints preserve colimits* – which for instance implies that if the functor $F$ is left adjoint, then $F(X+Y) \cong (FX) + (FY)$, and not only as objects, but also "as sums": If we apply the functor $F$ to the canonical inclusions $\iota_X\colon X \hookrightarrow X+Y$ and $\iota_Y\colon Y \hookrightarrow X+Y$ implicit to the sum $X+Y$, we obtain morphisms $F\iota_X\colon FX \to F(X+Y)$ and $F\iota_Y\colon FY \to F(X+Y)$. Cotupling these morphisms we obtain $[F\iota_X, F\iota_Y]\colon FX + FY \to F(X+Y)$, and the statement is that this morphism is an isomorphism. It is not hard to verify that this implies that the triple $(F(X+Y), F\iota_X, F\iota_Y)$ is a coproduct, and that $[F\iota_X, F\iota_Y]$ is an isomorphism of coproducts as defined in Lemma 1.1.4.

Now neither of these statements is strictly necessary for the results obtained in this thesis, but in the opinion of the author they provide relevant context. For instance, in Chapter 7, we will define an "exponential object", which is the right adjoint to the product functor $_ \times \mathcal{M}$ in the category of interest. From such a description, we not only know that this object is defined uniquely, but also that in this category, the product automatically distributes over the sum.

*Example*: For groups, we have the correspondence

$$FX \to_{\underline{\mathsf{Grp}}} G \quad \cong \quad X \to_{\underline{\mathsf{Set}}} UG.$$

This is an adjunction between the free group functor $F\colon \underline{\mathsf{Set}} \to \underline{\mathsf{Grp}}$ and the "forgetful functor" $U\colon \underline{\mathsf{Grp}} \to \underline{\mathsf{Set}}$, which assigns a group $G$ its underlying set $UG$. It is not hard to verify that this correspondence is natural.

*Example (for those familiar with universal algebra)*: If $\mathcal{V}$ is a variety of $\Omega$-algebras, we have the correspondence

$$F_{\mathcal{V}}X \to_{\mathcal{V}} \mathbb{A} \quad \cong \quad X \to_{\underline{\mathsf{Set}}} U\mathbb{A}$$

between the free algebra functor $X \mapsto F_{\mathcal{V}}X := T_{\Omega}(X)/\mathrm{Eq}(\mathcal{V})$ and the forgetful functor $U$.

*Example*: Let $\underline{\mathsf{Gr}}$ denote the category of (directed) graphs. The forgetful functor $U\colon \underline{\mathsf{Gr}} \to \underline{\mathsf{Set}}$ has a left adjoint, which sends $X$ to the *discrete graph* on $X$: Indeed, every function $X \to UG$ uniquely corresponds to a homomorphism from the discrete graph of $G$.

However, $U$ has also a *right adjoint*: If we send $X$ to the graph with all the edges, i.e., the clique with loops, then every function $UG \to X$ uniquely corresponds to a homomorphism into this graph.

*Example (for topologists)*: Analogously, the forgetful functor $U\colon \underline{\mathsf{Top}} \to \underline{\mathsf{Set}}$ from topological spaces to sets has both a left adjoint (assigning the discrete topology $\mathcal{P}(X)$) and a right adjoint (assigning the indiscrete topology $\{\emptyset, X\}$).

*Example*: In the category ${}_R\underline{\mathsf{Mod}}$ of $R$-modules and $R$-linear maps, we have the "tensor-hom" correspondence

$$M \otimes_R N \to_{{}_R\underline{\mathsf{Mod}}} K \quad \cong \quad M \to_{{}_R\underline{\mathsf{Mod}}} \mathrm{Hom}_R(N, K)$$

establishes an adjunction $_ \otimes_R N \dashv \mathrm{Hom}_R(N, _)$. This example is special, because it is an adjunction of functors between the same category – and in particular it does not fall into the "free-forget" pattern.

*Example*: Considering posets $(P, \leq)$ and $(Q, \leq)$ as categories as described previously, we note that a functor $F\colon P \to Q$ should map a morphism $p \underset{P}{\to} p'$, i.e., the arrow witnessing the inequality $p \leq p'$, to a morphism $Fp \underset{q}{\to} Fp'$, i.e., the arrow witnessing the inequality $Fp \leq Fp'$. Thus, functors are just maps $p \mapsto Fp$ which are monotone. Furthermore, an adjunction $F \dashv G$ between $F\colon Q \to P$ and $G\colon P \to Q$ is a correspondence

$$Fp \to_Q q \quad \cong \quad p \to_Q Gq,$$

whose existence precisely means that $Fp \leq q$ if and only if $p \leq Gq$. That is to say that $F$ and $G$ form a *(monotone) Galois connection* between $P$ and $Q$.

Following the same ideas, we would see that an *antitone* Galois connection is a pair of functors $F \dashv G$ where $F\colon P \to Q^{\mathrm{op}}$ and $G\colon Q^{\mathrm{op}} \to P$, since antitone (i.e., order-reversing) maps correspond to contravariant functors.

Basic facts about antitone Galois connections will be provided in Section A.1, as we will utilize the concept in the section about polymorphisms and preservation of relation pairs.

## 1.2 (Promise) Constraint Satisfaction Problems

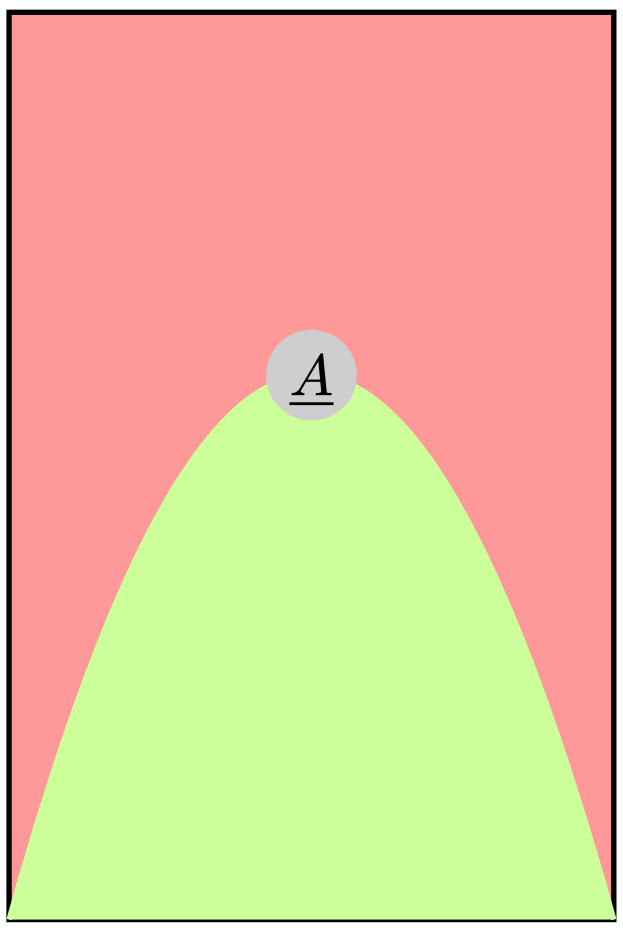

Figure 5: CSP($\underline{A}$) visualized as a lower cone in the homomorphism preorder of $\tau$-structures.
A program should answer "yes" if the input lies in the green area, and "no" if it lies in the red area.

While abstract and concrete minions can easily be defined and analyzed in isolation without many dependencies to other areas of mathematics, it is necessary to provide the context in which they rose to prominence. The intent of this section is to give a high-level, informal overview of the relevant background in constraint satisfaction problems, where the notion has been introduced.

To maintain the focus of this thesis, we will not formally introduce every concept defined here. For a comprehensive introduction to finite-domain constraint satisfaction problems, the reader is referred to the lecture notes [10], and the infinite-domain case is covered extensively in [2]. For an introduction to universal algebra, i.e., the notions of "algebra", "identity", and a precise elaboration on how clones relate to varieties, we refer to [11].

Although much of (pure) mathematics is concerned with determining whether a solution to a problem exists – and whether it is unique – at some point, existence guarantees are not enough: In order to determine such a solution, there must be a (constructive) method. Complexity theory is concerned with the question of how hard it is – measured asymptotically in the "size" of the input – to determine the desired solution.

We shall focus on the existence part, formally called a *decision problem*: That is, given some object $G$ as input, the problem is to determine whether $G$ satisfies some property – or more succinctly, to determine whether $G$ corresponds to some family $H$.

An important family of such problems are the so called *constraint satisfaction problems*, *CSP* for short, which can be formulated mathematically as *homomorphism problems* between relational structures. That is, fixing a finite signature $\tau$ and a $\tau$-structure $\underline{H}$, called the *CSP template*, the goal is to determine whether an input $\tau$-structure $\underline{G}$, called the *CSP instance*, has homomorphism into $\underline{H}$. In other words, we decide the class

$$\mathrm{CSP}(\underline{H}) := \{\underline{G}\ \ \tau\text{ -structure} \mid \underline{G} \to \underline{H}\}.$$

In Figure 5, we see this illustrated quite concisely: First, we draw all $\tau$-structures as a Hasse diagram ordered by homomorphicity. That means that if $\underline{A} \to \underline{B}$, we would draw $\underline{A}$ below $\underline{B}$ and connect them with a line. In that visualization, the class of CSP($\underline{A}$)-instances for which a program should answer "yes" is precisely the area "below $\underline{A}$", shaded in green. The complement region, shaded red, should encompass all the "no"-instances.

To explain why this is called a *constraint satisfaction problem*, it is helpful to view the input structure as a logical formula

$$\exists (x_g)_{g\in G}. \bigwedge_{\substack{R\in\tau \\ e\in R^{\underline{G}}}} R\left(x_{e_1}, ..., x_{e_{\mathrm{ar}\,R}}\right).$$

A homomorphism from $\underline{G}$ to $\underline{H}$ is then precisely an assignment of an element in $H$ such that this constraint formula is satisfied in $\underline{H}$.

Crucially, a wide variety of classical satisfiability problems can be formulated as a CSP:

- The *3-satisfiability problem*, denoted 3SAT, is concerned with deciding satisfiability of formulas in conjunctive normal form whose conjuncts consist of at three disjunctions ([12, Chapter 7]). This works by providing ternary relations

$$\begin{aligned} R_1(x, y, z) &:= x \vee y \vee z, \\ R_2(x, y, z) &:= \neg x \vee y \vee z, \\ R_3(x, y, z) &:= \neg x \vee \neg y \vee z, \\ R_4(x, y, z) &:= \neg x \vee \neg y \vee \neg z, \end{aligned}$$

  on the boolean domain for the template, and rewriting a given formula $\psi$ to only use these relations as atomic formulas to obtain an instance structure.
- The variants 1-in-3-3SAT, NAE-3SAT are defined as the problem of finding a boolean assignment which satisfies the formula, under the additional constraints that this assignment shall produce precisely one thruth value per conjunct, or shall exclude the truth values $(0, 0, 0)$ and $(1, 1, 1)$ from the conjuncts, respectively.
- The *positive* 3SAT, 1-in-3-3SAT, and NAE-3SAT refers to the variant where no negations are allowed in the clauses. For 3SAT, this amounts to only considering the relation $R_1$, as all other will never be populated in a CSP instance.
- The problem HORN-SAT is concerned with the satisfiability of formulas which are conjunctions of *horn clauses*, i.e., formula of the form $x_1 \wedge x_2 \wedge ... \wedge x_k \Rightarrow x_0$.
- The *3-coloring* problem is determined whether a given graph $\underline{G}$ is three-colorable, i.e., has a homomorphism to the three-color-graph $(\{\mathrm{red}, \mathrm{green}, \mathrm{blue}\}, \neq)$ ($\cong K_3$)

### 1.2.1 The CSP Dichotomy and the "Algebraic Approach"

In this section we want to introduce the notion of complexity of a CSP, and want to present the universal algebra-based approach to analyzing this complexity. On a high level, we will associate to CSP templates $\underline{H}$ their *polymorphism clone* $\mathrm{Pol}(\underline{H})$ – an inherently algebraic object – and if we know that $\mathrm{Pol}(\underline{H})$ and $\mathrm{Pol}(\underline{H'})$ are related by a special kind of map, then we can use this information to build an efficient[9] reduction from CSP($\underline{H'}$) to CSP($\underline{H}$), which roughly speaking gives us a way to solve CSP($\underline{H'}$) using a program which solves CSP($\underline{H}$). This implies that to solve the latter cannot be inherently easier than to solve the former. But let us start from the basics.

[9]Meaning "logarithmic-space" (L), but for our purposes "polynomial-time" (P) is sufficient.

❦

Associated to a decision problem is its complexity, formally communicated by membership in certain complexity classes. We define $P$ to be the class of all languages which are decided by a deterministic Turing machine with polynomial runtime, and NP to be the class of all languages which are decided by a *nondeterministic* Turing machine with polynomial runtime.

We say that a language $A$ is *polynomial-time reducible* to $B$ [12, Def. 7.29], written $A \leq_p B$, if there is a polynomial-time computable[10] function $f\colon \Sigma^* \to \Sigma^*$ such that $w \in A$ if and only if $f(w) \in B$.

More specifically, in the context of CSPs, we say that CSP($\underline{H}$) is *polynomial-time reducible* to CSP($\underline{H'}$) if there is a polynomial-time algorithm which translates an input $\underline{G}$ to the first CSP to an input $\underline{G'}$ to the second CSP in such a way that $G \to H$ if and only if $G' \to H'$.

This gives rise to the *polynomial-time reduction preorder* on languages (i.e., decision problems): A language $L$ is called *hard* with respect to a class $C$, or $C$*-hard*, if it is an upper bound to $C$ in this preorder, and $C$-complete if it is in $C$ and $C$-hard.

It is known that 3SAT, 1-in-3-3SAT, and NAE-3SAT as well as their positive variants are NP -complete [2, Example 1.2.2], while HORN-SAT is P-complete [2, Section 1.6.7].

Loosely speaking, the famous open P versus NP-problem adresses the question of whether every "hard" problem is actually "hard", or whether we just did not find a good algorithm. More precisely, we ask whether the inclusion P $\subseteq$ NP is strict. Most intriguingly, Ladner [13] proved that when P $\neq$ NP, then there also must be intermediate problems, i.e., problems which are neither "hard" (NP -complete) nor "easy" (in P).

In CSPs, these intermediate theorems cannot be found, as stated in the following theorem. Previously conjectured by Feder and Vardi in 1993 [14], two proofs have been announced in 2017 independently by Bulatov and by Zhuk.

**Theorem[11] 1.2.I** (Bulatov [15]; Zhuk [16] 2017): The decision problem CSP($\underline{H}$) is either in P or NP-complete.

Both solutions to this conjecture used well-developed universal algebra to characterize NP-hardness; the remaining duty is to provide a polynomial-time algorithm which solves the CSP in all the other cases.

A key tool in deriving hardness results (showing that CSP($\underline{H}$) is NP-hard) or tractability results (showing that CSP($\underline{H}$) is in $P$) is the following statement:

[10] i.e., a function computable by a deterministic Turing machine with polynomial runtime [12, Def. 7.28].

[11] *Editorial remark*: The statements and definitions of this section bear no logical connections to other statements and definitions in this thesis. The numbering is adapted to emphasize this "explanatory" role.

**Proposition 1.2.11** ([17, Theorems 13, 15]): Let $\underline{A}$ and $\underline{B}$ be $\sigma$- and $\tau$-structures, respectively. If there is clone homomorphism $\mathrm{Pol}(\underline{A}) \to \mathrm{Pol}(\underline{B})$ which preserves arities and composition, then

1. there is a primitive positive interpretation [2, Chapter 3] of $\underline{B}$ in $\underline{A}$
2. there is an efficient[12] reduction from $\mathrm{CSP}(\underline{B})$ to $\mathrm{CSP}(\underline{A})$.

Here $\mathrm{Pol}(\underline{A}) = (\underline{A}^n \to A)_{n>0}$ is the *polymorphism clone* of $\underline{A}$, where the attribute "clone" highlights that this set of functions contains the projections and is closed under composition.

*Proof sketch of the above proposition*: Without formally introducing all the relevant notions, let us sketch a proof of this proposition:

- A clone homomorphism $\mathrm{Pol}(\underline{A}) \to \mathrm{Pol}(\underline{B})$ allows us to make polymorphisms of $\underline{A}$ act on $B$.
- This means we can beget $B$ with the structure of a $\mathrm{Pol}(\underline{A})$-algebra $\mathbb{B}$, which satisfies all the identities which hold in the $\mathrm{Pol}(\underline{A})$-algebra $\mathbb{A}$.
- Birkhoff's theorem now guarantees that this new algebra $\mathbb{B}$ arises as a homomorphic image of a substructure of a power of $\mathbb{A}$.
- A classical duality, which relates polymorphisms with the relations they preserve, can now be applied to see that $\underline{B}$ has something called a *primitive positive interpretation* in $\underline{A}$.
- To such an interpretation $\underline{B} = I(\underline{A})$ we can associate an efficiently-computable *gadget replacement* $G$, which turns a $\tau$-structure into a $\sigma$-structure in such a way that $\underline{C} \to_\tau \underline{B} = I(\underline{A})$ if and only if $G(\underline{C}) \to_\sigma \underline{A}$.
- This procedure gives us an efficient way to answer $\mathrm{CSP}(\underline{B})$ given an oracle which answers $\mathrm{CSP}(\underline{A})$, viz. a reduction. □

There is a bit of an oversight in this description: Namely, it might quite easily be that $\mathrm{Pol}(\underline{A})$ does *not* have a clone homomorphism to $\mathrm{Pol}(\underline{B})$, but that $\underline{B}$ is homomorphically equivalent to a structure $\underline{B'}$ for which this is the case. In this scenario, where $\mathrm{Pol}(\underline{A}) \to \mathrm{Pol}(\underline{B'})$, we might be tempted to construct a mapping $h\colon \mathrm{Pol}(\underline{B'}) \to \mathrm{Pol}(\underline{B})$ as follows: Since the homomorphic equivalence is witnessed by structure homomorphisms $g\colon \underline{B} \to \underline{B'}$ and $f\colon \underline{B'} \to \underline{B}$, we can take a polymorphism $p \in \mathrm{Pol}(\underline{B'})$ of some arity $n$, and form the composite

$$\underline{B}^n \xrightarrow{g^{\times n}} \underline{B'}^n \xrightarrow{p} \underline{B'} \xrightarrow{f} \underline{B},$$

which, being a composition of homomorphisms, is in $\mathrm{Pol}(\underline{B})$. However, this map $h$ is not compatible with the clone structure[13];

for instance, if $f$ and $g$ are not inverse to each other, the identity over $B'$ need not be mapped to the identity over $B$. However, $h$ does preserve a subset of the information contained in the polymorphism clone: Namely, the *height-one identities*.

[12] As before, the reduction is in $L$, but for our purposes, the fact that it is in $\mathsf{P}$ suffices.

[13] i.e., commutes with composition and preserves projections, in particular the identity.

**Definition 1.2.i**: A map $h\colon \mathcal{A} \to \mathcal{B}$ between function clones with domain $A$ and $B$, respectively, is said to *preserve height-one identities*, if for all $f \in \mathcal{A}$ of arity $n$ and $g \in \mathcal{B}$ of arity $k$ and coordinate reassignments $\alpha\colon n \to l$ and $\beta\colon k \to l$, whenever

$$f\left(x_{\alpha_0}, ..., x_{\alpha_{n-1}}\right) = g\left(x_{\beta_0}, ..., x_{\beta_{k-1}}\right)$$

holds in $\mathcal{A}$ for all possible choices of tuples $\underline{x} \in A^l$, then so does

$$h(f)\left(x_{\alpha_0}, ..., x_{\alpha_{n-1}}\right) = h(g)\left(x_{\beta_0}, ..., x_{\beta_{k-1}}\right).$$

It is common in universal algebra to denote such an equation universally quantified in the input tuple as

$$f\left(x_{\alpha_0}, ..., x_{\alpha_{n-1}}\right) \approx g\left(x_{\beta_0}, ..., x_{\beta_{k-1}}\right),$$

where the $x_i$ are now interpreted as variable symbols instead of as specific elements. In other words, the left- and right-hand side expressions are now *terms*.

Crucially – and somewhat surprisingly – Barto, Opršal and Pinsker proved a weaker version of Birkhoff they dubbed “Linear Birkhoff”, which ensures precisely enough to make the above proof sketch work, albeit with a few modifications:

**Proposition 1.2.iii** (“Linear Birkhoff”, [18, Proposition 5.3]): If $\mathcal{A} \subseteq \mathcal{O}(A)$ and $\mathcal{B} \subseteq \mathcal{O}(B)$ are function clones over finite sets $A, B$ and there is a function $\mathcal{A} \to \mathcal{B}$ preserving arities and height-one identities (“height-one map"), then $\mathcal{B}$ is an expansion of a reflection of a finite power of $\mathcal{A}$, viewed as a $\mathcal{A}$-algebra over the domain $A$.

Here, a *reflection* of an algebra $\mathbb{A}$ along maps $f\colon A \to B$ and $g\colon B \to A$ is an algebra over $B$ whose operations are formed as follows: For each algebra operation $h\colon A^n \to A$ on h, we pre- and post-compose with $g$ and respectively $h$ to obtain the algebra operation $g^{\times n};h;f\colon B^n \to B$ on $B$. This operation captures homomorphic equivalence, and it is easy to see that it preserves height-one identities.

Applying duality results similar to the statement above, we can conclude

**Proposition 1.2.iv**: If there is a map $\mathrm{Pol}(\underline{A}) \to \mathrm{Pol}(\underline{B})$ preserving arities and height-one identities, then there is a logspace reduction $\mathrm{CSP}(\underline{B}) \to \mathrm{CSP}(\underline{A})$.

Crucially, this gives us two options – we can derive *hardness* of $\mathrm{CSP}(\underline{A})$ if we have a homomorphism $\mathrm{Pol}(\underline{A}) \to \mathrm{Pol}(\underline{B})$ where $\underline{B}$ is a known hard template $\underline{B}$, or we can derive *tractability* (being in P) of $\mathrm{CSP}(\underline{B})$ if we have a homomorphism $\mathrm{Pol}(\underline{A}) \to \mathrm{Pol}(\underline{B})$ from a known tractable template $\underline{A}$. For instance:

1. If there is a height-one map $\mathrm{Pol}(\underline{A}) \to \mathrm{Pol}(K_3)$, then $\mathrm{CSP}(\underline{A})$ must be NP-hard.
2. If there is a height-one map from $\mathrm{Pol}(*) \to \mathrm{Pol}(\underline{B})$, where $*$ is the discrete graph with a loop, then $\mathrm{CSP}(\underline{B})$ must be tractable.

The last key observation, motivating the definition of a minion, is the following.

**Proposition 1.2.v**: A map $h\colon \mathcal{A} \to \mathcal{B}$ preserves height-one identities if and only if it preserves the action of *minor operations*: A minor operation $\alpha\colon n \to k$ "acts" on $\mathcal{A}$ by associating to $f \in \mathcal{A}$ of arity $n$ the $k$-ary function

$$(f\alpha)(x_0, \dots, x_{k-1}) := f\big(x_{\alpha(0)}, \dots, x_{\alpha(n-1)}\big),$$

and this action is *preserved by* $h$ if $h(f\alpha) = h(f)\alpha$.

*Proof*: If $h$ preserves height-one identities, then interpreting equations $g = f\alpha$ in $\mathcal{A}$ as height-one identities gives us that $h(g) = h(f)\alpha$ must hold in $\mathcal{B}$, which is the desired homomorphicity.

For the reverse direction, assume $h$ preserves minor operations, and assume $f\alpha = g\beta$ holds in $\mathcal{A}$. Then $h(f)\alpha \overset{*}{=} h(f\alpha) = h(g\beta) \overset{*}{=} h(g\beta)$, where the equalities $(*)$ hold due to homomorphicity. □

The proof can deservedly be called "symbol-pushing", as it does not contain any original idea. However, the perspective shift it provides is interesting nevertheless: Instead of considering functions which preserve certain identities, we can now talk about functions *compatible with certain operations* – that is, we have a notion of *homomorphism* between something akin to an algebraic structure!

### 1.2.2 Promise CSPs and Minions

A generalization of the CSP is the *promise CSP*, which is given for any *pair* of $\tau$-structures $(\underline{A}, \underline{B})$ with $\underline{A} \to \underline{B}$. Its decision variant asks for a given $\tau$-structure $\underline{C}$ whether $\underline{C} \to \underline{A}$ – under the *promise* that if $\underline{C} \nrightarrow \underline{A}$, then $\underline{C} \nrightarrow \underline{B}$. As visualized in Figure 6, the promise amounts to excluding the region "below" $\underline{B}$, but "not below" $\underline{A}$.

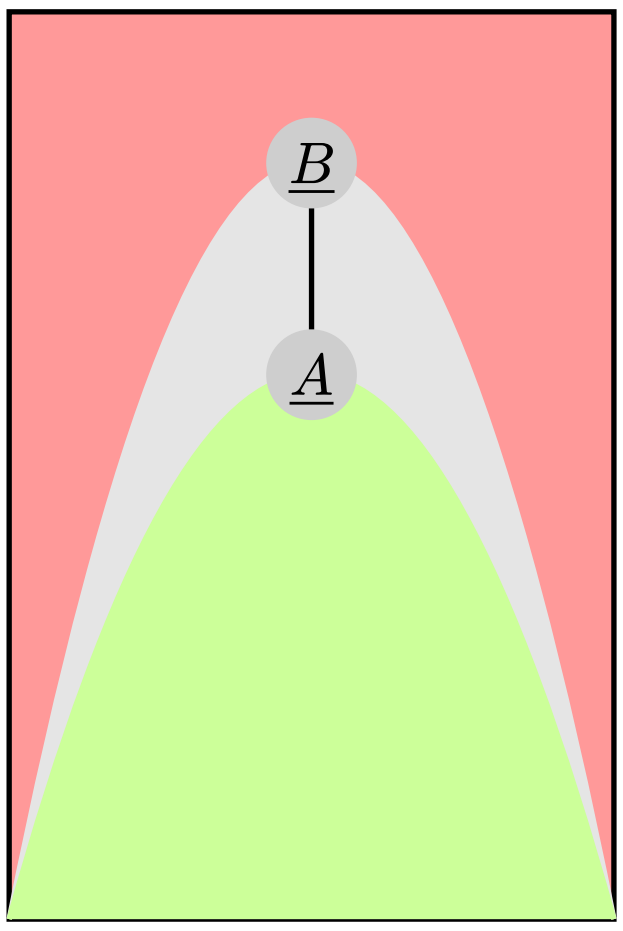

Figure 6: PCSP$(\underline{A}, \underline{B})$ visualized in the homomorphism preorder of $\tau$-structures.
A program should behave like in the CSP (Figure 5) for the red and green areas, but the answer in the gray in-between area does not matter, since the promise is that an input structure never lies in that region.

The algebraic approach has been extended by Bulín, Barto, Opršal and Krokhin in [19] in the form of the following results:

**Proposition 1.2.vi** ([19, Theorem 3.1]): If $\mathrm{Pol}(\underline{A}, \underline{B}) \to \mathrm{Pol}(\underline{A}', \underline{B}')$, then there is a log-space reduction from $\mathrm{PCSP}(\underline{A}', \underline{B}')$ to $\mathrm{PCSP}(\underline{A}, \underline{B})$.

Here, $\mathrm{Pol}(\underline{A}, \underline{B}) = (\underline{A}^n \to \underline{B})_{n \geq 1}$ is the object of all polymorphisms[14] *from* $\underline{A}$ *to* $\underline{B}$ – which crucially does not carry the structure of a clone, because the difference in domain and codomain obstructs composition. However, it is still possible to talk about height-one identities, because statements like $f(x, y) = g(y, y)$ are well-defined. More formally, $\mathrm{Pol}(\underline{A}, \underline{B})$ still carries an action of *minor operations*, because we are still able to permute, identify, and add dummy arguments.

As we will define in the next section, that is what makes $\mathrm{Pol}(\underline{A}, \underline{B})$ – just like $\mathrm{Pol}(\underline{A})$ – into a *minion*. We will henceforth call arity- and height-one preserving maps *minion homomorphisms*, although we will define them more generally in the next chapter.

❦

By applying reductions from gap label cover, a known NP-complete problem, we can establish the following hardness criterion:

**Proposition 1.2.vii** ([19, Proposition 5.14]): Let $\tau$ be a finite signature and $\underline{A}$, $\underline{B}$ be $\tau$-structures over finite domain. If there is a minion homomorphism from $\mathrm{Pol}(\underline{A}, \underline{B})$ to a minion $\mathcal{M}$ of bounded essential arity, then $\mathrm{PCSP}(\underline{A}, \underline{B})$ is NP-hard.

We will define what bounded essential arity means and discuss its implications in Section 2.6; for the time being, it is only relevant to keep in mind that minions of bounded essential arities are "hard".

One more thing deserves attention: While the naive conjecture would be that the promise CSP between two hard structures would also be hard, there is actually a counterexample:

**Proposition 1.2.viii** ([19, Example 7.15]): Let

$$\underline{T} := \left(2; \left\{ \begin{pmatrix}1\\0\\0\end{pmatrix}, \begin{pmatrix}0\\1\\0\end{pmatrix}, \begin{pmatrix}0\\0\\1\end{pmatrix} \right\}\right), \qquad \underline{H_2} := \left(2; \left\{ \begin{pmatrix}1\\0\\0\end{pmatrix}, \begin{pmatrix}0\\1\\0\end{pmatrix}, \begin{pmatrix}0\\0\\1\end{pmatrix}, \begin{pmatrix}0\\1\\1\end{pmatrix}, \begin{pmatrix}1\\0\\1\end{pmatrix}, \begin{pmatrix}1\\1\\0\end{pmatrix} \right\}\right)$$

be the templates for positive 1-in-3-3SAT and NAE-3SAT, respectively. Both $\mathrm{CSP}(\underline{T})$ and $\mathrm{CSP}(\underline{H_2})$ are NP -complete, however $\mathrm{PCSP}(\underline{T}, \underline{H_2})$ is in P.

This can be seen from something called an "affine integer programming relaxation", which algebraically can be viewed from a minion homomorphism from $\mathrm{Pol}(\underline{Z_{\mathrm{aff}}})$. Here, $\underline{Z_{\mathrm{aff}}}$ is the relational structure over the integers $\mathbb{Z}$ whose relations are all (solution sets of) affine equations [19, Definition 7.14, Theorem 7.19].

[14] These objects are sometimes called *weak polymorphisms*[20].

Even more remarkably, the following holds:

**Proposition 1.2.ix** ([19, Section 8]): Let $\tau$ be a finite signature and let $\underline{A}$ be a $\tau$-structure with finite domain. If $\mathrm{Pol}(\underline{A}) \to \mathrm{Pol}\left(\underline{T}, \underline{H_2}\right)$, then $\mathrm{CSP}(\underline{A})$ is NP-hard.

This implies that the tractability cannot be witnessed by a reduction from a tractable *finite-domain* promise CSP. Nonetheless, the finite-domain case will be the main focus of the thesis.

# 2 Minions

As discussed in the preliminaries, associated to a function $f\colon A^n \to B$ and a map $\alpha\colon n \to k$ is the associated *$\alpha$-minor of $f$*, a function $f\alpha\colon A^k \to B$, defined by

$$(f\alpha)(x_0, ..., x_{k-1}) := f\left(x_{\alpha(0)}, ..., x_{\alpha(n-1)}\right).$$

Viewing an input tuple $\underline{x} = (x_0, ..., x_{k-1}) \in A^k$ as a map $k \to A$, we can interpret the tuple $\left(x_{\alpha(0)}, ..., x_{\alpha(n-1)}\right)$ as the precomposition of $\underline{x}$ with $\alpha$, i.e., $\alpha;\underline{x}$, because we first map $i$ to $\alpha(i)$, and subsequently look up the $\alpha(i)$th component of $x$. This allows for much more concise notation:

$$f\alpha := f(\alpha;_).$$

Apart from being short, this notation also makes a few properties quite obvious:

1. $f\,\mathrm{id}_n = f(\mathrm{id}_n;_)$.
2. $f(\alpha;\beta) = f((\alpha;\beta);_) = f(\alpha;(\beta;_)) = (f\alpha)\beta$.
3. If $f\colon \underline{A}^n \to \underline{B}$ is a homomorphism of $\tau$-structures, we know that the precomposition $\alpha_\colon \underline{A}^k \to \underline{A}^n$ is a homomorphism because it is just a tupling of the projections $\pi^k_{\alpha(i)}\colon \underline{A}^k \to \underline{A}$, which are always homomorphisms. Therefore, $f(\alpha;_)$ is a composition of homomorphisms, so a homomorphism as well.
4. If $A = m$ is a natural number and we interpret elements of $A^n$ as minor operations $\beta\colon n \to m$, we know that $f(\beta) = f(\beta;\mathrm{id}_m) = (f\beta)(\mathrm{id}_m)$.

The first two points are a *functoriality* property, and motivate the definition of an abstract minion. The third point describes that even when the sets $A$ and $B$ are equipped with *extra structure*, the object $\mathrm{Pol}(\underline{A}, \underline{B})$ will be a minion. Finally, the last point is a key insight which we will exploit extensively in Chapter 4.

We should remark that the usual way to denote the $\alpha$-minor of a function $f$ is either $f^\alpha$ (e.g. [21]) or $f_\alpha$ (e.g. [10, Definition 9.2]). However, to emphasize the interplay with the coordinates, and as we will study minions as *abstract* objects, we chose to give the minor operation the status of a "first-class citizen" by placing it on the baseline. Furthermore, given the frequent ocurrence of minor operations in this thesis, this should aid in readability.

## 2.1 Abstract Minions

As the functoriality property strongly suggests, we can define minions as functors:

**Definition 2.1.1**: Let $\mathbb{F}_{>0}$ denote the category of all nonzero finite ordinals with arbitrary set-maps as homomorphisms. We will call the morphisms $\alpha\colon n \to k$ of $\mathbb{F}_{>0}$ *minor operations*, and the elements $i \in n$ a *coordinate* of $n$.

We define the *category of minions*, denoted $\underline{\mathsf{Min}}$, to be the functor category $[\mathbb{F}_{>0}, \underline{\mathsf{Set}}]$. An *abstract minion*, in this thesis just called *minion*, is an object of this functor category, usually denoted as $\mathcal{M}$, $\mathcal{L}$, and so on. For a minion $\mathcal{M}$, we call the evaluation $\mathcal{M}_n$ of this functor at the object $n \in \mathbb{F}_{>0}$ the *set of arity $n$* and its inhabitants *elements of arity $n$*.

Of course, a coordinate is just a number, making the definition feel a little pointless. However, that definition is better viewed as a "concept with an attitude" – that is, when we talk about a number as a "coordinate", it is clear that it will be used in the context of arity-$n$ elements of a minion. This word will carry more meaning in the later definitions.

Spelling out Definition 2.1.1, we see that a minion is a sequence of sets $(\mathcal{M}_n)_{n\geq 1}$ together with an assignment

$$\begin{array}{ccc} n & & \mathcal{M}_n \\ \big\downarrow{\scriptstyle \alpha} & \rightsquigarrow & \big\downarrow{\scriptstyle \mathcal{M}\alpha} \\ k & & \mathcal{M}_k \end{array}$$

where the induced function $\mathcal{M}\alpha$ maps $f \in \mathcal{M}_n$ to the minor $(\mathcal{M}\alpha)(f)$, which we will rather denote as $f\alpha$.

*Remark*: It is also possible to define minions with a more general notion of "arity" by allowing arbitrary finite sets. A minion would then be a functor $\underline{\mathsf{FinSet}}_{>\emptyset} \to \underline{\mathsf{Set}}$, where $\underline{\mathsf{FinSet}}_{>\emptyset}$ is the category of all nonempty finite sets and functions.

This is an equivalent[15] perspective with some advantages, which we will make use of in Section 3.1.

## 2.2 Polymorphism-, concrete-, and function minions

**Definition 2.2.1** (Polymorphism minion in arbitrary categories): Let $\mathcal{C}$ be a category with (finite) powers, i.e., for every object $A \in \mathcal{C}_0$ and natural number $n$, the $n$th power $A^n$ exists in $\mathcal{C}$. Then we define $\mathrm{Pol}_{\mathcal{C}}(A, B)$ to be the minion defined in each arity by

$$\mathrm{Pol}_{\mathcal{C}}\,(A, B)_n := A^n \to_{\mathcal{C}} B$$

where $\alpha\colon n \to k$ induces the map $f \mapsto \alpha^*;f$, with $\alpha^*$ being the induced morphism $A^k \to A^n$ obtained from tupling the projections $\pi^k_{\alpha_i}$ for each $i \in n$.

We call a minion $\mathcal{M}$ a *polymorphism minion* if it is of the form $\mathrm{Pol}_{\mathcal{C}}(A, B)$. Furthermore, we abbreviate $\mathrm{Pol}_{\mathcal{C}}(A, A)$ as $\mathrm{Pol}_{\mathcal{C}}(A)$.

Note that we use slightly different notation here: $\alpha^*$ is the so-called *pullback* of $\alpha$, which we previously denoted by $\alpha;_$. That we did not re-use the notation $\alpha;_$ is due to the fact that in abstract categories, it is not necessarily the case that the product $A^k$ has the interpretation of maps $k \to_{\mathcal{C}} A$, so "precomposing" with $\alpha\colon n \to k$ makes no formal sense.

We can summarize the definition visually as the composition of functors

[15] For a precise discussion see Section A.2

$$\begin{array}{ccccc} n & & A^k & & A^n \to_{\mathcal{C}} B \\ \Big\downarrow \alpha & \rightsquigarrow & \Big\downarrow \alpha^* & \rightsquigarrow & \Big\downarrow \alpha^*;__ . \circ \\ k & & A^n & & A^k \to_{\mathcal{C}} B \end{array}$$

This shows that polymorphism minions are in fact minions.

While in the study of (promise)-CSPs – and hence in this thesis – one is mainly concerned with polymorphism minions of *relational structures*, we should give these a separate name.

**Definition 2.2.2**:
- A *concrete minion* is a polymorphism minion of $\tau$-structures, i.e., of the form $\mathrm{Pol}_\tau(\underline{A}, \underline{B})$.
- A *function minion* is a subminion of $\mathcal{O}(A, B) := \mathrm{Pol}_{\underline{\mathsf{Set}}}(A, B)$.

Clearly, every concrete minion is a function minion; however, the reverse also holds, as we will see in Section 3.1.

However, even if one is interested in concrete minions only (and, say, their homomorphism order), abstract minions can have their use. For instance, we can use the following principle to use other polymorphism minions as intermediate steps when the goal is to find a minion homomorphism to a minion of bounded essential arity, a special property we will introduce later.

**Proposition 2.2.1** (implicit in [21], [22]): Let $F\colon \mathcal{C} \to \mathcal{D}$ be a functor which *preserves finite powers*, i.e., for every $A \in \mathcal{C}_0$ and $n \in \mathbb{N}$, the canonical morphism $F(A^n) \to_{\mathcal{D}} (FA)^n$ obtained by tupling the morphisms $F\pi_i\colon F(A^n) \to_{\mathcal{D}} FA$ is an isomorphism.

Then, for every $A, B \in \mathcal{C}_0$, the functor $F$ induces a minion homomorphism

$$\mathrm{Pol}_{\mathcal{C}}(A, B) \to \mathrm{Pol}_{\mathcal{D}}(FA, FB).$$

*Proof*: Let $\eta_{A,n} := (F\pi_1, ..., F\pi_n)$ denote the aforementioned isomorphism $F(A^n) \to (FA)^n$. For each $n$, define $h_n(f) := \eta_{A,n}^{-1};(Ff)$. This composite takes an $f \in \mathrm{Pol}_{\mathcal{C}}(A, B)_n$, i.e., a morphism $A^n \to_{\mathcal{C}} B$, and then precompose $Ff\colon F(A^n) \to B$ wtih $\eta_{A,n}^{-1}\colon (FA)^n \to F(A^n)$. Evidently, this gives us an element of $\mathrm{Pol}_{\mathcal{D}}(FA, FB)_n$.

To see that $h = (h_n)_n$ is a minion homomorphism, let $\alpha\colon n \to n'$. We need to show that for all $f \in \mathrm{Pol}_{\mathcal{C}}(A, B)$, we have $h_{n'}(f\alpha) = h_n(f)\alpha$. Indeed, we have

$$\begin{aligned} h_{n'}(f\alpha) &= h_{n'}(\alpha^*;f) && \text{Minion structure of } \mathrm{Pol}_{\mathcal{C}}(A, B) \\ &= \eta_{A,n'}^{-1};F(\alpha^*;f) && \text{Definition of } h_n \\ &= \eta_{A,n'}^{-1}; F\alpha^*;Ff && F \text{ is a functor} \\ &= \alpha^*; \eta_{A,n}^{-1};Ff && \text{Lemma 1.1.3} \\ &= \alpha^*;h_n(f) && \text{Definition of } h_n \\ &= h_n(f)\alpha && \text{Minion structure of } \mathrm{Pol}_{\mathcal{D}}(FA, FB), \end{aligned}$$

as desired. □

Although we will not use the above statement in this thesis, it is an important tool to keep in mind. For instance, in the (first half of the) paper [21], the authors used this principle implicitly to derive hardness results for $\mathrm{PCSP}(K_3, K_c)$ for certain $c$. For that, they used the chain of functors

- $\mathrm{Hom}(K_2, _)$ from graphs to $\mathbb{Z}_2$-complexes [23, Chapter 18]
- Geometric realization $|_|$, going from $\mathbb{Z}_2$-complexes to $\mathbb{Z}_2$-spaces [24, Chapter III]
- the fundamental group, going from $\mathbb{Z}_2$-spaces to groups [25, Section 2.5]

While not all of these functors preserve products, their composition does:

- $\mathrm{Hom}(K_2, _)$ preserves products *up to homotopy equivalence* [23, Proposition 18.17]
- Geometric realization preserves products [24, Theorem 14.3]
- the fundamental group functor preserves products [25, Theorem 2.6.2] and is homotopy invariant [25, Proposition 2.5.4].

In general, this gives us a minion homomorphism

$$\mathrm{Pol}_{\underline{\mathsf{Gr}}}(\underline{G}, \underline{H}) \to \mathrm{Pol}_{\underline{\mathsf{Grp}}}(\pi_1(|\mathrm{Hom}(K_2, \underline{G})|), \pi_1(|\mathrm{Hom}(K_2, \underline{G})|)),$$

where in the cases $(\underline{G}, \underline{H})$ analyzed in the paper, the target minion evaluated to $\mathrm{Pol}_{\underline{\mathsf{Grp}}}(\mathbb{Z})$, and furthermore, the map could be shown to target a sub-minion of $\mathrm{Pol}_{\underline{\mathsf{Grp}}}(\mathbb{Z})$ of bounded essential arity. Finally, the $\mathbb{Z}_2$-action has been used to show that the minion avoids a constant, meaning we satisfy all the premises for the hardness assumption presented in the preliminaries.

## 2.3 Constructions on minions

In this chapter, we will present obvious constructions known from algebraic structures or from sets. All of them follow either from well-known facts about functor categories (cf. [26, Part I]) or from universal algebra, generalized to the many-sorted[16] setting [27].

Repeating these definitions is necessary due to the fact that abstract minions do not have a base set *per se*, but rather a "base family of sets", indexed by positive numbers. This is a conscious choice of the author: While it is custom in the literature to write e.g. $\mathcal{O}(X) = \bigcup_n (X^n \to X)$ and then define an arity function $\mathrm{ar}\colon \mathcal{O}(X) \to \mathbb{N}_{>0}$ do designate the separate "sorts", i.e., arities, there are some drawbacks with this notation.

- First, size issues. While $\mathcal{M} = \bigcup_n \mathcal{M}_n$ makes sense if we define arities to be finite ordinals, it ceases to do so when we extend the concept of an arity to arbitrary finite sets. There is some use in doing so – as we will see in Section 3.1 – and we will arrive at a semantically equivalent definition of "minion"[17].
- Second, suggesting the wrong semantics. If $\mathcal{M} = \bigcup_n \mathcal{M}_n$, then a cartesian product is already defined on $\mathcal{M}$ – however, that is not the correct product, as we will see in a second. For another example, the notion of "cardinality" associated to a set – in this case the minion – is also of interest only in relation to the arity. We will explore this concept further in Chapter 6. For a last, and crucial, example, the notion of equivalence relation needs to be defined more restrictively, since we would not want to allow identifying elements of different arities. Otherwise, not every equivalence relation would lie below a congruence (to be defined shortly), making the phrase "congruence generated by" dangerously ill-defined.

[16] A minion can be viewed as a many-sorted algebra with one sort per natural number and minor operations as unaries.

Let us make this notion of "family of base sets" explicit by giving it a concise name.

**Definition 2.3.1** (Set sequence): We define the *category of set sequences* to be the functor category $[\mathbb{N}_{>0}, \underline{\mathsf{Set}}]$. In other words: A *set sequence* is a family $(X_n)_{n\in\mathbb{N}_{>0}}$, and a *morphism of set sequences* $h\colon X \to Y$ is a family $(h_n)_{n\in\mathbb{N}_{>0}}$ of maps $h_n\colon X_n \to Y_n$.

If $h \in X_n$, we call $n$ the *arity* of $h$ and sometimes denote it by $\operatorname{ar} h$. Furthermore, we will use $(s_i\colon k_i)_i$ as a shorthand notation to denote an indexed set with elements $s_i$ of arity $k_i$.

Furthermore, we will override the symbol $\emptyset$ to denote the set sequence which is empty in each arity, and the symbol $*$ to denote any set sequence consisting of a single element in each arity.

**Definition 2.3.2**: Let $X$ and $Y$ be set sequences and $h\colon X \to Y$ a morphism. We define

- $X \subseteq Y$ to mean $X_n \subseteq Y_n$ for all arities $n$, in which we call $X$ a *sub-set sequence of* $Y$,
- $h$ to be *surjective* (resp. *injective*, *bijective*) if $h$ is surjective (resp. injective, bijective) on all arities,
- $\operatorname{Im} h$ to be the sub-set sequence $(\operatorname{Im} h_n)_n \subseteq Y$, and
- a *partition* of $X$ to be a collection $(P_i)_{i\in I}$ of set sequences such that none of the $P_i$ are the empty set sequence, and for each arity $n$, the collection $\left(P_{i,n}\right)_{i\in I}$ consists of pairwise disjoint (but potentially empty) sets whose union is $X_n$.

Injectivity will at times be stylized as $h\colon X \hookrightarrow Y$ and surjectivity as $h\colon X \twoheadrightarrow Y$.

**Proposition 2.3.1**: Let $I$ be any set and $(\mathcal{M}_i)_{i\in I}$ a collection of minions. The product minion

$$\left(\prod_{i\in I} \mathcal{M}_i\right)_n := \prod_{i\in I} (\mathcal{M}_i)_n,$$

for which minor operations act componentwise, is a minion, which together with the projection homomorphisms $\pi_i\colon \left(\prod_j \mathcal{M}_j\right) \to \mathcal{M}_i$, defined in arity $n$ as the $i$th projection of the arity-$n$ product, satisfies the universal property of the product.

Dually, the coproduct minion

$$\left(\sum_{i\in I} \mathcal{M}_i\right)_n := \sum_{i\in I} (\mathcal{M}_i)_n,$$

with minor operations acting in each summand, together with the canonical injections $\iota_i$ defined arity-wise, satisfies the universal property of the coproduct.

[17]Category-theoretically, this is clear because the equivalence of functor categories $\mathbb{F}_{>0} \simeq \underline{\mathsf{FinSet}_{>\emptyset}}$ extends to an equivalence $[\mathbb{F}_{>0}, \underline{\mathsf{Set}}] \simeq [\underline{\mathsf{FinSet}_{>\emptyset}}, \underline{\mathsf{Set}}]$.

*Proof*: Defining the tupling of homomorphisms $f_i\colon \mathcal{L} \to \mathcal{M}_i$ arity-wise as $\left\langle (f_i)_n \right\rangle$, i.e., the tupling of the arity-$n$ components, inherits the universal property of the arity-$n$-products $\prod_i (\mathcal{M}_i)_n$ because composition of minion homomorphisms is defined in each arity and hence "commutes" with restriction to a fixed arity $n$. The only remaining nontrivial part is to show that the projections $\pi_i$ are actually minion homomorphisms. Indeed, if $f = (f_i)_i \in \left(\prod_i (\mathcal{M}_i)\right)_n = \prod_i \left((\mathcal{M}_i)_n\right)$ and $\alpha\colon n \to k$, we have $\pi_n(f\alpha) = \pi_n\left((f_i\alpha)_i\right) = f_n\alpha = \pi_n(f)\alpha$.

The argument for the coproduct is completely analogous. □

A subset of this argument shows that the arity-wise product or coproduct of set sequences satisfies the respective universal properties of set sequences. We shall thus define relations for set sequences as sub-set sequences of the product.

Observe that the categorical product $X \times Y$ of set sequences can be modelled by taking the product (and projections) in each arity. We will hence take $X \times Y$ to refer to the set sequence $(X_n \times Y_n)_n$.

**Definition 2.3.3**: Let $X$ and $Y$ be set sequences. We call a sub-set sequence $R \subseteq X \times Y$ a *relation*, and call a relation $R \subseteq X \times X$ *symmetric*, *reflexive*, and *transitive*, and an *equivalence relation* if it is so in every arity.

As we will see in the next proposition, minions share numerous properties with algebras (e.g. groups, rings, or modules):

**Definition 2.3.4**: A *subminion* of $\mathcal{M}$ is a sub-set sequence $\mathcal{L} \subseteq \mathcal{M}$ which is closed under minor operations, i.e., which inherits the minion structure from $\mathcal{M}$.

**Proposition 2.3.2**:

1. A bijective minion homomorphism $h\colon \mathcal{M} \to \mathcal{N}$ is an isomorphism.
2. The image of a minion homomorphism is a subminion.
3. An injective minion homomorphism is isomorphic to its image.
4. Every subminion arises as the image of an injection.

*Proof*: 1) Let $h^{-1} := \left((h_n)^{-1}\right)_n$ be the inverse to $h$ in each arity. Let furthermore $g \in \mathcal{N}_k$ and $\alpha\colon k \to k'$ be a minor operation. Since

$$h_{k'}^{-1}(g\alpha) = h_{k'}^{-1}((h_k(h_k^{-1}(g)))\alpha) = h_{k'}^{-1}(h_k(h_k^{-1}(g))\alpha) = (h_k^{-1}(g))\alpha,$$

we see that $(h_k^{-1})_k\colon \mathcal{N} \to \mathcal{M}$ is indeed a homomorphism of minions. Hence, $h$ has a homomorphic inverse, as desired. 2) Let $h\colon \mathcal{M} \to \mathcal{N}$ and $\alpha\colon n \to k$. If $g \in (\operatorname{Im} h)_n$, then $g = h_n(f)$ for some $f \in \mathcal{M}_{\operatorname{ar} g}$, hence $g\alpha = h_n(f)\alpha = h_k(f\alpha) \in (\operatorname{Im} h)_k$ by homomorphicity of $h$. 3) Let $h\colon \mathcal{M} \hookrightarrow \mathcal{N}$ be injective. Corestricting to the image $\operatorname{Im} h \leq \mathcal{N}$ yields a bijective minion homomorphism, hence an isomorphism. 4) If $\mathcal{N} \leq \mathcal{M}$, then the inclusion $\mathcal{N} \hookrightarrow \mathcal{M}$ is a minion homomorphism with image $\mathcal{N}$ by the previous insight. □

Since injective minion homomorphisms induce an isomorphism to their image, we shall sometimes refer to these as *embeddings*.

Following the theme of algebras, we have a notion of "quotient", i.e., surjective homomorphisms, which are classified by a certain notion of congruence.

**Definition 2.3.5** (Congruence): Let $\mathcal{M}$ be a minion. A *congruence* $\theta$ on $\mathcal{M}$ is a subminion of $\mathcal{M}^2$ which is an equivalence relation in each arity. We will denote by $\mathrm{Cong}(\mathcal{M})$ the poset of all congruences on $\mathcal{M}$, ordered by the subminion relation.

*Remark*: Since being a congruence is a conjunction of closure properties, we easily see that the congruences on a minion form a complete bounded lattice, where (meet-)completeness means we are closed under arbitrary intersections, and the bounds are given by the diagonal congruence $\Delta_{\mathcal{M}}$ (identifying an element only with itself) and the full congruence $\mathcal{M}^2$. Hence, arbitrary *joins* exist as well, because to a given set of congruences we can form the infimum of all the congruences above it, since the full congruence guarantees this set is always nonempty.

**Definition 2.3.6** (Quotient minion): Let $\mathcal{M}$ be a minion and $\theta$ a congruence on $\mathcal{M}$. We define a minion $\mathcal{M}/\theta$ to be $\mathcal{M}_n/\theta_n$ in arity $n$ and let $\alpha\colon n \to k$ act via $[f]_{\theta_n}\alpha := [f\alpha]_{\theta_k}$.

*Proof that $\mathcal{M}/\theta$ is a minion*: The requirement that $\theta$ is a subminion means that $(f, g) \in \theta_n$ implies $(f\alpha, g\alpha) \in \theta_k$ for every $\alpha\colon n \to k$. This ensures precisely that the map $[f] \mapsto [f\alpha]$ is well-defined. The other minion axioms are inherited from $\mathcal{M}$. □

**Proposition 2.3.3**: Let $\mathcal{M}, \mathcal{N}$ be minions and $\theta \leq \mathcal{M}^2$ be a congruence. There is a one-to-one correspondence between minion homomorphisms $h\colon \mathcal{M} \to \mathcal{N}$ which are constant on $\theta$, i.e., $h_n(f) = h_n(g)$ whenever $(f, g) \in \theta_n$, and minion homomorphisms $\mathcal{M}/\theta \to \mathcal{N}$.

*Proof*: To each function $h\colon \mathcal{M}/\theta \to \mathcal{N}$, we can associate the precomposition with the canonical quotient map $[_]_\theta\colon \mathcal{M} \to \mathcal{M}/\theta$, which is constant on $\theta$ because $[_]_\theta$ is. Conversely, given a homomorphism $g\colon \mathcal{M} \to \mathcal{N}$ constant on $\theta$, the map $[f]_\theta \mapsto g(f)$ is a well-defined minion homomorphism $\mathcal{M}/\theta \to \mathcal{N}$. These associations are mutually inverse, since $f \mapsto [f] \mapsto g(f)$ is the same as $g$, and if $g$ maps $f \mapsto [f] \mapsto h([f])$, then $[f] \mapsto g(f) = h([f])$ is the same as $h$. □

*Remark*: It is also not hard to see – albeit for this thesis not important – that this correspondence is natural in $\mathcal{N}$ and $\mathcal{M}$, in the sense that this correspondence commutes with postcomposing by a minion homomorphism $\mathcal{N} \to \mathcal{N}'$, and it also commutes with postcomposing by a minion homomorphism $l\colon \mathcal{M}' \to \mathcal{M}$ (where $\theta$ is to be replaced by its preimage $l^{-1}(\mathcal{M})$ under the corresponding map $\mathcal{M}'^2 \to \mathcal{M}^2$).

Let us close off with a structural observation regarding sums of minions. To put this observation into the proper context, we shall shift our attention to groups. The sum (coproduct) of two groups $G_1$ and $G_2$ is given by the *free (or "amalgamated") product* $G_1 * G_2$, which consists of words $x_1 x_2 ... x_k$ of arbitrary length $k \geq 0$, where odd indices are populated by elements of $G_1 \setminus \{e\}$ and even indices by elements of $G_2 \setminus \{e\}$, or vice versa. The group structure is given by concatenation $w_1, w_2 \mapsto w_1 w_2$ where we multiply the end of $w_1$ with the start of $w_2$ if both inhabit the same group $G_i$.[18] The canonical injections $\iota_i$ then identify elements from $G_i$ with their respective one-letter words. It is not hard to see that this structure satisfies the universal property of a coproduct: Indeed, giving maps allows us precisely to evaluate such "mixed words".

Notable about this "sum" is that e.g. $\mathbb{Z}/2\mathbb{Z} * \mathbb{Z}/2\mathbb{Z}$ is an *infinite group*, even though the summands were finite. Similar observations can be made in many categories of algebras, and just as well for clones: The coproduct of two abstract clones is given by the free clone generated by their (arity-wise) disjoint union, modulo the identities satisfied in each summand.

On the other hand, for relational structures like graphs, but also for topological spaces, the coproduct behaves much more like a disjoint union, and the sum cannot be much "larger than its parts". In this sense, relational structures can be seen as much more "tame" than algebras.

It is thus with delight that we observe the following.

**Proposition 2.3.4**: Let $\mathcal{M}$ be a minion.
1. A partition $(\mathcal{L}_i)_{i \in I}$ of $\mathcal{M}$ where all of the $\mathcal{L}_i$ are subminions of $\mathcal{M}$ gives rise to an isomorphism $h\colon \sum_i \mathcal{L}_i \xrightarrow{\sim} \mathcal{M}$ which satisfies $\iota_i;h = h_i$, where $h_i$ is the canonical inclusion $\mathcal{L}_i \hookrightarrow \mathcal{M}$.
2. Conversely, any isomorphism $\mathcal{M} \cong \sum_i \mathcal{L}_i$ induces a partition.

*Proof*: 1) Let $h\colon \sum_i \mathcal{L}_i \to \mathcal{M}$ be given by the cotupling of the $[h_i]_i$, where $h_i\colon \mathcal{L}_i \to \mathcal{M}$ is the canonical inclusion. This homomorphism is surjective, as every element of $\mathcal{M}$ lies in one of the $\mathcal{L}_i$. To see that it is injective, note that each of the $h_i$ is injective; hence, the only source for noninjectivity can be to identify two elements $f, g$ of the same arity, but of different components $\mathcal{L}_i \neq \mathcal{L}_j$. Since $(\mathcal{L}_i)_i$ is a partition, the intersection must be empty, so $h_i(f)$ can never equal $h_j(g)$. That $\iota_i;h = h_i$ follows from the definition of $h$ as the cotupling.
2) Conversely, if $h\colon \sum_i \mathcal{L}_i \to \mathcal{M}$ is an isomorphism, and $\iota_i\colon \mathcal{L}_i \to \sum_j \mathcal{L}_j$ are the injections associated to the coproduct, we see that the images $\mathrm{Im}(\iota_i;h)$ must be disjoint (by virtue of $h$ being injective) and must cover $\mathcal{M}$ completely (by virtue of $h$ being surjective). □

## 2.4 The projection minion $\mathcal{P}$ and free minions

Just as we can consider subalgebras *generated* by certain sets of elements – e.g. subgroups, subrings, subclones, etc. – we can do the same for minions. In fact, we can generalize the abstract concept of a "free algebra" over a set $S$. For instance, recall that the *free monoid* $S^*$ over $S$ is the set of *words* $w_1 ... w_n$ of arbitrary (possibly zero) length. Consequently, it is often dubbed the *word monoid*.

For another example, assume we want to describe algebras with operations $f_i$ of arities $n_i$ – this data is called the *signature* of an algebra. We can then form what is known as "terms" over $S$ inductively as follows:
- elements of $S$ are terms (called *variables* or *ground terms*)
- for every $n_i$-tuple $\left(t_0, ..., t_{n_i - 1}\right)$ of terms and $n$-ary operation $f_i$, the word $f_i\left(t_0, ..., t_{n_i-1}\right)$ is a term.

[18] More compactly, we could view this object as the free group over the disjoint union $G_1 + G_2$ of the base sets, modulo the identities satisfied in $G_1$ and $G_2$.

These terms together with the natural interpretation of $f_i$ as an $n_i$-ary operation $\left(t_0, \ldots, t_{n_i-1}\right) \mapsto f(t_0, \ldots, t_{n-1})$ form what is called the *term algebra* of $S$ over the aforementioned signature. This is the most general algebra one can build from the given set $S$ of generators, because if $\mathbb{A}$ is an algebra and $S \subseteq \mathbb{A}$, then the algebra generated by $S$ – namely the smallest sub-algebra of $\mathbb{A}$ containing $S$ – is precisely the image of the "evaluation homomorphism", which is the unique homomorphism from the term algebra over $S$ to $\mathbb{A}$ mapping the ground terms to their respective elements in $S$.

We shall apply the same idea to minions, with the slight difference that we do not have a single set of variable symbols, but variables sorted by arity. Conveniently, the notion of "term" we need is rather simple: An $n$-ary term is a sequence $f\alpha$ where $\alpha$ is a minor operation with codomain $n$ and $\operatorname{ar} f = \operatorname{dom} \alpha$.

**Definition 2.4.1** (free minion): Let $S = (S_k)_{k>0}$ be a set sequence. Letting

$$\langle S \rangle_n := \{ s\alpha \mid k > 0, s \in S_k, \alpha\colon k \to n \},$$

we define $\langle S \rangle$ to be the *free minion generated by* $S$ whose underlying set sequence is $\left(\langle S \rangle_n\right)_n$, and on which minor operations $\beta\colon n \to l$ act as

$$\begin{aligned} _\beta\colon \langle S \rangle_n &\to \langle S \rangle_l, \\ s\alpha &\mapsto (s\alpha)\beta := s(\alpha;\beta). \end{aligned}$$

By slight abuse of notation, we shall identify elements $s \in S_n$ with the element $s := s\,\mathrm{id}_n \in \langle S \rangle_n$. Furthermore, in analogy to set sequences, we introduce the shorthand notation $\langle f_i\colon n_i \rangle_i := \left\langle (f_i\colon n_i)_i \right\rangle$, .

In the above notation, $\langle f\colon 2 \rangle$ would be the free minion generated by one binary element, and $\langle f\colon 1, g\colon 2 \rangle$ describes a minion with two unaries, namely $f = f\,\mathrm{id}_1$ and the unique unary minor $g\alpha$ of $g$, where $\alpha$ is the unique map $2 \to 1$.

*Proof that the free minion is a minion*: Let $\alpha, \beta$ as in the definition, and $\gamma\colon l \to m$. Clearly, $(s\alpha)\,\mathrm{id}_n = s(\alpha;\mathrm{id}_n) = s\alpha$ and $((s\alpha)\beta)\gamma = (s(\alpha;\beta))\gamma = s(\alpha;\beta;\gamma) = (s\alpha)(\beta;\gamma)$ by the associativity of function composition. □

**Proposition 2.4.1**: $\langle S \rangle \cong \sum_{k>0} \sum_{s \in S_k} \langle s_k\colon k \rangle$.

*Proof*: Assigning to every $k > 0$ and $s \in S_k$ the subset $\langle s \rangle$, the smallest subminion of $\langle S \rangle$ containing $s$, gives rise to a partition: Indeed, every element of $\langle S \rangle$ is of the form $s\alpha$ for some $s$ and some $\alpha$, hence lies in $\langle s \rangle$. On the other hand, if $s\alpha = t\beta$, note that $s = t$ and $\alpha = \beta$, because we need to understand these items as formal words. Hence, we can invoke Proposition 2.3.4 to see that $\langle S \rangle$ decomposes as the sum of these one-element generated minions $\langle s \rangle$.

It remains to be shown that the "internally generated minion" $\langle s \rangle$ is isomorphic to the abstract free minion $\langle s\colon \operatorname{ar} s \rangle$. However, this follows trivially, because the elements of $\langle s \rangle$ are also just words of the form $s\alpha$, just as in the external description. This assignment is clearly an isomorphism. □

*Remark (for category theorists)*: What we have constructed here is precisely the left Kan extension of $S$, viewed as a functor from the discrete category $\mathbb{N}_{>0}$ to $\underline{\mathsf{Set}}$, along the inclusion $\iota\colon \mathbb{N}_{>0} \hookrightarrow \mathbb{F}_{>0}$. Indeed, using the coend formula for (pointwise) Kan extensions [28, Proposition 2.3.6], we have

$$(\mathrm{Lan}_\iota S)_n = \int^{k\in\mathbb{N}_{>0}} S_k \otimes \left(k \to_{\mathbb{F}_{>0}} n\right),$$

where the tensor product denotes the *copower* – in our case of $\underline{\mathsf{Set}}$ as a target category, this is just the cartesian product. Furthermore, since $\mathbb{N}_{>0}$ is discrete, the coend simplifies to a sum. Together, this yields

$$(\mathrm{Lan}_\iota S)_n = \sum_{k>0} S_k \times n^k,$$

where the minor operations act on the $n^k$-factor in each summand via postcomposition (or *pointwise*, if we view elements thereof as tuples).

**Proposition 2.4.2**: Let $S$ be a set sequence and $\mathcal{M}$ a minion. The restriction of homomorphisms $h\colon \langle S\rangle \to \mathcal{M}$ to $S \subseteq \langle S\rangle$ is a one-to-one correspondence between minion homomorphisms $\langle S\rangle \to_{\underline{\mathsf{Min}}} \mathcal{M}$ and maps of set sequences $S \to_{[\mathbb{N}_{>0},\, \underline{\mathsf{Set}}]} \mathcal{M}$.

In the proof, we will drop the subscripts referring to the arity, as it can be inferred from the context.

*Proof*: Let $h\colon S \to_{[\mathbb{N}_{>0},\, \underline{\mathsf{Set}}]} \mathcal{M}$ be a map of set sequences. To extend $s \mapsto h_{\mathrm{ar}\, s}(s)$ to a minion homomorphism would demand $h_{\mathrm{cod}\,\alpha}(s\alpha) = h_{\mathrm{ar}\, s}(s)\alpha$, which already determines an extension $\langle S\rangle \to \mathcal{M}$ uniquely. Indeed, this is a homomorphism: If $t = s\alpha$ and $\mathrm{cod}\,\beta = \mathrm{ar}\, t$, then this extension maps $h(t\beta) = h(s(\alpha\beta))$ to $h(s)(\alpha\beta) = (h(s)\alpha)\beta = h(s\alpha)\beta = h(t)\beta$. Evidently, restricting this extension to $S$ yields the same map $S \to \mathcal{M}$, and extending the restriction of $h\colon \langle S\rangle \to \mathcal{M}$ reconstructs the same homomorphism. □

*Remark*: We can also easily see that this correspondence is natural in $S$ and in $\mathcal{M}$. Indeed, if $a\colon S' \to S$, then we have a natural map $\langle S'\rangle \to \langle S\rangle$ which maps $s\alpha \mapsto a(s)\alpha$. Then, given a homomorphism $h\colon \langle S\rangle \to \mathcal{M}$, first restricting to $S$ yields $s \mapsto h(s)$, and precomposing with $a$ yields $s' \mapsto h(a(s'))$. On the other hand, first precomposing with $a$ and then restricting to $S'$ gives $s' \mapsto h(a(s'))$ as well. Similarly, the restriction is compatible with postcomposition $\mathcal{M} \to_{\underline{\mathsf{Min}}} \mathcal{N}$.

This insight means that the correspondence is formally an *adjunction* $\langle _ \rangle \dashv U$ (recall Definition 1.1.22), where $U\colon\ \underline{\mathsf{Min}} \to [\mathbb{N}_{>0}, \underline{\mathsf{Set}}]$ is the "forgetful functor", which assigns a minion its underlying set sequence.

The following is an important corollary.

**Corollary 2.4.3**: Let $\mathcal{M}$ be a minion and $\langle f\colon n\rangle$ a free minion generated by a single arity-$n$ element. There is a one-to-one correspondence between arity-$n$-elements and homomorphisms from $\langle f\colon n\rangle$ into $\mathcal{M}$, denoted visually as

$$\mathcal{M}_n \quad \cong \quad \langle f\colon n\rangle \to_{\underline{\mathsf{Min}}} \mathcal{M},$$

which is natural in $\mathcal{M}$ (i.e., compatible with postcomposition $\mathcal{M} \to \mathcal{M}'$).

There is an even simpler way to describe the one-element generated free minions.

**Definition 2.4.2** (Concrete projection minions): Let $\mathrm{Proj}_A$ denote the sub-minion of $\mathcal{O}(A)$ consisting of all the projections, i.e., $(\mathrm{Proj}_A)_k = \{\pi_i^k\colon A^k \to A \mid i \in k\}$, where $\pi_i^k$ projects onto the $i$-th coordinate. Note that in particular, $\mathrm{id}_A = \pi_0^1$.

**Lemma 2.4.4**: In $\mathrm{Proj}_A$, an argument operation $\alpha\colon n \to k$ we acts as $\pi_i^n\alpha = \pi_{\alpha(i)}^k$. In particular, $\pi_i^k = \mathrm{id}_A\,\eta$ when $\eta\colon 1 \to k$ maps $0 \mapsto i$.

*Proof*: Whenever $a \in A^n$, We have $\pi_i^n\alpha(a) = \pi_i^n(\alpha;a) = \pi_i^n\left(a_{\alpha(\bullet)}\right) = a_{\alpha(i)} = \pi_{\alpha(i)}^k(a)$. The second statement is the special case where $n = 1$. □

**Definition 2.4.3** (Abstract projection minion): Let $\mathcal{P}$ denote the minion induced by the inclusion functor $\iota\colon \mathbb{F}_{>0} \to \underline{\mathsf{Set}}$, i.e., $\mathcal{P}_n := n$ with minor operations acting by function evaluation.

**Proposition 2.4.5**:
1. Whenever $|A| > 2$, we have $\mathcal{P} \cong \mathrm{Proj}_A \cong \langle f\colon 1\rangle$.
2. For all $n \geq 1$, we have $\mathcal{P}^n \cong \langle f\colon n\rangle$.

*Proof*: 1) First, note how $|A| > 2$ precisely ensures that projections $\pi_i^n$ and $\pi_j^n$ differ whenever $i \neq j$. By Lemma 2.4.4, the family of maps $i \mapsto \pi_i^n$ gives rise to a homomorphism $\mathcal{P} \to \mathrm{Proj}_A$. Since it it is clearly surjective, and – by the starting consideration – injective, it must be an isomorphism. To see that $\mathcal{P} \cong \langle f\colon 1\rangle$, we note that elements $f\alpha \in \langle f\colon 1\rangle_n$ are uniquely determined by the value of $\alpha\colon 1 \to n$ at $0$, and thus consider the family of maps $\mathcal{P}_n \to \langle f\colon 1\rangle_n$ which send $i \mapsto f(0 \mapsto i)$. This gives rise to a (clearly bijective) homomorphism $h\colon \mathcal{P} \to \langle f\colon 1\rangle$, since for every $\alpha\colon n \to k$, the element $i\alpha = \alpha(i) \in \mathcal{P}_k$ gets mapped to

$$h(i\alpha) = f(0 \mapsto \alpha(i)) = f((0 \mapsto i);\alpha) = (f(0 \mapsto i))\alpha = h(i)\alpha.$$

2) By definition of the product, we have $(\mathcal{P}^n)_k = (\mathcal{P}_k)^n = k^n$. Therefore, let us define $h_k\colon (\mathcal{P}^n)$ to map $(\alpha_0, ..., \alpha_k) \mapsto f\alpha$, where we understand $\alpha$ to be the map $k \mapsto \alpha_k$ associated to the tuple. This gives rise to a minion homomorphism $h\colon \mathcal{P}^n \to \langle f\colon n\rangle$, since post-composition $\alpha \mapsto \alpha;\beta$ corresponds on the tuple's side to pointwise application of $\beta$. Since that assignment is bijective, the result follows. □

We can thus reformulate the corollary above:

**Corollary 2.4.6**: For every minion $\mathcal{M}$, we have

$$\mathcal{M}_n \;\cong\; \mathcal{P}^n \to_{\underline{\mathsf{Min}}} \mathcal{M}.$$

Similarly, we can reformulate Proposition 2.4.1:

**Corollary 2.4.7**: $\langle S\rangle \cong \sum_k \sum_{s\in S_k} \mathcal{P}^k$.

## 2.5 Locally finite minions

As alluded to before, due to the many-sorted nature of minions, the concept of a "finite" minion does not quite make sense: as soon as we have an element in arity $n$, we have elements in all arities $k$, just by applying minor operations $n \to k$. This means that for every nonempty minion, $\bigcup_k \mathcal{M}_k$ can never be finite. However, an interesting property is whether the minion is finite *in each arity*.

**Definition 2.5.1**: A minion $\mathcal{M}$ is called *locally finite* if $\mathcal{M}_n$ is a finite set for each arity $n$.

Locally finite minions are interesting in two regards:
1. Firstly, polymorphism minions of (pairs of) finite structures $(\underline{A}, \underline{B})$ are subminions of $\mathcal{O}(A, B)$, which contains only finitely many elements of each arity (to be more precise, $|B|^{|A|^n}$ in arity $n$).
2. Secondly, they enjoy particularly nice properties: As we will see, whether a homomorphism $\mathcal{M} \to \mathcal{N}$ between locally finite minions exists is determined purely by so-called *primitive positive sentences* ("pp-sentences") – to be defined briefly – in the sense that if all pp-setnences satisfied in $\mathcal{M}$ are also satisfied in $\mathcal{M}$, then a homomorphism exists.

The author is thankful to Jakub Opršal and Manuel Bodirsky, who both suggested the use of a compactness argument as a method to attack Proposition 2.5.4, and to Sebastian Meyer, who pointed out that one could consider homomorphisms "up to level $n$" as a useful intermediate step.

### 2.5.1 Compactness

We warm up with a weaker version of the statement we want to prove, to introduce the tool of compactness. For a complete, rigorous reference on the concepts we are going to introduce, the reader is referred to [29, Section 2.3].

Before stating the theorem, let us recall the basic topological concepts involved.

A *topology* on a set $X$ is a system of sets containing $\emptyset$, $X$, and closed under finite intersections and arbitrary unions. A pair $(X, \mathcal{T})$ of a set and a topology on that set is called a *topological space*. Elements of a topology are called *open sets*, and complements of open sets are *closed*.

A topological space is called *compact* if and only if one of the following equivalent conditions hold:

1. Every open cover has a finite subcover
2. Every system of closed sets which has the finite intersection property ("FIP") – i.e., for which every finite intersection is nonempty – has nonempty intersection in total.

For a sequence of topological spaces $(X_i, \mathcal{T}_i)_{i \in I}$ there is a natural topology on the cartesian product $\prod_{i \in I} X_i$ defined as the smallest topology making the projections continuous.

With this topology, $\left(\prod_i X_i, (\pi_i)_i\right)$ is a product in the category of topological spaces and continuous maps in the sense of Definition 1.1.14.

Most importantly, we can now state the following theorem by Tychonoff:

**Theorem 2.5.1** (Tychonoff): The product of compact topological spaces is compact.

For a proof the reader is referred to [29, section 3.4.2].

**Corollary 2.5.2**: Let $\mathcal{M}$ and $\mathcal{N}$ be locally finite minions. The space of function sequences

$$\prod_{n=1}^{N} \mathcal{N}_n^{\mathcal{M}_n},$$

where each of the (finite) factors $\mathcal{N}_n^{\mathcal{M}_n}$ is viewed as a discrete topological space, is compact.

*Remark*: Tychonoff's theorem is equivalent to the axiom of choice[19]. However, for the statements of this thesis, we could have also relied on the compactness theorem of first order logic, which is equivalent to a strictly weaker axiom called the ultrafilter lemma ("UF") – or, alternatively, the boolean prime ideal theorem ("BPI") [30, Theorem 2.2].

**Proposition 2.5.3**: Let $\mathcal{M}$ and $\mathcal{N}$ be a locally finite minions. The following are equivalent:

- For every $N$, there is a sequence of functions $(h_n)_{n=1}^{N} \in \prod_{n=1}^{N} \mathcal{N}_n^{\mathcal{M}_n}$ which satisfy $h_n(m\alpha) = h_n(m)\alpha$ for every operation $\alpha\colon k \to k'$ between ordinals $k, k' \leq n$
- There is a homomorphism $\mathcal{M} \to \mathcal{N}$.

*Proof*: The backward direction is trivial. For the forward direction, define

$$H_N := \left\{ (h_n)_{n=1}^{N} \in \prod_{n=1}^{N} \mathcal{N}_n^{\mathcal{M}_n} \,\middle|\, h_n(m\alpha) = h_n(m)\alpha \forall n, k, k' \leq N, m \in \mathcal{M}_n, \alpha\colon k \to k' \right\}$$

as in the theorem. Define further $C_N := \pi_{\{1,\dots,N\}}^{-1} H_N$, where

[19] See [30, Section 2.3, Problem 8].

$$\pi_{\{1,\dots,N\}}\colon \prod_{n\geq 1} \mathcal{N}_n^{\mathcal{M}_n} \twoheadrightarrow \prod_{n=1}^{N} \mathcal{N}_n^{\mathcal{M}_n},$$
$$(h_n)_n \mapsto (h_n)_{n=1}^{N}$$

is the projection onto the first $N$ coordinates. Since $H_n$ is a subset of a finite product of discrete topological spaces – hence itself a discrete topological space – it is closed, so $C_n$, the preimage under a continuous map[20] must be closed as well.

We furthermore note that $C_N \leq C_{N+1}$, because if a sequence $(h_n)_{n=1}^{N+1}$ satisfies the axioms of $H_{N+1}$, then the subsequent restriction $(h_n)_{n=1}^{N}$ satisfies the axioms of $H_N$.

Now note that the assumptions of the theorem imply that each $C_N$ is nonempty. Hence the system $\{C_N\}_N$ of closed sets has the FIP: Indeed, $\bigcap_{i=1}^{k} C_{n_i} = C_{\max_i\{n_i\}}$, which is nonempty. By compactness, $\bigcap_N C_N$ is nonempty, so there is a sequence of functions $(h_n)_n$ which is simultaneously in all $C_N$. Since such a sequence must satisfy $h_k(m\alpha) = h_n(m)\alpha$ for *all* $m \in \mathcal{M}_n$ and $\alpha\colon n \to k$, it must be a minion homomorphism. □

### 2.5.2 Primitive positive sentences

In the context of first-order logic, a *primitive positive formula* is a formula built using only conjunctions ($\wedge$), equality, and existential quantification. Crucially, satisfying a given primitive positive sentence[21] $\varphi$ is stable under homomorphism: If it is satisfied in a structure $\underline{A}$ and there is a homomorphism $\underline{A} \to \underline{B}$, it is also satisfied in $\underline{B}$.

For brevity, we will define primitive positive formulas for minions in such a way that we automatically demand a certain normal form; namely, it is well known that every primitive positive sentence is semantically equivalent[22] to one where existential quantifiers appear only in the front.

**Definition 2.5.2** (primitive positive sentence of minions): A *primitive positive sentence* or *pp-sentence* (of minions) $\varphi$ consists of

- a finite set sequence $B$ of „bound variables"
- a sequence of words $(b_i\alpha_i \approx b_i'\alpha_i')_{i\in I}$ over some finite index set $I$, where the domain of $\alpha_i$ (resp. $\alpha_i'$) is the arity of $b_i$ (resp. $b_i'$), the domain of $\alpha_i'$ has domain $\operatorname{ar} b_i'$, and the codomains of $\alpha_i$ and $\alpha_i'$ agree.

Such a sentence is visually expressed as

$$\varphi = \exists B.\ \bigwedge_{i\in I}(b_i\alpha_i \approx b_i'\alpha_i')$$

An *choice of witnesses* of the pp-sentence $\varphi$ in a minion $\mathcal{M}$ is a map $[\![_]\!]\colon B \to \mathcal{M}$ of indexed sets such that for all $i \in I$, we have $[\![b_i]\!]\alpha_i = [\![b_i']\!]\alpha_i'$. Note that the restrictions on arity, domain, and codomain of the symbols involved ensure precisely that the above equality is well-defined.

[20] The projections are continuous by construction/definition of the product topology.
[21] Recall that a *sentence* is a formula without free variables, i.e., every variable needs to be quantified over.
[22] Sentences $\varphi$ and $\varphi'$ are *semantically equivalent* if they hold in the same structures.

Instead of looking at primitive positive sentences as formulas, we can also view them as structures: To each primitive positive sentence $\varphi$, there is a *gadget minion* $\mathcal{M}_\varphi$ such that $\varphi$ holds in $\mathcal{M}$ if and only if $\mathcal{M}_\varphi \to \mathcal{M}$. This explains most clearly why homomorphisms preserve the primitive positive sentences of a structure (here: of a minion). To do so, let us introduce some notation.

**Definition 2.5.3**: Let $B$ be a set sequence and $b_i\alpha_i, f_i'\alpha_i'$ elements of $B^2$ where $i \in \{1, ..., r\}$. We define

$$\langle B \mid b_1\alpha_1 \approx f_1'\alpha_1', ..., b_r\alpha_r \approx f_r'\alpha_r'\rangle := \langle B\rangle / \left\langle (b_i\alpha_i, f_i'\alpha_i)_{i=1}^r \right\rangle_{\text{Cong}},$$

where the dividend refers the smallest congruence of $\langle B\rangle$ generated by the tuples $(b_i\alpha_i, f_i'\alpha_i')$.

**Observation 2.5.1**: By the universal properties of the free minion and of quotient minions, a homomorphism

$$h\colon \langle B \mid b_1\alpha_1 \approx f_1'\alpha_1', ..., b_r\alpha_r \approx f_r'\alpha_r'\rangle \to \mathcal{M}$$

corresponds precisely to a map $[\![_]\!]\colon B \to \mathcal{M}$ of set sequences which satisfies $[\![b_i]\!]\alpha_i = [\![f_i']\!]\alpha_i'$ for all $i \in \{1, ..., r\}$.

Hence, a primitive positive sentence

$$\varphi = \exists B.\ \bigwedge_{i=1}^r (b_i\alpha_i \approx b_i{}'\alpha_i{}')$$

is satisfied in $\mathcal{M}$ if and only if there is a homomorphism

$$h\colon \langle B \mid b_1\alpha_1 \approx b_1'\alpha_1', ..., b_r\alpha_r \approx b_r'\alpha_r'\rangle \to \mathcal{M}.$$

Let us give some examples for primitive positive sentences for minions. Most of these are motivated by certain types of operations whose existence in $\mathrm{Pol}(\underline{A})$ allows a certain algorithm to solve the CSP.

To do so, we will introduce some compact notation for argument operations $\alpha\colon n \to k$ by denoting them as a tuple $(\alpha_0, ..., \alpha_{n-1})$. Since this leaves the codomain implicit, we shall attach it as a subscript if necessary. So for instance, the tuple $(1, 2)_3$ is to be interpreted as a minor operation $2 \to 3$ mapping $0 \mapsto 1$, and $1 \mapsto 2$.

- Existence of a constant: For concrete functions, the statement is $f(x_0) = f(x_1)$, so the corresponding primitive positive formula is $\exists(f\colon 1).f(0)_2 \approx f(1)_2$.
- Existence of a binary symmetric function: $\exists(f\colon 2).f \approx f(10)$.
- Existence of an $n$-ary symmetric function: $\exists(f\colon n).f \approx f\sigma$, where $\sigma\colon n \to n$ is the right cyclic shift.
- Existence of a 4-ary "Siggers term": This is a term satisfying $f(a, r, e, a) = f(r, a, r, e)$, so as a primitive positive minion formula this is modelled as $\exists(f\colon n).f(0, 1, 2, 0) \approx f(1, 0, 1, 2)$.

- Existence of a *quasi majority term*: this is a term satisfying $f(x,x,y) = f(x,y,x) = f(y,x,x) = f(x,x,x)$ [10, Exercise 59], which corresponds to $\exists(f\colon 3).f(0,0,1) \approx f(0,1,0) \approx f(1,0,0) \approx f(0,0,0)$.
- Existence of a *quasi Maltsev term*: This is a term satisfying $f(x,x,y) = f(y,x,x) = f(y,y,y)$, i.e., $\exists(f\colon 3).f(0,0,1) \approx f(1,0,0) \approx f(1,1,1)$.

### 2.5.3 Classifying the existence of a homomorphism

We are now equipped to prove a full characterization of when there is a minion homomorphism $\mathcal{M} \to \mathcal{N}$. It is a a generalization of [19, Thm. 4.12]where they proved a similar thing for concrete minions $\mathrm{Pol}(\underline{A}, \underline{B})$ where $\underline{A}, \underline{B}$ are relational structures with finite domain.[23]

To understand the second point in the following proposition, recall that the full transformation monoid $\mathcal{T}_n$ is the set $n^n$ of all self-functions on $n$ endowed with the composition as monoid structure. A $\mathcal{T}_n$-set is now a set $X$ with a monoid action $\mathcal{T}_n \curvearrowright X$, i.e., a monoid homomorphism $\mathcal{T}_n \to X^X$, or a map $X \times \mathcal{T}_n \to X$ satisfying the compatibility property $(x.\alpha).\beta = x.(\alpha;\beta)$, i.e., precisely these conditions to make words $x\alpha\beta$ well-defined. A homomorphism $X \to Y$ of $\mathcal{T}_n$-sets is then a map $h$ which commutes with the action of each $\alpha \in \mathcal{T}_n$. In particular, the arity-$n$ set $\mathcal{M}_n$ of a minion $\mathcal{M}$ has the natural structure of a $\mathcal{T}_n$-set, and a minion homomorphism $h\colon \mathcal{M} \to \mathcal{N}$ restricts to a homomorphism $h_n\colon \mathcal{M}_n \to \mathcal{N}_n$ of $\mathcal{T}_n$-sets.

**Proposition 2.5.4**: Let $\mathcal{M}, \mathcal{N}$ be locally finite minions. The following are equivalent:

1. there is a minion homomorphism $\mathcal{M} \to \mathcal{N}$
2. for each $N \geq 1$, there is a sequence of functions $(h_n)_n\colon \prod_n \mathcal{N}_n^{\mathcal{M}_n}$ which is compatible with the action of minor operations $\alpha\colon k \to k'$ forall $k, k' \leq N$
3. for each $N \geq 1$, there is a homomorphism of $T_N$-sets $\mathcal{M}_n \to \mathcal{N}_n$
4. for each $N \geq 1$, there is a sequence of functions $(h_n)_{n=1}^N\colon \prod_{n=1}^N (\mathcal{N}_n)^{\mathcal{M}_n}$ which is compatible with the action of minor operations $\alpha\colon k \to k'$ forall $k, k' \leq N$
5. $\mathcal{N}$ satisfies all the primitive positive sentences which are satisfied in $\mathcal{M}$.

*Proof*: The cases where one of $\mathcal{M}$ or $\mathcal{N}$ is empty are trivial, so assume both minions are nonempty. $1 \Leftrightarrow 2$ has been established in Proposition 2.5.3. $2 \Rightarrow 3$ is trivial. For $3 \Rightarrow 4$, let $h_N\colon \mathcal{M}_N \to \mathcal{N}_N$ be a morphism of $\mathcal{T}_N$-sets. Then, for each $n \leq N$, we define $h_n(m) := h_N(m\iota_n)\rho_n$, where $\iota_n\colon n \hookrightarrow N$ is the canonical inclusion and $\rho_n$ some post-inverse, say $\min(n-1, _)$. Then for every $\alpha\colon n \to n'$, the sequence $(h_n)_{n=1}^N$ satisfies

$$\begin{aligned}
h_{n'}(m\alpha) &= h_N(m\alpha\iota_{n'})\rho_{n'} && \text{Definition of } h_{n'} \\
&= h_N(m\iota_n\rho_n\alpha\iota_{n'})\rho_{n'} && \iota_n;\rho_n = \mathrm{id}_n \\
&= h_N(m\iota_n)(\rho_n\alpha\iota_{n'})\rho_{n'} && \text{By assumption} \\
&= h_N(m\iota_n)\rho_n\alpha && \iota_{n'};\rho_{n'} = \mathrm{id}_{n'} \\
&= h_n(m)\alpha. && \text{Definition of } h_n
\end{aligned}$$

[23]Strictly speaking, we will not provide a full generalization, as the theorem in [19] shows some more equivalent conditions to $\mathcal{M} \to \mathcal{N}$ which make reference to the "promise template" $(\underline{A}, \underline{B})$. Of course, these don't carry over to abstract minions.

For $4 \Rightarrow 2$ we can simply extend $(h_n)_{n=1}^{N}$ to an infinite sequence using arbitrary maps[24] $h_i\colon \mathcal{M}_i \to \mathcal{N}_i$ for $i > N$, as we place no conditions on these higher $h_i$ in our desired sequence.

To see $4 \Rightarrow 5$, let $\varphi = \exists B. \bigwedge_{i\in I} f_i\alpha_i \approx f_i{}'\alpha_i{}'$ be a pp-formula satisfied in $\mathcal{M}$, as witnessed by the choice of witnesses $[\![_]\!]\colon B \to \mathcal{M}$. Let $N := \max\big(\{\mathrm{ar}\ f_i\}_i \cup \{\mathrm{ar}\ f_i{}'\}_i \cup \{\mathrm{cod}\ \alpha_i\}_i \cup \{\mathrm{cod}\ \alpha_i{}'\}_i\big)$. We claim that the composite $[\![_]\!] := h([\![_]\!])$, which maps a bound variable $f \in B_n$ to $h_{n([\![f]\!])}$, witnesses $\mathcal{N} \models \varphi$. Indeed, for every $i \in I$ we have

$$\begin{aligned}
[\![f_i]\!]'\alpha_i &= h_{\mathrm{ar}\ f_i}([\![f_i]\!])\alpha_i && \text{Definition of } [\![_]\!]' \\
&= h_{\mathrm{cod}\ \alpha_i}([\![f_i]\!]\alpha_i) && \text{By assumption} \\
&= h_{\mathrm{cod}\ \alpha_i{}'}([\![f_i{}']\!]\alpha_i{}') && \mathcal{M} \models \varphi \\
&= h_{\mathrm{ar}\ f_i{}'}([\![f_i{}']\!])\alpha_i{}' && \text{By assumption} \\
&= [\![f_i{}']\!]'\alpha_i{}'. && \text{Definition of } [\![_]\!]'
\end{aligned}$$

For the converse $5 \Rightarrow 4$, let $N$ be arbitrary. We claim that $\mathcal{T}_N$-equivariance is encoded by the sentence

$$\varphi := \bigwedge_{f\in\mathcal{M}_N} \bigwedge_{\alpha\colon N\to N} s_f\alpha \approx s_{f\alpha},$$

where we introduced a bound variable $s_f$ for each $f \in \mathcal{M}_N$, i.e., $B := \{s_f\colon N \mid f \in \mathcal{M}_N\}$. Indeed, the choice of witnesses $[\![s_f]\!] := f$ clearly satisfies $[\![s_f]\!]\alpha = [\![s_{f\alpha}]\!]$, so preservation of pp-formulas implies existence of a choice of witnesses $[\![_]\!]'\colon B \to \mathcal{N}$ such that $[\![s_f]\!]'\alpha = [\![s_{f\alpha}]\!]'$. This means that the map $h\colon \mathcal{M}_N \to \mathcal{N}_N$ sending $f \mapsto [\![s_f]\!]'$ is compatible with the $\mathcal{T}_N$-action, as desired. □

## 2.6 Essential coordinates and bounded essential arity

In this section, we will investigate the notion of *inessential coordinates*, and how this concept interacts with homomorphisms and argument operations. These facts will be relevant throughout this thesis. For instance, the projection $\pi_0^2\colon A^2 \to A$ onto the first component does not depend on the second coordinate. in general, we say that for a $n$-ary function $f\colon A^n \to B$, a coordinate $i \in n$ is *inessential* if for all elements $a_0, \ldots, a_{n-1}, a_n$ of $A$ we have

$$f(a_0, \ldots, a_{i-1}, a_i, a_{i+1}, \ldots, a_{n-1}) = f(a_0, \ldots, a_{i-1}, a_n, a_{i+1}, \ldots, a_{n-1}).$$

To generalize this notion to abstract minions, note that the left hand side is $f\iota$, where $\iota\colon n \to n+1$ is the canonical inclusion, and the right hand side is $f\alpha$ where $\alpha\colon n \to n+1$ is the minor operation sending $i$ to $n$ and leaving $n \setminus \{i\}$ as-is.

The concept of (in)essential coordinates is well-known in the context of function clones and function minions (cf. [19, Definition 5.13]). In the context of abstract minions, the notion has been known under slightly different names, and is defined in much greater generality. We will elaborate on the details in Section 6.3. While the results presented in this chapter follow from known general results about endofunctors on subcategories of Set, we will present and derive them in isolation, as we are e.g. not concerned with infinitary behavior.

[24] which exist due to the nonemptiness assumption

First off however, we will introduce some notation which enables for a more concise and effective presentation.

**Definition 2.6.1**:
- $\iota\colon n \to n+1$ is the canonical inclusion (whose domain will be clear from the context)
- We may denote a minor operation $\alpha\colon n \to k$ as

$$\alpha = (\alpha_0, ..., \alpha_{n-1})_k$$

  where $\alpha_i \in k$, omitting the codomain subscript if it can be inferred from the context.
- For $\alpha\colon n \to k$, we let $\alpha^\wedge\colon n+1 \to k+1$ denote the extension of $\alpha$ which sends the new domain element $n$ to the new codomain element $k$.
- For every natural number $n$, $S \subseteq n$, $k > \max\{n \setminus S\}$ and $j \in k$, we let ${}_n(S \mapsto j)_k : n \to k$ denote the function which takes the constant value $j$ on $S$ and is the identity anywhere else.
- The shorthand notation ${}_n(i \mapsto j)_k$ means ${}_n(\{i\} \mapsto j)_k$.

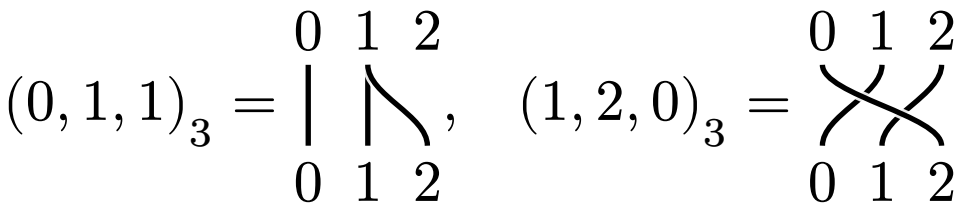

Figure 7: Two minor operations $3 \to 3$ encoded as tuples (read from bottom to top).

The following is trivial, but very useful:

**Lemma 2.6.1**: For every $\alpha\colon n \to k$, we have $\alpha\iota = \iota\alpha^\wedge$.

It is instructive for the reader to verify that given an arity for $\alpha$, the meaning of $\iota$ is unambiguous, both on the left- and on the right hand side.

**Definition 2.6.2**: Let $f \in \mathcal{M}_n$. An element $i \in n$ is called an *inessential coordinate of f* if

$$f\iota = f\,{}_n(i \mapsto n)_{n+1},$$

and an *essential coordinate* otherwise. Accordingly, set

$$\begin{aligned} \operatorname{Iness} f &:= \{i \in n \mid i \text{ inessential in } f\} \subset n \\ \operatorname{Ess} f &:= n \setminus \operatorname{Iness} f. \end{aligned}$$

We call $f$ *nondegenerate* if $\operatorname{Iness} f = \emptyset$, *degenerate* if $\operatorname{Iness} f \neq \emptyset$, and *constant* if $\operatorname{Ess} f = \emptyset$.

**Proposition 2.6.2**: Let $f \in \mathcal{M}_n$ with $n \geq 2$. The following are equivalent:
1. $i \in \operatorname{Iness} f$
2. $f$ is fixed by a non-invertible endomorphism $\varepsilon\colon n \to n$ such that $i \notin \operatorname{Im} \varepsilon$
3. $f = g\alpha$ for some $g$ and $\alpha$ with $i \notin \operatorname{Im} \alpha$.
4. $f = g\eta$ for some $g$ and injective $\eta$ with $i \notin \operatorname{Im} \eta$.

In particular, $f$ is degenerate if and only if it is fixed by a non-invertible endomorphism.

*Proof*: $1 \Rightarrow 2$: We know that $f\iota = f\,{}_n(i \mapsto n)_{n+1}$. Let $j \in n \setminus \{i\}$, which exists because $n \geq 2$. Then the map $\rho := {}_{n+1}(n \mapsto j)_n$ is a retraction to $\iota$, i.e., satisfies $\iota; \rho = \mathrm{id}_n$. On the other hand, ${}_n(i \mapsto n)_{n+1}\,; \rho = {}_n(i \mapsto j)_n$ is an endomorphism of $n$ avoiding $i$ in its image.

$2 \Rightarrow 3$: Let $f = f\varepsilon$ for $\varepsilon\colon n \to n$ with $i \notin \operatorname{Im} \varepsilon$. We can factorize $\varepsilon$ as $\rho; \eta$ where $\rho\colon n \twoheadrightarrow |\operatorname{Im} \varepsilon|$ is a surjection and $\eta\colon\ |\operatorname{Im} \varepsilon| \hookrightarrow n$. Since the inclusion $\eta$ must be proper, we have $f = f\varepsilon = (f\rho)\eta$. Since $\rho$ is surjective, we have that $\operatorname{Im} \eta = \operatorname{Im} (\rho; \eta) = \operatorname{Im} \varepsilon$, so $i \notin \operatorname{Im} \eta$. Therefore, we can pick $g := (f\rho)$ to satisfy 3.

$3 \Rightarrow 4$: If $f = g\alpha$ with $i \notin \operatorname{Im} \alpha$ and we factorize $\alpha$ as $\beta;\eta$ where $\beta$ is a surjection and $\eta$ an injection, then certainly $i \notin \operatorname{Im} \eta$ and $f = (g\beta)\eta$.

$4 \Rightarrow 1$: If $f = g\eta$ and $i \notin \operatorname{Im} \eta$, we can conclude that $\eta\iota = \eta\,{}_n(i \mapsto n)_{n+1}$, hence $f\iota = f\,{}_n(i \mapsto n)_{n+1}$.

For the corollary, let $f$ be degenerate. Choosing an $i \in \operatorname{Iness} f$ we get that $f$ must be fixed by some endomorphism with $i \notin \operatorname{Im} \varepsilon$, implying $\varepsilon$ is not surjective. On the other hand, if $f = f\varepsilon$ with $\varepsilon$ not surjective, then choosing an $i \notin \operatorname{Im} \varepsilon$ means that $i \in \operatorname{Iness} f$, witnessing its degeneracy. □

### 2.6.1 Essentiality and minor operations

Let $f \in \mathcal{M}_n$. If we know which coordinates of $f$ are essential, what can we say about the coordinates of a minor $f\alpha$?

In the concrete case, i.e., if $f$ is a function $A^n \to B$, we can make a few intuitive observations:
1. If $\pi\colon n \to n$ is a permutation, then $f\pi(\underline{x}) = f\left(x_{\pi_0}, ..., x_{\pi_{n-1}}\right)$, so $i$ is essential in $f\pi$ if and only if $\pi_i$ is essential in $f$.
2. If $\iota\colon n \hookrightarrow n + k$ is the canonical inclusion, we know that $f\iota(\underline{x}) = f(x_0, ..., x_{n-1})$, so it certainly does not depend on the last $k$ coordinates. Generalizing this, $f\alpha$ cannot depend on any $i$ outside of $\operatorname{Im} \alpha$.
3. If $f(x_0, x_1, x_2)$ depends on $x_1$ and $x_2$, and we identify $x_1$ and $x_2$, then the minor $(x_0, x_1) \mapsto f(x_0, x_1, x_1)$ should still depend on $x_1$.
4. If $f(x_0, x_1)$ depends on both coordinates, but we restrict to the diagonal via the minor $x_0 \mapsto f(x_0, x_0)$, it might perfectly well happen that it is constant: An example is the mod-2 addition $\oplus\colon 2^2 \to 2$, which sends $(0, 0)$ and $(1, 1)$ both to $0$. Theferore, its minor $\oplus\, \alpha$ with $\alpha\colon 2 \to 1$ has only inessential coordinates, while $\oplus$ had only essential coordinates.

Let us prove that these statements hold in the abstract setting as well, where the definition of "essentiality" is encoded using minor identities.

**Lemma 2.6.3**: Let $f \in \mathcal{M}_n$. If $\alpha\colon n \to k$ is injective, then
1. $\operatorname{Ess} f\alpha = \alpha[\operatorname{Ess} f]$
2. $\alpha[\operatorname{Iness} f] \subseteq \operatorname{Iness} f\alpha$

*Proof*: Claim 1, "$\supseteq$": Assume $i \in \operatorname{Ess} f$, i.e., $f\iota \neq f(i \mapsto k+1)$. Since $\alpha$ and thus $\alpha^\wedge$ are injective, we have $f\alpha\iota = f\iota\alpha^\wedge \neq f(i \mapsto k+1)\alpha^\wedge = f\alpha(\alpha(i) \mapsto k+1)$, so $\alpha(i) \in \operatorname{Ess} f\alpha$, as desired.

For "$\subseteq$", let $j \in \operatorname{Ess} f\alpha$, i.e., $f\alpha(j \mapsto k+1) \neq f\alpha$ Then in particular $j$ must lie in the image of $\alpha$, because otherwise $\alpha(j \mapsto k+1) = \alpha\iota$. Since $\alpha$ is injective, denote the unique element in $\alpha^{-1}(j)$ by $i$. Thus we have $f\iota\alpha^\wedge = f\alpha\iota \neq f\alpha(\alpha(i) \mapsto k+1) = f(i \mapsto k+1)\alpha^\wedge$. Cancelling the injection $\alpha^\wedge$ implies $f\iota \neq f(i \mapsto k+1)$, which means that $i \in \operatorname{Ess} f$, hence $j = \alpha(i) \in \alpha[\operatorname{Ess} f]$.

Claim 2 follows because

$$\alpha[\operatorname{Iness} f] = \alpha\left[(\operatorname{Ess} f)^{\complement}\right] \subseteq \alpha[\operatorname{Ess} f]^{\complement} = (\operatorname{Ess} f\alpha)^{\complement} = \operatorname{Iness} f\alpha.$$

□

We should note however that this result does not hold *locally*: Indeed, it is possible to construct a ternary map $f\colon 3^3 \to 3$ which has only essential coordinates, but for which $f(x_0, x_1, x_1)$ is constant. That means, while $0 \in \operatorname{Ess} f$ and $\alpha = (0, 1, 1)\colon 3 \to 2$ was injective "on 0", the induced coordinate $\alpha(0)$ is *not* essential in $f\alpha$. Such a map can be obtained e.g. as

$$f(x, y, z) := \begin{cases} 0 & \text{if } |\{x, y, z\}| < 3 \\ \max(x, y, z) & \text{otherwise} \end{cases},$$

which is clearly constant on all 2-minors, but has only essential coordinates.

**Lemma 2.6.4**: Let $f \in \mathcal{M}_n$, $\alpha\colon n \to k$, and $j \in k$. If $\alpha^{-1}(j) \subseteq \operatorname{Iness} f$, then $j \in \operatorname{Iness} f\alpha$. In other words, if $j$ is inessential and $f$ exclusively identifies $j$ with other coordinates which are inessential, then $\alpha(j)$ is inessential in $f\alpha$.

*Proof*: To show that $j \in \operatorname{Iness} f\alpha$ we have to verify the equality

$$f\alpha \underbrace{{}_{k}(j \mapsto k+1)_{k+1}}_{=:\gamma} = f\alpha\iota.$$

To do so, consider the operation

$$\beta := {}_{n}(\alpha^{-1}(j) \mapsto n+1)_{n+1}\,.$$

On the one hand, it satisfies $j \notin \operatorname{Im} \beta\alpha^\wedge$, and hence $(\beta\alpha^\wedge)\gamma^\wedge = (\beta\alpha^\wedge)\iota$. On the other hand, since every element in $\alpha^{-1}(j)$ is inessential, we can see inductively that $f\iota = f\beta$ by fixing an enumeration $(a_i)_i$ of $\alpha^{-1}(j)$ and then factoring $\beta = \beta_1;\ldots;\beta_{|\alpha^{-1}(j)|}$ $\beta_0 = {}_{n}(a_0 \mapsto n)_{n+1}$ and $\beta_i = {}_{n+1}(a_i \mapsto n)_{n+1}$. We thus see

$$f\alpha\gamma\iota = f\iota\alpha^\wedge\gamma^\wedge = f\beta\alpha^\wedge\gamma^\wedge = f\beta\alpha^\wedge\iota = f\iota\alpha^\wedge\iota = f\alpha\iota\iota$$

which shows the desired statement. □

Note that the proof works even if $\alpha^{-1}(j) = \emptyset$: in that case $\beta$ just degenerates to $\iota$.

**Corollary 2.6.5**: For $f \in \mathcal{M}_n$ and $\gamma\colon n \to l$, we have $|\operatorname{Ess} f\gamma| \leq |\operatorname{Ess} f|$.

*Proof*: Let $j \in \operatorname{Ess} f\gamma$. By Proposition 2.6.2, we know that $j \in \operatorname{Im} \gamma$. Now due to Lemma 2.6.4, the preimage $\gamma^{-1}(j)$ cannot consist of only inessential coordinates.

Therefore, all the essential coordinates must lie in the image of the restriction $\gamma\,|_{\operatorname{Ess} f}$, which has cardinality at most $\operatorname{Ess} f$. □

### 2.6.2 The constant minion

If a concrete minion $\operatorname{Pol}(\underline{A}, \underline{B})$ has a constant element $f$ with constant value $b \in \underline{B}$, then $\operatorname{PCSP}(\underline{A}, \underline{B})$ is trivial: $\{b\}$ must induce a one-element substructure of $\underline{B}$ in which all relation symbols of the signature are either unpopulated or populated by a constant tuple $(b, ..., b)$, i.e., a *loop*. Deciding the PCSP is now trivial, because for each instance $\underline{I}$, we answer *yes* if every relation populated in $\underline{I}$ is populated in the substructure induced by $\{b\}$ and *no* otherwise.

The predicate[25] "f is a constant" can be formulated in terms of minor conditions: An element $f$ is constant if and only if it has no essential coordinates. The property "$\mathcal{M}$ has a constant element of some arity" can be encoded as "$\mathcal{M}$ admits a homomorphism from" a special kind of minion, which we will call *the constant minion*. This minion is the terminal object in the category of minions, and hence represents the greatest element of the homomorphism order.

**Definition 2.6.3** (constant minion): Let $S$ be any set with one element (a "singleton"). We define a minion $*$ by $*_n := S$, i.e. with the same set in each arity, where a minor operation $\alpha\colon n \to k$ sends the unique element $s \in *_n$ to the unique element $s \in *_k$.

**Proposition 2.6.6**: The minion $*$ satisfies the universal property of a *terminal object* in the category of minions, i.e., for every minion $\mathcal{M}$, there is a unique minion homomorphism from $\mathcal{M}$ to $*$.

*Proof*: There is only one possible candidate, which is clearly a minion homomorphism. □

*Remark*: This property determines $*$ uniquely up to isomorphism: Indeed, if both $\mathcal{M}$ and $\mathcal{N}$ were terminal minions, then there had to be unique minion homomorphisms $h\colon \mathcal{M} \to \mathcal{N}$ and $h'\colon \mathcal{N} \to \mathcal{M}$. Now, the composite $h;h'\colon \mathcal{M} \to \mathcal{M}$ must agree with the identity $\mathrm{id}_{\mathcal{M}}$ due to the uniqueness of morphisms $\mathcal{M} \to \mathcal{M}$. The uniqueness part of the universal property now enforces that these must agree, so $h;h' = \mathrm{id}_{\mathcal{M}}$. Analogously, $h';h = \mathrm{id}_{\mathcal{N}}$, so $h$ and $h'$ are isomorphisms.

This justifies referring to the minions $*$, whose definition depended on a choice of singleton set $S$, as *the* constant minion.

[25] More precisely, a family of predicates $\mathrm{const}_k$ for each arity $k \geq 1$.

**Proposition 2.6.7**: Let $f \in \mathcal{M}_n$. The following statements are equivalent:

1. $f$ is a constant.
2. For all $k$ and $\alpha, \beta \colon n \to k$, we have $f\alpha = f\beta$.
3. There exists a homomorphism $* \to \mathcal{M}$ sending the unique $(\operatorname{ar} f)$-ary constant element to $f$.

Furthermore,

- if $n \geq 2$, this is equivalent to $f$ being fixed by all endomorphisms $n \to n$, and
- if $n = 1$, this is equivalent to $f$ being degenerate.

*Proof*: $2 \Rightarrow 1$ is clear, because for each $j$ we have $f\iota = f\,{}_n(j \mapsto n)_{n+1}$.

$1 \Rightarrow 2$: If $\operatorname{ar} f = 1$ then this is clear: The assumption implies $f\iota = f(0 \mapsto 1)$, and given $\alpha, \beta \colon 1 \to k$, we can apply $(\alpha(0), \beta(0))_k$ to the initial equation and arrive at $f\alpha = f\beta$. If the arity is strictly higher, we can invoke Proposition 2.6.2 inductively to see that $f = u\eta$ where $u$ is a unary and $\eta \colon 1 \hookrightarrow \operatorname{ar} f$ an injection. by Lemma 2.6.3 we know that

$$\eta[\operatorname{Iness} u] \subseteq \operatorname{Iness} u\eta = \operatorname{Iness} f = \emptyset,$$

hence $u$ is a constant as well. Applying $1 \Rightarrow 2$ to $u$, we see that $f\alpha = u\eta\alpha = u\eta\beta = f\beta$ must hold for all $\alpha$ and $\beta$.

$2 \Rightarrow 3$: If $*_n = \{c_n\}$, then the map $c_k = c_n\alpha \mapsto f\alpha$ is well-defined regardless of choice of $\alpha \colon n \to k$ precisely because $f\alpha = f\beta$ whenever $\alpha, \beta \colon n \to k$.

$3 \Rightarrow 2$: Let $*_n = \{c_n\}$ and $h \colon * \to \mathcal{M}$ be the asserted homomorphism. Since $c_n\alpha = c_n\beta$ for all $\alpha, \beta \colon n \to k$, we have $f\alpha = h(c_n)\alpha = h(c_n\alpha) = h(c_n\beta) = h(c_n)\beta = f\beta$.

The second to last point follows because if $f$ is fixed by all endomorphisms, then for every $i \in n$ it is also fixed under ${}_n(i \mapsto \min(n \setminus \{i\}))_n$, hence $i \in \operatorname{Iness} f$.

Finally, last point follows because degenerate unaries having one coordinate implies them having only inessential coordinates, so they are automatically constants. □

**Corollary 2.6.8**: For every minion $\mathcal{M}$, there is the same number of constants in each arity.

*Proof*: First note that if $f$ is a constant, $f\alpha$ is a constant as well, regardless of the arities involved.

If $C_n$ denotes the constants of $\mathcal{M}$ in arity $n$, then applying the unique minor operation $\alpha \colon n \to 1$ gives a map $a \colon C_n \to C_1$. This map is well-defined because of the opening remark. Claim: This map is a bijection, as witnessed by the inverse map $c \mapsto c\beta$, where $\beta$ is an arbitrary operation $\beta \colon 1 \to n$. Indeed, we have $f\alpha\beta = f\,\mathrm{id}_n = f$ because $f \in C_n$ is a constant, and on the other hand $c\beta\alpha = c\,\mathrm{id}_1 = c$, because $c \in C_1$ is a constant. □

### 2.6.3 Essentiality and minion homomorphisms

Here, we want to prove some small properties of essential coordinates in regard to minion homomorphisms. Since inessentiality of a coordinate is expressible as a minor condition, it is intuitively clear that minion homomorphisms $h$ "preserve inessential coordinates" in the sense that $i \in \operatorname{Iness} f$ implies $i \in \operatorname{Iness} h(f)$. The remaining statements provide fine-grained statements about the converse, i.e., when *essentiality* is preserved. These insights will become relevant in Chapter 3.

**Lemma 2.6.9**: Let $h\colon \mathcal{M} \to \mathcal{N}$ be a minion homomorphism, $f \in \mathcal{M}_n$. We have the following statements:
1. $\operatorname{Iness} f \subseteq \operatorname{Iness} h(f)$, i.e., homomorphisms preserve inessentiality of coordinates.
2. If $h$ is injective in arity $n$, then $\operatorname{Ess} f \subseteq \operatorname{Ess} h(f)$.
3. If $h$ is injective in arity $n+1$, then $\operatorname{Ess} f = \operatorname{Ess} h(f)$.
4. If $n \geq 2$ and $h$ is injective in arity $n$, then $\operatorname{Ess} f = \operatorname{Ess} h(f)$.

*Proof*: Throughout, let $\gamma_j := {}_n(j \mapsto n+1)_{n+1}$ and let $\iota\colon n \to n+1$ be the canonical inclusion. 1) Let $j \in \operatorname{Iness} f$, i.e., $f\iota = f\gamma_j$. Then $h_n(f)\gamma_j = h_{n+1}(f\gamma_j) = h_{n+1}(f\iota) = h_n(f)\iota$, so $j \in \operatorname{Iness} h(f)$.

2) Let $j \in \operatorname{Ess} f$, i.e., $f\iota \neq f\gamma_j$. Injectivity of $h_n$ yields

$$h_n(f)\iota = h_n(f\iota) \neq h_n(f\gamma_j) = h_n(f)\gamma_j,$$

which means that $j \in \operatorname{Ess} h(f)$.

3) Let $j \in \operatorname{Iness} h_n(f)$, i.e., $h_n(f)\iota = h_n(f)\gamma_j$. This implies $h_{n+1}(f\iota) = h_{n+1}(f\gamma_j)$, and by injectivity of $h_{n+1}$ we conclude $f\iota = f\gamma_j$, i.e., $j \in \operatorname{Iness} f$. Combined with the first statement, we conclude that $\operatorname{Iness} f = \operatorname{Iness} h(f)$, whence $\operatorname{Ess} f = \operatorname{Ess} h(f)$.

4) Let $\varepsilon_j := {}_n(j \mapsto \min(n \setminus j))_n$. Since $\operatorname{ar} f \geq 2$, we have that $j \in \operatorname{Iness} f$ if and only if $f = f\varepsilon_j$ which by injectivity of $h_n$ is equivalent to $h(f) = h(f\varepsilon_j) = h(f)\varepsilon_j$, i.e., $j \in \operatorname{Iness} h(f)$. We reach the same conclusion as in the last point. □

### 2.6.4 Bounded essential arity and finite generatedness

As we have seen before, we can define minions by generators and relations, just as for algebras. This gives us a notion of being "finitely generated", which is an important notion of smallness. This notion can be characterized by its growth, as we will expand on in Chapter 6.

In this section, we will establish a few easy technical results about this notion, which will be needed later.

**Definition 2.6.4**: A minion $\mathcal{M}$ is said to have
- *essential arity at most $N$* if all elements $f \in \mathcal{M}$ have at most $N$ essential coordinates, and
- *bounded essential arity* if there is an $N$ such that it has essential arity at most $N$.

**Lemma 2.6.10**: If $\mathcal{M}$ has essential arity at most $N$ and $h\colon \mathcal{M} \twoheadrightarrow \mathcal{N}$ is a surjective minion homomorphism, then $\mathcal{N}$ has essential arity at most $N$ as well.

*Proof*: If $n > N$, then any $f \in \mathcal{M}_n$ must have by assumption at least $N - n$ inessential coordinates. Since homomorphisms preserve inessential coordinates (Lemma 2.6.9), we know that $h(f)$ must have at least $N - n$ inessential coordinates as well. □

**Definition 2.6.5**: A minion $\mathcal{M}$ is called *finitely generated* if there is a finite set $S$ of generators, i.e., a finite[26] (indexed) set $S$ such that the smallest subminion $\langle S \rangle$ containing $S$ is $M$.

*Remark*: With this language, we can rephrase Observation 2.5.1 by saying "gadget minions $\mathcal{M}_\varphi$ of primitive positive sentences $\varphi$ are precisely finitely generated minions":

On one hand, gadget minions are quotients of the free minion generated by all of the (finitely many) symbols involved.

On the other hand, to each finitely generated minion we can pick a finite set of generators $S$ and enforce all of the minor conditions which hold between them. These are finitely many, because a minor condition $f\alpha \approx g\beta$ can easily be seen to be equivalent to a minor condition where both $\alpha$ and $\beta$ have codomain at most $\operatorname{ar} f + \operatorname{ar} g$.

Recall that we can explicitly describe the generated subminion $\langle S \rangle$ as

$$\langle S \rangle_n = \{ s\alpha \mid s \in S_k, \alpha\colon k \to n \}.$$

**Proposition 2.6.11**: A minion is finitely generated if and only if it is locally finite and has bounded essential arity.

*Proof*: Assume first that $\mathcal{M}$ is finitely generrated by an (indexed) set $S$.

Then it is a quotient of the free minion $\langle S \rangle \cong \sum_k S_k \cdot \mathcal{P}^k$, which is a sum of finitely many locally finite minions – hence locally finite.

Regarding bounded essential arity, note that an element of the form $s\alpha$ can have at most as many essential coordinates as $s$ due to Corollary 2.6.5, so the essential arity of $\mathcal{M}$ is bounded by $\max\{ |\operatorname{Ess} s| \mid s \in S \}$.

For the converse, if $\mathcal{M}$ has essential arity bounded by $N$ and is locally finite, then the elements up to arity $N$ form a finite set of generators. □

**Proposition 2.6.12**: Let $\mathcal{M}$ have essential arity at most $M$ and $\mathcal{N}$ have essential arity at most $N$. Then $\mathcal{M} \times \mathcal{N}$ has essential arity at most $M + N$.

*Proof*: Let $(f, g) \in (\mathcal{M} \times \mathcal{N})_n$. We have

$$i \in \operatorname{Iness}(f, g) \Leftrightarrow (f, g)\iota = (f, g)\sigma \Leftrightarrow (f\iota, g\iota) = (f\sigma, g\sigma) \Leftrightarrow i \in \operatorname{Iness} f \wedge i \in \operatorname{Iness} g,$$

where $\sigma = {}_i(n \mapsto n+1)_n$ is the map from Definition 2.6.2. Therefore, $\operatorname{Iness}(f, g) = \operatorname{Iness} f \cap \operatorname{Iness} g$, which dually implies $\operatorname{Ess}(f, g) = \operatorname{Ess} f \cup \operatorname{Ess} g$ which contains at most $N + M$ elements. □

**Corollary 2.6.13**: The sum and product of two finitely generated minions is also finitely generated.

[26] Note that we actually mean *finite* instead of only *locally finite* here: $S = (S_n)_n$ is finite if and only if it is locally finite and of "finite support" (i.e., nonempty only in finitely many arities).

*Proof*: Since being locally finite and – by the preceding proposition – having bounded essential arity is stable under products, by Proposition 2.6.11 being finitely generated is as well. □

**Proposition 2.6.14**: Bounded essential arity is preserved by finite sums. Consequently, the sum of two finitely generated minions is also finitely generated.

*Proof*: If $\mathcal{M}$ has essential arity at most $M$ and $\mathcal{N}$ has essential arity at most $N$, then any element $f$ of $\mathcal{M} + \mathcal{N}$ either lies in the $\mathcal{M}$- or the $\mathcal{N}$-component, so has essential arity at most $\max\{M, N\}$.

As before, since local finiteness is preserved under sums, the claim for finite generatedness follows. □

**Proposition 2.6.15**: Bounded essential arity and being finitely generated are heredirary properties, i.e., they are stable under taking subminions.

*Proof*: The first statement is clear because having bounded essential arity is an elementwise condition. For the second statement, use Proposition 2.6.11 by noting that local finiteness is hereditary. □

**Proposition 2.6.16**: if $\mathcal{M}$ is locally finite and $\mathcal{N}$ has bounded essential arity, then a homomorphism $h\colon \mathcal{M} \to \mathcal{N}$ factors through a minion $\mathcal{K}$ which is finitely generated.

*Proof*: Take $\mathcal{K} := \operatorname{Im} h$: Since the image of a homomorphism from a locally finite minion must be locally finite and subminions remain to be of bounded essential arity, $\operatorname{Im} h$ must be finitely generated. Thus, $h$ factors as usual into the corestriction $\tilde{h}$ onto its image followed by the inclusion $\iota\colon \operatorname{Im} h \to \mathcal{N}$. □

For a last remark, recall the theorem we quoted in the introduction:

**Proposition** ([19, Proposition 5.14]): Let $\tau$ be a finite signature and $\underline{A}$, $\underline{B}$ be $\tau$-structures over finite domain. If there is a minion homomorphism from $\operatorname{Pol}(\underline{A}, \underline{B})$ to a minion $\mathcal{M}$ of bounded essential arity, then $\operatorname{PCSP}(\underline{A}, \underline{B})$ is NP-hard.

The above considerations imply that without loss of generality, we can take the "hardness witness" $\mathcal{M}$ to be finitely generated.

# 3 Representability: Cayley Theorems for Minions

In this section, we shall prove a "Cayley theorem for minions", by which we mean that every minion is isomorphic to a concrete minion. The results demonstrated here are known to experts[27], but to the knowledge of the author, proofs have not yet been published anywhere.

Furthermore, we should point out that there is a rather simple "abstract cayley theorem", which states the following:

**Proposition** (cf. [31]): For every minion $\mathcal{M}$, we have

$$\mathcal{M} \cong \mathrm{Pol}_{\underline{\mathsf{Min}}}(\mathcal{P}, \mathcal{M}).$$

*Proof sketch*: By Corollary 2.4.6, we know that the elements of $\mathcal{M}_n$ correspond to minion homomorphisms $\mathcal{P}^n \to \mathcal{M}$. In the backward direction, this amounts to evaluating at the "rainbow tuple" $\mathrm{id}_n = (0, ..., n-1)$. It is straightforward to verify that this correspondence is compatible with minor operations. □

This is a conceptually elegant statement, which establishes that every minion is isomorphic to a polymorphism minion[28], however it does not show that every minion is isomorphic to a *concrete minion*, i.e., a polymorphism minion of relational structures.

We will open with the discussion of a classical duality to establish that the notions of "concrete minion" and of "function minion" (recall Definition 2.2.2) coincide. This means that for the desired cayley theorem, it is sufficient to prove that every minion is isomorphic to a function minion.

## 3.1 Preservation of relation pairs

The goal of this section is to establish that every subminion of $\mathcal{O}(A, B)$ can be defined as the polymorphism minion $\mathrm{Pol}(\underline{A}, \underline{B})$ where $\underline{A}$ and $\underline{B}$ are suitable structures over base sets $A$ and $B$, respectively. To do so requires the notion of preservation, which provides a duality – or, more formally, a Galois connection – between subminions of $\mathcal{O}(A, B)$ and certain sets of relation pairs.

We follow results from the work of Pöschel and Lehtonen [32], who generalized the notion of preservation to the many-sorted setting.[29] Crucially, we can apply this to our setting, because the polymorphism minion $\mathrm{Pol}(\underline{A}, \underline{B})$ can be viewed as a many-sorted algebra with domain $A + B$, where $A$ is of sort 1 and $B$ of sort 2, and for which all operations map from a power of sort 1 to sort 2.

For easier, self-contained presentation, we will not only consider functions $A^n \to B$, but also "$X$-ary" functions $A^X \to B$ where $X$ is an arbitrary finite set. This simplifies some arguments – crucially, the canonical inclusion $l \hookrightarrow A^{A^l}$ of $i$ to the $i$-th projection can be viewed as an $A^l$-ary relation.

---

[27] At the very least, Libor Barto and Manuel Bodirsky are aware of this.

[28] Recall that Definition 2.2.1 defines "polymorphism minion" in an abstract way, since the notion of "polymorphism" can exist in contexts unrelated to relational structures (such as in this case).

[29] Earlier, Pippenger [33] developed a Galois connection between function minions and relation pairs which is almost suitable for our purposes. However, he restricts to the boolean domain, while Pöschel and Lehtonen generalized his ideas to the full arbitrary many-sorted setting over arbitrary finite domains.

While considering functions with input set $A^X$ as “operations” may seem a little unusual, there is an intuitive interpretation: If for instance $X = \{\text{apple}, \text{orange}\}$, then an $X$-ary function would take a named input tuple $(\text{apple}\colon a_1, \text{orange}\colon a_2)$ and return an element of $B$ as output. From a programming language perspective, this is the difference of a function $f$ taking *positional arguments* as opposed to a function taking *named arguments*[30].

```python
def f(a: int, b: int, /) -> int
    return a + b

print(f(3, 5))
```

```python
def g(*, a: int, b: int) -> int
    return a + b

print(g(a=3, b=5))
```

Listing 1: Two functions using unamed an named arguments respectively, implemented in the programming language python. Both invocations output the value 8.

The relational counterpart of this generalization, the concept of an “$X$-ary relation”, would be a subset of $A^X$. This is known in database theory as the “named perspective” [34, Definitions 2.1 and 2.2]. We will show in Section A.2 that generalizing the notion of coordinate to arbitrary finite sets yields an equivalent notion of minion.

We will now introduce the central property of preservation, which is just spelling out what we need for something to be a homomorphism. To do so, we will need to talk about sequences of tuples, i.e., doubly indexed entities. To emphasize the distinction between the kinds of indices, we will place one as a subscript and another as a superscript. Since no notion of multiplication exists in this context, this should be perfectly unambiguous.

**Definition 3.1.1** (preservation): Let $A, B, X, Y$ be sets with $X, Y$ finite. We call $f\colon A^X \to Y$ an *$X$-ary function* from $A$ to $B$, and a pair $(R, S)$ where $R \subseteq A^Y$ and $S \subseteq B^Y$ a *$Y$-ary relation pair* on $(A, B)$. We say that $f$ *preserves* $(R, S)$, written

$$f \triangleright_{X,Y} (R, S),$$

if for any $X$-tuple $(r_x)_{x \in X}$ of elements $r_x = r_x^\bullet$ in $R$, we have

$$\left(f(r_\bullet^y)\right)_{y \in Y} \in S.$$

If the arities are clear from the context or do not require emphasis, we shall omit the subscript and just write $f \triangleright (R, S)$. Furthermore, if

$$\mathcal{F} = (\mathcal{F}_X)_{X \in \underline{\mathsf{FinSet}}_{>\emptyset}}, \qquad \mathcal{R} = (\mathcal{R}_Y)_{Y \in \underline{\mathsf{FinSet}}_{>\emptyset}}$$

are a collection of sets of $X$-ary functions and a collection of sets of $Y$-ary relation pairs, respectively, we write $\mathcal{F} \triangleright \mathcal{R}$ if and only if $f \triangleright_{X,Y} (R, S)$ for all sets $X$, $Y$, functions $f \in \mathcal{F}_X$ and relation pairs $(R, S) \in \mathcal{R}_Y$.

[30] In programming also known as “keyword arguments”.

$$\underbrace{f}_{A^X \to B} \;\triangleright_{X,Y} (\underbrace{R}_{\subseteq A^Y}, \underbrace{S}_{\subseteq B^Y})$$

Figure 8: The notation used to denote "$f$ preserves the relation pair $(R, S)$".

This definition benefits from some visualization, for which we shall consider $X = n$ and $Y = k$. It is common to view elements of the relations, i.e., $R \subseteq A^n$ and $S \subseteq B^n$, as column vectors. Given $n$ column vectors $(r_i)_{i \in n}$ of height $k$ with values in $A$, we can apply the function in each row to obtain a new column vector of height $k$, now with values in $B$. If for every pair of $n$ column vectors taken from $S$, the resulting column vector is an element of $R$, then the condition is fulfilled. Visually, this looks as follows:

$$\begin{array}{cccc} r_0^0 & r_1^0 & \cdots & r_{n-1}^0 \\ r_0^1 & r_1^1 & \cdots & r_{n-1}^1 \\ \vdots & \vdots & \ddots & \vdots \\ r_0^{k-1} & r_1^{k-1} & \cdots & r_{n-1}^{k-1} \\ \| & \| & & \| \\ \underbrace{r_0^\bullet}_{\in R} & \underbrace{r_1^\bullet}_{\in R} & & \underbrace{r_{n-1}^\bullet}_{\in R} \end{array} \quad \mapsto \quad \begin{array}{c} f(r_\bullet^0) \in S \\ f(r_\bullet^1) \in S \\ \vdots \\ f(r_\bullet^{k-1}) \in S \end{array}$$

Figure 9: $f \triangleright_{n,k} (R, S)$ means: If the columns are in $R$, the application of $f$ to the rows lies in $S$.

Associated to this notion of preservation are operators Inv and Pol, whose definition should not surprise any reader familiar with the notion of a Galois connection.

**Definition 3.1.2** (Inv, Pol): Let $A, B$ be sets, $\mathcal{F}$ an indexed set of functions and $\mathcal{P}$ an indexed set of relation pairs as in the last definition. We define

- $\mathrm{Inv}(\mathcal{F})_Y := \{(R, S) \in A^Y \times B^Y \mid \mathcal{F} \triangleright (R, S)\}$
- $\mathrm{Pol}(\mathcal{P})_X := (\{f\colon A^X \to B \mid f \triangleright \mathcal{P}\}),$

and by $\mathrm{Inv}(\mathcal{F})$ resp. $\mathrm{Pol}(\mathcal{P})$ we mean the respective $\underline{\mathsf{FinSet}}_{>\emptyset}$-indexed sets. As before, we will abbreviate $\mathrm{Inv}(\{f\}) =: \mathrm{Inv}(f)$ and $\mathrm{Pol}(\{(R, S)\}) =: \mathrm{Pol}((R, S))$.

To formulate some basic properties of preservation, in particular to minor operations, recall that given a set $B$, a function $\gamma\colon Y \to Y'$ induces the *precomposition* (or *pullback*) $\gamma^*\colon B^{Y'} \to B^Y$, which sends $\left(b^{y'}\right)_{y' \in Y'}$ to $\left(b^{\gamma(y)}\right)_{y \in Y'}$, which we can write more concisely as $b_{\gamma(\bullet)}$. In particular, if $\beta\colon Y' \to Y$, the image $\beta^*[R]$ consists of elements of the form $r'_{\beta(\bullet)}$ where $r' \in R$, and if $\gamma\colon Y \to Y'$, then the preimage $(\gamma^*)^{-1}[R]$ consists of elements $r \in B^Y$ such that $r_{\gamma(\bullet)} \in R$.

Furthermore, recall the notation $f[S]$ for the image of a set under a function and $f^{-1}[T]$ for the preimage, respectively.

Not all of the properties we prove are strictly necessary for the proof of Corollary 3.1.4, however they are simple to prove and all interesting on their own.

**Proposition 3.1.1**: Let $\alpha\colon X \to X'$, $\beta\colon Y' \to Y$, and $\gamma\colon Y \to Y'$.
1. If $f \rhd_{X,Y} (R,S)$, then $f\alpha \rhd_{X',Y} (R,S)$
2. If $f \rhd_{X,Y} (R,S)$, then $f \rhd_{X,Y'} (\beta^*[R], \beta^*[S])$, where $\beta^*$ is the induced map $B^Y \to B^{Y'}$
3. If $f \rhd_{X,Y} (R,S)$, then $f \rhd_{X,Y'} \left((\gamma^*)^{-1}[R], (\gamma^*)^{-1}[S]\right)$
4. If $f \rhd (R,S)$ and $S \subseteq S'$, then $f \rhd (R,S')$
5. If $f \rhd (R,S)$ and $R \subseteq R'$, then $f \rhd (R',S)$
6. If $f \rhd (R,S_i)$ for all $i \in I$, then $f \rhd (R, \cap_i S_i)$
7. If $f \rhd_{X,Y_i}$ (R_i, $S_i$) for all $i$ in some finite index set $I$, then $f \rhd_{X,\sum_i Y_i} \left(\prod_i R_i, \prod_i S_i\right)$.

Before delving into the proof, a few remarks are in order.

First, note that the compatibility with minor operations implies that it is sufficient to consider arities which are finite ordinals: Indeed, if $X$ is enumerated by $\xi\colon n \to X$ and $Y$ is enumerated by $v\colon k \to Y$, then

$$\begin{aligned} & f \rhd_{X,Y} (R,S) \\ \Leftrightarrow\ & f\xi^{-1} \rhd_{n,Y} (R,S) \\ \Leftrightarrow\ & f\xi^{-1} \rhd_{n,k} \left((v^{-1})^*[R], (v^{-1})^*[S]\right) \end{aligned}$$

because both $\xi$ and $v$ are bijective.

Second, let us illustrate the subtle – but substantial – difference between the operations $\beta$, used as $\beta^*[_]$ in point 2, and $\gamma$, used as $(\gamma^*)^{-1}[_]$, in point 3. Formally, we are taking the image and preimage of the pullback maps, respectively.

As mentioned before, we will view elements of a relation $R$ as column vectors.

- Let $\beta_1\colon 2 \to 3$ be the inclusion. Then the operation $\beta^*[_]$ takes a ternary relation $R \subseteq B^3$ and produces a binary relation

$$\beta_1^*[R] = \{\beta_1^*(r) \mid r \in R\} = \left\{ \begin{pmatrix} r^{\beta_1(0)} \\ r^{\beta_1(1)} \end{pmatrix} \middle| \begin{pmatrix} r^0 \\ r^1 \\ r^2 \end{pmatrix} \in R \right\} = \left\{ \begin{pmatrix} r^0 \\ r^1 \end{pmatrix} \middle| \begin{pmatrix} r^0 \\ r^1 \\ r^2 \end{pmatrix} \in R \right\},$$

  i.e., injective, monotone operations correspond to the deletion of rows.
- Let $(\beta_1)_2\colon 3 \to 2$ map 2 to 0 but leave $\{0,1\}$ as-is. Then a binary relation $R \subseteq B^2$ induces a ternary relation

$$\beta_2^*[R] = \{\beta_2^*(r) \mid r \in R\} = \left\{ \begin{pmatrix} r^{\beta_2(0)} \\ r^{\beta_2(1)} \\ r^{\beta_2(2)} \end{pmatrix} \middle| \begin{pmatrix} r^0 \\ r^1 \end{pmatrix} \in R \right\} = \left\{ \begin{pmatrix} r^0 \\ r^1 \\ r^0 \end{pmatrix} \middle| \begin{pmatrix} r^0 \\ r^1 \end{pmatrix} \in R \right\},$$

  which corresponds to duplicating the first row and appending it to the end.
- Now, let $\gamma_1 \coloneqq \beta_1$, and consider what $(\gamma^*)^{-1}[_]$ does. Since we are taking the preimage, the direction is reversed: Given a binary relation $R \subseteq B^2$, we induce the relation

$$(\gamma_1^*)^{-1}[R] = \{r \mid \gamma_1^*(r) \in R\} = \left\{ \begin{pmatrix} r^0 \\ r^1 \\ r^2 \end{pmatrix} \in B^3 \middle| \begin{pmatrix} r^{\gamma_1(0)} \\ r^{\gamma_1(1)} \end{pmatrix} \in R \right\} = \left\{ \begin{pmatrix} r^0 \\ r^1 \\ r^2 \end{pmatrix} \in B^3 \middle| \begin{pmatrix} r^0 \\ r^1 \end{pmatrix} \in R \right\}.$$

This amounts to adding a row at the bottom which may be filled with arbitrary elements.

- Analogously, let $\gamma_2 := \beta_2$ and consider a ternary relation $R \subseteq B^3$. Then

$$(\gamma_2^*)^{-1}[R] = \{r \mid \gamma_2^*(r) \in R\} = \left\{ \begin{pmatrix} r^0 \\ r^1 \end{pmatrix} \in B^2 \,\middle|\, \begin{pmatrix} r^{\gamma_2(0)} \\ r^{\gamma_2(1)} \\ r^{\gamma_2(2)} \end{pmatrix} \in R \right\} = \left\{ \begin{pmatrix} r^0 \\ r^1 \end{pmatrix} \in B^2 \,\middle|\, \begin{pmatrix} r^0 \\ r^1 \\ r^0 \end{pmatrix} \in R \right\}.$$

  This amounts to selecting precisely these tuples in the relation where the first and last coordinate agree, and omitting the redundant row.
- In the case where $\beta$ and $\gamma$ are bijective, the associated operations on the relations correspond in both cases to the permutation of rows.

So to summarize, the contravariant action $\beta \rightsquigarrow \beta^*[_]$ corresponds to the usual row operations such as deletion or duplication, whereas the covariant action $\gamma \rightsquigarrow (\gamma^*)^{-1}[_]$ corresponds to expansion by arbitrary columns and restriction of the tuples according to agreement on some columns.

These operations are well-known in the context of database theory, although under different names ("selection", "join"); the interested reader is referred to [35, Section 4.4].

*Remark (for category theorists)*: Curiously enough, we have thus discovered that a "relational minion", i.e., a set of relation pairs for each arity, is canonically a *minion*, but also what we might call a *co-minion*, i.e., a functor $\underline{\mathsf{FinSet}_{>\emptyset}}^{\mathsf{op}} \to \underline{\mathsf{Set}}$. Furthermore, we note that the co- and contravariant action agree on invertible operations, because when $\gamma$ is bijective, we have $(\gamma^*)^{-1} = (\gamma^{-1})^*$ by functoriality of the pullback, and because the preimage of an invertible function is the same as the image of its inverse, we know that the covariant action of $\gamma$ agrees with the contravariant action of $\gamma^{-1}$.

We can summarize this information by stating that $\mathrm{Inv}(\mathcal{F})$ is a functor

$$\left(\underline{\mathsf{FinSet}_{>\emptyset}}^{\mathsf{op}} + \underline{\mathsf{FinSet}_{>\emptyset}}\right)/\{\beta^{\mathsf{op}} \sim (\beta^{-1}) \mid \beta \colon X \to X' \text{ iso}\} \to \underline{\mathsf{Set}}.$$

To simplify the wording of the proof, we are going to refer to a sequence $(r_x)_{x \in X} \in R^X$ for a $Y$-ary relation $R \subseteq A^Y$ as *an* $(X,Y)$*-matrix* $r^{\bullet}_{\bullet}$ *with columns in* $R$, meaning that $r^{\bullet}_{\bullet} \in A^{X \times Y}$ is a doubly-indexed sequence of elements of $A$ such that for every $x \in X$ we have $r^{\bullet}_x \in R$. We then have $f \triangleright_{X,Y} (R,S)$ if and only if for every $(X,Y)$-matrix $r^{\bullet}_{\bullet}$ with columns in $R$, we have $(f(r^y_{\bullet}))^{y \in Y} \in S$.

*Proof of Proposition 3.1.1*: 1) Let $f \triangleright (R,S)$, $\alpha \colon X \to X'$, and $r^{\bullet}_{\bullet}$ a $(X',Y)$-matrix with columns in $R$. Since every column $r^{\bullet}_{\alpha(x)}$ of the $(X,Y)$-matrix $r^{\bullet}_{\alpha(\bullet)}$ lies in $R$, we know that

$$(f\alpha(r^y_{\bullet}))^y = \left(f\left(r^y_{\alpha(\bullet)}\right)\right)^y \in S.$$

2) Let $f \triangleright (R,S)$, $\beta \colon Y' \to Y$, and $q^{\bullet}_{\bullet}$ an $(X,Y')$-matrix with columns in $\beta^*[R]$. By definition of the pullback and the image, this means that every $q^{\bullet}_x$ is of the form $r^{\beta(\bullet)}_x$ for some $r^{\bullet}_x \in R$. Therefore,

$$\left(f\left(q^{y'}_{\bullet}\right)\right)^{y'} = \left(f\left(r^{\beta(y')}_{\bullet}\right)\right)^{y'} = \beta^*((f(r^y_{\bullet}))^y) \in \beta^*[S].$$

3) Let $f \triangleright (R, S)$, $\gamma\colon Y \to Y'$, and $q_\bullet^\bullet$ an $(X, Y')$-matrix with columns in $(\gamma^*)^{-1}[R]$. By definition of the preimage and the pullback, this means that $q_x^{\gamma(\bullet)} \in R$ for every $x \in X$. We conclude that

$$\gamma^*\left(\left(f\left(q_\bullet^{y'}\right)\right)^{y'}\right) = \left(f\left(q_\bullet^{\gamma(y)}\right)\right)^{y} \in S,$$

which means that indeed $\left(f\left(q_\bullet^{y'}\right)\right)^{y'} \in (\gamma^*)^{-1}[S]$.

Points 4) and 5) are clear because they are a weakening of the implication and a strengthening of the assumptions, respectively, and point 6) follows straightforwardly from the definition of intersection.

To see point 7), let $r_\bullet^\bullet$ be an $\left(X, \sum_i Y_i\right)$-indexed tuple with $(r_x^y)^{y \in \sum_i Y_i} \in \prod_i R_i$ for every $x$ in $X$. Projecting to a component $R_j$ gives us an $(X, Y_j)$-indexed tuple such that $(r_x^y)^{y \in Y_j} \in R_j$ for every $x$. The assumption $f \triangleright (R_j, S_j)$ implies now

$$(f(r_\bullet^y))^{y \in R_j} \in S_j.$$

Since this holds for every $j$, this combines to

$$(f(r_\bullet^y))^{y \in \sum_i R_i} \in \prod_i S_i,$$

as desired. □

We can strengthen the last point to an equivalence.

**Corollary 3.1.2**: Let $A, B$ be arbitrary sets and let $(R_i, S_i)_{i \in I}$ be a finite set of relation pairs on $(A, B)$. Then

$$\mathrm{Pol}\left(\{(R_i, S_i)\}_{i \in I}\right) = \mathrm{Pol}\left(\prod_{i \in I} R_i, \prod_{i \in I} S_i\right)$$

*Proof*: The forward inclusion follows from point 7 of Proposition 3.1.1. For $\supseteq$, assume

$$f \triangleright \left(\prod_{i \in I} R_i, \prod_{i \in I} S_i\right).$$

Since the inclusion $\iota_{X_j}\colon X_j \hookrightarrow \sum_i X_i$ pulls back to the projection $\pi_{X_j}\colon B^{\sum_i X_i} \twoheadrightarrow B^{X_j}$, we conclude from point 2 of Proposition 3.1.1 that

$$f \triangleright \left(\iota_{X_j}^*\left[\prod_{i \in I} R_i\right], \iota_{X_j}^*\left[\prod_{i \in I} S_i\right]\right) = \left(\pi_{X_j}\left[\prod_{i \in I} R_i\right], \pi_{X_j}\left[\prod_{i \in I} S_i\right]\right) = (R_j, S_j),$$

as desired. □

**Definition 3.1.3**: If $A$ is finite, let $p_l := \{\pi_j \mid j \in l\} \subseteq A^{A^l}$ be the set of $l$-ary projections, each viewed as an $A^l$-ary relation on $A$ – i.e, $(\pi_j)^{\underline{a}} = \pi_j(\underline{a}) = a_j$.

**Proposition 3.1.3**: Let $A$ be finite, $B$ arbitrary, and $f\colon A^n \to B$. We have

$$f \triangleright_{n,A^l} (p_l, S) \Leftrightarrow \forall n \xrightarrow{\alpha} l\colon f\alpha \in S,$$

where we view $f\alpha\colon A^l \to A$ as an $A^l$-ary relation on $A$, i.e., $(f\alpha)^{\underline{a}} = (f\alpha)(\underline{a})$ as in the previous Definition 3.1.3.

*Proof*: For this proof, first note that there is a one-to-one correspondence to elements $(r_i)$ of $(p_l)^n$ and maps $\alpha\colon n \to l$, because every such $r_i$ must be associated to a unique projection $\pi_j$. With this insight, it becomes straightforward to verify the claim from the definitions:

$$\begin{aligned} & f \triangleright_{n,A^l} (p_l, S) \\ \Leftrightarrow \forall (r_i)_i \in (p_l)^n\colon \quad & \left(f\left((r_i^{\underline{a}})_{i\in n}\right)\right)^{\underline{a}\in A^l} \in S \\ \Leftrightarrow \forall \alpha\colon n \to l\colon \quad & \left(f\left(\left(\pi_{\alpha(i)}^{\underline{a}}\right)_{i\in n}\right)\right)^{\underline{a}\in A^l} \in S \end{aligned}$$

and then it suffices to observe that

$$\left(f\left(\left(\pi_{\alpha(i)}^{\underline{a}}\right)_{i\in n}\right)\right)^{\underline{a}\in A^l} = \left(f\left(\left(a_{\alpha(i)}\right)_{i\in n}\right)\right)^{\underline{a}\in A^l} = ((f\alpha)(\underline{a}))^{\underline{a}\in A^l} = f\alpha. \qquad \square$$

We are now equipped to prove the main result of this section, which asserts that every subminion of $\mathcal{O}(A, B)$ is of the form $\mathrm{Pol}(\underline{A}, \underline{B})$ for some relational structures $\underline{A}$ and $\underline{B}$.

**Corollary 3.1.4**: Let $\mathcal{M}$ be a subminion of $\mathcal{O}(A, B)$ where $A$ is finite. Then

$$\mathcal{M} = \mathrm{Pol}(\{(p_l, \mathcal{M}_l) \mid l \in \mathbb{N}_{>0}\}).$$

Formulated differently, $\mathcal{M} = \mathrm{Pol}(\underline{A}, \underline{B})$, where $\underline{A} := \left(A; (p_l)_{l\in L}\right)$ and $\underline{B} := \left(B; (\mathcal{M}_l)_{l\in L}\right)$ are structures over the (countable) relational signature $\tau := \left(l\colon A^l\right)_{l\in\mathbb{N}}$.

*Proof*: To see $\subseteq$, let $f \in \mathcal{M}_n$. Since $\mathcal{M}$ is a minion, for each $l$ and $\alpha\colon n \to l$, we have $f\alpha \in \mathcal{M}_l$, and so by Proposition 3.1.3, $f \triangleright (p_l, \mathcal{M}_l)$. For the reverse direction, let $f\colon A^n \to B$ preserve $(p_l, \mathcal{M}_l)$ for all $l$. Then certainly it preserves $(p_n, \mathcal{M}_n)$, so applying Proposition 3.1.3 in the other direction we must have $f = f\,\mathrm{id}_n \in \mathcal{M}_n$. $\square$

The following follows from basic facts about so-called *Galois connections*, to which some background can be found in Section A.1.

**Lemma 3.1.5**:
1. $\mathrm{Pol} \circ \mathrm{Inv}$ is extensive (in the lattice of indexed subsets of $\mathcal{O}(A, B)$)
2. $\mathrm{Pol} \circ \mathrm{Inv}$ is monotone (in the lattice of indexed subsets of $\mathcal{O}(A, B)$)
3. $\mathrm{Pol} \circ \mathrm{Inv} \circ \mathrm{Pol} = \mathrm{Pol}$

in particular:

4. $\mathrm{Pol} \circ \mathrm{Inv}$ is a closure operator.

*Proof*: These follow from Section A.1, since Pol and Inv form the Galois connection induced by the relation $\triangleright$. □

**Corollary 3.1.6** ("Pol-Inv = Min"): Let $A$ be finite and $\mathcal{S} := (S_n)_{n>0}$ be an indexed subset of $\mathcal{O}(A, B)$, i.e., for every $n$ we have $S_n \subseteq A^n \to B$.

Then $\mathrm{Pol}(\mathrm{Inv}(\mathcal{S})) = \mathrm{Min}(\mathcal{S})$, where $\mathrm{Min}(\mathcal{S})$ is the minion generated by $\mathcal{S}$.[31]

*Proof*: First, we will prove the fact where $\mathcal{S} = \mathcal{M}$ is a minion. For this let $f \in \mathrm{Pol}(\mathrm{Inv}(\mathcal{M}))_n$, which means that $f \triangleright \mathrm{Inv}(\mathcal{M})$.

*Claim:* $(p_n, \mathcal{M}_n) \in \mathrm{Inv}(\mathcal{M})$. Indeed, because $\mathcal{M}$ is a minion, for every $g \in \mathcal{M}_k$, all $n$-ary minors $g\alpha$ where $\alpha\colon k \to n$ must be in $\mathcal{M}_n$. This implies in one direction that $g \triangleright (p_n, \mathcal{M}_n)$, and because we imposed no restrictions on $g$, we conclude $(p_n, \mathcal{M}_n) \in \mathrm{Inv}(\mathcal{M})$, as desired.

Therefore – by assumption – $f \triangleright (p_n, \mathcal{M}_n)$, which allows us to apply Proposition 3.1.3 to conclude that for all $\alpha\colon n \to n$ we have $f\alpha \in \mathcal{M}_n$ – so in particular $f = f\,\mathrm{id}_n \in \mathcal{M}_n$.

Finally, to show that the statement extends to the case where $\mathcal{S}$ is not a minion, we can argue that $\mathrm{Pol}(\mathrm{Inv}(\mathcal{S})) = \mathrm{Min}(\mathcal{S})$ with the following circular chain of inclusions:

$$\mathrm{Pol}(\mathrm{Inv}(\mathcal{S})) \overset{1}{\subseteq} \mathrm{Pol}(\mathrm{Inv}(\mathrm{Min}(\mathcal{S}))) \overset{!}{=} \mathrm{Min}(\mathcal{S}) \overset{3}{\subseteq} \mathrm{Min}(\mathrm{Pol}(\mathrm{Inv}(\mathcal{S}))) \overset{4}{=} \mathrm{Pol}(\mathrm{Inv}(\mathcal{S})),$$

where in $(1)$ we have used that $\mathrm{Pol} \circ \mathrm{Inv}$ is monotone (Lemma 3.1.5) and Min is extensive, i.e., $\mathcal{S} \subseteq \mathrm{Min}(\mathcal{S})$, in $(!)$ we have used the statement we just proved because $\mathrm{Min}(\mathcal{S})$ is a minion, in $(3)$ we used that Min is monotone and $\mathrm{Pol} \circ \mathrm{Inv}$ is extensive (Lemma 3.1.5), and in $(4)$ we used that $\mathrm{Pol}(\ldots)$ is always a minion (Proposition 3.1.1). □

### 3.1.1 Finite relatedness

Again, let $A$ and $B$ be finite sets. Having phrased Corollary 3.1.4, a natural question is to ask which function minions are defined by *finitely many relations*. In other words: Which subminions are of the form $\mathrm{Pol}(\underline{A}, \underline{B})$ of finite-domain structures $\underline{A}, \underline{B}$ with finite signature? After all, these are the minions of interest when analyzing finite-domain promise CSPs.

Let us introduce a few notions which help to ease the discussion.

[31] i.e., the smallest subminion of $\mathcal{O}(A, B)$ containing $\mathcal{S}$. This operation, Min, is also a closure operator: The generated minion contains its generators (extensivity), removing generators cannot enlarge the generated minion (monotonicity), and the minion generated by a set of generators which already form a minion is clearly that same minion (idempotence).

**Definition 3.1.4**: Let $\mathcal{M} \leq \mathcal{O}(A, B)$ be a subminion.
1. A set $\mathcal{S}$ of relation pairs on $A$ and $B$ is called a *relational basis* of $\mathcal{M}$ if $\mathcal{M} = \mathrm{Pol}(\mathcal{S})$
2. $\mathcal{M}$ is called *finitely related* if it has a finite relational basis.

Could it be that all sub-minions of $\mathcal{O}(A, B)$ are finitely related? By the following cardinality argument, the answer is no.

**Proposition 3.1.7**: There exist subminions of $\mathcal{O}(2, 2)$ which are not finitely related.

*Proof*: We will argue that there are countably many finitely related minions, but uncountably many subminions of $\mathcal{O}(2)$.

For the first part, note that any finite set of relations must have a finite definition, which can be written down as a finite word over some fixed alphabet $\Sigma$. Since the set $\Sigma^*$ of all words over a fixed, finite alphabet is countable, we have a partial surjection from $\Sigma^*$ to finitely related subminions of $\mathcal{O}(2, 2)$ (and in fact of $\mathcal{O}(A, B)$ for any finite set $A, B$).

That there are continuum many boolean minions has been established by Sparks in [36]. □

But when is a minion $\mathcal{M} \leq \mathcal{O}(A, B)$ finitely related? Brakensiek and Guruswami provided an abstract characterization in [20]. Paraphrased slightly to fit the narrative of this thesis, the statement is:

**Proposition 3.1.8** (cf. [20, Theorem 6.5]): Let $A, B$ be finite sets. A minion $\mathcal{M} \leq \mathcal{O}(A, B)$ is finitely related if and only if there exists an $r$ such that for all $f \in \mathcal{O}(A, B)_k$, we have $f \in \mathcal{M}_k$ if and only if all $r$-minors $f\alpha \in \mathcal{M}_r$.

Let us introduce some language for this central property.

**Definition 3.1.5**: We call $\mathcal{M} \leq \mathcal{N}$ an *$r$-strong subminion (of $\mathcal{N}$)* if for all $k$ and $f \in \mathcal{N}_k$ we have $f \in \mathcal{M}_k$ if and only if $f\alpha \in \mathcal{M}_r$ for all $\alpha\colon k \to r$.

*Remark*: This is not an intrinsic property of $\mathcal{M}$: If $\mathcal{M} \leq \mathcal{N}$ is $r$-strong in $\mathcal{N}$ and $\mathcal{M}' \leq \mathcal{N}$ is a different subminion isomorphic to $\mathcal{M}$, then there is a priori no reason why $\mathcal{M}'$ should be an $r$-strong subminion of $\mathcal{N}$ as well.

*Remark (For category theorists)*: Viewing subobjects as monomorphisms, we might adapt the definition to say that $\iota\colon \mathcal{M} \hookrightarrow \mathcal{N}$ is $r$-strong if all elements $f \in \mathcal{N}$ - viewed as homomorphisms $\mathcal{P}^{\mathrm{ar}\, f} \to \mathcal{N}$ – factor through $\iota$ if and only if all of their $r$-minors factor through $\iota$.

It is not hard to verify that this is equivalent to $\iota$ being $\mathrm{Res}_r$-*cartesian*, where $\mathrm{Res}_r\colon\ \underline{\mathsf{Min}} \to {}_{\mathcal{T}_r}\underline{\mathsf{Set}}$ is the restriction functor assigning a minion its underlying arity-$r$ set, acted upon by the full transformation monoid $\mathcal{T}_r = (r^r, \circ)$.

Let us use this nomenclature to formulate a more specific statement, which implies Proposition 3.1.8.

**Proposition 3.1.9**: Let $\mathcal{M} \leq \mathcal{O}(A, B)$ be a subminion.
1. If $\mathcal{M}$ is an $r$-strong subminion, then $\mathcal{M} = \operatorname{Pol}(p_r, \mathcal{M}_r)$
2. If there is a finite set $\{(R_i, S_i)\}_{i \in I}$ of relation pairs such that $\mathcal{M} = \operatorname{Pol}\big(\{(R_i, S_i)\}_i\big)$, then $\mathcal{M}$ is an $\max_i |R_i|$-strong subminion.

*Proof*: The first point is a direct consequence of Proposition 3.1.3.

For the second point, assume $\mathcal{M} = \operatorname{Pol}\big(\{(R_i, S_i)\}_i\big)$. We will show that any $f \notin \mathcal{M}$ has some minor of arity at most $r$ which is also not in $\mathcal{M}$. Such an $f$ would have to violate preservation of some pair $(R_i, S_i)$ for some $i$. Picking an enumeration $\underline{r}_0, \dots, \underline{r}_{|R_i|-1}$ of $R_i$ – which is finite because it is a subset of $A^{\operatorname{ar} R}$ – this means that there exists a sequence $\underline{r}_{\rho(0)}, \dots, \underline{r}_{\rho(\operatorname{ar} f - 1)}$ of elements of $R_i$ with $\rho\colon \operatorname{ar} f \to |R_i|$ such that $f\big(\underline{r}_{\rho(\bullet)}\big)$, evaluated componentwise, is not in $S_i$. However, $f\big(\underline{r}_{\rho(\bullet)}\big) = f\rho(\underline{r}_\bullet)$, and so we conclude that the $|R_i|$-minor $f\rho$ is also not in $\mathcal{M}$, and neither is any arity-$r$ lift $f\rho\iota$, since $f\rho$ is a minor of it. □

We can use this criterion to verify that a well-known boolean function clone is not finitely related as a subminion of $\mathcal{O}(2)$.

**Proposition 3.1.10**: The function minion

$$\mathcal{B}_\infty := \operatorname{Pol}\big(\{\operatorname{NAZ}(k)\}_{k \in \mathbb{N}} \cup \{\{(0)\}, \{(1)\}\}\big) \leq \mathcal{O}(2)$$

is not finitely related, where $\operatorname{NAZ}(k) := 2^k \setminus \{(0, \dots, 0)\}$ ("not all zero").

*Proof*: According to Proposition 3.1.9, it suffices to check that the inclusion is not $r$-strong for any $r$.

For that, let us define $m_k$ as the boolean $k$-ary function checking if the input tuple contains at least two ones, where $k \geq 2$. Evaluating at the identity matrix shows that $m_k \not\triangleright \operatorname{NAZ}(k)$. However, if $\alpha\colon k \to l$ is noninjective, then $m_k\alpha \triangleright \operatorname{NAZ}(r)$ for all $r$: If $\underline{c}_0, \dots, \underline{c}_{l-1} \in \operatorname{NAZ}(r)$, then the $r$-column vector

$$m_k\alpha\big(\underline{c}_0, \dots, \underline{c}_{l-1}\big) = m_k\Big(\underline{c}_{\alpha(0)}, \dots, \underline{c}_{\alpha(k-1)}\Big)$$

must have a one at some row. Indeed, if $j \in l$ has $|\alpha^{-1}(j)| > 2$, then the $\underline{c}_j$ column must occur at least twice in the evaluation of $m_k$, so any row $i$ at which $\underline{c}_j$ is nonzero evaluates to one. Lastly, note that since $m_k(b, \dots, b) = b$[32], we have that all minors $m_k\alpha$ preserve $\{(0)\}$ as well as $\{(1)\}$.

Now let $r \geq 1$. The above considerations demonstrated that for any $s > r$, the function $m_s$ does not lie in the subminion, while all proper $r$-minors $m_s\alpha$ – with $\alpha$ necessarily noninjective – do. Thus, $r$-strength is violated. □

As we will see briefly, this is equivalent to the well-known fact that this minion is not finitely related as a clone, i.e., there is no finite set of relations $\{R_i\}_{i=1}^n$ such that $\mathcal{B}_\infty = \operatorname{Pol}\big(\{R_i\}_{i=1}^n\big)$.

[32] This requires $k \geq 2$.

## 3.2 Preservation of relations

Throughout this section, let $A$ be a finite set.

The results of the previous section allow for a special case when $A = B$, i.e., when we consider functions whose domain and codomain agree. Indeed, if we consider not $\mathrm{Pol}(\underline{A}, \underline{B})$, but $\mathrm{Pol}(\underline{A}, \underline{A}) =: \mathrm{Pol}(\underline{A})$, i.e., $n$-ary homomorphisms from a structure to itself, we will note that this (indexed) set of operations preserving the relations of $\underline{A}$ is not only a minion, but also closed under composition, and must contain the projections. Such sets of operations are called *function clones*, also called *operation clones*, on $A$.

Since a large part of the algebraic theory of the CSP on finite structures is devoted to understanding the homomorphism order of all clones of the form $\mathrm{Pol}(\underline{A})$ where $\underline{A}$ is a finite relational structure over finite signature, a similar question as in the last section arises: Is every sub-clone of $\mathcal{O}(A)$ of the form $\mathrm{Pol}(\underline{A})$ for some relational structure on $A$?

Unsurprisingly, the tools from the previous section can be sharpened to answer this question positively, namely in Proposition 3.2.3.

**Definition 3.2.1**: A relation pair of the form $(R, R)$ is called *diagonal*. We introduce the shorthand $f \triangleright R :\Leftrightarrow f \triangleright (R, R)$. Furthermore, $\mathrm{Pol}\big(\{R_i\}_i\big) := \mathrm{Pol}\big(\{(R_i, R_i)\}_i\big)$.

**Definition 3.2.2** (Function clones): Let $A$ be a set. A subminion $\mathcal{M} \le \mathcal{O}(A)$ is called a *sub-clone of* $\mathcal{O}(A)$ if it contains the identity $\mathrm{id}_A$ and it is closed under composition, i.e., for all $f \in \mathcal{M}_n$ and $g_0, \dots, g_{n-1} \in \mathcal{M}_l$, the composite

$$f \bigcirc (g_0, \dots, g_{n-1}) \colon A^l \to A, \qquad \underline{x} \mapsto f(g_0(\underline{x}), \dots, g_{n-1}(\underline{x}))$$

is contained in $\mathcal{M}_l$.

Sub-clones of some $\mathcal{O}(A)$ are called *function clones*.

**Proposition 3.2.1**: If $R$ is a relation on a set $A$, the minion $\mathrm{Pol}(R) = \mathrm{Pol}((R, R))$ is always a sub-clone of $\mathcal{O}(A)$. Furthermore, if $\underline{A}$ is a relational structure with domain $A$, then $\mathrm{Pol}(\underline{A})$ is a clone as well.

*Proof*: Let $f \in \mathrm{Pol}(R)_n$ and $g_1, \dots, g_n \in \mathrm{Pol}(R)_k$. If $\underline{r_1}, \dots, \underline{r_k} \in R$, then

$$(f \bigcirc (g_1, \dots, g_n))\big(\underline{r_1}, \dots, \underline{r_k}\big) = f\left(\underbrace{g_0\big(\underline{r_1}, \dots, \underline{r_k}\big)}_{\in R}, \dots, \underbrace{g_{n-1}\big(\underline{r_1}, \dots, \underline{r_k}\big)}_{\in R}\right)$$

must lie in $R$, because both all the $g_i$ and $f$ preserve $R$. This extends to multiple diagonal relation pairs because $\mathrm{Pol}\big(\{R_i\}_i\big) = \bigcap_i \mathrm{Pol}(R_i)$, and being a sub-clone is closed under intersection. □

We shall finish the chapter by reproving (in Proposition 3.2.3) the well-known theorem that every function clone is of the form $\mathrm{Pol}\big(\{R_i\}_i\big)$ for some set of relations $\{R_i\}_i$. Let us illustrate the strategy with the following simple example: Let $\mathcal{C}$ be the (abstract) clone generated by a single unary function $u\colon A \to A$, i.e., the clone which contains all powers $u^n$ and minors thereof. Since $u^n$ is unary, all minors must factor through some surjection, so all functions are of the form $f(x_0, \ldots, x_{n-1}) = u^n(x_i)$ for some $n$ and $i$ depending on $f$.

By Corollary 3.1.4, we know that $\mathcal{C} = \mathrm{Pol}(\Sigma)$, where $\Sigma$ is an indexed set of relation pairs $(R, S)$. Our goal is to replace this relation pair by a "diagonal" relation pair $(R, R)$, which is preserved by the same polymorphisms.

To do so, let us assume that $\mathcal{C} \triangleright (R, S)$, and let us see how much this assumption can be strengthened to a diagonal relation pair. For instance, tracking an $\underline{r} \in R$, we observe[33]:

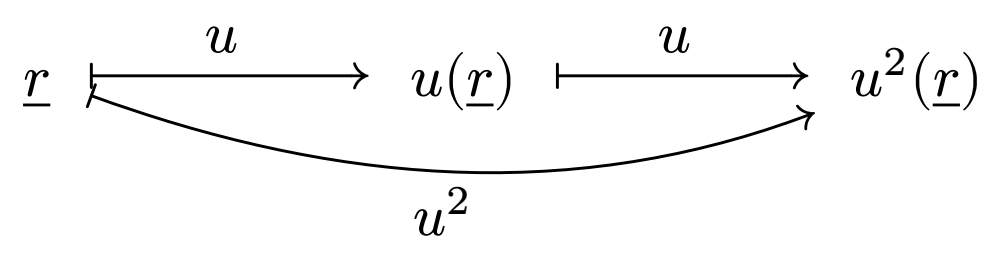

Since $u \triangleright (R, S)$, we can deduce $u(\underline{r}) \in S$, and because $u^2 \triangleright (R, S)$, we can deduce that $u^2(\underline{r}) = u(u(\underline{r})) \in S$. This suggests that $u$ preserves the subset of $S$ given by the image of $u$.

**Lemma 3.2.2**: Let $\mathcal{C}$ be a sub-clone of $\mathcal{O}(A)$, and let $R, S$ be relations on $A$ of the same arity. Then in the following set of predicates, $1 \Rightarrow 2$ and $2 \Leftrightarrow 3$:

1. $\mathcal{C} \triangleright (R, S)$
2. $\mathcal{C} \triangleright (R, \mathcal{C}[R])$
3. $\mathcal{C} \triangleright \mathcal{C}[R]$

where $\mathcal{C}[R] := \{f[R^n] \mid n > 0, f \in \mathcal{C}_n\}$ .

*Proof*: For $1 \Rightarrow 2$, assume $\mathcal{C} \triangleright (R, S)$. Since $\mathcal{C}$ contains the identity, we know that $R \subseteq S$.

*Claim:* $\mathcal{C}[R] = \bigcap\{S \supseteq R \mid \mathcal{C} \triangleright (R, S)\} =: \bigcap \Sigma$. To see $\supseteq$, note that $\mathcal{C}[R]$ contains $R$ because $\mathcal{C}$ contains the identity. Furthermore, if $f \in \mathcal{C}_n$ and $\underline{r_0}, \ldots, \underline{r_{n-1}} \in R$, then $f\big(\underline{r_0}, \ldots, \underline{r_{n-1}}\big) \in \mathcal{C}[R]$ by definition. Therefore, $\mathcal{C}[R]$ is in $\Sigma$, so in particular it contains $\bigcap \Sigma$. To see $\subseteq$, let $S$ be any other set which contains $R$ and for which $\mathcal{C} \triangleright (R, S)$.

Now, if $\mathcal{C} \triangleright (R, S)$, we conclude that $S \in \Sigma$, so by Proposition 3.1.1 we know that $\mathcal{C} \triangleright (R, \bigcap \Sigma) = (R, \mathcal{C}[R])$. This proves the first implication.

For $2 \Rightarrow 3$, note that since $\mathcal{C} \triangleright (\mathcal{C}[R], A^{\text{ar R}})$, we can invoke $1 \Rightarrow 2$ to conclude $\mathcal{C} \triangleright (\mathcal{C}[R], \mathcal{C}[\mathcal{C}[R]])$. To bridge the gap, we claim that $\mathcal{C}[\mathcal{C}[R]] = \mathcal{C}[R]$: Indeed, the left hand side consists of elements of the form $\underline{r} = f\big(\underline{r_0}, \ldots, \underline{r_{F-1}}\big)$ with $\text{ar } f =: F$, where to each $r_i$ there exist a $g_i$ and $\underline{r_{i,j}}$ for $j \in \text{ar } g_i =: G_i$ such that $r_i = g_i\big(\underline{r_{i,0}}, \ldots, \underline{r_{i,G_i-1}}\big)$. With that, we can identify our initial tuple as

$$\underline{r} = (f \bigcirc (g_0\iota_0, \ldots, g_{F-1}\iota_{F-1}))\Big(\underline{r_{0,0}}, \ldots, \underline{r_{0,G_0}}, \ldots, \underline{r_{F-1,0}} \ldots, \underline{r_{F-1,G_{F-1}}}\Big),$$

where $\iota_k\colon G_k \to \Sigma_i G_i$ is the inclusion of $G_k = \{0, \ldots, k-1\}$ into the $k$th summand of the codomain. Since $\mathcal{C}$ is closed under composition, this shows that $\underline{r} \in \mathcal{C}[R]$, proving the desired equality.

[33] Recall the usual abuse of notation $u(\underline{r}) = u\big((r_i)_{i\in \text{ar } R}\big) := (u(r_i))_{i\in \text{ar } R}$. In other words, we evaluate at a sequence of column vectors by applying $u$ row-wise.

The converse 3 $\Rightarrow$ 2 is just a relaxation (Proposition 3.1.1) because $R \subseteq \mathcal{C}[R]$. □

**Proposition 3.2.3**: Let $\mathcal{C}$ be a sub-clone of $\mathcal{O}(A)$. Then there exist relations $R_i$ such that $\mathcal{C} = \mathrm{Pol}(R_i)$.

*Proof*: Let $\mathcal{C} \leq \mathcal{O}(A)$ be a sub-clone. By Corollary 3.1.6, we know that $\mathcal{C} = \mathrm{Pol}(\mathrm{Inv}(\mathcal{C}))$, so let us define $\Sigma := \mathrm{Inv}(\mathcal{C})$. Our goal is to show that $\mathcal{C} = \mathrm{Pol}(\Sigma')$, where $\Sigma' := \{(\mathcal{C}[R], \mathcal{C}[R]) \mid (R, S) \in \Sigma\}$. The previous Lemma 3.2.2 shows that $\Sigma' \subseteq \Sigma$, so certainly $\mathrm{Pol}(\Sigma) \subseteq \mathrm{Pol}(\Sigma')$. The converse inclusion remains to be shown, so let $f \in \mathrm{Pol}(\Sigma')$. That means that for all $(R, S) \in \Sigma$, we have $f \triangleright \mathcal{C}[R]$. Since $\mathcal{C}$ is a clone, $R \subseteq \mathcal{C}[R]$, so this implies that $f \triangleright (R, \mathcal{C}[R])$. However, due to the universal property of $\mathcal{C}[R]$ established in the proof of Lemma 3.2.2, we know that $S \supseteq \mathcal{C}[R]$, so we can relax the right hand side (Proposition 3.1.1) to obtain $f \triangleright (R, S)$. □

We can now augment the earlier discussion about finitely related subminions by proving that for function clones, the two possible notions of "finite relatedness" coincide.

**Lemma 3.2.4**: Let $A$ be a finite set and $R, S$ be relations on $A$ of the same arity. Then $\mathrm{Pol}(R, S) = \mathrm{Pol}(\mathcal{C}[R])$, where $\mathcal{C} = \mathrm{Pol}(R, S)$.

*Proof*: The forward inclusion has been shown in Lemma 3.2.2. For the backward inclusion, let $f \triangleright \mathcal{C}[R]$, or equivalently, $f \triangleright (R, \mathcal{C}[R])$. Since $\mathcal{C}[R] = \bigcap\{S' \supseteq R \mid \mathcal{C} \triangleright (R, S')\}$ and $\mathcal{C} \triangleright (R, S)$, we have that $\mathcal{C}[R] \subseteq S$, so we can relax the preservation to $f \triangleright (R, S)$, as desired. □

**Proposition 3.2.5**: Let $A$ be a finite set and $\mathcal{C} \leq \mathcal{O}(A)$ be a function clone. Then the following are equivalent:

1. There exists a finite set $\Sigma$ of relation pairs on $(A, A)$ such that $\mathcal{C} = \mathrm{Pol}(\Sigma)$.
2. There exists a finite set $\Sigma'$ of relations on $A$ such that $\mathcal{C} = \mathrm{Pol}(\Sigma')$.

*Proof*: The backwards implication is trivial. For the forwards implication, let $\mathcal{C}$ be of the form $\mathrm{Pol}(\Sigma)$ for some finite set $\Sigma$ of relation pairs. Taking a product (Corollary 3.1.2) ensures that $\mathcal{C}$ is even of the form $\mathrm{Pol}(R, S)$ for a single relation pair $(R, S)$ on $(A, A)$. But by the previous lemma, this is the same as $\mathrm{Pol}(\mathcal{C}[R])$. □

In particular, this shows that Proposition 3.1.10 not only demonstrated that $\mathcal{B}_\infty \leq \mathcal{O}(2)$ is finitely related as a function minion – i.e., not of the form $\mathrm{Pol}(R, S)$ – but also not finitely related as a function clone – i.e., not of the form $\mathrm{Pol}(R)$.

There is extensive literature around the question when a function clone is finitely related. Often, one talks instead about *finitely related algebras*, or *algebras of finite type*. For an entrypoint to the discussion, the reader is referred to [37].

## 3.3 Representability of minions as $\mathrm{Pol}(A, B)$

In this section we are concerned with the question: Is every abstract minion $\mathcal{M}$ isomorphic to a concrete minion, i.e., a minion of the form $\mathrm{Pol}(\underline{A}, \underline{B})$?

For the first answer to this question, we will introduce a notion of "height-one terms": given a minion $\mathcal{M}$ and a variable set $X$: if $x_0, \ldots, x_{n-1}$ are elements of $X$ and $f$ is an element of $\mathcal{M}_n$, we would like something like $f(x_0, \ldots, x_{n-1})$ to be a "term" associated to $f$ with the respective variables. For a set of such terms to function as terms of minions, we want to have the compatibility conditions $f\alpha(x_0, \ldots, x_{n-1}) = f\left(x_{\alpha(0)}, \ldots, x_{\alpha(k-1)}\right)$ for a given $f \in \mathcal{M}_n$ and $\alpha\colon k \to n$. This is the motivation of Definition 3.3.1. We can then represent an element $f$ of a minion as the function sending an $n$-tuple of elements in $X$ to the height-one term $f(x_0, \ldots, x_{n-1})$.

**Definition 3.3.1** (height-one terms, tensor product): Let $\mathcal{M}$ be a minion and $X$ be a set. We define

$$\mathcal{M} \otimes X := \left(\coprod_{n \in \mathbb{F}_{>0}} \mathcal{M}_n \times X^n\right) / \sim,$$

where $\sim$ is the smallest equivalence relation containing

$$(f\alpha, \underline{x}) \sim (f, \alpha; \underline{x})$$

for all $f \in \mathcal{M}_n$, $\alpha\colon n \to k$, and $\underline{x} \in X^k$. We denote the equivalence class associated to a 2-tuple $(f, \underline{x})$ by $[f, \underline{x}]$.[34]

Since the notation $\alpha; \underline{x}$ will strike many readers as a little awkward, a justification is in order. While $\underline{x}$ can be viewed as a tuple $X^k$, we interpret it as a map $k \to X$ between sets, and so an argument operation $\alpha\colon n \to k$ can be *precomposed* with $\underline{x}$ to give a new tuple $\alpha; \underline{x}$. In coordinates, the resulting tuple is $\left(x_{\alpha(0)}, \ldots, x_{\alpha(k-1)}\right)$.

This "minor operation-then-tuple" notation has the advantage that it closely mirrors the usual context in which $\otimes$ is used: Namely, the tensor product $M \otimes_R S$ combines a right $R$-module $M$ and a left $R$-module $S$ by demanding – among other things – that $(mr) \otimes n = m \otimes (rn)$.

What seems to be just a visual suggestion however allows us to completely trivialize the proof of the following proposition.

**Proposition 3.3.1**: Given a map $m\colon X \to Y$, the induced map

$$\begin{aligned} \mathrm{id}_{\mathcal{M}} \otimes m\colon \mathcal{M} \otimes X &\to \mathcal{M} \otimes Y, \\ [f, \underline{x}] &\mapsto [f, \underline{x}; m], \end{aligned}$$

where $\underline{x}; m = (m(x_0), \ldots, m(x_{k-1}))$ is the post-composition of $x$ with $m$, is well-defined. Furthermore, this procedure turns $M \otimes _$ into a functor $\underline{\mathrm{Set}} \to \underline{\mathrm{Set}}$.

*Proof*: For well-definedness, note that $[f, \alpha; \underline{x}]$ gets mapped to $[f, (\alpha; \underline{x}); m]$ whereas $[f\alpha, \underline{x}]$ gets mapped to $[f\alpha, \underline{x}; m] = [f, \alpha; (\underline{x}; m)]$ which is the same. Since these pairs form a generator of the equivalence relation we quotiented out, well-definedness follows.

[34] This is in analogy to the coordinates in homogeneous space, e.g. $(x, y) \in \mathbb{RP}^1$. However, we use square brackets to emphasize that we are talking about an equivalence class.

For functoriality, just note that $[f, \underline{x}; \mathrm{id}_X] = [f, \underline{x}]$, so the identity on $X$ induces the identity on $\mathcal{M} \otimes X$, and since $[f, (\underline{x}; m); m'] = [f, \underline{x}; (m; m')]$, we are compatible with composition. □

**Proposition 3.3.2**: Let $\mathcal{M}$ be a minion and $X$ be a set. The map[35]

$$\begin{aligned} \eta_X \colon \mathcal{M} &\to \mathcal{O}(X, \mathcal{M} \otimes X), \\ f &\mapsto [f, _], \end{aligned}$$

where $[f, _](\underline{x}) := [f, \underline{x}]$, is a minion homomorphism.

*Proof*: Let $\alpha \colon n \to k$ and $f \in \mathcal{M}_n$. Recall that in $\mathcal{O}(X, \mathcal{M} \otimes X)$, the minor operation $\alpha$ takes an $n$-ary function $p \colon X^n \to \mathcal{M} \otimes X$, and sends it to the associated function which maps $\underline{x} \in X^k$ to first the "pulled-back" tuple $\alpha; \underline{x} = \left(x_{\alpha(0)}, ..., x_{\alpha(n-1)}\right)$, and then applies the polymorphism $p$ – i.e., $p\alpha = p(\alpha; _)$.

With that in mind, we can verify that

$$\eta_X(f)\alpha = \eta_X(f)(\alpha; _) = [f, \alpha; _] = [f\alpha, _] = \eta_X(f\alpha),$$

and we are done. □

*Remark*: It is not hard to verify that we have a one-to-one correspondence between minion homomorphisms $\mathcal{M} \to \mathrm{Pol}(X, Y)$ and maps $\mathcal{M} \otimes X \to Y$: Indeed, given such a homomorphism $h$ and a term $[f, \underline{x}]$, we can evaluate $h(f)$ at the $\mathrm{ar}\ f$-tuple $\underline{x}$. On the other hand, given such an evaluation map $\varepsilon \colon \mathcal{M} \otimes X \to Y$, we can consider the "curried" function $\varepsilon([f, _])$, which induces a polymorphism $X^{\mathrm{ar}\ f} \to Y$.

In fact, if instead of sets $X$ and $Y$ we consider relational structures $\underline{A}$ and $\underline{B}$ over a signature $\tau$, there is a way to endow $\mathcal{M} \otimes A$ with the structure of a $\tau$-algebra in the following way: if $R \in \tau$ is a $n$-ary relation symbol and $R^{\underline{A}} \subseteq A^n$ its interpretation in $\underline{A}$, then each of the coordinate projections $\pi_i \colon R^{\underline{A}} \to A$ induces a map $M \otimes R^{\underline{A}} \to M \otimes A$, puzzling together to a map $M \otimes R^{\underline{A}} \to (\mathcal{M} \otimes A)^n$. The image of this map defines an interpretation of $R$ in $\mathcal{M} \otimes A$. More concretely, we send every height-one term $\left[f, e_0, ..., e_{\mathrm{ar}\ f}\right]$, where all $e_i$ are tuples in $R^{\underline{A}}$, to the tuple of height-one-terms $\left[f, e_{0,i}, ..., e_{\mathrm{ar}\ f,i}\right]_i$. In the case where $\underline{A}$ has finite domain, we can apply the bijection $\mathcal{M} \otimes |A| \cong \mathcal{M}_{|A|}$ and arrive at what is known as the "free structure" associated to the structure $\underline{A}$ and the minion $\mathcal{M}$ (see [19, section 4.2]).

*Remark (for category theorists)*: The preceding remarks can be put in context by noting that $\mathcal{M} \otimes X$ is nothing but the coend $\int^n \mathcal{M}_n \times X^n$: Indeed, the usual manipulations give us

$$\begin{aligned} \mathcal{M} \otimes X \to_{\underline{\mathsf{Set}}} Y = \left(\int^{n \in \mathbb{F}_{>0}} \mathcal{M}_n \times X^n\right) \to_{\underline{\mathsf{Set}}} Y \quad &= \int_{n \in \mathbb{F}_{>0}} (\mathcal{M}_n \times X^n \to_{\underline{\mathsf{Set}}} Y) \\ = \int_{n \in \mathbb{F}_{>0}} (\mathcal{M}_n \to_{\underline{\mathsf{Set}}} (X^n \to_{\underline{\mathsf{Set}}} Y)) &= \mathcal{M}_\bullet \to_{[\mathbb{F}_{>0},\ \underline{\mathsf{Set}}]} (X^\bullet \to_{\underline{\mathsf{Set}}} Y) \\ &= \mathcal{M} \to_{[\mathbb{F}_{>0},\ \underline{\mathsf{Set}}]} \mathrm{Pol}_{\underline{\mathsf{Set}}}(X, Y). \end{aligned}$$

[35]More explicitly, the sequence of maps $\mathcal{M}_n \to \mathcal{O}(X, \mathcal{M} \otimes X)_n$ for each $n$

As every step is natural in $\mathcal{M}$ and $Y$, this forms an adjunction $_ \otimes X \dashv \mathrm{Pol}(X, _)$. Clearly, these manipulations do not only work for $X, Y \in \underline{\mathsf{Set}}$, but also for objects in functor categories $[\mathcal{C}, \underline{\mathsf{Set}}]$, as these categories remain to be (co-)complete and also posess an internal hom[36]. This encompasses two special cases:

First, relational structures. To sketch the analogy, note that a graph $G$ may be described by sets $G_0$ of vertices and $G_1$ of edges, together with source- and target-maps $e_0, e_1\colon G_1 \rightrightarrows G_0$. Indeed, this allows us to view graphs as functors $\mathcal{C} \to \underline{\mathsf{Set}}$ where $\mathcal{C}$ is the small category given by the diagram $e_0, e_1\colon 1 \rightrightarrows 0$. While on first sight, this category describes more general structures – namely *multigraphs*[37], i.e., graphs with potentially more than one edge between the same source and target – the act of collapsing multiple edges gives us a reflector[38] $R\colon [\mathcal{C}, \underline{\mathsf{Set}}] \to \underline{\mathsf{Graph}}$. In fact, the canonical map $G \to RG$ for a multigraph $G$ is always post-invertible by choosing an arbitrary edge in the preimage. In particular, this means that $G$ and $RG$ are always homomorphically equivalent, and by extension also the minions $\mathrm{Pol}(G, G')$ and $\mathrm{Pol}(RG, RG')$. Thus, for the purposes of studying promise CSPs and polymorphism minions, the distinction between graphs and multigraphs does not matter. The same considerations apply to relational structures over other signatures $\tau$, giving us an adjunction $_ \otimes \underline{A} \dashv \mathrm{Pol}_\tau(\underline{A}, _)$.

Second, minions themselves, i.e., the case $\mathcal{C} = \mathbb{F}_{>0}$: There we we can summarize the adjunction as $_ \otimes \mathcal{N} \dashv \mathrm{Pol}_{\underline{\mathsf{Min}}}(\mathcal{N}, _)$. Taking the projection minion $\mathcal{P}$ as a designated element now endows $\underline{\mathsf{Min}}$ with a monoidal structure. In fact, our adjunction shows that this structure is *right closed* with respect to Pol (see [38, sec. 4.5]). But what are the monoids $(\mathcal{M}, \mu, \eta)$ over this monoidal structure? As per our adjunction, the map $\mu\colon \mathcal{M} \otimes \mathcal{M} \to \mathcal{M}$ corresponds to a map $\mathcal{M} \to \mathrm{Pol}_{\underline{\mathsf{Min}}}(\mathcal{M}, \mathcal{M})$: This means that to each $n$-ary $f \in \mathcal{M}_n$, we associate an $n$-ary polymorphism $\mathcal{M}^n \to \mathcal{M}$. This can be interpreted as the *composition* of $f$ with $n$ elements $g_i$ of potentially different, but fixed, arity $k$. If we interpret the unit $\eta\colon \mathcal{P} \to \mathcal{M}$, which designates a special unary element, as choosing an *identity*, we can observe that we have precisely described an *abstract clone*[39].

This product of minions is well-known in category theory, especially in the context of clubs and operads. The interested reader is referred to [39], where the relevant product is introduced as $T \circ S$ in Equation 3.2.

**Lemma 3.3.3**:
1. $\mathcal{M}_n \cong \mathcal{M} \otimes n$
2. The inclusion $\iota\colon n \hookrightarrow \omega$ gives rise to an injection $\mathrm{id}_{\mathcal{M}} \otimes \iota\colon \mathcal{M} \otimes n \hookrightarrow \mathcal{M} \otimes \omega$.

*Proof*: The idea for the first part is rewriting $[f, \alpha]$ to $[f\alpha, \mathrm{id}_n]$, i.e., establishing a normal form. More formally, assiging $[f, \alpha] \mapsto f\alpha$ gives a well-defined map $\mathcal{M} \otimes n \to \mathcal{M}_n$, whose counterpart $\mathcal{M}_n \ni g \mapsto [g, \mathrm{id}_n] \in \mathcal{M} \otimes n$ is easily checked to be an inverse.

For the second part, note that the map $m := \max(n, _)\colon \omega \twoheadrightarrow n$ is a post-inverse to the inclusion $\iota\colon n \hookrightarrow \omega$, i.e., we have $\iota; \omega = \mathrm{id}_n$. By functoriality, the induced maps

[36]The existence of the internal hom is less trivial to see (but also well-known); see Chapter 7 where we explicitly give the construction for the case of minions, i.e., $\mathcal{C} = \mathbb{F}_{>0}$.

[37]In representation theory also known as *quivers*

[38]i.e., left adjoint to the inclusion functor

[39]For a self-contained definition of abstract clones, see [10, Section 8.6]. See also [40, Prop. 3.4] for a reference to this claim.

$$\mathcal{M} \otimes n \xrightarrow{\mathrm{id}_{\mathcal{M}} \otimes \iota} \mathcal{M} \otimes \omega \xrightarrow{\mathrm{id}_{\mathcal{M}} \otimes m} \mathcal{M} \otimes n$$

must also compose to the identity; in particular, $\mathrm{id}_{\mathcal{M}} \otimes \iota$ must be injective. □

**Proposition 3.3.4**: The minion homomorphism $\eta_\omega$: $\mathcal{M}$ to cal(O)(ω, M tens ω) is injective. Hence, every abstract minion is isomorphic to a function minion.

*Proof*: Let $f, g \in \mathcal{M}_n$ satisfy $[f, _] = [g, _]$ as maps $\omega^n \to \mathcal{M} \otimes \omega$. In particular, evaluating at the inclusion $\iota = (0, ..., n-1)\colon n \hookrightarrow \omega$ we get $[f, \iota] = [g, \iota]$.

Observe that the composition

$$\begin{aligned} \mathcal{M}_n &\hookrightarrow \mathcal{M} \otimes n \hookrightarrow \mathcal{M} \otimes \omega, \\ f &\mapsto [f, \mathrm{id}_n] \\ &\quad [f, \alpha] \quad \mapsto [f, \alpha; \iota] \end{aligned}$$

of the two maps in the previous Lemma 3.3.3 is injective and sends $f$ resp. $g$ to $[f, \iota]$ resp. $(g\colon \iota)$. By the same lemma, these two maps must be injective, so $[f, \iota] = [g, \iota]$ implies $f = g$, as desired. □

In the case where every element of the minion $\mathcal{M}$ has essential arity at most $N$, we can show that the homomorphism $\eta_N\colon \mathcal{M} \to \mathcal{O}(N, \mathcal{M} \otimes N)$ is injective. To do so, we first need a few observations about essential coordinates and injective maps.

**Lemma 3.3.5**: Let $h\colon \mathcal{M} \to \mathcal{N}$ be a minion homomorphism and $n > 0$. If $h_n$ is injective, then $h_{n'}$ is injective for all $n' \leq n$.

*Proof*: Let $f, g \in \mathcal{M}_{n'}$ for some $n' \leq n$ with $h(f) = h(g)$. If $\iota\colon n' \hookrightarrow n$ is the inclusion, then we have $h(f\iota) = h(g\iota)$, and thus $f\iota = g\iota$ by injectivity in $n$. By operating with some retraction $\rho\colon n \twoheadrightarrow n'$ we conclude $f = f\iota\rho = g\iota\rho = g$. □

**Proposition 3.3.6**: Let $\mathcal{M}$ be a minion with essential arity at most $N$ and $h\colon \mathcal{M} \to \mathcal{N}$ a minion homomorphism. If $h_n$ is injective for all $n \leq \max\{N, 2\}$, then $h$ is injective.

*Proof*: *Claim*: For every $f$ in $\mathcal{M}_n$, we have $\operatorname{Ess} f = \operatorname{Ess} h(f)$. If $f$ is a constant, then $h(f)$ is as well, so assume it is not. If $n = 1$, then injectivity of $h$ in arity 2 allows us to use Lemma 2.6.9, and if $2 \leq n \leq N$, we can apply the same lemma by injectivity of $h_n$. For $n > N$, we know that $f \in \mathcal{M}_n$ can have at most $N$ essential coordinates, so it must factorize as $\tilde{f}\eta$ with $\operatorname{Iness} \tilde{f} = \emptyset$ and $\eta$ an injection. This implies that $\operatorname{ar} \tilde{f} = |\operatorname{Ess} \tilde{f}| \leq N$, whereby

$$\begin{aligned}
\operatorname{Ess} f &= \operatorname{Ess} \tilde{f}\eta && \text{Definition} \\
&= \eta\left[\operatorname{Ess} \tilde{f}\right] && \text{Lemma 2.6.3} \\
&= \eta\left[\operatorname{Ess} h\left(\tilde{f}\right)\right] && \operatorname{ar} \tilde{f} \leq N \\
&= \operatorname{Ess} h\left(\tilde{f}\right)\eta && \text{Lemma 2.6.3} \\
&= \operatorname{Ess} h\left(\tilde{f}\eta\right) = \operatorname{Ess} h(f), && \text{Homomorphicity, Definition}
\end{aligned}$$

which proves the claim.

Let now $h(f) = h(g)$ for some $f, g \in \mathcal{M}_n$. The preceding observation forces $\operatorname{Ess} f = \operatorname{Ess} h(f) = \operatorname{Ess} h(g) = \operatorname{Ess} g$, which we will denote by $E$. If $E = \emptyset$, we will instead set $E = \{0\}$, and note that this choice also satisfies the property that $n \setminus E =: I$ consists of solely inessential coordinates.

Now, we claim that both $f$ and $g$ factor through an enumeration $\eta\colon |E| \hookrightarrow n$ of $E$: Indeed, let $\rho\colon n \twoheadrightarrow |E|$ send $j \mapsto \eta^{-1}(j)$ if $j \in \operatorname{Im}\eta$ and $j \mapsto 0$ otherwise, $\rho;\eta\colon n \to n$ is the identity on $\operatorname{Im}\eta = E$ and sends every other element to $0$. Clearly, $\eta;\rho = \mathrm{id}_{|E|}$. Since $n \geq 2$, and because $n \setminus E \subseteq \operatorname{Iness} f = \operatorname{Iness} g$, we can conclude that $f = f\rho\eta$ and $g = g\rho\eta$. We now know that $h(f\rho) = h(f)\rho = h(g)\rho = h(g\rho)$, which by injectivity in arity $|E|$ implies that $f\rho = g\rho$, hence $f = f\rho\eta = g\rho\eta = g$. □

**Proposition 3.3.7**: If $\mathcal{M}$ has essential arity $\leq N$ and $N \geq 2$, then

$$\eta_N\colon \mathcal{M} \to \mathcal{O}(N, \mathcal{M} \otimes N)$$

is injective.

*Proof*: We will show that $\eta_N$ is injective in arity $N$, which by the preceding considerations is sufficient.

Indeed, let Let $f, g \in \mathcal{M}_N$. If $[f, _] = [g, _]$, then certainly $[f, \mathrm{id}_N] = [g, \mathrm{id}_N]$. Since the evaluation at the identity is a bijection from $\mathcal{M} \otimes N \to \mathcal{M}_N$ (Lemma 3.3.3), we know that this can only be the case if $f = g$. □

**Corollary 3.3.8**: Every finitely generated minion is finitely representable.

*Proof*: Recall that being finitely generated is equivalent to being locally finite and of bounded essential arity. The preceding Proposition 3.3.7 guarantees existence of an injection $\eta_n\colon \mathcal{M} \hookrightarrow \mathcal{O}(n, \mathcal{M}_n)$, where $\mathcal{M}_n$ is finite due to local finiteness of $\mathcal{M}$. Hence, this injection witnesses finite representability. □

A nice application of this insight is the following:

**Corollary 3.3.9**: To each primitive positive sentence $\varphi$ of minions, there are finite relational structures $\underline{A}$ and $\underline{B}$ such that $\mathcal{M}$ satisfies $\varphi$ if and only if $\operatorname{Pol}(\underline{A}, \underline{B}) \to \mathcal{M}$.

*Proof*: The gadget minion $\mathcal{M}_\varphi$ associated to $\varphi$ is finitely generated, hence finitely representable, hence of the form $\mathrm{Pol}(\underline{A}, \underline{B})$ for some finite relational structures $\underline{A}$ and $\underline{B}$. $\square$

# 4 Finite representability

Of main interest in this thesis are minions which are of the form $\mathrm{Pol}(\underline{A}, \underline{B})$, where $\underline{A}$ and $\underline{B}$ are relational structures (with potentially infinite signature) over finite domain. These will be called *finitely representable*. By Corollary 3.1.4, this is equivalent to asking which minions are isomorphic to a sub-minion of $\mathcal{O}(n, k)$. In this chapter, we will present a simple characterization of this property, which – roughly speaking – exploits the fact that minion homomorphisms $\mathcal{M} \to \mathcal{O}(n, k)$ are determined by their value in arity $n$. We will demonstrate how to construct embeddings $\mathcal{O}(n, k) \hookrightarrow \mathcal{O}(n', k')$ when $n \leq n'$ and $k \leq k'$.

The material presented here lays the groundwork for analyzing the homomorphism order of finitely representable minions.

## 4.1 The correspondence $_(n) \dashv \mathcal{O}(n, _)$

A central insight for function minions on finite sets is that for any function

$$f \in \mathcal{O}(n, k)_l = n^l \to k,$$

we can interpret its inputs $\alpha \in n^l$ as minor operations $\alpha\colon l \to n$. This allows us to see that

$$f(\alpha_1, \ldots, \alpha_l) = (f\alpha)(1, \ldots, n),$$

or, written more abstractly,

$$f(\alpha) = f(\alpha; \mathrm{id}_n) = (f\alpha)(\mathrm{id}_n),$$

where the last equality holds by definition of the minor action on function minions.

**Definition 4.1.1**: Let $n, k$ be natural numbers and $\mathcal{M}$ an abstract minion. Let the maps

$$\mathcal{M} \to_{\underline{\mathsf{Min}}} \mathcal{O}(n, k) \quad \underset{\tilde{\square}}{\overset{\overline{\square}}{\rightleftarrows}} \quad \mathcal{M}_n \to_{\underline{\mathsf{Set}}} k$$

be defined as follows: for every function $\varphi\colon \mathcal{M}_n \to_{\underline{\mathsf{Set}}} k$, define

$$\tilde{\varphi}\colon \mathcal{M} \to_{\underline{\mathsf{Min}}} \mathcal{O}(n, k),$$

to be in each degree $l$ the map

$$\begin{aligned} \tilde{\varphi}_l\colon \mathcal{M}_l &\to_{\underline{\mathsf{Min}}} \mathcal{O}(n, k)_l, \\ f &\mapsto (\alpha \mapsto \varphi(f\alpha)) \end{aligned}$$

and for every minion homomorphism $h = (h_l)_l\colon \mathcal{M} \to_{\underline{\mathsf{Min}}} \mathcal{O}(n, k)$, set

$$\begin{aligned} \overline{h}\colon \mathcal{M}_n &\to k, \\ f &\mapsto \underbrace{h_n(f)}_{\in n^n \to k}(\mathrm{id}_n). \end{aligned}$$

**Lemma 4.1.1**: The map $\tilde{\square}$ from Definition 4.1.1 is well-defined, i.e., for every $\varphi\colon \mathcal{M}_n \to_{\underline{\mathsf{Set}}} k$, the induced map $\tilde{\varphi}$ is a minion homomorphism.

*Proof*: Let $\beta\colon l \to l'$ be a minor operation and $f \in M_l$. Then, for every $\alpha' \in n^{l'}$, we observe

$$\begin{aligned}\tilde{\varphi}_{l'}(f\beta)(\alpha') = \varphi((f\beta)\alpha') = \varphi(f(\beta;\alpha')) &= \tilde{\varphi}_l(f)(\beta;\alpha') \\ &= \tilde{\varphi}_l(f)\left(\alpha'_{\beta(1)}, \ldots, \alpha'_{\beta(l)}\right) \\ &= (\tilde{\varphi}_l(f)\beta)(\alpha').\end{aligned}$$

Therefore, $\tilde{\varphi}_{l'}(f\beta) = \tilde{\varphi}_l(f)\beta$, proving homomorphicity. □

**Proposition 4.1.2**: The maps fromy Definition 4.1.1 are mutually inverse.

*Proof*: Let $\varphi\colon \mathcal{M}_n \to k$. This induces $\tilde{\varphi}\colon \mathcal{M} \to \mathcal{O}(n,k)$, which induces $\overline{\tilde{\varphi}}\colon \mathcal{M}_n \to k$. Given a $f \in \mathcal{M}_n$, we unravel the definition:

$$\begin{aligned}\overline{\tilde{\varphi}}(f) &= \tilde{\varphi}_l(f)(\mathrm{id}_n) && \text{Definition of } \overline{\square} \\ &= (\alpha \mapsto \varphi(f\alpha))(\mathrm{id}_n) && \text{Definition of } \tilde{\square} \\ &= \varphi(f\,\mathrm{id}_n) && \text{Evaluation} \\ &= \varphi(f),\end{aligned}$$

completing the first direction.

For the converse, let $h = (h_l)_l\colon \mathcal{M} \to \mathcal{O}(n,k)$ be a minion hom. We need to prove that in every degree $l$, the induced minion hom $\tilde{\overline{h}}_l = h_l$. Indeed, let $f \in \mathcal{M}_l$. Then both $h_l(f)$ and $\tilde{\overline{h}}_l(f)$ are elements of $\mathcal{O}(n,k)_l = n^l \to k$, so we shall verify their equality for every $\alpha \in n^l$. We observe

$$\begin{aligned}\left(\tilde{\overline{h}}\right)_l(f)(\alpha) &= \overline{h}(f\alpha) && \text{Definition of } \tilde{\square} \\ &= h_n(f\alpha)(\mathrm{id}_n) && \text{Definition of } \overline{\square} \\ &= (h_l(f)\alpha)(\mathrm{id}_n) && \text{Homomorphicity of } h \\ &= h_l(f)(\alpha;\mathrm{id}_n) && \text{Definition of minor operations in } \mathcal{O}(n,k) \\ &= h_l(f)(\alpha),\end{aligned}$$

and conclude that indeed $\tilde{\overline{h}} = h$. □

In order to demonstrate, we will quickly verify that this correspondence is natural.

**Proposition 4.1.3**: The correspondence from Definition 4.1.1 is natural in $\mathcal{M}$ and $k$.

*Proof*: We will show that the evaluation $\varphi \mapsto \overline{\varphi}$ commutes with precomposition by a minion homomorphism $h'\colon \mathcal{M}' \to \mathcal{M}$ and with postcomposition by a set-map $r\colon k \to k'$. Indeed, we observe

$$\overline{h \circ h'}(f) = (h \circ h')_n(f)(\mathrm{id}_n) = h_n(h'_n(f))(\mathrm{id}_n) = \overline{h}(h'_n(f)) = \overline{h} \circ h'_n(f),$$

and

$$\overline{r_* \circ h}(f) = r(h_n(f))(\mathrm{id}_n) = r\big(\overline{h}(f)\big) = r \circ \overline{h}(f),$$

as desired. Naturality of the inverse $\widetilde{\square}$ follows automatically because $\overline{\square}$ is a natural isomorphism. □

**Corollary 4.1.4**: The operations defined in Definition 4.1.1 give rise to a natural one-to-one correspondence between minion homomorphisms $\mathcal{M} \to \mathcal{O}(n, k)$ and functions $\mathcal{M}_n \to k$, i.e., an adjunction between functors

$$\underline{\mathsf{Set}} \quad (_)_n \left( \dashv \right) \mathcal{O}(n, _). \quad \underline{\mathsf{Min}}$$

*Proof*: This follows by definition of an adjunction (Definition 1.1.22) from Proposition 4.1.2 and Proposition 4.1.3. □

### 4.1.1 Free structures

Lastly, let us remark that this adjunction is a special case of a more general situation: Whenever $\underline{A}, \underline{B}$ are $\tau$-structures, we have a correspondence

$$(\mathcal{M} \to_{\underline{\mathsf{Min}}} \mathrm{Pol}(\underline{A}, \underline{B})) \quad \cong \quad (F_{\mathcal{M}}(\underline{A}) \to_\tau \underline{B}),$$

where $F_{\mathcal{M}}(\underline{A})$ is known as the *free structure* of $\mathcal{M}$ generated by $\underline{A}$ (see [19, section 4.1]).

One way of describing the construction goes as follows: Assuming without loss of generality that the domain of $\underline{A}$ is $n$, We set the domain of $F_{\mathcal{M}}(\underline{A})$ to be $\mathcal{M}_n$. Now, for every $r$-ary relation $R \in \tau$, assume we enumerated the corresponding interpretation $R^{\underline{A}}$ in $\underline{A}$ as $\varepsilon\colon\ s := |R^{\underline{A}}| \hookrightarrow A^r = n^r$. Postcomposing with the $i$th projections gives functions $\varepsilon_i\colon\ s \to n$, which we can then interpret as minor operations to get an induced map $\mathcal{M}_s \to \mathcal{M}_n$. Tupling these maps together for every $i$ gives us a map $\mathcal{M}_s \to (\mathcal{M}_n)^r = F_{\mathcal{M}}(\underline{A})^r$, the image of which will be our interpretation of $R$ in our free structure. It is not obvious, but also not hard to see that taking $F_{\mathcal{M}}(\underline{A})$ to be $\mathcal{M}_n$ with these relations makes the adjunction work. For a full proof, the reader is referred to [19, Lemma 4.4].

The adjunction we elaborated on is of course the special case $\tau = \emptyset$, i.e., no structure.

## 4.2 Characterizing finite representability

**Definition 4.2.1**: Let $n, k$ be natural numbers. A minion $\mathcal{M}$ is defined to be
- … *$(n, k)$-representable* if there is an embedding $\mathcal{M} \hookrightarrow \mathcal{O}(n, k)$.
- …*$n$-representable* if it is $(n, k)$-representable for some $k$.
- …*finitely representable* if it is $(n, k)$-representable for some $n$ and $k$.

Note that since $\mathcal{O}(A, B) \cong \mathcal{O}(|A|, |B|)$, we can also make sense of $(A, B)$-representability for arbitrary finite sets $A$ and $B$.

**Lemma 4.2.1**: Let $\mathcal{M}$ be a minion. If $\mathcal{M}$ is $(n, k)$-representable, it is also $(n, k')$ representable for any $k' \geq k$.

*Proof*: Any embedding $h\colon \mathcal{M} \hookrightarrow \mathcal{O}(n, k)$ can be postcomposed with the canonical inclusion $\mathcal{O}(n, k) \subseteq \mathcal{O}(n, k')$. □

**Proposition 4.2.2**: Let $\mathcal{M}$ be a minion such that $|\mathcal{M}_n|$ is finite. The following are equivalent for a map $\varphi\colon \mathcal{M}_n \to k$:
1. $\varphi$ induces an injection $\mathcal{M} \to \mathcal{O}(n, k)$ under the correspondence of Definition 4.1.1.
2. For all $l \geq 1$, if $f \neq g \in \mathcal{M}_l$, there must be a minor operation $\alpha\colon l \to n$ such that $\varphi(f\alpha) \neq \varphi(g\alpha)$.

In particular, the following are equivalent:
1. $\mathcal{M}$ is $(n, k)$-representable.
2. There is a map $\varphi\colon \mathcal{M}_n \to k$ such that for all $l \geq 1$, if $f \neq g \in \mathcal{M}_l$, there must be a minor operation $\alpha\colon l \to n$ such that $\varphi(f\alpha) \neq \varphi(g\alpha)$.

*Proof*: The induced minion homomorphism $\tilde{\varphi}$ is given in arity $l$ by $f \mapsto \varphi(f_)\colon n^l \to k$, so it is injective in arity $l$ if and only if

$$\forall f, g \in \mathcal{M}_l\colon\ f \neq g \Rightarrow \varphi(f_) \neq \varphi(g_).$$

Since two functions are distinct if and only they differ at a point of the domain, this stamenent is equivalent to

$$\forall f, g \in \mathcal{M}_l\colon\ f \neq g \Rightarrow \exists \alpha \in n^l\colon\ \varphi(f\alpha) \neq \varphi(g\alpha).$$

For the second part, an injection $h\colon \mathcal{M} \to \mathcal{O}(n, l)$ gives rise to a map $\overline{h}$ which satisfies the criterion of the first part, and conversely, if such a $\varphi$ exists, the induced homomorphism $\tilde{\varphi}$ is an injection as per the preceding discussion. □

**Definition 4.2.2**: If $\varphi\colon \mathcal{M}_n \to k$ induces an embedding as per the charatccterization above, we say that $\varphi$ *witnesses $(n, k)$-representability* of $\mathcal{M}$.

**Proposition 4.2.3**: Let $\mathcal{M}$ be a minion such that $|\mathcal{M}_n|$ is finite. The following are equivalent:
1. $\mathcal{M}$ is $n$-representable.
2. $\mathcal{M}$ is $(n, |\mathcal{M}_n|)$-representable.
3. The identity $\mathrm{id}_{\mathcal{M}_n}\colon \mathcal{M}_n \to \mathcal{M}_n$ induces an embedding $\mathcal{M} \hookrightarrow \mathcal{O}(n, \mathcal{M}_n)$.
4. For all $k \geq 1$, if $f \neq g \in \mathcal{M}_k$, then there must be $\alpha\colon k \to n$ such that $f\alpha \neq g\alpha$ (*"elements can be distinguished in arity n"*).

*Proof*: The directions $3 \Rightarrow 2 \Rightarrow 1$ are trivial. For $1 \Rightarrow 3$, let $\mathcal{M}$ be $n, k$-representable. We claim that if $\varphi\colon \mathcal{M}_n \to k$ an embedding – i.e., for all $l$ and $f \neq g \in \mathcal{M}_l$, there must be an $\alpha\colon l \to n$ such that $\varphi(f\alpha) \neq \varphi(g\alpha)$ – then $\mathrm{id}_n\colon \mathcal{M}_n \to n$ induces an embedding as well: Indeed, if $f \neq g \in \mathcal{M}_l$, then $\varphi(f\alpha) \neq \varphi(g\alpha)$ most certainly implies that $f\alpha \neq g\alpha$. Statement 4 is just a rewording for the criterion that $\mathrm{id}_n$ induces an embedding, hence equivalent to 3. □

**Lemma 4.2.4**: If $\mathcal{M}$ is $n$-representable, it is $n'$-representable for any $n' \geq n$.

*Proof*: Let $\mathcal{M}$ be $n$-representable and $f \neq g \in \mathcal{M}_l$. By Proposition 4.2.3, we need to find a minor operation $\alpha'\colon l \to n'$ such that $f\alpha' \neq g\alpha'$. Indeed, invoking the same characterization for $n$-representability ensures the existence of a minor op $\alpha\colon l \to n$ such that $f\alpha \neq g\alpha$ in arity $n$. Since the canonical injection $\iota\colon n \hookrightarrow n'$ is invertible, $f\alpha\iota \neq g\alpha\iota$. Thus $\alpha;\iota\colon l \to n'$ witnesses $f \neq g$, and we are done. □

**Proposition 4.2.5**: If $\mathcal{M}$ is $(n, k)$-representable, then it is $(n+1, k)$-representable.

*Proof*: Let $\varphi\colon \mathcal{M}_n \to k$ be the witness for $n$-representability and choose a minor operation $\alpha_{f,g}\colon l \to n$ for each $f \neq g \in \mathcal{M}_l$ such that $\varphi(f\alpha_{f,g}) \neq \varphi(g\alpha_{f,g})$.

To witness $n+1$-representability, We define a map $\psi\colon \mathcal{M}_{n+1} \to k$ as follows: Let $\iota$ denote the inclusion $n \subseteq n+1$, and let $c \in n+1$ arbitrary. We then set

$$\psi(f) := \begin{cases} \varphi(e) & \exists e \in \mathcal{M}_n\colon f = e\iota \\ c & \text{else.} \end{cases}$$

Note that this is well-defined, because injectivity of $\iota$ forces $e$ to be unique: If $e'$ were another element with $f = e'\iota$, then $e = e\iota\rho = e'\iota\rho = e'$ where $\rho\colon n+1 \to n$ is any (post-)inverse to $\iota$.

Now, if $f \neq g \in \mathcal{M}_l$, we see that

$$\psi(f\alpha_{f,g}\iota) = \varphi(f\alpha_{f,g}) \neq \varphi(g\alpha_{f,g}) = \psi(g\alpha_{f,g}\iota),$$

so $\alpha_{f,g}\iota$ witnesses the inequality in degree $n+1$. □

**Corollary 4.2.6**: $\mathcal{O}(n, k) \hookrightarrow \mathcal{O}(n', k')$ whenever $n \leq n'$ and $k \leq k'$. Hence, $(n, k)$-representable minions are also $(n', k')$-representable. In particular, $\mathcal{O}(n, k) \hookrightarrow \mathcal{O}(\max(n, k))$.

*Proof*: We can compose embeddings

$$\mathcal{O}(n,k) \hookrightarrow \mathcal{O}(n+1,k) \hookrightarrow \ldots \hookrightarrow \mathcal{O}(n',k) \hookrightarrow \mathcal{O}(n',k+1) \hookrightarrow \ldots \hookrightarrow \mathcal{O}(n',k') \qquad \square$$

**Corollary 4.2.7**: For any finite set of finitely representable minions $\{\mathcal{M}_1, \ldots, \mathcal{M}_f\}$, there are $n$ and $k$ such that all embed into the minion $\mathcal{O}(n,k)$.

*Proof*: If $\mathcal{M}_i$ is $(n_i, k_i)$-representable, then by the previous Corollary 4.2.6, all of them are $\left(\max_{i=1}^{f}(n_i), \max_{i=1}^{f}(k_i)\right)$-representable. $\square$

# 5 The homomorphism order of finitely representable minions

In this chapter, we shall define and analyze the main interesting homomorphism orders of minions, most prominently the homomorphism order of finitely representable minions. We want to demonstrate how much can be said about these homomorphism orders by

1. taking abstract constructions, like products and sums,
2. showing that they preserve certain properties such as local finiteness, and
3. translating these notions to the homomorphism orders.

Using this approach, we will come to the conclusion that some of the homomorphism orders are not only bounded lattices but also *distributive* lattices. In Chapter 7, we will continue with this approach by constructing an exponential object, and deduce furthermore that the orders in question are even *bi-Heyting algebras*.

## 5.1 Definitions

**Definition 5.1.1**: We define

- $\mathfrak{M} := \underline{\mathsf{Min}} /\leftrightarrow$ ,
- $\underline{\mathsf{LFin}}$ to be the category of locally finite minions and $\mathfrak{LF} := \underline{\mathsf{LFin}} /\leftrightarrow$ ,
- $\underline{\mathsf{FinRep}}$ to be the category of finitely representable minions and $\mathfrak{F} := \underline{\mathsf{FinRep}} /\leftrightarrow$ ,
- $\underline{\mathsf{Rep}}_{A,B}$ to be the category of $(A, B)$-representable minions and $\mathfrak{F}_{A,B} := \underline{\mathsf{Rep}}_{A,B}/\leftrightarrow$ ,
- $\underline{\mathsf{Rep}}_{A}$ to be the category of $(A, A)$-representable minions and $\mathfrak{F}_{A} := \mathfrak{F}_{A,A} = \underline{\mathsf{Rep}}_{A,A}/\leftrightarrow$ , and finally
- $\underline{\mathsf{FinGen}}$ to be the category of finitely generated minions and $\mathfrak{FG} := \underline{\mathsf{FinGen}} /\leftrightarrow$ .

**Proposition 5.1.1**: Let $A$ and $B$ be sets. The canonical inclusions

$$\underline{\mathsf{Rep}}_{A,B} \hookrightarrow \underline{\mathsf{FinRep}} \hookrightarrow \underline{\mathsf{LFin}} \hookrightarrow \underline{\mathsf{Min}}$$

induce poset embeddings

$$\mathfrak{F}_{A,B} \hookrightarrow \mathfrak{F} \hookrightarrow \mathfrak{LF} \hookrightarrow \mathfrak{M}.$$

Similarly, the inclusion $\underline{\mathsf{FinGen}} \hookrightarrow \underline{\mathsf{FinRep}}$ provided by Corollary 3.3.8 induces a poset embedding $\mathfrak{FG} \hookrightarrow \mathfrak{F}$.

*Proof*: This follows from Lemma 1.1.1 since inclusions of categories are full. □

In the above proposition, the inclusions of the homomorphism orders seem to be the identity: After all, they map $[\mathcal{M}]$ to $[\mathcal{M}]$. However, that is not precisely true, because we are talking about different equivalence classes: For instance, if $\mathcal{M} \leq \mathcal{O}(A)$, the inclusion $\mathfrak{F}_A \hookrightarrow \mathfrak{F}$ would map

$$\underbrace{\{\mathcal{N} \leq \mathcal{O}(A) \mid \mathcal{N} \leftrightarrow \mathcal{M}\}}_{=[\mathcal{M}]_{\mathfrak{F}_A}} \mapsto \underbrace{\{\mathcal{N} \text{ finitely representable} \mid \mathcal{N} \leftrightarrow \mathcal{M}\}}_{=[\mathcal{M}]_{\mathfrak{F}}}.$$

We can strengthen some of the embeddings to embeddings of lattices with the following observation:

**Lemma 5.1.2**: Let $\iota\colon P \to L$ be a strong poset homomorphism and $\{p_i\}_{i\in I}$ a set of elements of $P$.

1. If $\{\iota(p_i)\}_i$ has an infimum (resp. supremum) in $L$ which is of the form $\iota(p)$, then $p$ is the infimum (resp. supremum) of $\{p_i\}_i$ in $P$.
2. If $L$ has a top (resp. bottom) element which is contained in the image of $\iota$, then its preimage under $\iota$ is the top (resp. bottom) element of $P$.

*Proof*: For the first point, assume that $\iota(p) = \bigwedge_i \iota(p_i)$. Since in particular for any $i$, we have $\iota(p) \leq \iota(p_i)$, strong homomorphicity of $\iota$ allows us to conclude $p \leq p_i$, so $p$ is indeed a common lower bound. For maximality, let $p'$ be another common lower bound to $\{p_i\}_i$. Applying $\iota$, we see that $\iota(p')$ must be a common lower bound to $\{p_i\}_i$, so by $\iota(p)$ being an infimum we have $\iota(p') \leq \iota(p)$. Strong homomorphicity once again entails $p' \leq p$, so $p$ is indeed the *greatest* lower bound.

The proof for suprema works in the same way, and the second statement regarding top and bottom elements is the special case of an empty infimum resp. supremum. □

Note that this lemma implies in particular that if infima exist in $L$ and happen to be in the image, we know that the infima exist in $P$. In particular, if we have an embedding like $\mathfrak{F} \hookrightarrow \mathfrak{L}\mathfrak{F}$, and we knew a) that $\mathfrak{L}\mathfrak{F}$ was a bounded lattice, and b) that the embedding contains all (finite) infima / suprema, we can in particular conclude that $\mathfrak{F}$ must be a bounded lattice as well, and that this embedding is a lattice homomorphism. This is precisely what we are going to do, especially since in the abstract setting of (locally finite) minions, the lattice structure can be shown with little effort.

Indeed, the following proposition is a precise way to say that the existence of products in a category translates to the existence of meets (infima) in the hom-order, and dually, that the existence of coproducts (i.e., disjoint unions) translates to the existence of joins (suprema) in the hom-order.

**Proposition 5.1.3**: Let $I$ an index set, $\mathcal{C}$ be a category, and $X_i \in \mathcal{C}_0$ objects for every $i \in I$.

1. If the coproduct $\coprod_{i\in I} X_i$ exists, then the supremum $\bigvee_{i\in I}[X_i]$ exists in $\mathcal{C}/\leftrightarrow$ and we have
$$\bigvee_{i\in I}[X_i] = \left[\coprod_{i\in I} X_i\right].$$
2. If a initial object $O$ exists, then $\mathcal{C}/\leftrightarrow$ has a least element $[O]$.
3. If the product $\prod_{i\in I} X_i$ exists, then the infimum $\bigwedge_{i\in I}[X_i]$ exists in $\mathcal{C}/\leftrightarrow$ and we have
$$\bigwedge_{i\in I}[X_i] = \left[\prod_{i\in I} X_i\right].$$
4. If an terminal object $T$ exists, then $\mathcal{C}/\leftrightarrow$ has a greatest element $[T]$.

*Proof*: For the supremum, note that the canonical embeddings $X_i \hookrightarrow \coprod_{i\in I} X_i$ witness $[X_i] \leq \left[\coprod_{i\in I} X_i\right]$, proving that the coproduct is indeed a common upper bound. To see that it is the least such, let $[Y] \in \mathcal{C}/\leftrightarrow$ be any other upper bound, i.e., $[X_i] \leq [Y]$ for all $i$.

These inequalities must be witnessed by some homomorphisms $X_i \to Y$, which combine via the cotupling to a morphism $\coprod_{i \in I} X_i \to Y$. This witnesses $\left[\coprod_{i \in I} X_i\right] \leq [Y]$ as desired.

Since an initial object can be considered the empty coproduct and a least element is the supremum over an empty set (the least upper bound of nothing), the next claim follows directly.

The proofs for infima and terminal objects work analogously. □

**Corollary 5.1.4**:
- $\mathfrak{M}$ is a complete bounded lattice.
- $\mathfrak{LF}$ is a bounded lattice.

*Proof*: Since we know from Proposition 2.3.1 that we can form arbitrary products and coproducts of abstract minions, Proposition 5.1.3 applies.

For locally finite minions, arbitrary products and coproducts might not exist, because it might cease to remain finite in each arity. However, *finite* products and coproducts of locally finite minions certainly respect this cardinality constraint, so finite infima and suprema exist in $\mathfrak{LF}$, as well as a least and greatest element. □

In each of the orders, we have an obvious atom (i.e., unique minimal element above the minimum).

**Observation 5.1.1**: Since the unaries of a minion correspond to minion homomorphisms $\mathcal{P} \to \mathcal{M}$, every nonempty minion lies above $\mathcal{P}$ in the homomorphism order. Since $\mathcal{P}$ is in all of the respective classes (finitely generated, locally finite, $n, k$-representable for every $n, k \neq 1$), this means that $\mathfrak{M}$, $\mathfrak{LF}$, $\mathfrak{F}$, $\mathfrak{F}_n$ for $n \geq 2$ and $\mathfrak{F}_{n,k}$ for $n, k \geq 2$ have $[\mathcal{P}]$ as an atom.

That the above proofs work so effortlessly is because we are in the abstract domain: Since the category of sets allows for lots of constructions like products or coproducts, so does the category of minions, being a Set-valued functor category. In Chapter 7, we will follow the same theme with another construction, called the exponential.

First, let us see which of these properties hold in the hom-order of finitely representable minions, $\mathfrak{F}$.

## 5.2 Closure under product and sum

Now onto the question of whether sums and products of finitely representable minions remain to be representable. Indeed, there are simple ways to compose the witnesses of representability in each case.

**Proposition 5.2.1**: Let $\mathcal{M}_1, ..., \mathcal{M}_f$ be a finite set of minions such that $\mathcal{M}_i$ is $(n, k_i)$-representable. Then $\mathcal{M} := \coprod_i \mathcal{M}_i$ is $\left(n, \sum_i k_i\right)$-representable.

*Proof*: Let $\varphi_i\colon \mathcal{M}_{i,n} \to k_i$ be the respective functions which witness the $(n, k_i)$-representability. Then define

$$\varphi := \underbrace{\coprod_i \mathcal{M}_{i,n}}_{=\mathcal{M}_n} \to \sum_i k_i,$$

$$\mathcal{M}_{i,n} \ni f \mapsto \sum_{j=1}^{i-1} k_i + \varphi_i(f).$$ [40]

Now let $f \neq g \in \mathcal{M}_l$. Again, there are different possibilities: Either $f$ and $g$ lie in the same sumand $\mathcal{M}_i$. Then there must be an $\alpha\colon l \to n$ such that $\varphi_i(f\alpha) \neq \varphi_i(g\alpha)$, because $\varphi_i$ witnesses the $(n, k_i)$-representability of $\mathcal{M}_i$. This implies that

$$\varphi(f\alpha) = \sum_{j=1}^{i-1} k_i + \varphi_i(f) \neq \sum_{j=1}^{i-1} k_i + \varphi_i(g) = \varphi(g\alpha).$$

Otherwise, if $f$ and $g$ lie in different summands, so do $f\alpha$ and $g\alpha$, and because $\varphi$ sends elements of different summands to elements of different summands, we conclude $\varphi(f\alpha) \neq \varphi(g\alpha)$ as well. □

**Proposition 5.2.2**: Let $\mathcal{M}_1, ..., \mathcal{M}_f$ be a finite set of minions such that $\mathcal{M}_i$ is $(n, k_i)$-representable. Then $\mathcal{M} := \prod_i \mathcal{M}_i$ is $\left(n, \prod_i k_i\right)$-representable.

*Proof*: Let $\varphi_i\colon (\mathcal{M}_i)_n \to k_i$ be the maps witnessing $(n, k_i)$-representability of the minion $\mathcal{M}_i$, and construct from it

$$\varphi := (\varphi_1, ..., \varphi_f)\colon \mathcal{M}_n = \prod_i (\mathcal{M}_i)_n \to \prod_i k_i.$$

We claim that this witnesses $\left(n, \prod_i k_i\right)$-representability: Indeed, let $(f_i)_i \neq (g_i)_i \in \mathcal{M}_l$. Since these tuples disagree, there must be a coordinate $j \in \{1, ..., f\}$ such that $f_j \neq g_j$. By $(n, k_j)$-representability of $\mathcal{M}_j$, there must exist an $\alpha\colon l \to n$ such that $\varphi_j(f_j\alpha) \neq \varphi_j(g_j\alpha)$. Since minor operations act pointwise, this implies that $\varphi\big((f_i)_i\alpha\big) = \varphi\big((f_i\alpha)_i\big)$ and $\varphi\big((g_i)_i\alpha\big) = \varphi\big((g_i\alpha)_i\big)$ must disagree in the $j$th component, and in particular are not the same. Therefore, $\varphi$ witnesses $\left(n, \prod_i k_i\right)$-representability. □

**Corollary 5.2.3**: Products and sums of finitely many finitely representable minions are finitely representable. In particular, $\mathfrak{F}$ is a bounded lattice, and the canonical inclusion $\mathfrak{F} \hookrightarrow \mathfrak{L}\mathfrak{F}$ is a homomorphism of bounded lattices.

Notably, it is not a priori clear whether quotients of finitely representable minions are still finitely representable. The reason for the difficulty is that a quotient of a, say, 2-representable minion, may identify too many elements of arity 2, but keep too many elements in a higher arity like $N$, so that distinguishing $f \neq g \in (\mathcal{M}/\theta)_n$ may not be possible anymore in arity 2. Trying to prove any higher representability would require an approximation as to "how high one should go", which is information that needs to depend on the congruence. What this information needs to be is not obvious.

[40] Abstractly, this is just the juxtaposition $\sum_i \varphi_i$ of the morphisms $\varphi_i$ which is guaranteed to exist due to the universal property of the coproduct. Our explicit construction then uses that the coproduct of the ordinals $k_i$ is isomorphic to the sum ordinal $\sum_i k_i$.

**Open Question 5.1**: Are quotients of finitely representable minions finitely representable?

Another question that remains unclear is the following:

**Open Question 5.2**: For finite sets $A, B$, is $\mathfrak{F}_{A,B}$ a lattice?

Indeed, although the product of $(A, B)$-representable minions need not be $(A, B)$-representable, a homomorphy equivalence class satisfying the property of a meet might still exist, it might just be represented by some other minion.

## 5.3 The coatom

As we saw in the beginning, the empty minion is the unique minimal element, and $[\mathcal{P}]$ is the unique element above it. If we now look at the top instead of the bottom of these orders, we see that the unique maximal element is const. But which elements are maximally below it[41]? Does such an element even exist, or do we have an infinite ascending chain whose supremum is the maximum? As it turns out, in most of these cases there is precisely one such element, i.e., we have a *coatom*.

To see what it is, we need to gather a few insights about the role of unaries in minions.

### 5.3.1 Unary decomposition: The coatom of $\mathfrak{F}_{n,k}$

In a minion $\mathcal{M}$, every element $f \in \mathcal{M}_n$ is associated to a unique unary $f\tau$, because there is only one minor operation $\tau\colon n \twoheadrightarrow 1$. If $f$ were an element of a function minion $\mathcal{O}(A, B)$, then $f\tau$ would be the function $x \mapsto f(x, x, ..., x)$.

What we will present here is not new, but often exposed using pp-contsructions; see [41, Sections 1.6,1.7].

Throughout this section, let $\tau$ denote the unique constant map $n \to 1$. Note that while this notation is ambiguous with respect to $n$, the arity can always be inferred from the context, since in an expression such as $f\tau$, we necessarily have $n = \operatorname{ar} f$.

**Definition 5.3.1**: Let $\mathcal{M}$ be a minion and $S \subseteq \mathcal{M}_1$ a set of unaries. We define the (indexed) set

$$\mathcal{M}{\downarrow_S} := (\{f \in \mathcal{M}_n \mid f\tau \in S\})_n$$

of all elements whose unary minor is $u$ to be the *restriction of $\mathcal{M}$ onto the unary $u$*. As usual, we abbreviate $\mathcal{M}{\downarrow_{\{u\}}} =: \mathcal{M}{\downarrow_u}$.

[41] An element maximally below another element $y$ is called a "co-cover" of $y$.

**Proposition 5.3.1**: $\mathcal{M}{\downarrow_S}$ is a subminion of $\mathcal{M}$.

*Proof*: If $f\tau = u$, then also $(f\alpha)\tau = f\tau = u$, as there is only one map ar $f \to 1$, which implies the special case $\mathcal{M}{\downarrow_u} \leq \mathcal{M}$ for every $u$. Thus, $\mathcal{M}{\downarrow_S}$ is a union of subminions, so a subminion itself. □

**Proposition 5.3.2**: Let $\mathcal{L}$ be a minion with a single unary and $\mathcal{M}, \mathcal{N}$ arbitrary minions. Then every homomorphism $\mathcal{L} \to \mathcal{M} + \mathcal{N}$ must factor through one of the components, i.e., $\mathcal{L} \to \mathcal{M} + \mathcal{N} \cong (\mathcal{L} \to \mathcal{M}) + (\mathcal{L} \to \mathcal{N})$.

*Proof*: Let $u \in \mathcal{L}_1$ be the unique unary and $h\colon \mathcal{L} \to \mathcal{M} + \mathcal{N}$ be a homomorphism. If $h(u) \in \mathcal{M}$, then for all $f \in \mathcal{M}_n$ we have $h(f)\alpha = h(f\alpha) = h(u) \in \mathcal{M}$, so in particular $h(f)\alpha \notin \mathcal{N}$. Since $\mathcal{N}$ is a subminion of $\mathcal{M} + \mathcal{N}$, we know by contraposition that $h(f) \notin \mathcal{N}$, but this means that $h(f) \in \mathcal{M}$. Therefore $\operatorname{Im} h \subseteq \mathcal{M}$. The analogous argumentation works for $\mathcal{N}$.□

Note that the argument actually works for arbitrary sums! Following the language of lattice theory, we might say that minions with one unary are a *completely join-prime elements* in the lattice that is the homomorphism order of minions.

Another statement very similar in nature is very similar statement is the following.

**Lemma 5.3.3**: Let $\mathcal{L}$ be a minion with a single unary $u$ and $\mathcal{M}$ arbitrary. Then $h\colon \mathcal{L} \to \mathcal{M}$ has image in $\mathcal{M}{\downarrow_{h(u)}}$

*Proof*: Letting $\tau\colon n \twoheadrightarrow 1$ as usual, we see that $h(f)\tau = h(f\tau) = h(u)$ □

**Proposition 5.3.4**: For every minion $\mathcal{M}$, we have $\mathcal{M} \cong \sum_{u \in \mathcal{M}_1} \mathcal{M}{\downarrow_u}$.

*Proof*: Let $h\colon \sum_{u \in \mathcal{M}_1} \mathcal{M}{\downarrow_u} \to \mathcal{M}$ be the canonical inclusion of the subminions which on each component sends $m \mapsto m$. This map is an injection on each of the summands $\mathcal{M}{\downarrow_u}$. Since in each degree, the set $\{(\mathcal{M}{\downarrow_u})n \mid u \in \mathcal{M}_1\}$ is a partition of $\mathcal{M}_n$, we know that the subminions $\mathcal{M}{\downarrow_u}$ have trivial intersection. Therefore, the combined homomorphism $h$ must be injective as well. Since we always have $f \in \mathcal{M}{\downarrow_{f\tau}}$, it is also surjective, so a bijective minion homomorphism. By Proposition 2.3.2, $h$ is an isomorphism. □

For the following considerations, let $n, k > 1$. The reason for that is the following observation:

**Observation 5.3.1**:
- $\mathcal{O}(1, k)$ contains each arity precisely $k$ constant operations, hence $\mathcal{O}(1, k) \cong k \cdot *$.
- $\mathcal{O}(n, 1)$ contains precisely one constant in each arity, hence $\mathcal{O}(n, 1) \cong *$.

In particular, this implies that the homomorphism orders $\mathfrak{F}_{1,k}$ and $\mathfrak{F}_{n,1}$ are both isomorphic to the two-element poset $0 \leq 1$ (signifying $[\emptyset]$ and $[*]$, respectively).

**Proposition 5.3.5**: If $\mathcal{M} \leq \mathcal{O}(n,k)$ does not contain a constant, it is a sub-minion of $\mathcal{O}(n,k)\downarrow_{\mathrm{nc}}$, where $\mathrm{nc} := \{u\colon n \to k \mid u \text{ is not constant}\}$.

*Proof*: Arguing contrapositively, assume $\mathcal{M} \not\leq \mathcal{O}(n,k)\downarrow_{\mathrm{nc}}$. By Proposition 5.3.4[42], this implies that there is some $f \in \mathcal{M}_n$ such that $f\tau \notin \mathrm{nc}$, i.e., $f\tau$ is constant. But $\mathcal{M}$ is a subminion, so $f\tau \in \mathcal{M}$ – meaning it contains a constant. □

**Corollary 5.3.6** (cf. [41, Theorem 1.38]): $[\mathcal{O}(n,k)\downarrow_{\mathrm{nc}}]$ is a coatom in $\mathfrak{F}_{n,k}$ (unless $n = 1$ or $k = 1$).

*Proof*: Otherwise, assume $[\mathcal{M}] < [*]$ in $\mathfrak{F}_{n,k}$, where we choose a representative $\mathcal{M}$ which is a subminion of $\mathcal{O}(n,k)$. Since the inequality means that $* \nrightarrow \mathcal{M}$, it follows from Proposition 5.3.5 that $[\mathcal{M}] \leq [\mathcal{O}(n,k)\downarrow_{\mathrm{nc}}]$. □

### 5.3.2 Idempotents: The coatom of $\mathfrak{F}$

If $n = k$, i.e., if we consider $\mathcal{O}(n)$, there is a special unary element: The *identity*. Functions $f$ whose associated unary $f\tau$ is the identity, i.e., for which $f(x,...,x) = x$, are called *idempotent*. we will demonstrate in this section how all of the coatom-representatives $[\mathcal{O}(n,k)\downarrow_{\mathrm{nc}}]$ are actually all homomorphically equivalent to the minion[43] of idempotents in $\mathcal{O}(2)$.

**Definition 5.3.2**: We call $\mathcal{I}_n := \mathcal{O}(n)\downarrow_{\mathrm{id}_n}$ the *minion of idempotents over $n$*.

**Definition 5.3.3** (Reflection): Let $l\colon A' \to A$, $r\colon B \to B'$. Define the *reflection* along $(l,r)$

$$R_{l,r}\colon \mathcal{O}(A,B) \to \mathcal{O}(A',B')$$

in each arity to be

$$\begin{aligned} \mathcal{O}(A,B)_i &\to \mathcal{O}(A',B')_i \\ f &\mapsto \underbrace{f(_;l);r}_{=\underline{x}\mapsto r(f(\underline{x};l))} . \end{aligned}$$

**Lemma 5.3.7**: Let $l\colon A' \to A$, $r\colon B \to B'$. Then $R_{l,r}$ is a minion homomorphism. Furthermore, if $l'\colon A'' \to A'$, $r'\colon B' \to B''$, then $R_{l',r'} \circ R_{l,r} = R_{l';l,\ r;r'}$.

*Proof*: Let $\alpha\colon i \to i'$ be a minor operation, and $f\colon A^i \to B$. Then

[42] and because nonconstant elements in $\mathcal{O}(n,k)$ *exist*, as guaranteed by $n,k > 1$

[43] this is actually closed under composition, so a *clone*

$$\begin{aligned}
R_{l,r}(f\alpha) &= (f\alpha)(_;l);r && \text{Definition of } R_{l,r} \\
&= f((\alpha;_);l);r && \text{Minion structure of } \mathcal{O}(A,B) \\
&= f(\alpha;(_;l));r && \text{Associativity} \\
&= (f(_;l);r)\alpha && \text{Minion structure of } \mathcal{O}(A,B) \\
&= R_{l,r}(f)\alpha && \text{Definition of } R_{l,r}
\end{aligned}$$

For the second conclusion, note that

$$R_{l',r'}\big(R_{l,r}(f)\big) = R_{l',r'}(f(_;l);r) = f(_;l';l);r;r' = R_{l';l,\ r;r'},$$

as desired. □

We can visualize what the reflection does using a commutative diagram as depicted below:

$$\begin{array}{ccccc}
A''^i & \xrightarrow{\;l'\;} & A'^i & \xrightarrow{\;l\;} & A^i \\
{\scriptstyle R_{l',r'}(R_{l,r}(f))}\big\downarrow & & \big\downarrow{\scriptstyle R_{l,r}(f)} & & \big\downarrow{\scriptstyle f} \\
B'' & \xleftarrow[\;r'\;]{} & B' & \xleftarrow[\;r\;]{} & B
\end{array}$$

Figure 10: How reflections behave with respect to composition.

**Lemma 5.3.8**: If $l\colon A' \to A$ and $r\colon B \to B'$, then $R_{l,r}[\mathcal{O}(A,B)\!\downarrow_u] \subseteq \mathcal{O}(A',B')\!\downarrow_{lur}$. This means that $R_{l,r}$ witnesses a homomorphism $\mathcal{O}(A,B)\!\downarrow_u \to \mathcal{O}(A',B')\!\downarrow_{lur}$.

*Proof*: This is clear by homomorphicity of $R_{l,r}$: If $f\tau = u$, then

$$R_{l,r}(f)\tau = R_{l,r}(f\tau) = R_{l,r}(u) = lur.$$

□

**Proposition 5.3.9**: For all $n, k > 1$ and $u\colon n \to k$, we have $\mathcal{O}(n,k)\!\downarrow_u \leftrightarrow \mathcal{I}_{|\operatorname{Im} u|} = \mathcal{O}\!\downarrow_{|\operatorname{Im} u|}$.

*Proof*: Let $l := |\operatorname{Im} u|$. Then $u$ factors as $u = e;m$ for some surjection $e\colon n \twoheadrightarrow l$ and injection $m\colon l \hookrightarrow k$. For the backward direction, note how $u = em = e\,\mathrm{id}_l\, m$, so $R_{e,m}$ witnesses a homomorphism $\mathcal{I}_l \to \mathcal{O}(n,k)\!\downarrow_u$.

For the converse, we can choose a pre-inverse $e'\colon l \to n$ and post-inverse $m'\colon k \to l$. This ensures that $\mathrm{id}_l = e'emm' = e'um'$, hence $R_{e',m'}$ witnesses $\mathcal{O}(n,k)\!\downarrow_u \to \mathcal{I}_l$. □

**Proposition 5.3.10**: For all $n > 1$, we have $\mathcal{I}_2 \leftrightarrow \mathcal{I}_n$.

*Proof*: For $\to$, let $\iota\colon 2 \hookrightarrow n$ be the canonical inclusion. Choosing a retraction $\rho\colon n \twoheadrightarrow 2$, we observe that $\mathrm{id}_2 = \iota;\rho = \iota\,\mathrm{id}_n\,\rho$, hence $R_{\iota,\rho}$ serves as a witness for $\mathcal{I}_n \to \mathcal{I}_2$.

For the other direction, we shall show that there is a homomorphism $\mathcal{I}_2 \to \mathcal{O}(n)$ which maps $\mathrm{id}_2$ to $\mathrm{id}_n$. Under the correspondence of Corollary 4.1.4, homomorphisms $\mathcal{I}_2 \to \mathcal{O}(n)$ are induced by functions $\varphi\colon (\mathcal{I}_2)_n \to n$ via the assignment $\tilde{\varphi}\colon f \mapsto \varphi(f__)$. So if $f = \mathrm{id}_2$, the corresponding unary element of $\mathcal{O}(n)$ is $\eta \mapsto \varphi(\mathrm{id}_2\, \eta)$. Since the lift of the identity along $\eta\colon 1 \to n$ is the projection $\pi_{\eta(0)}$, we see that $\tilde{\varphi}(\mathrm{id}_2) = i \mapsto \varphi(\pi_i)$. This is the identity precisely if we choose $\varphi$ such that $\varphi(\pi_i) = i$ for all $i \in n$, which is certainly possible, as the $\pi_i\colon 2^n \to 2$ are all distinct. □

**Corollary 5.3.11** (cf. [41, Theorem 1.38]): Whenever $n, k > 1$, we have $\mathcal{O}(n,k)\!\downarrow_{\mathrm{nc}} \leftrightarrow \mathcal{I}_2$.

*Proof*: By the unary decomposition Proposition 5.3.4, we have

$$\begin{aligned}
\mathcal{O}(n,k)\!\downarrow_{\mathrm{nc}} &= \sum_{i>1}^{k} \sum_{\substack{u\in \mathrm{nc} \\ |\operatorname{Im} u| = i}} \mathcal{O}(n,k)\!\downarrow_u && \text{Definition} \\
&\leftrightarrow \sum_{i>1}^{k} \sum_{\substack{u\in \mathrm{nc} \\ |\operatorname{Im} u| = i}} \mathcal{I}_i && \text{Proposition 5.3.9} \\
&\leftrightarrow \sum_{i>1}^{k} \sum_{\substack{u\in \mathrm{nc} \\ |\operatorname{Im} u| = i}} \mathcal{I}_2 && \text{Proposition 5.3.10} \\
&\leftrightarrow \mathcal{I}_2.
\end{aligned}$$

□

### 5.3.3 Cores: A canonical representative of the coatom of $\mathfrak{F}$

While we have seen that $\mathcal{O}(n,k)_{\mathrm{nc}}$ is homomorphically equivalent to $\mathcal{I}_2$ – arguably a "smaller" minion – we will demonstrate soon that there is an even smaller representative. But what do ewe mean by "small"? The relevant notion here is that of a *core*, a concept well known from graphs [42, section 1], [10, section 2.4].

The basic idea is that homomorphic equivalences lead to endomorphisms, whose image will always be a homomorphically equivalent sub-minion. If we have an object which does not admit any non-surjective endomorphisms, i.e., where every endomorphism is an automorphism[44].

**Definition 5.3.4**: A *core* in a category $\mathcal{C}$ is an object for which every endomorphism is an automorphism. We say that an object $X \in \mathcal{C}_0$ *has a core* if it is homomorphically equivalent to a core.

**Proposition 5.3.12**: Any two homomorphically equivalent cores are isomorphic.

[44] These statements are equivalent for finite graphs, but not in general.

*Proof*: Let $C$ and $D$ be cores wtih $f\colon C \to D$ and $g\colon D \to C$. Since $f;g \in \mathrm{End}(C)$ and $g;f \in \mathrm{End}(D)$ are invertible by assumption, the compositions $fg(fg)^{-1} = \mathrm{id}_C$ and $(gf)^{-1}gf = \mathrm{id}_d$ imply that $f$ has a post- as well as a pre-inverse, and hence must be an isomorphism.□

Since it is unique up to isomorhism, it is justified to associate "the" core to an object.

**Definition 5.3.5**: If $X \in \mathcal{C}_0$ has a core, we will denote it by $\mathrm{core}(X)$ .

**Proposition 5.3.13**: If a locally finite minion $\mathcal{M}$ has only finitely many endomorphisms, it has a core.

For the proof, recall the abuse of terminology in which we refer to minion homomorphisms as "surjective" or "injective" if they are so in every arity, and recall furthermore that the image is not a set, but an indexed set – specifically, a subminion – so $\operatorname{Im} f \subseteq \operatorname{Im} g$ is also understood to be a statement upheld in each arity.

*Proof*: Note first that the subminions of $\mathcal{M}$ are naturally ordered by inclusion. Since $\mathrm{End}(\mathcal{M})$ is finite, there must be an $e\colon \mathcal{M} \to \mathcal{M}$ with minimal image. Let the surjection $h\colon \mathcal{M} \to \operatorname{Im} e$ denote its corestriction onto the image, i.e., $e = h;\iota$ where $\iota$ is the inclusion $\operatorname{Im} e \subseteq \mathcal{M}$. Now if $f \in \mathrm{End}(\operatorname{Im} e)$, we know that $e;f;\iota$ is another endomorphism of $\mathcal{M}$ with image contained in $\operatorname{Im} e$. Thus, by minimality, they must agree, forcing $\operatorname{Im} f = \operatorname{Im} e$, i.e., $f\colon \operatorname{Im} e \to \operatorname{Im} e$ must have been surjective. Since $\mathcal{M}$ is locally finite, $f$ must already be bijective. This implies that it is an isomorphism by Proposition 2.3.2, whereby we have demonstrated that $\mathrm{End}(\operatorname{Im} e) = \mathrm{Aut}(\operatorname{Im} e)$. □

Let us give some examples of cores and non-cores.

- Clearly, the constant minion the minion $*$ is a core because $\mathrm{End}_{\underline{\mathsf{Min}}}(*) = \{\mathrm{id}_*\} = \mathrm{Aut}_{\underline{\mathsf{Min}}(*)}$.
- If $\mathcal{M}$ contains a constant $u$ but also a non-constant, then the map $m \mapsto u\iota_{\mathrm{ar}\ m}$ (where $\iota_n\colon 1 \to n$ is the canonical inclusion) is a non-surjective, hence non-invertible, minion endomorphism. Therefore, $\mathcal{M}$ is not a core.
- The disjoint union $\mathcal{M} + \mathcal{M}$ is never a core (unless $\mathcal{M}$ is empty) as witnessed by the projection onto the first component.
- If $\underline{A}, \underline{B}$ are relational structures such that $\mathrm{Aut}(B)$ acts transitively via postcomposition on $\underline{A} \to \underline{B}$, then $\mathrm{Pol}(\underline{A}, \underline{B})$ is not a core, becase all the unary components are isomorphic: Indeed, if $u' = u;h$ for some automorphism $h \in \mathrm{Aut}(B)$, we see that $\mathrm{Pol}(\underline{A}, \underline{B})\downarrow_{u'} = R_{\mathrm{id}_{\underline{B}},h}\,\mathrm{Pol}(\underline{A}, \underline{B})\downarrow_u$, and reflections of polymorphism minions along isomorphisms are isomorphisms.
- A special case of the above is the scenario where $\underline{A}$ is a relational structure with more than one element which is a core. In that case, the unary decomposition gives us

$$\mathcal{M} := \mathrm{Pol}(\underline{A}) = \mathrm{Pol}(\underline{A})\downarrow_{\mathrm{nc}} \twoheadrightarrow \mathrm{Pol}(\underline{A})\downarrow_{\mathrm{id}_n} \subseteq \mathrm{Pol}(\underline{A}),$$

  so $\mathrm{Pol}(\underline{A})$ is not a core if $\underline{A}$ has any non-identity endomorphisms. A structure *without* any non-identity endomorphisms is called *rigid* [42, section 4.2], also called a *rigid core* to emphasize that rigid structures must in particular be cores.
- The projection minion $\mathcal{P}$ is a core (even rigid), because homomorphisms $\mathcal{P} \to \mathcal{P}$ correspond to a choice of unary element in $\mathcal{P}$, of which there is only one.

We now turn our attention to finitely representable minions, and show that they have only finitely many endomorphisms. However, we will be a bit more fine-grained in our analysis to provide a useful criterion of when such a minion is already core. Recall from the characterization of homomorphism existence that $\mathcal{T}_n$ is the full transformation monoid, and denote by $\mathcal{T}_n$ -set the category of all $\mathcal{T}_n$-sets $\mathcal{T}_n \curvearrowright X$. For brevity, we will denote endomorphism monoids in this category as $\mathrm{End}_{\mathcal{T}_n}(X)$.

**Proposition 5.3.14**: Let $\mathcal{M}$ be $(n, k)$-representable. Then the restriction map

$$\begin{aligned} \mathrm{res}_n\colon\ \mathrm{End}_{\underline{\mathsf{Min}}}(\mathcal{M}) &\to \mathrm{End}_{\mathcal{T}_n}(\mathcal{M}_n), \\ e = (e_l)_l &\mapsto e_n \end{aligned}$$

is a well-defined injective monoid homomorphism.

*Proof*: First note that it is clearly a homomorphism, as the composition $f;g$ of endomorphisms in $\mathcal{M}$ is given by $(f;g)_l := f_l;g_l\colon \mathcal{M}_l \to \mathcal{M}_l$ for each arity $l$. To see injectivity, note that the map

$$\begin{aligned} \mathrm{End}_{\underline{\mathsf{Min}}}(\mathcal{M}) &\to (\mathcal{M}_n \to k), \\ e &\mapsto (f \mapsto e_n(f)(\mathrm{id}_n)), \end{aligned}$$

being a restriction of the bijection ⬚ defined in Definition 4.1.1, must be injective. Since the map does not require any information other than the arity-$n$-part of the endomorphism, it factors through $\mathrm{res}_n$, which we can conclude to be injective as well. □

This result begs the question:

**Open Question 5.3**: When is $\mathrm{res}_n$ surjective, i.e. when does an endomorphism of $\mathcal{T}_n \curvearrowright \mathcal{M}_n$ extend to an endomorphism of the full minion (assuming $n$-representability)?

**Corollary 5.3.15** (cf. [43, Thm. 3.4]): Finitely representable minions have finitely many endomorphisms. Hence, every finitely representable minion has a core.

*Proof*: This follows from Proposition 5.3.14 because $\mathrm{End}_{\mathcal{T}_n}(\mathcal{M}_n) \subseteq (\mathcal{M}_n)^{\mathcal{M}_n}$ must be finite. This allows us to apply Proposition 5.3.13. □

**Corollary 5.3.16** (cf. [43, Lem. 3.5]): Let $\mathcal{M}$ be $n$-representable. If $\mathcal{T}_n \curvearrowright \mathcal{M}_n$ is a core, then $\mathcal{M}$ is as well.

*Proof*: This follows from the observation that submonoids of finite groups are subgroups: Indeed, if $M \leq G$ is a submonoid, then every $m \in M$ must have finite order in $G$, so $m^k = 1$ for some $k$. This implies that the inverse, $m^{k-1}$, is in $M$.

Now let $\mathcal{M}$ be $n$-representable, which guarantees the restriction $\mathrm{End}_{\underline{\mathsf{Min}}}(\mathcal{M}) \to \mathrm{End}_{\mathcal{T}_n}(\mathcal{M}_n)$ to be injective. If $\mathcal{T}_n \curvearrowright \mathcal{M}_n$ is a core, then $\mathrm{End}_{\mathcal{T}_n}(\mathcal{M}_n) = \mathrm{Aut}_{\mathcal{T}_n}(\mathcal{M}_n)$ is a group, so by the preceding argument, $\mathrm{End}_{\underline{\mathsf{Min}}}(\mathcal{M})$ is as well. This proves $\mathrm{End}(\mathcal{M}) = \mathrm{Aut}(\mathcal{M})$. □

**Lemma 5.3.17** (cf. [44, Corollary 1.2]): Every element $m$ in a finite monoid $M$ has an idempotent power.

*Proof*: Since the set $\{m^n \mid n > 1\}$ is finite, there must be $1 \leq c < c'$ such that $m^c = m^{c'}$. Writing this as $m^{c+d} = m^c$, we see that this implies $m^e = m^{e+kd}$ for any $e \geq c$ and any $k \in \mathbb{N}$. Thus, $\left(m^{dc}\right)^2 = m^{dc+dc} = m^{dc}$, as desired. □

In the case of finitely many endomorphisms, we therefore can always choose the retraction on a core substructure to be idempotent.

In order to determine cores of $\mathcal{T}_n$-sets, it is useful to have a small generating set for the full transformation monoid. In the following lemma, we determine three generators, which for easier comprehension are depicted as wiring diagrams to be read from bottom to top.

**Lemma 5.3.18**:

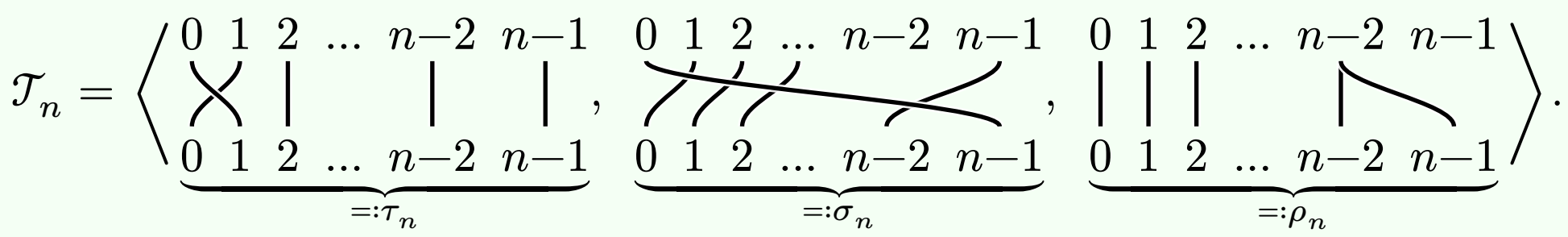

*Proof*: It is well-known that a transposition and an $n$-cycle generate the symmetric group $S_n$, so $\langle \sigma_n, \tau_n, \rho_n \rangle$ is guaranteed to contain all the injective (equivalently: bijective) maps. For the remaining cases let $\varepsilon \in \mathcal{T}_n$ be noninjective, and observe that $\varepsilon$ can be written as a sequence $\varepsilon_1;\ldots;\varepsilon_l$ where each $\varepsilon_k$ identifies precisely two values.

We claim that such a morphism must always be of the form $\alpha;\rho_n;\beta$, where $\alpha$ and $\beta$ are permutations. Indeed, assume $\varepsilon_k$ identifies $i < j \in n$. Let $\alpha$ be the map which swaps $i$ with $n-2$ and $j$ with $n-1$, and define

$$\beta\colon \begin{cases} i & \mapsto \varepsilon(n-2) \\ j & \mapsto \varepsilon(n-1) \\ n-2 \mapsto \varepsilon(i) & (\ldots = \varepsilon(j)) \\ n-1 \mapsto r & \\ s & \mapsto \varepsilon(s) \qquad \text{for all } s \notin \{i, j, n-2, n-1\} \end{cases}$$

where $r$ is the unique element in $n \setminus \operatorname{Im} \varepsilon$. To see that $\beta$ is a well-defined function, let us analyze the overlaps. For the overlaps of $i$, note that $i \neq j$ and $i \neq s$ by definition. Furthermore, since $i < j$, we also exclude $n-1$ as a possible value. The only remaining possible overlap is thus $i = n-2$, but in that case $\varepsilon(n-2) = \varepsilon(i)$, so we assign the same function value. Note that $\beta$ is surjective. The composition $\alpha;\tau_n;\beta$ can now be seen to map

$$\begin{aligned} i &\mapsto n-2 \mapsto n-2 \mapsto \varepsilon(i), \\ j &\mapsto n-1 \mapsto n-2 \mapsto \varepsilon(i), \\ n-2 &\mapsto i \qquad\ \ \mapsto i \qquad\ \ \mapsto \varepsilon(n-2), \\ n-1 &\mapsto j \qquad\ \ \mapsto j \qquad\ \ \mapsto \varepsilon(n-1), \\ s &\mapsto s \qquad\ \ \mapsto s \qquad\ \ \mapsto \varepsilon(s), \end{aligned}$$

where $s$ ranges through $n \setminus \{i, j, n-2, n-1\}$ as above, so we have indeed described $\varepsilon_k$.□

### 5.3.4 The core of $\mathcal{J}_2$

To determine the core of $\mathcal{J}_2$, we shall first determine a core of the $\mathcal{T}_2$-set $(\mathcal{J}_2)_2$, i.e., the binary idempotent operations.

have their values at $(0,0)$ and $(1,1)$ determined by idempotence, so they can freely choose their value at $(0,1)$ and $(1,0)$. The four possibilities are $\wedge$, $\vee$, $\pi_0$, and $\pi_1$.

**Proposition 5.3.19**: The $\mathcal{T}_2$-set structure of $(\mathcal{J}_2)_2$ restricts to the subset $\{\pi_0, \pi_1, \wedge\}$. This sub-$\mathcal{T}_2$-set is a core of $(\mathcal{J}_2)_2$, as witnessed by the idempotent endomorphism $e\colon (\mathcal{J}_2)_2 \to (\mathcal{J}_2)_2$ which sends $\vee$ to $\wedge$ and is the identity otherwise.

Furthermore, this endomorphism is in the image of $\mathrm{res}_2\colon\ \mathrm{End}_{\underline{\mathsf{Min}}}(\mathcal{J}_2) \to \mathrm{End}_{\mathcal{T}_2}\big((\mathcal{J}_2)_2\big)$.

*Proof*: To visualize the $\mathcal{T}_2$-set structure of $(\mathcal{J}_2)_2$, we note that $\sigma := \sigma_2 = \tau_2 = (1,0)$, $\rho := \rho_2 = (0,0)$, so a set of generators is given by $\sigma$ and $\rho$. Endomorphisms of $(\mathcal{J}_2)_2$ are thus precisely these maps which commute with the action of both $\sigma$ and $\rho$. We observe that

- $\sigma$ swaps $\pi_0$ and $\pi_1$, but leaves $\wedge$ and $\vee$ fixed, as they are symmetric, and
- for all $f \in (\mathcal{J}_2)_2$, we have $(f\rho)(x_0, x_1) = f(x_0, x_0) = x_0$ by idempotence, so $f\rho = \pi_0$.

These relations can be visualized using a labelled digraph just as for group actions, where we have a $\sigma$-labelled edge $(f, g)$ whenever $g = f\sigma$ and analogously for $\rho$. Endomorphisms of this digraph then correspond by construction to endomorphisms of $\mathcal{T}_2 \curvearrowright (\mathcal{J}_2)_2$.

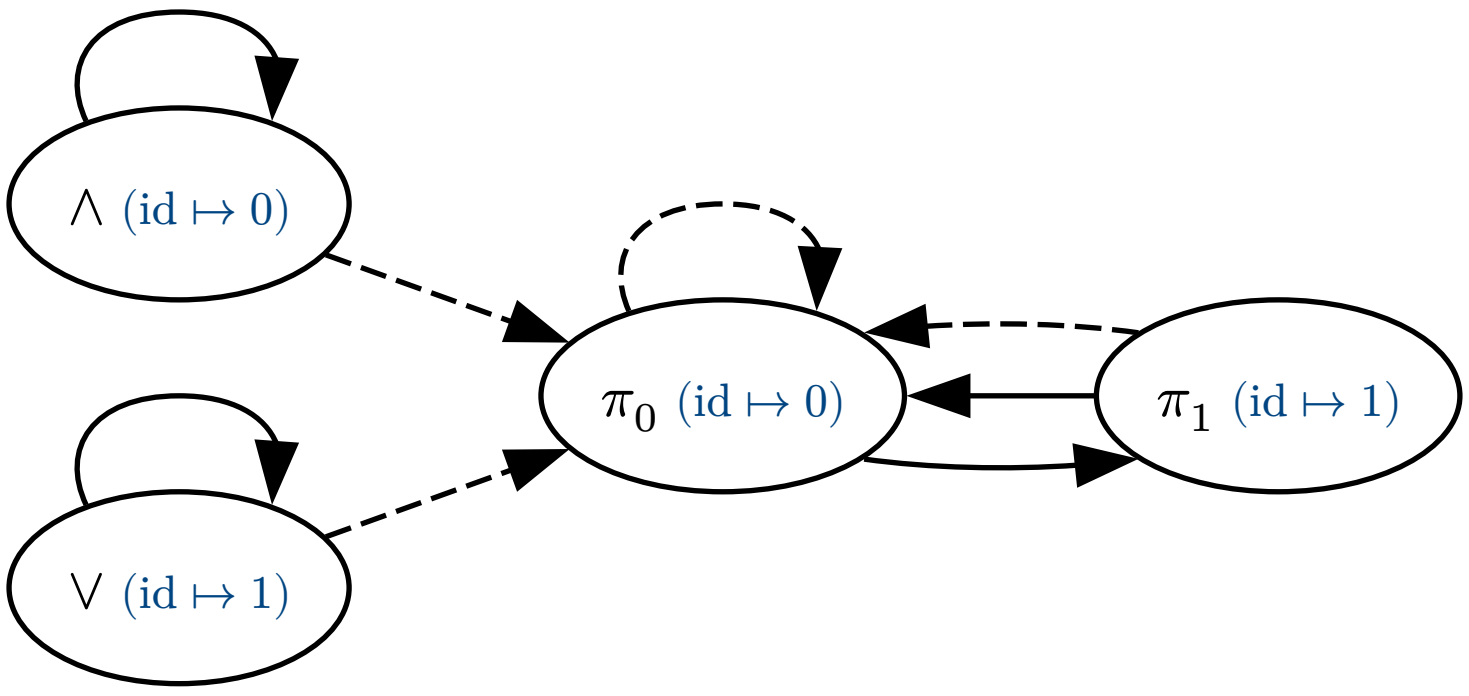

Figure 11: The action of the generators $\sigma$ (solid) and $\tau$ (dashed) of $\mathcal{T}_2$ on $(\mathcal{J}_2)_2$

We see that sending $\wedge$ to $\vee$ is indeed a nontrivial idempotent endomorphism onto $\{\pi_0, \pi_1, \wedge\}$, which induces a $\mathcal{T}_2$-subset because there are now edges leaving this subset. Evidently, this is a core, because every identification of two of these elements would create either a $\sigma$- or a $\rho$-loop which does not exist, and an identification of three of these elements would demand an element with a $\sigma$- and a $\rho$-loop, which also does not exist.

To see that this endomorphism $e$ lifts to an endomorphism of $\mathcal{I}_2$, note that that the induced homomorphism $\mathcal{I}_2 \to \mathcal{O}(2)$, which sends $f$ to $\alpha \mapsto e(f\alpha)(\mathrm{id}_2)$, maps the identity $f = \mathrm{id}_2$ to $\alpha \mapsto e(\mathrm{id}_2\,\alpha)(\mathrm{id}_2)$. Identifying unary minor operations $2^1$ with one-tuples, we see that we map

$$(i) \mapsto e(\mathrm{id}_2(i))(\mathrm{id}_2) = e(\pi_i)(\mathrm{id}_2) = \pi_i(\mathrm{id}_2) = i,$$

where we have used that the lift of the identity along the operation $(i)\colon 1 \to 2$ is the $i$-th projection, and that $e$ fixes projections. By Lemma 5.3.3, the lifted homomorphism $\mathcal{I}_2 \to \mathcal{O}(2)$ must have image in $\mathcal{O}(2){\downarrow_{\mathrm{id}_2}}$, so it is[45] an endomorphism. □

We are now ready to identify $\mathrm{core}(\mathcal{I}_2) = \mathrm{Im}\, h$, where

$$h(f)(\alpha) = e(f\alpha)(\mathrm{id}_2)$$

is the lift of $e$ to $\mathrm{End}_{\mathsf{Min}}(\mathcal{M})$. Since $h$ is idempotent, elements of the image are precisely the fixed points, i.e., those $f \in (\mathcal{I}_2)_n$ for which $h(f) = f$. This happens if for all $\alpha\colon n \to 2$, we have

$$f(\alpha) = e(f\alpha)(\mathrm{id}_2) = \begin{cases} 1 \text{ if } f\alpha = \pi_1 \\ 0 \text{ else} \end{cases},$$

or in other words, $f \in \mathrm{Im}\, h$ if and only if

$$f(\alpha) = 1 \Leftrightarrow f\alpha = \pi_1.$$

The right hand side is an equality of binary idempotent functions, so let us evaluate at all 2-tuples to see what the constraint $f\alpha = \pi_1$ entails.

- $(f\alpha)(0,0) = \pi_1(0,0)$ means $(f\alpha)(0,0) = 0$, but we already know that because of idempotence of $f$ and thus $f\alpha$.
- Evaluating at $(1,1)$ yields the same insight.
- Evaluating at the identity $(0,1)$ yields $(f\alpha)(\mathrm{id}_2) = f(\alpha) = 1$, which is the statement of the left hand side.
- Evaluating at the swap $\sigma_1 = (1,0)$ yields $(f\alpha)(\sigma) = f(\alpha\sigma) = \pi_1(\sigma) = 0$.

This analysis simplifies the characterization to

$$f \in \mathrm{Im}\, h \quad \Leftrightarrow \quad \forall \alpha \in 2^{\mathrm{ar}\, f}\colon (f(\alpha) = 1 \Rightarrow f(\alpha\sigma) = 0).$$

Let us further investigate what the right hand side means in practice. If for instance $\mathrm{ar}\, f = 4$ and we know that $f(0,0,1,1) = 1$, then the requirement demands that $f(1,1,0,0) = 1$. Putting this into a more familiar schematical form, we have

$$\begin{aligned} & f(\ \ \alpha_0, \ldots, \ \ \alpha_{n-1}) = 1 \\ \Rightarrow\ & f(\neg\alpha_0, \ldots, \neg\alpha_{n-1}) = 0. \end{aligned}$$

This condition is a preservation property: If the arguments are filled with column tuples out of the relation $\{(0,1),(1,0)\}$ – resulting in complementary rows – then the function values shall

[45] or, to be pedantic, "corestricts to"

form everything but the forbidden column tuple $(1, 1)$, i.e., lie in $\{(0, 0), (0, 1), (1, 0)\}$. We have thus proven the following result, which is due to Libor Barto:

**Proposition 5.3.20** (cf. [45, Thm. 3.8]):

$$\mathrm{core}(\mathcal{I}_2) = \mathrm{Pol}\Big( 0 \longrightarrow 1 \;,\; \circlearrowright 0 \text{ --- } 1 \Big).$$

## 5.4 Distributivity

**Proposition 5.4.1**: Let $\mathcal{L}, \mathcal{M}, \mathcal{N}$ be minions. We have
1. $\mathcal{L} \times (\mathcal{M} + \mathcal{N}) \cong (\mathcal{L} \times \mathcal{M}) + (\mathcal{L} \times \mathcal{N})$ and
2. $\mathcal{L} + (\mathcal{M} \times \mathcal{N}) \leftrightarrow (\mathcal{L} + \mathcal{M}) \times (\mathcal{L} + \mathcal{N})$.

*Proof*: Throughout this proof, we will use the universal properties of the product and coproduct (sum), i.e., Proposition 2.3.1.

*Statement 1:* Let

$$h\colon (\mathcal{L} \times \mathcal{M}) + (\mathcal{L} \times \mathcal{N}) \to \mathcal{L} \times (\mathcal{M} + \mathcal{N})$$

be the canonical homomorphism induced by the injections $(\mathrm{id}_{\mathcal{L}} \times \iota_1)\colon (\mathcal{L} \times \mathcal{M}) \to \mathcal{L} \times (\mathcal{M} + \mathcal{N})$ and $(\mathrm{id}_{\mathcal{L}} \times \iota_2)\colon (\mathcal{L} \times \mathcal{N}) \to \mathcal{L} \times (\mathcal{M} + \mathcal{N})$. Since $\mathrm{id}_{\mathcal{L}}$, $\iota_1$ and $\iota_2$ are minion homomorphisms, this must be a minion homomorphism. Furthermore, in each degree $n$ we have that the map

$$h\colon (\mathcal{L}_n \times \mathcal{M}_n) + (\mathcal{L}_n \times \mathcal{N}_n) \to \mathcal{L}_n \times (\mathcal{M}_n + \mathcal{N}_n)$$

is a bijection between sets, because each tuple $(l, r)$ can be uniquely assigned to one of the summands on the right hand side depending on whether $r \in \mathcal{M}_n$ or $r \in \mathcal{N}_n$. Therefore we have a minion homomorphism which is bijective in every degree, and thus an isomorphism of minions by Proposition 2.3.2.

*Statement 2:* The forward direction follows from first tupling the injections $\mathcal{L} \hookrightarrow \mathcal{L} + \mathcal{M}$ and $\mathcal{L} \hookrightarrow \mathcal{L} + \mathcal{N}$ to a minion homomorphism $\mathcal{L} \to (\mathcal{L} + \mathcal{M}) \times (\mathcal{L} + \mathcal{N})$, and for the second summand combining the projection-injection composites $\mathcal{M} \times \mathcal{N} \to \mathcal{L} + \mathcal{M}$ and $\mathcal{M} \times \mathcal{N} \to \mathcal{L} + \mathcal{N}$ to a homomorphism $\mathcal{M} \times \mathcal{N} \to (\mathcal{L} + \mathcal{M}) \times (\mathcal{L} + \mathcal{N})$. Cotupling these gives a desired homomorphism $\mathcal{L} + (\mathcal{M} \times \mathcal{N}) \to (\mathcal{L} + \mathcal{M}) \times (\mathcal{L} + \mathcal{N})$.

For the backward direction, note that we know $(\mathcal{L} + \mathcal{M}) \times (\mathcal{L} + \mathcal{N}) \cong \mathcal{L}^2 + \mathcal{M} \times \mathcal{L} + \mathcal{L} \times \mathcal{N} + \mathcal{M} \times \mathcal{N}$ by the distributivity we established earlier. The first three of these summands have a homomorphism to $\mathcal{L}$, and the last one is the second summand. Thus all of the summands map to $\mathcal{L}$ or $\mathcal{M} \times \mathcal{N}$, hence into their sum $\mathcal{L} + \mathcal{M} \times \mathcal{N}$, as desired. □

## 5.5 An embedding of $(\mathcal{P}(\mathbb{N}), \subseteq)$

While minions are a relatively recent notion, *(function) clones*, i.e., composition-closed subsets of $\mathcal{O}(A)$ which contain the identity, have undergone extensive research. A well-known fact is that there are countably many sub-clones of $\mathcal{O}(2)$, i.e., "boolean clones". In fact, they have been completely classified by Post [46], and the poset of all boolean clones modulo containment is nowadays known as "Post's lattice". However, there are much more *subminions*, as being closed under minor operations only is a weaker demand. Indeed, sparks demonstrated [36] that there even on two elements, we have uncountably many minions.

A drawback of this result is that it only makes statements about the number of *distinct* subminions, while for the purposes of this thesis we would be interested in the number of subminions *modulo (minion-)homomorphic equivalence*. Indeed, as demonstrated by Bodirsky and Vucaj [47], when collapsing Post's lattice modulo minion homomorphisms, the resulting poset is substantially smaller[46].

In an unpublished manuscript [48], Kazda and Moore demonstrated how to embed the poset of subsets of natural numbers into the lattice $\mathfrak{F}_2$. We will present a slightly adapted version of their proof which is more in line with the tools we presented so far, however the essence of the argument remains unchanged. Namely, they consider minions generated by sets of near-unanimity[47] functions.

**Definition 5.5.1**: For $s \in \mathbb{N}$, let $n_s \colon 2^s \to 2$ denote the function returning 1 if and only if at most one of the inputs is zero, i.e.,

$$n_s(\underline{a}) := [|\text{supp } \underline{a}| \geq s - 1]$$

using the Iverson bracket notation[48].

Note that all of these are idempotent, i.e., $n_s(0, ..., 0) = 0$ and $n_s(1, ..., 1) = 1$, and, more importantly, the definition is symmetric. Furthermore, if $s \geq 3$, they satisfy the *$s$-quasi-near-unanimity condition* ("$s$-QNU")

$$n_s(x_0, x_1, ..., x_1) \approx n_s(x_1, x_0, x_1, ..., x_1) \approx ... \approx n_s(x_1, ..., x_1, x_0) \approx n_s(x_1, ..., x_1),$$

which is clear for $(x_0, x_1) \in \{(0,0), (1,1)\}$ by idempotence, and for $(x_0, x_1) = (0, 1)$ this is quite literally the definition. However, for $(x_0, x_1) = (1, 0)$, the statements evaluate to

$$n_s(1, 0, ..., 0) = n_s(0, 1, 0, ..., 0) = ... = n_s(0, ..., 0, 1) = n_s(0, ..., 0),$$

which is certainly false for $s = 2$, as $f(0, 1) = 1$, but $f(0, 0) = 0$. Note that as a minor condition, this means

$$n_s(0, 1, ..., 1)_2 = n_s(1, 0, 1, ..., 1)_2 = ... = n_s(1, ..., 1, 0)_2 = n_s(1, ..., 1)_2,$$

which we can denote more concisely by stating $n_s \chi^s_s = n_s \chi^s_{s \setminus i}$ for all $i \in s$, where $\chi^s_T \colon s \to 2$ denotes the characteristic function of a set $T \subseteq s$.

Note furthermore that when we take a proper (i.e., surjective, non-bijective) minor $n_s \alpha$ like $n_4(0, 1, 2, 2)_3$, i.e., the function mapping $(x_0, x_1, x_2)$ to $n_4(x_0, x_1, x_2, x_2)$, the 3-QNU will not be satisfied: Indeed, $n_4\alpha(1, 1, 0) = n_4(1, 1, 0, 0) = 0$, whereas $n_4\alpha(1, 1, 1) = n_4(1, 1, 1, 1) = 1$.

Let us generalize this statement properly.

[46] …albeit still countable, as witnessed by two infinite descending chains.

[47] A a polymorphism $f$ is called a *quasi-near unanimity* polymorphism if and only if it satisfies the identity $f(y, x, ..., x) \approx f(x, y, x, ..., x) \approx ... \approx f(x, ..., x, y) \approx f(x, ..., x)$. This can be formulated as a minor condition for elements $f$ of abstract minions. A *near-unanimity* polymorphism additionally satisfies $f(x, ..., x) \approx x$.

[48] $[P] = 1$ if and only if $P$ is true.

**Lemma 5.5.1**: Let $s, t \geq 3$ and $\alpha\colon s \to t$ be a non-injective operation. Then $n_s\alpha$ does not satisfy the condition $t$-QNU.

*Proof*: Let $j \in s$ be a coordinate on which $\alpha$ is not injective, i.e., $|\alpha^{-1}(j)| > 1$. We claim that $n_s\alpha\chi_t^t \neq n_s\alpha\chi_{t\setminus j}^t$. We first note that $\alpha\chi_t^t = \chi_{\alpha^{-1}(t)}^s = \chi_s^s$, and $\alpha\chi_{t\setminus j}^t = \chi_{\alpha^{-1}(t\setminus j)}^s = \chi_{s\setminus\alpha^{-1}(j)}^s$. Evaluating at the identity tuple $\mathrm{id}_2 = (0,1)$, we see that $\chi_s$ is a tuple with only ones, whereas $\chi_{s\setminus\alpha^{-1}(j)}$ is a tuple with $|\alpha^{-1}(j)|$ many zeroes. Therefore, $n_s\alpha\chi_t(\mathrm{id}_2) = 1$, whereas $n_s\alpha\chi_{t\setminus j}(\mathrm{id}_2) = 0$. □

Another important observation is the following, illustrating the “rigidity” of symmetric elements.

**Lemma 5.5.2**: A symmetric element $f \in \mathcal{M}_n$ of arity $n \geq 2$ is either nondegenerate or constant.

*Proof*: If every coordinate is essential, we are done. If $i \in \operatorname{Iness} f$, then for every $j \in s$ there exists a permutation $\sigma \in S_n$ such that $\sigma(i) = j$. Since injections preserve inessentiality (Lemma 2.6.3), we have $j = \sigma(i) \in \operatorname{Iness} f\sigma = \operatorname{Iness} f$. Thus, all coordinates are inessential, hence $f$ is a constant. □

It is clear that this argument generalizes to the situation where $f$ is fixed by any transitive subgroup $G \leq S_{\mathrm{ar}\, f}$, for instance if $f$ is *cyclic*.

The approach of Kazda and Moore is now to consider for each subset $S \subseteq \mathbb{N}$ the minion generated by $\{n_s\}_{s\in S}$. The previous observations will help us in seeing that for different sets $S$ and $T$, these minions are already homomorphically *inequivalent*.

Recall that an embedding of posets is a map which preserves not only the order relation $\leq$, but also its negation.

**Theorem 5.5.3** ([48, Thm. 5.2]): The assignment

$$S \mapsto \mathcal{N}_S := \langle n_s \mid s \in S \rangle$$

induces an embedding of posets (“order embedding”) $(\mathcal{P}(\{3,4,\ldots\}), \subseteq) \hookrightarrow \mathfrak{F}_2$.

*Proof*: The claim amounts to showing that $\mathcal{M}_S \to \mathcal{M}_T$ if and only if $S \subseteq T$. The backward direction is clear, because $S \subseteq T$ implies the inclusion $\mathcal{M}_S \subseteq \mathcal{M}_T$.

For the forward direction, let $h\colon \mathcal{M}_S \to \mathcal{M}_T$ and $s \in S$. Recalling that the generated minion $\mathcal{M}_T$ contains in arity $k$ only elements of the form $n_t\alpha$ where $t \in T$ and $\alpha\colon t \to k$. Therefore, $h(n_s) = n_t\alpha$ for some $t \in T$ and $\alpha\colon t \to s$.

*Claim*: $\alpha$ must be surjective. Indeed, if $j \notin \operatorname{Im}\alpha$, then $j$ would be an inessential coordinate. However, since $n_s$ is symmetric, so must be $h(n_s)$, so the existence of a inessential coordinate implies by Lemma 5.5.2 that it must be constant. However, $n_t\alpha$ can never be a constant, because the unary minor $n_t(0,\ldots,0)_1$ of $n_t$ is the identity, which is not a constant.

*Claim*: $\alpha$ must be injective. If $\alpha$ were not injective, then by Lemma 5.5.1 the minor $n_t\alpha$ would not satisfy $s$-QNU, however $n_s$ and thus $h(n_s)$ does.

Since $\alpha\colon t \to s$ needs to be bijective, we conclude that $s = t \in T$. Therefore, $S \subseteq T$, as desired. □

Having a poset as rich as the power set of natural numbers in $\mathfrak{F}_2$ – and by the canonical embeddings also in $\mathfrak{F}$ and $\mathfrak{LF}$ – we can deduce many properties of these orders.

**Corollary 5.5.4**: Every local property of the power set of natural numbers also holds in $\mathfrak{F}_2, \mathfrak{F}$, for instance:

1. There exists an order embedding of every countable order $(L, \leq)$, so in particular
2. …there exist infinite ascending and descending chains, and
3. …there exist dense intervals.
4. There exists an order embedding of the continuum $(\mathbb{R}, \leq)$.
5. There exists an uncountable antichain, i.e., a set for which distinct elements are incomparable.

All of these facts are well known statements about the lattice of subsets of natural numbers. For a readable exposition of this lattice, the author recommends [49]. Let us summarize the arguments given there:

*Proof sketches*: We shall show that $(P(\mathbb{N}), \subseteq)$ satisfies the specified properties.

1) Any countable order $(L, \leq)$ order-embeds into $(\mathcal{P}(L), \subseteq)$ via $l \mapsto l^{\downarrow}$, and any bijection $L \to \mathbb{N}$ gives rise to an order isomomrphism $(\mathcal{P}(L), \subseteq) \cong (\mathcal{P}(\mathbb{N}), \subseteq)$ by element-wise evaluation. 2) is clear because $(\mathbb{Q}, \leq)$ has infinite ascending and descending chains, e.g. $\{1/n \mid n \in \mathbb{N}_{>0}\}$ and its negation. 3) is clear because order embeddings map intervals which are not gaps to intervals which are not gaps. For 4), note that $(\mathcal{P}(\mathbb{N}), \subseteq)$ is complete, as we can take arbitrary infima (given by intersections) and suprema (given by unions). Therefore, the image of the embedding $(\mathbb{Q}, \leq) \hookrightarrow \mathcal{P}(\mathbb{N})$ can be extended to its completion, which can be explicitly described as the set consisting of all suprema of $L$ where $(L, R)$ is a Dedekind cut of the image of $(\mathbb{Q}, \leq)$ under the embedding.

To see 5) Let $(B, \leq)$ denote the countable, infinite (rooted) binary tree, ordered with its root $0$ as a minimum. We claim that the set $P$ of (maximal) paths through that tree is an antichain in $(\mathcal{P}(B), \subseteq)$. Indeed, if $0 = b_0 < b_1 < \ldots$ and $0 = b'_0 < b'_1 < \ldots$ are distinct, they must start diverging at some index $k$, hence $b_{k+1} \notin \{b'_i\}_i$ and. Thus, they are incomparable elements in $(\mathcal{P}(B), \subseteq)$, so the collection of all paths is an antichain of uncountable cardinality (because each path is encoded by an infinite sequence of bits). Finally, a vertex enumeration $B \to \mathbb{N}$ witnesses that $(\mathcal{P}(\mathbb{N}), \subseteq)$ has such an antichain as well. □

This tells us that there are some regions in which this lattice is dense. But does this embedding already provide us with some knowledge of *gaps*? Unfortunately not, because if $h\colon (P, \leq) \to (Q, \leq)$ is an embedding of posets and $[p, p']$ is a gap in $P$, then $[h(p), h(p')]$ need not be a gap in $Q$, because intermediate elements might be absent from $\operatorname{Im} h$. We will give gaps some more detailed attention in Chapter 7, specifically Section 7.6.1.

In the meantime, it is both perfectly conceivable that the interval spanned by the atom $[\mathcal{P}]$ and the coatom $[\mathcal{I}_2]$ is dense, but also that this embedding is more or less the only dense "region" in $\mathfrak{F}_2$ or $\mathfrak{F}$.

# 6 Growth

As mentioned in Section 2.5, due to the multi-sorted nature of minions, we do not have a notion of "cardinality"; rather, for a given minion $\mathcal{M}$, we can talk about the cardinality of $\mathcal{M}_n$ for each positive number $n$. Since there are injections $\mathcal{M}_n \hookrightarrow \mathcal{M}_{n+1}$, this cardinality sequence is monotonically increasing. It seems therefore appropriate to call this the *growth sequence* associated to a minion.

The central questions of this chapter are:
1. What information can we recontsruct from its minion by its growth sequence?
2. For which sequences of natural numbers is there a minion whith that sequence as growth sequence?

Naturally, we will assume the minions studied in this chapter to be locally finite, for otherwise the growth sequence will not retain much information.

While the observations in this chapter are not deep, one should note that question 2 has been thoroughly investigated by Randall Dougherty in [50]. He obtained a complete characterization of when a growth sequence $\alpha\colon \omega \to \omega$ is "realized" by a functor $\underline{\mathsf{FinSet}} \to \underline{\mathsf{FinSet}}$ (i.e., he did not only investigate minions, but more generally endofunctors of the category of finite sets *including the emptyset*), and for that used comprehensive combinatorical insights and case distinctions. One corollary of his results is a "large-small"-dichotomy, i.e., a precise way to say whether a minion is large or small. We will extract the core argument and present the corollary on its own in Corollary 6.2.3. Luckily, the language of minions allows us to get a little more intuition for what the involved quantities mean. The interested reader is encouraged to consult the referenced paper [50] in detail.

We shall begin with the basic definitions and a few examples to illustrate some growth sequences of minions, as summarized in Table 1.

For the below considerations, recall that sums and products of minions are reflected in sums and products of their arity-$n$ cardinalities, because these constructions are defined arity-wise (Proposition 2.3.1). Also, the term "cardinality" is to be understood to be dependent on the arity $n$.

*Examples*:
- The constant minion $*$ has one element in each arity, namely the constant one.
- The projection minion $\mathcal{P}$ has $n$ elements of arity $n$, because it is modelled by the identity functor.
- More generally, the free minion $\langle s\colon k\rangle$ with one $k$-ary operation $s$ has $n^k$ elements in arity $n$ because it is isomorphic to $\mathcal{P}^n$.
- $\mathrm{Pol}(K_3) \cong \mathcal{P} + \mathcal{P} + \mathcal{P}$ has cardinality $3n$.
- The powerset minion $\mathcal{P}(_) := (n \mapsto \mathcal{P}(n))$, where an operation $\alpha\colon n \to k$ assigns a subset $S \subseteq n$ its image $\alpha[S]$, has $2^n$ elements of arity $n$. Note that this minion contains a constant in each arity, namely the empty set, hence is homomorphically equivalent to the constant minion.
- Its subminion $\mathcal{P}_{>0}$ of all nonempty sets has arity $2^n - 1$ (and no constants anymore).
- For a finite field $\mathbb{K}$, the minion $\mathrm{Pol}_{\underline{\mathbb{K}-\mathsf{Vec}}}(\mathbb{K})$ of polymorphisms of $\mathbb{K}$ viewed as a $\mathbb{K}$-vector space has precisely $|\mathbb{K}^n \to_{\underline{\mathbb{K}-\mathsf{Vec}}} \mathbb{K}| = |\mathbb{K}^n|$ elements of arity $n$.
- The double-powerset minion $\mathcal{P}(\mathcal{P}(_))$ has cardinality $2^{2^n}$.
- The function minion $\mathcal{O}(A, B)$ has cardinality $|B|^{|A|^n}$.

| $\mathcal{M}$ | $\lvert\mathcal{M}_n\rvert$ |
|---|---|
| $* \cong \mathcal{O}(A,1)$ | $1$ |
| $\langle s\colon k\rangle$ | $n^k$ |
| $\langle s\colon 1\rangle \cong \mathcal{P} \cong \langle \mathrm{id}_A\rangle_{\mathcal{O}(A,A)}$ | $n$ |
| $\mathrm{Pol}(K_3)$ | $3n$ |
| $\langle S\rangle$ | $\Sigma_{k=1}^{\infty} \lvert S_k\rvert \cdot n^k$ |
| $\mathcal{P}(_) \cong * + \mathcal{P}_{>0}(_)$ | $2^n$ |
| $\mathrm{Pol}_{\mathbb{K}-\underline{\mathsf{Vec}}}(\mathbb{K})$ | $\lvert\mathbb{K}\rvert^n$ |
| $\mathcal{P}\mathcal{P}(_)$ | $2^{2^n}$ |
| $\mathcal{O}(A,B)$ | $\lvert B\rvert^{\lvert A\rvert^n}$ |

Table 1: Examples of minions and their growth sequences.

## 6.1 Growth and essential growth

Recalling the definition of essential coordinates (Definition 2.6.2), we can introduce an important variation in counting the arity-$n$ elements of a minion: Namely, we would like to omit elements with inessential coordinates, as these can be understood to be "shadows" of lower arity elements, to which we applied an injective minor operation. This will be the key combinatorial insight to understand the growth of a minion.

**Definition 6.1.1**: Let $\mathcal{M}$ be a locally finite minion. We define functions $\alpha_{\mathcal{M}}, \gamma_{\mathcal{M}}\colon \omega \to \omega$ called *growth sequence* and *essential growth sequence* by

$$\alpha_{\mathcal{M}}(n) := \begin{cases} \lvert\mathcal{M}_n\rvert & n \geq 1 \\ \lvert\{c \in \mathcal{M}_1 \mid c \text{ constant}\}\rvert & n = 0 \end{cases},$$

$$\gamma_{\mathcal{M}}(n) := \begin{cases} \lvert\{f \in \mathcal{M}_n \mid \mathrm{Iness}\, f = \emptyset\}\rvert & n \geq 1 \\ \lvert\{c \in \mathcal{M}_1 \mid c \text{ constant}\}\rvert & n = 0 \end{cases}.$$

If clear from the context, we will denote this function just by $\alpha$.

**Proposition 6.1.1** (cf. [50, Lem. 2.6]): Let $\mathcal{M}$ be a locally finite minion. Then we have

$$\alpha_{\mathcal{M}}(n) = \sum_{k=0}^{n} \binom{n}{k} \gamma_{\mathcal{M}}(k).$$

*Proof*: The case $n = 0$ is clear, so assume $n \geq 1$. We have

$$\begin{aligned}\mathcal{M}_n &= \sum_{S\subseteq n}\{f\in\mathcal{M}_n \mid \operatorname{Ess} f = S\}\\ &= \sum_{k=0}^{n}\sum_{S\in\binom{n}{k}}\{f\in\mathcal{M}_n \mid \operatorname{Ess} f = S\}\\ &\overset{(!)}{\cong} \sum_{k=0}^{n}\sum_{S\in\binom{n}{k}}\underbrace{\{f\in\mathcal{M}_n \mid \operatorname{Ess} f = k\}}_{\cong\gamma_{\mathcal{M}}(k)\quad(*)}\\ &\cong \sum_{k=0}^{n}\binom{n}{k}\gamma_{\mathcal{M}}(k).\end{aligned}$$

Here, we have used in (!) that a bijection $\pi$ – in his case one which transports elements in $S\subseteq n$ to elements in $k=\{0,\ldots,k-1\}\subseteq n$ – commutes with the essentiality of coordinates. Furthermore, in $(*)$ we have used that 1. if $k=0$, then $\{f\in\mathcal{M}_n \mid \operatorname{Ess} f=\emptyset\}$ are all the arity-$n$ constants, which is in bijection to the unary constants (Corollary 2.6.8), i.e., $\gamma_{\mathcal{M}}(0)$, and 2. if $k\geq 1$, then the elements $f\in\mathcal{M}_n$ with $\operatorname{Ess} f = k$ must all (uniquely) factor through the inclusion $\iota\colon k\hookrightarrow n$; furthermore, by Lemma 2.6.3, we know that $\operatorname{Ess}(\tilde{f}\iota) = \iota[\operatorname{Ess}\tilde{f}]$, so the elements $f\in\mathcal{M}_n$ with $\operatorname{Ess} f = k$ cannot arise from those $\tilde{f}\in\mathcal{M}_k$ which posess an inessential coordinate. Thus, in that case the set is in bijection to $\gamma_{\mathcal{M}}(k)$ also. □

This rather simple argument allows us to calculate the numbers of elements purely by counting its nondegenerate elements, i.e., elements without inessential coordinates. Let us consider two examples of quotients of $\mathcal{P}^2\cong\langle f\colon 2\rangle$ to illustrate this.

*Example (1)*: Let $\mathcal{C}_2 := \langle f\colon 2 \mid f(0,1)\approx f\rangle$ be the most general minion with one binary symmetric operation, i.e., the gadget associated to the minor condition $f(x,y)=f(y,x)$. This minion has no constants, one nondegenerate unary element induced by our generator, and one binary one – in particular, the only potential second binary nondegenerate element would be the "swapped sibling" $f\sigma$ of $f$, which in this case is $f$. All higher-arity elements are degenerate as they are minors of the arity-2 generator. To summarize, $\gamma = (0,1,1,0,0,\ldots)$, from which we can deduce

$$\alpha(n) = \binom{n}{1}\cdot 1 + \binom{n}{2}\cdot 1 = n + \frac{n(n-1)}{2} = \frac{1}{2}n + \frac{1}{2}n^2.$$

This example shows that $\alpha$ can be a polynomial with non-integer coefficients.

*Example (2)*: Let $\mathcal{M} := \langle f\colon 2 \mid f(0,0)\approx f(1,1)\rangle$ be the most general minion with a (non-symmetric) binary operation which is constant on the diagonal, i.e., the gadget associated to the minor condition $f(x,x)=f(y,y)$. There, we have one constant, no nondegenerate unary, but two nondegenerate binary elements, because $f$ is not symmetric. In summary, $\gamma=(1,0,2,0,0,\ldots)$, which implies

$$\alpha(n) = \binom{n}{0}\cdot 1 + \binom{n}{2}\cdot 2 = 1 + n(n-1) = 1 - n + n^2.$$

*Example (3)*: Adding a copy of $\mathcal{P}$ to the previous example, i.e., considering the minion $\mathcal{P}+\mathcal{M}$, we get a minion with $\alpha(n)=1+n^2$, which is the same as the growth of the minion $*+\mathcal{P}$, to which $\mathcal{P}+\mathcal{M}$ is not isomorphic.

Using the well-known combinatorial tool of binomial inversion, we can derive a way to calculate $\alpha$ from $\gamma$:

**Proposition 6.1.2**: For functions $\alpha\colon \omega \to \omega$ and $\gamma\colon \omega \to \omega$, the following are equivalent:

$$\alpha(n) = \sum_{k=0}^{n} \binom{n}{k} \gamma(k)$$

$$\gamma(k) = \sum_{j=0}^{k} (-1)^{j-k} \binom{k}{j} \alpha(k).$$

For a proof of this statement see e.g. [51, Theorem 9.8] or [52, Section 3.6].

There is a slightly more compact way to write this inverse:

**Proposition 6.1.3**: We have

$$\sum_{j=0}^{k} (-1)^{j-k} \binom{k}{j} \alpha(k) = (\Delta^k \alpha)(0),$$

where $\Delta^k$ is the $k$-fold iteration of the *forward difference operator*

$$\Delta\colon \mathbb{Z}^\omega \to \mathbb{Z}^\omega, \quad (\Delta\alpha)(n) := \alpha(n+1) - \alpha(n).$$

*Proof*: Note that as a ($\mathbb{Z}$-)linear operator on $\mathbb{Z}^\omega$, we have $\Delta = L - \mathrm{id}$, where $L$ is the left shift operator $(L\alpha)(n) := \alpha(n+1)$. Since $L$ commutes with id, i.e., $L \circ \mathrm{id} = \mathrm{id} \circ L$, the subring of the endomorphism ring of $\mathbb{Z}^\omega$ spanned by $L$ and id is commutative, so we can apply binomial expansion to obtain

$$\Delta^k = (L - \mathrm{id})^k = \sum_{j=0}^{k} \binom{k}{j} L^j (-\mathrm{id})^{k-j} = \sum_{j=0}^{k} (-1)^{k-j} \binom{k}{j} L^j.$$

This implies that

$$(\Delta^k \alpha)(0) = \sum_{j=0}^{k} (-1)^{k-j} \binom{k}{j} (L^j \alpha)(0) = \sum_{j=0}^{k} (-1)^{k-j} \binom{k}{j} \alpha(j),$$

as desired. □

This implies:

**Corollary 6.1.4**:

$$\gamma(k) = \sum_{j=0}^{k} (-1)^{k-j} \binom{k}{j} \alpha(j) = (\Delta^k \alpha)(0).$$

This allows us to immediately know about the number of nondegenerates in each degree in the cases where it is easy to determine how many elements the minion has in general. For instance:

*Example*: Let $\mathcal{M} = \mathcal{O}(A, B)$ for some finite sets $A, B$. Then $\alpha(n) = B^{A^n}$, and hence

$$\gamma(k) = \sum_{j=0}^{k} (-1)^{k-j} \binom{k}{j} B^{A^j}.$$

## 6.2 The dichotomy

**Proposition 6.2.1** (cf. [50, Prop.1.4]): Let $\mathcal{M}$ be a locally finite minion which is not finitely generated[49]. Then $\limsup_n \sqrt[n]{\alpha(n)} \geq 2$.

*Proof*: Since $\mathcal{M}$ must have generators in arbitrarily large degrees, we know that $\gamma(m) > 0$ for arbitrarily large $m$. This has ramifications in arity $2m$: Indeed,

$$\alpha(2m) = \sum_k \binom{2m}{k} \gamma(k) \geq \binom{2m}{m},$$

and since $\lim_m \sqrt[2m]{\binom{2m}{m}} = 2$ (this follows e.g. from Stirling's formula), we conclude that 2 must be a lower bound for $\limsup_n \alpha(n)$. □

**Proposition 6.2.2**: Finitely generated minions have growth $\alpha$ bounded by a polynomial function of degree $N$, where $N$ is the maximum arity of a set of generators.

*Proof*: Let $S = (S_n)_n$ be a set of generators for $S$. Then $\mathcal{M}$ is a quotient of the free minion $\langle S \rangle = (\{s\alpha \mid k > 0, s \in S_k, \alpha\colon n \to k\})_n$, which is of the form

$$\sum_{k>0} |S_k| \; P^k,$$

so for each $n$, there is a surjection from $\left(\sum_{k>0} |S_k| \; P^k\right)_n = \sum_{k>0} |S_k| \; n^k$, which is a polynomial function in $n$ of arity $\max\{n \mid S_n > \emptyset\}$, as desried. □

For the main classification, define the following notions:

**Definition 6.2.1**: A sequence $\alpha\colon \omega \to \omega$ is called *subexponential* if $\limsup_n \sqrt[n]{\alpha(n)} \leq 1$.

One easily observes that (finite) sums and products of subexponential functions are subexponential, and since $\limsup_n \sqrt[n]{c} = 1$ for all $c \geq 1$ and $\limsup_n \sqrt[n]{n} = 1$, we deduce that all polynomial functions are subexponential, as the naming suggests.

[49] Equivalently: Not of bounded essential arity (Proposition 2.6.11).

**Corollary 6.2.3**: For a locally finite minion $\mathcal{M}$ with growth sequence $\alpha(n)$, precisely one of the following holds:

1. $\alpha(n)$ is subexponential, i.e., $\limsup_n \sqrt[n]{\alpha(n)} \leq 1$ or
2. $\limsup_n \sqrt[n]{\alpha(n)} \geq 2$

Furthermore, statement one is fulfilled if and only if the minion is finitely generated.

*Proof*: If the minion is finitely generated, then from Proposition 6.2.2 we know that $\alpha(n)$ is bounded above by a polynomial, and by the above remarks we know that polynomials $p(n)$ satisfy $\limsup_n \sqrt[n]{p(n)} = 1$, so condition 1 holds. If the minion is not finitely generated, then Proposition 6.2.1 shows we satisfy condition 2. □

Note that every polynomial growth function $\alpha(n) = \sum_k a_k n^k$ with nonnegative integer coefficients is the growth function of a (finitely generated) free minion, namely $\langle A \rangle$ where $A_k = a_k$.

## 6.3 Bibliographic remarks

On first sight, the notion of a minion seems rather recent: The word was coined in the 2018 paper [19] in the context of a "function minion", further weakening the more restrictive notion of a clonoid, which is a function minion closed under the action of an algebra over the codomain. Predating the word, the idea that only height-one identities are relevant had been established in 2015 in [18].

However, in the abstract context, a minion is just a functor from the category $\mathbb{F}_{>0}$ of finite nonzero ordinals to the category of sets. These objects – or rather, objects quite related to it, namely endofunctors on subcategories of $\underline{\mathsf{Set}}$ – have been studied extensively since the early days of category theory. Unsurprisingly, the material presented so far has substantial overlap with preexisting results.

Unfortunately, due to the mismatch in approach, notation and language, the connections are a bit obscured; however, if we pay close attention to the existing literature, we can see that the new perspective can enrich our understanding of the – previously quite technical – results.

### 6.3.1 From minions to endofunctors

As previously mentioned in Section 3.1, we can use the "named" perspective to extend a minion from a functor

$$\begin{aligned} \mathbb{F}_{>0} &\to \underline{\mathsf{Set}}, \\ n &\mapsto \mathcal{M}_n \end{aligned}$$

to a functor

$$\begin{aligned} \underline{\mathsf{FinSet}} &\to \underline{\mathsf{Set}}, \\ X &\mapsto \mathcal{M} \otimes X, \end{aligned}$$

which in particular sends the empty set to the empty set. In other words, a minion viewed this way has no "nullary" elements. Since we do not depend on the demand that $X$ is a finite set, we can analogously define extend a minion to a functor $\underline{\mathsf{Set}} \to \underline{\mathsf{Set}}$, i.e., an *endofunctor* on the category $\underline{\mathsf{Set}}$.

> *Remark*: This construction can be characterized by a universal property, and is known in category theory as a *Kan extension* [5, Chapter 6]: Specifically, extending a minion $\mathcal{M}$ this way constitutes the left Kan extension of $\mathcal{M}\colon \mathbb{F}_{>0} \to \underline{\mathsf{Set}}$ along the inclusion of categories $\iota\colon \mathbb{F}_{>0} \hookrightarrow \underline{\mathsf{Set}}$, written $\mathrm{Lan}_\iota\, \mathcal{M}$.
>
> While such an extension does not need to exist in general, the left Kan extension of $F\colon \mathcal{C} \to \mathcal{D}$ along any other functor $\mathcal{C}$ exists whenever $\mathcal{C}$ is small – i.e., its class of objects is a set – and $\mathcal{D}$ is cocomplete – i.e., it has all colimits [5, Thm. 6.2.1]. Roughly speaking, cocompleteness ensures that a "quotient of sums"-construction always exists, which of course is the case in the category of sets.

There is a second way to obtain an endofunctor of a minion: If the minion is locally finite, we can view it as a functor

$$\mathbb{F}_{>0} \to \underline{\mathsf{FinSet}},$$

and stop our extension at finite sets. This yields an endofunctor on $\underline{\mathsf{FinSet}}$.

More generally, if every arity $\mathcal{M}_n$ of $\mathcal{M}$ is bounded by an ordinal $\alpha$, then we can extend to an endofunctor of $\underline{\mathsf{Set}}_{<\alpha}$, the category of sets smaller than $\alpha$.

Modulo a technicality relating to the empty set, $\underline{\mathsf{FinSet}}$ can be viewed pretty much exactly as locally finite minions. More precisely: The category $\underline{\mathsf{LFin}}$ of locally finite minions and the category of $\underline{\mathsf{FinSet}}$-endofunctors fixing the empty set are equivalent[50] (cf. Proposition A.2.1).

> *Remark*: If such a functor does not fix the empty set, we can think of it as a minion with "designated nullaries": The only new minor operations to consider are the empty functions $()_n\colon 0 \to n$, which lift a nullary to an $n$-ary element – necessarily constant. However, the lack of functions $n \to 0$ implies that a constant need not arise from a nullary. We can thus view nullary elements, if they exist, as arbitrary labels for constants.

However, there are generally "more" endofunctors on $\underline{\mathsf{Set}}$ than minions, i.e., the aforementioned extension $\mathrm{Lan}_\iota\colon\ \underline{\mathsf{Min}} \hookrightarrow [\underline{\mathsf{Set}}, \underline{\mathsf{Set}}]$ is not an equivalence. The reason is that minions as we defined them do not have any proper infinitary operations: The "arity $\omega$"-set $\mathcal{M} \otimes \omega$ of the extended minion only contains equivalence classes of pairs $[f, \eta]$, where $f$ is of some finite arity and $\eta$ without loss of generality an embedding. Extending the notion of an essential coordinate to this new arity $\omega$ in the natural way, we see that $\mathrm{Ess}[f, \eta]$ must be finite.

We can describe this phenomenon by calling a $\underline{\mathsf{Set}}$-endofunctor $F$ *finitary* if for all sets $X$, the set $F_X$ consists exclusively of minors of finite-arity sets $F_Y$, i.e., for every $f \in F_X$ there exists a finite set $Y$, a $g \in F_Y$, and an operation $\eta\colon Y \to X$ such that $f = g\eta$ [53]. Replacing finiteness by another strict upper cardinality bound $\alpha$ gives rise to the notion of *$\alpha$-accessible* endofunctors of $\underline{\mathsf{Set}}$ [54]. Finitary ones are thus called *$\omega$-accessible*.

### 6.3.2 Research on endofunctors

Starting around 1970, V. Trnková and V. Koubek initiated a body of research around properties of endofunctors on $\underline{\mathsf{Set}}$ – hereafter just called $\underline{\mathsf{Set}}$ *(endo)functors* – and more generally $\mathcal{C}$-endofunctors for some subcategories $\mathcal{C}$ of $\underline{\mathsf{Set}}$ [55], [56], [57], [58]. When do endofunctors preserve certain properties or constructions? Does every such functor have a least or greatest fixed point? Under which cardinality constraints – i.e., conditions on the function $|X| \mapsto |F_X|$ associated to such a

[50] A functor $F\colon \mathcal{C} \to \mathcal{D}$ is called an equivalence if it has an inverse functor $G\colon \mathcal{D} \to \mathcal{C}$ up to natural isomorphism. Categories $\mathcal{C}$ and $\mathcal{D}$ are called *equivalent* if an equivalence exists.

functor – is it guaranteed to exist? Is it unique? Do weakly terminal objects exist? These are just few of the questions considered in this domain. A good entrypoint is a fairly recent paper of the aforementioned authors [59], which also has some comments on the historical development.

In [55, Prop. 2.1], V. Trnková established, among many considerations, the following key result: If $A_1, A_2 \subseteq X$ have nonempty intersection, then any $f \in F_X$ which is simultaneously induced by an element from $F_{A_1}$ and an element from $F_{A_2}$ (in the sense that they are minors of the respective subset inclusions) must already be induced by an element from $F_{A_1 \cap A_2}$. Thus, keeping in mind technical difficulties with constants, to any element $f \in F_X$ we can associate a *filter*

$$\mathcal{F}_X(f) := \begin{cases} \{A \subseteq X \mid \exists g \in F_A, \alpha\colon A \to X\colon f = g\alpha\} & \text{if } f \text{ is not constant} \\ \mathcal{P}(X) & \text{otherwise.} \end{cases}$$

In the case where $F$ is the finitary endofunctor arising from a minion, this is precisely the filter

$$(\operatorname{Ess} f)^{\uparrow} = \{A \subseteq X \mid \operatorname{Ess} f \subseteq A\}.$$

This notion has been further developed in [56], [57], [58].

Both L. Barto [45] and R. Dougherty [50] presented some descriptive results about the behavior of these filters: Roughly speaking, all the compatible ways to assign an element $f \in F_X$ and its proper[51] minors $f\alpha$ a filter on $\operatorname{cod}\alpha$ have been dubbed "($\kappa$-) filter structures" in [45] and "control structures" in [50], where the focus was put on endofunctors on *finite* sets, i.e., locally finite minions.

Dougherty noted that assigning each element $f$ the family $(\operatorname{Ess} f\alpha)_\alpha$ ranging over surjective minor operations $\alpha$ gives rise to a natural transformation – i.e., minion homomorphism – to the functor $Q$ of all such "control structures", and that applying this natural transformation to $Q$ itself yields the identity [50, Prop. 3.3]. He used these insights to derive combinatorial constraints on what we called the *growth function* of a FinSet-endofunctor. These were quite elaborate, and for instance required a nontrivial fact regarding primitive groups on 8 elements, whose proof alone spans four pages [50, Section 5]. Even phrasing the main statement is a bit involved, so the author of this text opted to prove a corollary [50, Prop. 1.4] in isolation.

L. Barto did a more detailed analysis of these structures in the more general setting of Set -endofunctors. For that, he introduced $\kappa$-filter structures, $\kappa$ being a cardinal. The definition of the $\kappa$-filter structure of an element $f \in F_X$ is analogous to that of a control structure, however instead of considering the filters of arbitrary proper minors $f\alpha$, it considers only those minors with $|\operatorname{Im}\alpha| < \kappa$. The functor containing all abstract $\kappa$-filter structures is called $T^\kappa$. He then notices that omitting the constant component of $T^3$ gives us a functor $W$ which is a *rigid core* and *weakly terminal* in the category of all faithful connected set functors [45, Thm. 3.8], which correspond to locally finite minions with a single, non-constant unary. This demonstrates that $\mathfrak{LF}$ has a coatom, whose core representative we constructed in Section 5.3.4 using a more "top-down" approach.

In a most recent[52] publication [43], L. Barto and M. Kapytka analyzed multisorted boolean clones, and also introduced the notion of a minion core [43, Section 3]. They prove that function minions on finite sets have a core [43, Thm. 3.4] and derive a statement analogous to Corollary 5.3.16 [43, Lem. 3.5].

[51] In the sense that $\alpha$ is a surjection

[52] A preprint of which appeared after submission of this thesis

# 7 The exponential minion: Dualities and gaps

In this section we will construct an abstract minion $\mathcal{N}^{\mathcal{M}}$ called the "exponential", characterized by the property that there is a (natural) one-to-one correspondence between $\mathcal{L} \times \mathcal{M} \to \mathcal{N}$ and $\mathcal{L} \to \mathcal{N}^{\mathcal{M}}$.

While the existence of an exponential in the abstract setting is not surprising and well known, it is less obvious that the exponential of finitely representable minions remains to be finitely representable, and that the exponential of finitely generated minions remains to be finitely generated. These are important properties, which ensure that the respective homomorphism orders form a *Heyting algebra*.

## 7.1 Exponential objects

We start with an abstract definition, which puts the construction of the exponential minion in the right theoretical context, and justifies the naming. The definition below is well-understood by category theorists, and defines our object uniquely up to isomorphism. We will later give an explicit construction, and show that it satisfies the following axioms of an exponential object.

**Definition 7.1.1**: Let $\mathcal{M}$ be a minion (or, more general, an object in a category with products). The $\mathcal{M}$-fold exponential $_^{\mathcal{M}}$ is defined to be the right adjoint functor $\underline{\mathsf{Min}} \to \underline{\mathsf{Min}}$ to the functor $_ \times \mathcal{M}$ as defined in Definition 1.1.22.

This means the following: The exponential is a functor $\underline{\mathsf{Min}} \to \underline{\mathsf{Min}}$ – i.e., a minion $\mathcal{N}^{\mathcal{M}}$ forevery minion $\mathcal{N}$, together with an induced homomorphism $h_*\colon \mathcal{N}^{\mathcal{M}} \to \mathcal{N}'^{\mathcal{M}}$ for every minion homomorphism $h\colon \mathcal{N} \to \mathcal{N}'$ compatible with identities and composition – such that for every minion $\mathcal{L}$ and $\mathcal{N}$, we have mutually inverse correspondences

$$(\mathcal{L} \times \mathcal{M} \to_{\underline{\mathsf{Min}}} \mathcal{N}) \underset{\Psi_{\mathcal{L},\mathcal{N}}}{\overset{\Phi_{\mathcal{L},\mathcal{N}}}{\rightleftarrows}} (\mathcal{L} \to_{\underline{\mathsf{Min}}} \mathcal{N}^{\mathcal{M}}),$$

Where $\Psi$ – and thus $\Phi$ – is natural in $\mathcal{M}$ and $\mathcal{N}$ in the sense of Definition 1.1.22.

Spelling out the naturality conditions reads as follows:

1. For every $h\colon \mathcal{N} \to \mathcal{N}'$, the (via postcomposition) induced homomorphisms $(\mathcal{L} \times \mathcal{M} \to \mathcal{N}) \to (\mathcal{L} \times \mathcal{M} \to \mathcal{N}')$ and $(\mathcal{L} \to \mathcal{N}^{\mathcal{M}}) \to \left(\mathcal{L} \to \mathcal{N}'^{\mathcal{M}}\right)$ make the following diagrams commute:

$$\begin{array}{ccc} (\mathcal{L} \times \mathcal{M} \to_{\underline{\mathsf{Min}}} \mathcal{N}) & \xrightarrow{\Phi_{\mathcal{L},\mathcal{N}}} & (\mathcal{L} \to_{\underline{\mathsf{Min}}} \mathcal{N}^{\mathcal{M}}) \\ \downarrow & & \downarrow \\ (\mathcal{L} \times \mathcal{M} \to_{\underline{\mathsf{Min}}} \mathcal{N}') & \xrightarrow{\Phi_{\mathcal{L},\mathcal{N}'}} & \left(\mathcal{L} \to_{\underline{\mathsf{Min}}} \mathcal{N}'^{\mathcal{M}}\right) \end{array}$$

2. For every $g\colon \mathcal{L} \to \mathcal{L}'$, the (via precomposition) induced homomorphisms $(\mathcal{L}' \times \mathcal{M} \to \mathcal{N}) \to (\mathcal{L} \times \mathcal{M} \to \mathcal{N})$ and $\left(\mathcal{L}' \to \mathcal{N}^{\mathcal{M}}\right) \to \left(\mathcal{L} \to \mathcal{N}^{\mathcal{M}}\right)$ make the following diagrams commute:

$$\begin{array}{ccc} (\mathcal{L}' \times \mathcal{M} \to_{\underline{\mathsf{Min}}} \mathcal{N}) & \xrightarrow{\Phi_{\mathcal{L}',\mathcal{N}}} & (\mathcal{L}' \to_{\underline{\mathsf{Min}}} \mathcal{N}^{\mathcal{M}}) \\ \downarrow & & \downarrow \\ (\mathcal{L} \times \mathcal{M} \to_{\underline{\mathsf{Min}}} \mathcal{N}) & \xrightarrow{\Phi_{\mathcal{L},\mathcal{N}}} & (\mathcal{L} \to_{\underline{\mathsf{Min}}} \mathcal{N}^{\mathcal{M}}) \end{array}$$

*Example*: If in the above Definition 7.1.1, we replace the category Min of minions with the category Gr of graphs, we recover the usual notion of the exponential graph $H^G$: Indeed, notice that vertices of a graph correspond to homomorphisms from the discrete, one-element graph $*$ – and then we see that

$$V(H^G) \cong (* \to_{\underline{\mathsf{Gr}}} H^G) \cong (* \times G \to_{\underline{\mathsf{Gr}}} H).$$

Since the product of a discrete graph with any graph is discrete, we conclude that every map out of $* \times G$ must be a homomorphism, so $* \times G \to_{\underline{\mathsf{Gr}}} H$ is the set of all functions $G \to_{\underline{\mathsf{Set}}} H$. To find the edges, notice that these are characterized by homomorphisms from $\overrightarrow{P_1}$, the path of length one, into our graph:

$$E(H^G) \cong \left(\overrightarrow{P_1} \to_{\underline{\mathsf{Gr}}} H^G\right) \cong \left(\overrightarrow{P_1} \times G \to_{\underline{\mathsf{Gr}}} H\right),$$

which we can describe as pairs $(f, g)$ of functions from $G$ to $H$ where for every $v \in G$ there is an edge from $f(v)$ to $g(v)$.

*Example*: For a radically simpler example, consider the category Set of sets and functions: There, the exponential object $Y^X$ will be the set of all functions $X \to Y$.

*Example*: Now consider an arbitrary category with arbitrary products – in particular with the *nullary product* $*$, which is characterized by the property that every other object $X$ shall have a unique morphism into it. Such an object is called a *terminal object*. The nomenclature is justified because we have $* \times X \cong X$: Indeed, the projection onto $X$ can be reversed by tupling the identity $X \to X$ with the unique map $X \to *$ to obtain a unique map $X \to * \times X$.[53]. This means that we have natural bijective correspondences

$$(X \to Y) \cong (* \times X \to Y) \cong (* \to Y^X),$$

and so we see that

*morphisms $X \to Y$ correspond to morphisms from the terminal object into the exponential object.*

In our particular instances, this means the following:
- In Set, the terminal object is a one-element "singleton" set. This implies that functions $X \to Y$ correspond to arbitrary elements in the exponential set.
- In Gr however, the terminal object is the graph *with a loop*, so graph homomorphisms $G \to H$ do not correspond to arbitrary elements in $H^G$, but *loops*
- In Min, the terminal minion is the constant minion. So we should expect the minion homomorphisms $\mathcal{M} \to \mathcal{N}$ to correspond to constants in $\mathcal{N}^{\mathcal{M}}$.

[53] The complete argument for why this arrow is an isomorphism is a bit longer, but not hard to do.

## 7.2 The construction

The fact that $\underline{\mathsf{Set}}$-valued functor categories, in category theory known as *(co-)presheaves*, admit exponential objects, is well-known in the area of topos theory [26, Proposition I.6.1]. There is also a more general way to derive that functor categories whose domain category is "cartesian closed" (i.e., admits exponentials), expressing the exponential as a special kind of limit (more specifically, an *end*) [60].

The following is an explicit construction for which we will verify in Proposition 7.2.5 that, in the context of the right correspondence, it indeed behaves like an "exponential object" of minions.

**Definition 7.2.1**: Let $\mathcal{M}$ and $\mathcal{N}$ be minions. The *exponential minion* $\mathcal{N}^{\mathcal{M}}$ is defined in degree $n$ as

$$\left(\mathcal{N}^{\mathcal{M}}\right)_n := (\mathcal{P}^n \times \mathcal{M} \to_{\underline{\mathsf{Min}}} \mathcal{N}),$$

where minor operations act on the exponent, i.e., $\alpha\colon n \to n'$ acts as

$$(\mathcal{P}^n \times \mathcal{M} \to \mathcal{N}) \to \left(\mathcal{P}^{n'} \times \mathcal{M} \to \mathcal{N}\right),$$
$$H \mapsto H(\alpha;_,_).$$

*Proof that the exponential minion is a minion*: We can either see that $\mathcal{N}^{\mathcal{M}}$ is a composition of functors

1. $\mathcal{P}^{\bullet}\colon \mathbb{F}_{>0} \to \underline{\mathsf{Min}}^{\mathsf{op}}$
2. $_ \times \mathcal{M}\colon\ \underline{\mathsf{Min}} \to \underline{\mathsf{Min}}$
3. $_ \to \mathcal{N}\colon \underline{\mathsf{Min}}^{\mathsf{op}} \to \underline{\mathsf{Set}}$,

forming altogether a functor $\mathbb{F}_{>0} \to \underline{\mathsf{Set}}$.

Alternatively, we can see the functoriality explicitly, because clearly $\alpha = \mathrm{id}_n$ acts as the identity, and $\alpha\colon n \to n'$ sends

$$H \mapsto H(\alpha;_,_),$$

so subsequently applying $\beta\colon n' \to n''$ gives $\ldots \mapsto H(\alpha;(\beta;_),_)$, which is the same as acting with the composite $\alpha;\beta$. □

Furthermore, we can use pre- and postcomposition to obtain minion homomorphisms between exponentials.

**Definition 7.2.2** (pushforward, pullback): Let $\mathcal{M}$ and $\mathcal{N}$ be minions and $l\colon \mathcal{M}' \to \mathcal{M}$, $r\colon \mathcal{N} \to \mathcal{N}'$ be minion homomorphisms. the map

$$r_*\colon \mathcal{N}^{\mathcal{M}} \to \mathcal{N}'^{\mathcal{M}}, \quad (H_n)_n \mapsto (r(H_n(_,_)))_n$$

is called the *pushforward* of $r$. the map

$$l^*\colon \mathcal{N}^{\mathcal{M}} \to \mathcal{N}^{\mathcal{M}'}, \quad (H_n)_n \mapsto (H_n(_,l(_)))_n$$

is called the *pullback* of $l$.

**Lemma 7.2.1**: For $r\colon \mathcal{N} \to \mathcal{N}'$ and $l\colon \mathcal{M}' \to \mathcal{M}$ as above, the pushforward $r_*$ and the pullback $l^*$ are minion homomorphisms.

*Proof*: Let $\alpha\colon n \to n'$. For every $H \in \left(\mathcal{N}^{\mathcal{M}}\right)_n = \mathcal{P}^n \times \mathcal{M} \to \mathcal{N}$ and $r, l$ as stated, we observe

$$r_*(H\alpha) = r_*(H(\alpha;_,_)) = r(H(\alpha;_,_)) = (r_*(H))(\alpha;_,_) = (r_*(H))\alpha,$$

as well as

$$l^*(H\alpha) = l^*(H(\alpha;_,_)) = H(\alpha;_,l(_)) = (l^*(H))\alpha.$$

□

The above proof seems like nothing has been happening: Indeed, the minor operations, pushforward, and pullback all act on the first argument, second argument, and on the output, respectively, so all of these "actions" do not interfere.

The exponentiation will play a role later on this thesis: The process of taking $\mathcal{M}$ and considering $\mathcal{N}^{\mathcal{M}}$ can be viewed as some kind of dualization.

**Proposition 7.2.2**: The assignment

$$\mathcal{M}' \xrightarrow{l} \mathcal{M} \quad \rightsquigarrow \quad \mathcal{N}^{\mathcal{M}} \xrightarrow{l^*} \mathcal{N}^{\mathcal{M}'}$$

turns $\mathcal{N}^{\bullet}$, the exponentiation of a fixed base $\mathcal{N}$, into a (contravariant) functor $\underline{\mathsf{Min}}^{\mathsf{op}} \to \underline{\mathsf{Min}}$. The assignment

$$\mathcal{N} \xrightarrow{r} \mathcal{N}' \quad \rightsquigarrow \quad \mathcal{N}^{\mathcal{M}} \xrightarrow{r_*} \mathcal{N}'^{\mathcal{M}}$$

turns $\bullet^{\mathcal{M}}$, exponentiation with a fixed exponent $\mathcal{M}$, into a (covariant) functor $\underline{\mathsf{Min}} \to \underline{\mathsf{Min}}$.

*Proof*: That identities act as identities is clear; to see that composition behaves well, consider first

$$\mathcal{M}'' \xrightarrow{l'} \mathcal{M}' \xrightarrow{l} \mathcal{M}.$$

We see that for every $H \in \left(\mathcal{N}^{\mathcal{M}}\right)_n$, we have

$$l'^*(l^*(H)) = l'^*(H(_, l(_))) = H(_, l(l'(_))) = H(_, (l';l)(_)) = (l';l)^*(H),$$

so contravariant functoriality holds for the pullback. For the pushforward, consider

$$\mathcal{N} \xrightarrow{r} \mathcal{N}' \xrightarrow{r'} \mathcal{N}'',$$

and note that

$$r'_*(r_*(H)) = r'_*(r(H(_,_))) = r'(r(H(_,_))) = (r;r')(H(_,_)) = (r;r')_*(H),$$

which verifies functoriality of the pushforward. □

To show that this object behaves like an exponential, we need to provide the "currying" and "decurrying morphisms".

**Definition 7.2.3** ((de-)currying): Let $\mathcal{M}$ be a minion. For each minion $\mathcal{L}$ and $N$, define the following functions:

$$\begin{aligned}
\Phi_{\mathcal{L},\mathcal{N}} &\colon (\mathcal{L} \times \mathcal{M} \to_{\underline{\mathsf{Min}}} \mathcal{N}) \to \left(\mathcal{L} \to_{\underline{\mathsf{Min}}} \mathcal{N}^{\mathcal{M}}\right), \\
&\quad (B_n)_n \mapsto \Big(l \mapsto \underbrace{((\alpha, m) \mapsto B_k(l\alpha, m))_k}_{P^n \times \mathcal{M} \to_{\underline{\mathsf{Min}}} \mathcal{N}}\Big)_n \\
\Psi_{\mathcal{L},\mathcal{N}} &\colon (\mathcal{L} \times \mathcal{M} \to_{\underline{\mathsf{Min}}} \mathcal{N}) \leftarrow \left(\mathcal{L} \to_{\underline{\mathsf{Min}}} \mathcal{N}^{\mathcal{M}}\right), \\
&((l, m) \mapsto A_n(l)(\mathrm{id}_n, m))_n \leftarrow\!\shortmid (A_n)_n.
\end{aligned}$$

Using some slightly ambiguous shorthand notation, we could have written the above more concisely as $\Phi(B)_n := l \mapsto B(l_, _)$ and $\Psi(A)_n := A_n(_)(\mathrm{id}_n, _)$.

**Proposition 7.2.3** (naturality of (de-)currying): The family $\left(\Phi_{\mathcal{L},\mathcal{N}}\right)_{\mathcal{L},\mathcal{N}}$ is natural both in $\mathcal{L}$ and in $\mathcal{N}$.

*Proof*: *Naturality in $\mathcal{L}$:* Let $h\colon \mathcal{L} \to \mathcal{L}'$. We need to show that the diagram

$$\begin{array}{ccc}
(\mathcal{L}' \times \mathcal{M} \to_{\underline{\mathsf{Min}}} \mathcal{N}) & \xrightarrow{\Phi_{\mathcal{L}',\mathcal{N}}} & \left(\mathcal{L}' \to_{\underline{\mathsf{Min}}} \mathcal{N}^{\mathcal{M}}\right) \\
\big\downarrow & & \big\downarrow \\
(\mathcal{L} \times \mathcal{M} \to_{\underline{\mathsf{Min}}} \mathcal{N}) & \xrightarrow{\Phi_{\mathcal{L},\mathcal{N}}} & \left(\mathcal{L} \to_{\underline{\mathsf{Min}}} \mathcal{N}^{\mathcal{M}}\right)
\end{array}$$

commutes, where the vertical arrows represent the appropriate precomposition with $h$ – i.e. on the left we have the map $\mathcal{L}' \times \mathcal{M} \to \mathcal{N} \ni B \mapsto B(h(_), _)$ and on the right we have $A \mapsto A(h(_))$. We will denote both of these maps by $h^*$. Commutativity of this diagram now means that for every binary minion homomorphism $B\colon \mathcal{L}' \times \mathcal{M} \to \mathcal{N}$, we have $h^*\left(\Phi_{\mathcal{L}',\mathcal{N}}(B)\right) = \Phi_{\mathcal{L},\mathcal{N}}(h^*(B))$. Indeed, we have

$$\begin{aligned}
h^*\big(\Phi_{\mathcal{L}',\mathcal{N}}(B)\big) &= h^*\Big(\big(l \mapsto B_n(l__, __)\big)_n\Big) && \text{Definition } \Phi_{\mathcal{L}',\mathcal{N}} \\
&= \big(l \mapsto B_n(h(l)__, __)\big)_n && \text{Definition right } h^* \\
&= \Phi_{\mathcal{L},\mathcal{N}}(B(h(__), __)) && \text{Definition } \Phi_{\mathcal{L},\mathcal{N}} \\
&= \Phi_{\mathcal{L},\mathcal{N}}(h^*(B)). && \text{Definition left } h^*
\end{aligned}$$

Analogously, for naturality in $\mathcal{N}$, let $h\colon \mathcal{N} \to \mathcal{N}'$. We now need to show that the diagram

$$\begin{array}{ccc}
(\mathcal{L} \times \mathcal{M} \to_{\underline{\mathsf{Min}}} \mathcal{N}) & \xrightarrow{\Phi_{\mathcal{L},\mathcal{N}}} & (\mathcal{L} \to_{\underline{\mathsf{Min}}} \mathcal{N}^{\mathcal{M}}) \\
\downarrow & & \downarrow \\
(\mathcal{L} \times \mathcal{M} \to_{\underline{\mathsf{Min}}} \mathcal{N}') & \xrightarrow{\Phi_{\mathcal{L},\mathcal{N}'}} & \big(\mathcal{L} \to_{\underline{\mathsf{Min}}} \mathcal{N}'^{\mathcal{M}}\big)
\end{array}$$

commutes, where the vertical arrows are postcomposition with $h$, which we denote by $h_*$ in both cases. Letting $B$ as above, we see

$$\begin{aligned}
h_*\big(\Phi_{\mathcal{L}',\mathcal{N}}(B)\big) &= h_*\Big(\big(l \mapsto B_n(l__, __)\big)_n\Big) && \text{Definition } \Phi_{\mathcal{L},\mathcal{N}} \\
&= \big(l \mapsto h_n(B_n(l__, __))\big)_n && \text{Definition right } h_* \\
&= \Phi_{\mathcal{L},\mathcal{N}}\Big(\big(h_n(B(__, __))\big)_n\Big) && \text{Definition } \Phi_{\mathcal{L},\mathcal{N}'} \\
&= \Phi_{\mathcal{L},\mathcal{N}}(h_*(B)). && \text{Definition left } h_*
\end{aligned}$$

□

**Proposition 7.2.4**: For each $\mathcal{L}$ and $\mathcal{M}$, the maps $\Phi_{\mathcal{L},\mathcal{M}}$ and $\Psi_{\mathcal{L},\mathcal{M}}$ are mutually inverse.

*Proof*: Let $B\colon \mathcal{L} \times \mathcal{M} \to \mathcal{N}$ a minion homomorphism. Let $C \coloneqq \Phi_{\mathcal{L},\mathcal{N}}$ be the curried minion homomorphism, i.e.,

$$C_n(l) = \big(B_k(l__, __)\big)_k\colon \mathcal{P}^n \times \mathcal{M} \to \mathcal{N}$$

for every $l \in \mathcal{L}_n$. Then

$$\big(\Psi_{\mathcal{L},\mathcal{N}}(C)\big)_n = (l, m) \mapsto C_n(l)(\mathrm{id}_n, m) = (l, m) \mapsto B_n(l\,\mathrm{id}_n, m) = B_n,$$

so indeed $\Psi_{\mathcal{L},\mathcal{N}}(C) = B$. For the other equality, let $C\colon \mathcal{L} \to \mathcal{N}^{\mathcal{M}}$, and let $B \coloneqq \Psi_{\mathcal{L},\mathcal{N}}(C)$, i.e.,

$$B_n = (l, m) \mapsto C_n(l)(\mathrm{id}_n, m)\colon \mathcal{L}_n \times \mathcal{M}_n \to \mathcal{N}_n.$$

We now see that for every $l \in \mathcal{L}_n$, we have

$$\begin{aligned}\left(\Phi_{\mathcal{L},\mathcal{N}}(B)\right)_n(l) &= (\alpha, m) \mapsto B_k(l\alpha, m) \\ &= (\alpha, m) \mapsto C_k(l\alpha)(\mathrm{id}_n, m), \\ &\overset{*}{=} (\alpha, m) \mapsto (C_n(l)\alpha)(\mathrm{id}_n, m) \\ &= (\alpha, m) \mapsto C_n(l)(\mathrm{id}_n\,\alpha, m) \\ &= C_n(l)\end{aligned}$$

where in $*$ we have used that $C$ is a minion homomorphism. □

We have thus verified that our construction is actually "the" exponential, as defined by the adjointness property:

**Proposition 7.2.5**: The construction from Definition 7.2.1 satisfies the axioms of an exponential object.

Unsurprisingly, we have a distributivity property:[54]

**Proposition 7.2.6**: We have $\mathcal{N}^{\sum_i \mathcal{M}_i} \cong \prod_i \mathcal{N}^{\mathcal{M}_i}$.

*Proof*: Given $H = H(_,_)\colon \mathcal{P}^k \times \sum_{i\in I} \mathcal{M}_i \to \mathcal{N}$, we can form restrictions

$$H_i(_,_)\colon \mathcal{P}^k \times \mathcal{M}_i \to \mathcal{N}$$

under the canonical inclusions to obtain a family $(H_i)_{i\in I}$ of maps, i.e., an arity-$k$ element of the product minion $\prod_i \mathcal{N}^{\mathcal{M}_i}$. This assignment is clearly a bijection, and it is a minion homomorphism, because this restriction does not interfere with the action on the first argument. □

**Corollary 7.2.7**: If $X$ is a discrete minion, i.e., a sum of constants, then $\mathcal{N}^X \cong \mathcal{N}^X$, in the sense that the $X$-fold iterated product and the exponential agree.

## 7.3 Preservation of local finiteness

It is a priori not clear why the exponential minion of two locally finite minions should be locally finite again. In fact, we will give a counterexample later in this section. However, we can prove that this is the case when the exponent is locally finite and the basis is *finitely representable*:

**Proposition 7.3.1**: Let $\mathcal{M}$ be locally finite and $\mathcal{N}$ be finitely representable. Then $\mathcal{N}^{\mathcal{M}}$ is locally finite.

[54] This holds for arbitrary exponentials due to abstract category-theoretical reasons, namely because $\mathcal{N}^\bullet$ is left adjoint to its opposite, which implies that it turns colimits into the respective limits.

*Proof*: Let $n$ be a natural number and, without loss of generality, assume $\mathcal{N}$ is a subminion of some $\mathcal{O}(l, k)$. We have the correspondences[55]

$$\begin{aligned}
\left(\mathcal{N}^{\mathcal{M}}\right)_n &= \mathcal{P}^n \times \mathcal{M} \to_{\underline{\mathsf{Min}}} \mathcal{N} && \text{Definition} \\
&\subseteq \mathcal{P}^n \times \mathcal{M} \to_{\underline{\mathsf{Min}}} \mathcal{O}(l,k) \\
&\cong (\mathcal{P}^n \times \mathcal{M})_l \to_{\underline{\mathsf{Set}}} k && \text{Corollary 4.1.4}, \\
&= l^n \times \mathcal{M}_l \to_{\underline{\mathsf{Set}}} k && \text{Definition } \mathcal{P} \\
&\cong k^{l^n \times \mathcal{M}_l},
\end{aligned}$$

i.e., for every element of arity $n$, we can assign a unique map from $(\mathcal{P}^n \times \mathcal{M})_l$ to $k$, and of such maps there are only finitely many, namely $k^{l^n \times |\mathcal{M}_l|}$. □

### 7.3.1 Counterexample

In this section we will prove that there is a locally finite minion $\Omega$ which has infinitely many endomorphisms. Since endomorphisms correspond to constants in the self-exponential $\Omega^\Omega$, This implies that $\Omega^\Omega$ cannot be locally finite as well. Note that due to Proposition 7.3.1, for this to work we need $\Omega$ to not be finitely representable.

The object we use is well-known in topos theory: It is the so-called *subobject classifier* of the category of minions, a highly abstract object.

**Definition 7.3.1** (Cosieve): Let $n \in \mathbb{N}_{>0}$. An *n-cosieve* [26, Section II.2] is a set $K$ of maps $n \to k$ for arbitrary finite codomain ordinals $k$ which is closed under postcomposition with arbitrary maps, i.e., for every $\alpha\colon n \to k \in K$ and $\beta\colon k \to k'$, the composite $\alpha;\beta$ must be in $K$ as well.

**Definition 7.3.2** (The minion $\Omega$): Define the minion $\Omega = (\Omega_n)_n$ by

$$\Omega_n := \{K \mid K \text{ is an } n\text{-cosieve}\},$$

where a minor operation $\alpha\colon n \to n'$ acts via

$$K\alpha := \{\beta\colon n' \to k \mid k \in \mathbb{N}_{>0}, \alpha;\beta \in K\}.$$

**Proposition 7.3.2**: The object $\Omega$ is a well-defined minion.

*Proof*: The only nontrivial thing to show is that when $K$ is an $n$-cosieve and $\alpha\colon n \to n'$, the set $K\alpha$ is indeed a cosieve. Indeed, let $\beta\colon n' \to k \in K\alpha$ and $\gamma\colon k \to k'$ arbitrary. By definition of $\beta \in K\alpha$, we have $\alpha;\beta \in K$, and by closure under postcomposition we conclude $\alpha;\beta;\gamma \in K$. This however means that $\beta;\gamma \in K\alpha$, so $K\alpha$ is indeed closed under postcomposition.

[55] Note that when we write $A \cong B$ where $A$ and $B$ are sets without additional structure, that just means that the two sets are in bijection, because isomorphisms in the categoy of sets are bijections. For our purposes this is just a convenient way to assert that the cardinalities are equal without having to surround the sets with vertical lines every time (i.e., $|A| = |B|$).

That $K\,\mathrm{id}_n = K$ and $(K\alpha)\beta = K(\alpha;\beta)$ follows just as easily. □

**Definition 7.3.3** (truth): For each $n \in \mathbb{N}_{>0}$, we define the *full $n$-cosieve* $\mathtt{true}_n := \{\beta\colon n \to k \mid k \in \mathbb{N}_{>0}\}$.

One easily verifies that for each $\alpha\colon n \to n'$, we have $\mathtt{true}_n\alpha = \mathtt{true}_{n'}$, which in particular means that $(\mathtt{true}_n)_n$ is a constant subminion of $\Omega$.

For the next observation, we need the following lemma.

**Lemma 7.3.3**: Let $K$ be an $n$-cosieve and $\alpha\colon n \to n'$. We have $K\alpha = \mathtt{true}_{n'}$ if and only if $\alpha \in K$.

*Proof*: The full cosieve $\mathtt{true}_{n'}$ is generated by the identity $\mathrm{id}_{n'}$, so $K\alpha = \mathtt{true}_{n'}$ if and only if $\mathrm{id}_{n'} \in K\alpha$. Untangling the definitions, this is the case if and only if $\alpha;\mathrm{id}_{n'} \in K$, however the left hand side is just $\alpha$. □

**Proposition 7.3.4**: The minion $\Omega$ with designated elements $(\mathtt{true}_n)_n$ satisfies the following universal property: For every minion $\mathcal{M}$ and subminion $\mathcal{N} \leq \mathcal{M}$, there is a unique minion homomorphism $\chi_{\mathcal{N}}\colon \mathcal{M} \to \Omega$ which sends $f \in \mathcal{M}_n$ to $\mathtt{true}_n$ if any only if $f \in \mathcal{N}$.

*Proof*: For each subminion $\mathcal{N} \leq \mathcal{M}$, we will construct this characteristic function $\chi := \chi_{\mathcal{N}}$ by assigning each $f \in \mathcal{M}_n$ the set $\chi(f)$ of all minor operations $\beta\colon n \to k$ for which $f\beta \in \mathcal{N}_k$. This set is clearly closed under postcomposition, so it is an $n$-cosieve. We need to verify that $\chi$ is a homomorphism of minions: Indeed, if $\alpha\colon n \to n'$, then

$$\chi(f)\alpha = \{\beta\colon n' \to k \mid \alpha;\beta \in \chi(f)\} = \{\beta\colon n' \to k \mid f(\alpha;\beta) = (f\alpha)\beta \in \mathcal{N}\} = \chi(f\alpha).$$

Finally, we note that $\chi(f)$ is the full cosieve $\mathtt{true}_n = \langle \mathrm{id}_n \rangle$ if and only if $f = f\,\mathrm{id}_n \in \mathcal{N}$. This means that the map $\mathcal{N} \mapsto \chi_{\mathcal{N}}$ from subminions of $\mathcal{M}$ to homomorphisms $\mathcal{M} \to \Omega$ is injective, because different subminions $\mathcal{N}$ and $\mathcal{N}'$ will give rise to functions $\chi_{\mathcal{N}}$ and $\chi_{\mathcal{N}}$ which send different elements to full cosieves, hence must be different.

To see surjectivity of the correspondence $\mathcal{N} \mapsto \chi_{\mathcal{N}}$, note that for every homomorphism $\psi\colon \mathcal{M} \to \Omega$, the preimage minion $\mathcal{N} := \psi^{-1}\big((\mathtt{true}_n)_n\big)$ of all the elements being sent to a full cosieve satisfies $\chi_{\mathcal{N}} = \psi$: Indeed, Let $f \in \mathcal{M}_n$. We then know by definition of $\mathcal{N}$ that for every $\alpha\colon n \to n'$ with $f\alpha \in \mathcal{N}$, we have $\psi(f\alpha) = \mathtt{true}_{n'}$, and by homomorphicity $\psi(f)\alpha = \mathtt{true}_{n'}$. Furthermore, whenever $f\alpha \notin \mathcal{N}$, we know that $\psi(f\alpha) = \psi(f)\alpha \neq \mathtt{true}_{n'}$. Therefore, $\psi(f)$ must be an $n$-cosieve which satisfies $\psi(f)\alpha = \mathtt{true}_{n'}$ if and only if $f\alpha \in \mathcal{N}$. By lemma Lemma 7.3.3, this means that $\alpha \in \psi(f)$ if and only if $f\alpha \in \mathcal{N}$, which is the definition of $\alpha \in \chi_{\mathcal{N}}(f)$. Therefore, $\psi = \chi_{\mathcal{N}}$. □

**Proposition 7.3.5**: For each $n \in \mathbb{N}_{>0}$, there are only finitely many $n$-cosieves – i.e., $\Omega$ is locally finite.

*Proof*: Let $n \in \mathbb{N}_{>0}$. For a set $G$ of maps $n \to _$, we define the generated cosieve $\langle G \rangle$ in the usual way as the smallest $n$-cosieve containing $G$. Explicitly, $\langle G \rangle = \{\alpha;\beta \mid \alpha \in G, \beta\colon k \to k'\}$. From that it is evident that $\langle G \rangle = \cup_{\beta \in G} \langle \beta \rangle$. We say that generating sets $G$ and $G'$ by demanding $\langle G \rangle = \langle G' \rangle$.

Claim: Every generating set is equivalent to a generating set whose every generator has codomain at most $n$. Indeed, if $\beta\colon n \to k$, then we can factor $\beta = \tilde{\beta};\eta$ where $\tilde{\beta}\colon n \twoheadrightarrow i$ is surjective and $\eta\colon i \hookrightarrow n$ injective, where we set $i := |\mathrm{Im}\ \beta|$. Therefore, $\beta \in \left\langle \tilde{\beta} \right\rangle$. Since the map $\eta$ is an injection, it has a post-inverse $\rho$, whence $\tilde{\beta} = \left(\tilde{\beta};\eta\right);\rho = \beta\rho$, showing that $\tilde{\beta} \in \langle \beta \rangle$. Therefore, $\beta$ and $\tilde{\beta}$ are equivalent. We now see that $\langle G \rangle = \cup_{\beta \in G} \langle \beta \rangle = \cup_{\beta \in G} \tilde{\beta} = \left\langle \left\{ \tilde{\beta} \,\middle|\, \beta \in G \right\} \right\rangle$, so a choice of $\tilde{\beta}$ for each $\beta \in G$ gives an equivalent generating set all of whose elements have codomain $G$.

Since every $n$-cosieve is generated by a subset of maps $n \to k$ with $k$ at most $n$ – of which there are only finitely many – there can only be finitely many $n$-cosieves. □

**Proposition 7.3.6**: The minion $\Omega$ is not finitely generated. In particular, $\Omega$ has infinitely many subminions.

*Proof*: For an arbitrary $n \geq 2$, consider the $n$-cosieve $\langle \mathrm{const}_n \rangle$ with $\mathrm{const}_n\colon n \to 1$ the constant map. We claim that every $i \in n$ is essential in $\langle \mathrm{const}_n \rangle$: Indeed, $\langle \mathrm{const}_n \rangle \iota$ does not contain the characteristic function $\chi_{i,n+1}\colon n+1 \to 2$, because $\iota;\chi_{i,n+1} = \chi_{i,n}$, which does not factor through $\mathrm{const}_n$, as it does not identify $i$ and any element $\neq i$. On the other hand, ${}_i(n \mapsto n+1)_n\ \chi_{i,n+1} = \mathrm{const}_n$, so $\chi_{i,n+1} \in \langle \mathrm{const}_n \rangle\ {}_i(n \mapsto n+1)_n$ . Therefore, these two cosieves are always distinct. Thus $\langle \mathrm{const}_n \rangle$ is nondegenerate for every $n \geq 2$.

For the last statement, note that the subminions generated by an element with $n$ essential coordinates must have essential arity at most $n$, so the sub-minions generated by the cosieve $\langle \mathrm{const}_n \rangle$ must all be distinct. □

**Corollary 7.3.7**: $\Omega$ has infinitely many endomorphisms.

*Proof*: Invoking the classifying property (Proposition 7.3.4) of $\Omega$ with respect to itself, we see that to each subminion $M \leq \Omega$ there must correspond a unique endomorphism $\Omega \to \Omega$, hence by Proposition 7.3.6 there are infinitely many such endomorphisms. □

**Corollary 7.3.8**: The self-exponential $\Omega^\Omega$ is not locally finite. In particular, there exist locally finite minions $\mathcal{M}$ and $\mathcal{N}$ such that $\mathcal{N}^\mathcal{M}$ is not locally finite.

*Proof*: Constants in $\Omega^\Omega$ – i.e., minion homomorphisms $* \to \Omega^\Omega$ – correspond to endomorphisms of $\Omega$, of which there are infinitely many by the last Corollary 7.3.7. This violates local finiteness. □

## 7.4 Preservation of finite representability

**Proposition 7.4.1**: If $\mathcal{M}$ is locally finite and if $\mathcal{N}$ is $l$-representable, then $\mathcal{N}^{\mathcal{M}}$ is $l$-representable.

*Proof*: Since $\mathcal{M}$ is locally finite and $\mathcal{N}$ is finitely representable, we can deduce from Proposition 7.3.1 that $\mathcal{N}^{\mathcal{M}}$ must be locally finite as well. This enables the characterization from Proposition 4.2.3: To see that $\mathcal{N}^{\mathcal{M}}$ is $l$-representable it suffices to show that for all $H \neq H' \in \left(\mathcal{N}^{\mathcal{M}}\right)_n$ there is a minor operation $\alpha\colon n \to l$ such that $H\alpha \neq H'\alpha$.

So let $H \neq H'$ be of arity $n$, i.e., they are minion homomorphisms $P^n \times \mathcal{M} \to \mathcal{N}$, and let $\beta\colon n \to k$ and $m \in \mathcal{M}_k$ be witnesses to their inequality, i.e., $H_k(\beta, m) \neq H'_k(\beta, m) \in \mathcal{N}_k$. By $l$-representability of $\mathcal{N}$, there must be an $\alpha\colon k \to l$ such that $H_k(\beta, m)\alpha \neq H'_k(\beta, m)\alpha$, which by homomorphicity of $H$ and $H'$ means that

$$\begin{aligned}
(H(\beta;\alpha))_l(\mathrm{id}_l, m\alpha) &= H_l(\beta;\alpha, m\alpha) && \text{Minion structure of } \mathcal{N}^{\mathcal{M}} \\
&= (H_k(\beta, m))\alpha && \text{Homomorphicity of } H \\
&\neq (H'_k(\beta, m))\alpha && \text{Definition of } \alpha \\
&= H'_l(\beta;\alpha, m\alpha) && \text{Homomorphicity of } H \\
&= (H'(\beta;\alpha))_l(\mathrm{id}_l, m\alpha). && \text{Minion structure of } \mathcal{N}^{\mathcal{M}}
\end{aligned}$$

Therefore, $H\alpha\beta \neq H'\alpha\beta$, as they differ when evaluated $\mathrm{id}_l\colon l \to l$ and $m\alpha \in \mathcal{M}_l$. □

For reasons which will become appearent in Section 7.6, $\mathcal{P}^{\mathcal{M}}$ can be seen as a sort of "dual" to the minion $\mathcal{M}$ relative to $\mathcal{P}$. Indeed, it is not hard to see the following:

**Lemma 7.4.2**: Let $\mathcal{M}$ and $\mathcal{N}$ be minions. Then $\mathcal{M} \to \mathcal{N}^{\mathcal{N}^{\mathcal{M}}}$.

Note that here we took $\to$ to denote the existence of a homomorphism, not the set of all homomorphisms.

*Proof*: We have

$$\left(\mathcal{N}^{\mathcal{M}} \to \mathcal{N}^{\mathcal{M}}\right) \Rightarrow \left(\mathcal{N}^{\mathcal{M}} \times \mathcal{M} \to \mathcal{N}\right) \Rightarrow \left(\mathcal{M} \times \mathcal{N}^{\mathcal{M}} \to \mathcal{N}\right) \Rightarrow \left(\mathcal{M} \to \mathcal{N}^{\mathcal{N}^{\mathcal{M}}}\right).$$

□

Quite noteworthy about the previous Proposition 7.4.1 is that the exponent need not be finitely representable for the exponential to be. If $\mathcal{M}$ is locally finite but not finitely representable and $\mathcal{N}$ is $n$-representable, then this implies that the homomorphism $\mathcal{M} \to \mathcal{N}^{\mathcal{N}^{\mathcal{M}}}$ cannot be an injection, because being finitely representable is a hereditary property.[56]

Nevertheless, this begs some questions:

[56] i.e., minions which are (isomorphic to) subminions of finitely representable minions are finitely representable themselves.

**Open Question 7.1**:
1. Given a locally finite minion $\mathcal{M}$ and a finitely representable minion $\mathcal{N}$, under which conditions does $\mathcal{N}^{\mathcal{N}^{\mathcal{M}}} \to \mathcal{M}$?
2. Given a locally finite minion $\mathcal{M}$, does there always exist a finitely representable minion $\mathcal{N}$ such that $\mathcal{N}^{\mathcal{N}^{\mathcal{M}}} \to \mathcal{M}$?

A positive answer to point 2 of this question would imply a positive answer to the following as well:

**Open Question 7.2**: Is every locally finite minion homomorphically equivalent to a finitely representable minion?

*Proof – assuming a positive answer to point 2 of Open Question 7.1*: Let $\mathcal{M}$ be locally finite. By assumption there is a finitely representable minion $\mathcal{N}$ such that $\mathcal{M}$ is homomorphically equivalent to $\mathcal{N}^{\mathcal{N}^{\mathcal{M}}}$, since a homomorphism $\mathcal{M} \to \mathcal{N}^{\mathcal{N}^{\mathcal{M}}}$ always exists. Since $\mathcal{N}^{\mathcal{M}}$ is the exponential of a finitely representable minion with locally finite exponent, it must be locally finite by Proposition 7.4.1, and applying the same proposition again to $\mathcal{N}^{\mathcal{N}^{\mathcal{M}}}$ proves that this minion is finitely representable. □

### 7.4.1 Relational basis

While our construction produces a finitely representable minion, we only live on one side of the inv-pol duality in the sense that we give an abstract description of a minion, but do not specify a relational basis, or how to compute one. In practice, minions of interest will mostly be given in the form $\mathrm{Pol}(\underline{A}, \underline{B})$, where $\underline{A}$ and $\underline{B}$ are structuers of a given – often finite – relational signature. The following would be helpful to answer:

**Open Question 7.3**: Let $\underline{A}, \underline{B}$ be relational structures and let $n := |A|$. If $\mathcal{M}$ is a minion with $\mathcal{M}_n < \infty$, what does a relational basis of $\mathrm{Pol}(\underline{A}, \underline{B})^{\mathcal{M}}$ look like?

Even if such an explicit description would exist, it would not necessarily have to produce a finite relational basis even if we started with one. A natural follow-up question is thus:

**Open Question 7.4**: Does the exponential minion preserve being finitely related? In particular, given a (finite) relational basis for $\mathcal{M}$ and $\mathcal{N}$, is there an algorithm to produce a relational basis of $\mathcal{N}^{\mathcal{M}}$?

## 7.5 Preservation of finite generatedness

In this section we will show that exponentiating minions which are both finitely generated – i.e., locally finite and of bounded essential arity – yields a finitely generated minion. For that, we will use the growth characterization of being finitely generated (Corollary 6.2.3).

**Proposition 7.5.1**: Let $S = (S_n)_n$ be a finite (indexed) set and $\mathcal{N}$ finitely generated. Then $\mathcal{N}^{\langle S \rangle}$ is finitely generated.

*Proof*: Let $S = (f_i \colon n_i)_{i \in I}$, which shall be taken to mean that $f_i \in S_{n_i}$. By the decomposition of free minions, we know that $S \cong \sum_{\substack{k>0 \\ s \in S_k}} \mathcal{P}^k$. By distributivity (Prop. 7.2.6), this means that the exponential decomposes as

$$\mathcal{N}^{\langle S \rangle} \cong \mathcal{N}^{\sum_{\substack{k>0 \\ s \in S_k}} \mathcal{P}^k} \cong \prod_{\substack{k>0 \\ s \in S_k}} \mathcal{N}^{\mathcal{P}^k}.$$

Focusing on a fixed arity $n$, we see that

$$\begin{aligned}
\left(\mathcal{N}^{\langle S \rangle}\right)_n &\cong \prod_{\substack{k>0 \\ s \in S_k}} \left(\mathcal{N}^{\mathcal{P}^k}\right)_n && \text{See above; Proposition 2.3.1} \\
&\cong \prod_{\substack{k>0 \\ s \in S_k}} \left(\mathcal{P}^n \times \mathcal{P}^k \to_{\underline{\mathsf{Min}}} \mathcal{N}\right) && \text{Definition 7.2.1} \\
&\cong \prod_{\substack{k>0 \\ s \in S_k}} \left(\mathcal{P}^{n+k} \to_{\underline{\mathsf{Min}}} \mathcal{N}\right) && \text{Proposition 7.2.6/Corollary 7.2.7} \\
&\cong \prod_{\substack{k>0 \\ s \in S_k}} \mathcal{N}_{n+k} && \text{Corollary 2.4.6.}
\end{aligned}$$

Therefore, $\alpha_{\mathcal{N}^{\langle S \rangle}}(n)$ is a product of $\sum_k S_k$ many (ergo, finitely many) subexponential sequences, so it is subexponential itself. By Corollary 6.2.3, we can conclude that $\mathcal{N}^{\langle S \rangle}$ must be finitely generated. □

**Proposition 7.5.2**: Let $h \colon \mathcal{M}' \to \mathcal{M}$ be a surjection of minions. Then the pullback $h^* \colon \mathcal{N}^{\mathcal{M}} \to \mathcal{N}^{\mathcal{M}'}$ is an injection.

*Proof*: Let $H \neq H' \in \left(\mathcal{N}^{\mathcal{M}}\right)_n = (\mathcal{P}^n \times \mathcal{M} \to_{\underline{\mathsf{Min}}} \mathcal{N})$. We will show that $h^*(H) \neq h^*(H')$, i.e., $H(_, h(_)) \neq H'(_, h(_))$. Assume the inequality $H \neq H'$ is witnessed in arity $k$ by a pair $(\alpha, m) \in k^n \times \mathcal{M}_k =$, i.e., we have $H_k(\alpha, m) \neq H'_k(\alpha, m)$. By surjectivity of $h \colon \mathcal{M}' \to \mathcal{M}$, there must be an $m' \in \mathcal{M}'_k$ in the preimage of $m$ under $h$. Therefore, $H_k(\alpha, h(m')) \neq H'_k(\alpha, h(m'))$, so the pair $(\alpha, h(m'))$ serves as a witness for the inequality of $h^*(H)$ and $h^*(H')$, as desired. □

**Proposition 7.5.3**: Let $\mathcal{M}$ and $\mathcal{N}$ be finitely generated. Then $\mathcal{N}^{\mathcal{M}}$ is finitely generated.

*Proof*: If $\mathcal{M}$ is finitely generated, then there is a finite (indexed) set $\Sigma$ of generators, i.e., a surjective minion homomorphism $h \colon \langle \Sigma \rangle \to \mathcal{M}$. By Proposition 7.5.2, we know that the induced pullback map $h^* \colon \mathcal{N}^{\mathcal{M}} \to \mathcal{N}^{\langle \Sigma \rangle}$ is an injection. From Proposition 7.5.1 we know that $\mathcal{N}^{\langle\Sigma\rangle}$ is finitely generated, so $\mathcal{N}^{\mathcal{M}}$ is a subminion of a finitely generated minion, hence finitely generated itself (Proposition 2.6.15). □

## 7.6 Dualization

**Definition 7.6.1**: A bounded lattice $(L, \wedge, \vee, 0, 1)$ together with a function $\rightarrow: L \times L \rightarrow L$, denoted as an infix operator, is called a *heyting algebra* if for all $x, y, z \in L$ we have the equivalence

$$(x \wedge y) \leq z \Leftrightarrow x \leq (y \rightarrow z).$$

The function $\rightarrow$ is called *relative pseudocomplement*. We will denote a Heyting algebra by $(P, \wedge, \vee, \rightarrow, 0, 1)$.

We will prove that this is additional structure is intrinsic in the sense that any two relative pseudocomplements must agree. To see tha this follows from the universal property, we shall introduce a lemma about posets that despite its triviality is remarkably useful.

**Lemma 7.6.1** (Yoneda lemma for posets): Let $(P, \leq)$ be a poset and $y, z \in P$. If $y^{\downarrow} = z^{\downarrow}$, i.e., $x \leq y$ if and only if $x \leq z$, then $y = z$.

*Proof*: By the assumptions, $y \leq y$ implies $y \leq z$, and $z \leq z$ implies $z \leq y$. □

**Proposition 7.6.2**: Any two Heyting algebra structures on a bounded lattice must agree.

*Proof*: If $\rightarrow'$ is any second Heyting algebra structure, then

$$x \leq (y \rightarrow z) \Leftrightarrow x \wedge y \leq z \Leftrightarrow x \leq (y \rightarrow' z),$$

so by Lemma 7.6.1, $y \rightarrow z$ and $y \rightarrow' z$ must agree. □

*Remark (for category theorists)*: Of course, this is just the yoneda lemma for posets in disguise: We have actually shown that $y \rightarrow z$ and $y \rightarrow' z$ have the same downward closures ($(y \rightarrow z)^{\downarrow}$ resp. $(y \rightarrow' z)^{\downarrow}$) in the poset, i.e., their hom functors $_ \rightarrow (y \rightarrow z)$ and $_ \rightarrow (y \rightarrow' z)$ are naturally isomorphic.

Since the Heyting algebra structure on a lattice is uniquely determined, we will conflate the notions of "Heyting algebra" and "lattice which admits a Heyting algebra structure" and just say that a lattice "is a Heyting algebra".

**Proposition 7.6.3**: A bounded lattice $(L, \wedge, \vee, 0, 1)$ is a Heyting algebra if and only if for every $b, c \in L$, the set $\{a \in L \mid a \wedge b \leq c\}$ has a maximum.

*Proof*: If $L$ is a Heyting algebra with relative pseudo complement $\rightarrow$, then $b \rightarrow c$ is the maximum of $\{a \in L \mid a \wedge b \leq c\}$: indeed, if $a \wedge b \leq c$, then $a \leq b \rightarrow c$ by the defining prop-erty of $\to$. Conversely, if $\{a \in L \mid a \wedge b \leq c\}$ has a maximum $m$, then $a \leq m$ if and only if $a$ is in that set, which means $a \wedge b \leq c$. Hence $m$ satisfies the desired universal property of the relative pseudocomplement. □

A lattice which is a Heyting algebra is often called a *Heyting lattice* (see e.g. [61, Def. 1.5.2]).

*Remark*: Instead of using a universal property, heyting algebras can be axiomatized using an equatonal theory using the symbols $\wedge, \vee, \to, 0, 1$ (see [61, Prop. 2.1.4]); this implies that Heyting algebras form a variety in the sense of universal algebra when we consider homomorphisms to be homomorphisms of bounded lattices which also preserve the relative pseudocomplement.

**Proposition 7.6.4**: Let $\mathcal{C}$ be a class of minions (i.e., a full subcategory of Min) such that
- The homomorphism order $\mathcal{C}/\leftrightarrow$ is a bounded lattice with $[\mathcal{M}] \wedge [\mathcal{N}] = [\mathcal{M} \times \mathcal{N}]$
- $\mathcal{C}$ is closed under taking exponentials.

Then $\mathcal{C}/\leftrightarrow$ is a Heyting algebra.

*Proof*: Let $\mathcal{L}, \mathcal{M}, \mathcal{K} \in \mathcal{C}$. Then $\mathcal{K}^{\mathcal{M}} \in \mathcal{C}$. Claim: $[\mathcal{K}^{\mathcal{M}}] \in \mathcal{C}/\leftrightarrow$ is a relative pseudocomplement. Indeed, we have $[\mathcal{L}] \wedge [\mathcal{M}] = [\mathcal{L} \times \mathcal{M}] \leq [\mathcal{K}]$ if and only if there is a homomorphism $h\colon \mathcal{L} \times \mathcal{M} \to \mathcal{K}$, which happens if and only if there is a homomorphism $\mathcal{L} \to \mathcal{K}^{\mathcal{M}}$, i.e., $[\mathcal{L}] \leq [\mathcal{K}^{\mathcal{M}}]$. □

**Corollary 7.6.5**: The following bounded lattices are Heyting algebras:
1. $\mathfrak{M}$, the homomorphism order of all (abstract) minions.
2. $\mathfrak{F}$, the homomorphism order of all finitely representable minions.
3. $\mathfrak{FG}$, the homomorphism order of all finitely generated minions.

*Proof*: In each case we verify that the premises of Proposition 7.6.4 are fulfilled. Case 1. is clearly true because Min has products and exponentials. For case 2, recall that the class of locally finite minions is closed under products due to Proposition 5.2.2, and under exponentials due to Proposition 7.4.1. For case 3, the product of two locally finite minions is locally finite by construction (Proposition 2.3.1), and closure under exponentials is Proposition 7.5.3.

That all of the three homomorphism orders form bounded lattices has been established in Corollary 5.1.4. □

For our discussion of density and gaps, We also need the dual notion of a Heyting lattice / lattice with a Heyting structure / Heyting algebra.

**Definition 7.6.2**: A *co-Heyting algebra* is a bounded lattice with an operation $\setminus\colon L^2 \to L$, denoted as an infix operator, which turns its opposite lattice[57] into a Heyting algebra.

Spelt out, it is a lattice with a function $a \setminus b$ such that $a \setminus b \leq c$ if and only if $a \leq b \vee c$.

[57] This is the lattice with the opposite order $a \overset{\text{op}}{\leq} b \Leftrightarrow b \leq a$. It is easy to see that if an order is a lattice (resp. bounded lattice), then its opposite order is also a lattice (resp. bounded lattice).

**Definition 7.6.3**: For minions $\mathcal{L}, \mathcal{M}$, define the *subtraction* $\mathcal{L} \setminus \mathcal{M}$ as

$$\sum_{\substack{u \in \mathcal{L}_1 \\ \mathcal{L}|_u \not\to \mathcal{M}}} \mathcal{L}|_u.$$

**Proposition 7.6.6**: Let $\mathcal{L}, \mathcal{M}, \mathcal{N}$ be minions. Then $\mathcal{L} \setminus \mathcal{M} \to_{\underline{\mathsf{Min}}} \mathcal{N}$ if and only if $\mathcal{L} \to_{\underline{\mathsf{Min}}} \mathcal{M} + \mathcal{N}$.

*Proof*: Let $M := \{u \in \mathcal{L}_1 \mid \mathcal{L}|_u \to \mathcal{M}\}$. We observe

$$\begin{aligned}
\mathcal{L} = \left(\sum_{u \in \mathcal{L}_1} \mathcal{L}|_u\right) &\to \mathcal{M} + \mathcal{N} && \text{Proposition 5.3.4} \\
\Leftrightarrow \forall u \in \mathcal{L}_1 \colon \mathcal{L}|_u &\to \mathcal{M} + \mathcal{N} && \text{Proposition 2.3.1} \\
\Leftrightarrow \forall u \in \mathcal{L}_1 \colon (\mathcal{L}|_u &\to \mathcal{M}) \vee (\mathcal{L}|_u \to \mathcal{N}) && \text{Proposition 5.3.2} \\
\overset{(*)}{\Leftrightarrow} \forall u \in \mathcal{L}_1 \setminus M \colon (\mathcal{L}|_u &\to \mathcal{N}) \\
\Leftrightarrow \mathcal{L} \setminus \mathcal{M} = \left(\sum_{u \in \mathcal{L}_1 \setminus M} \mathcal{L}|_u\right) &\to \mathcal{N} && \text{Proposition 2.3.1},
\end{aligned}$$

where the forward direction of $(*)$ follows because for all the $u \in \mathcal{L}_1 \setminus M$ by definition cannot satisfy $\mathcal{L}|_u \to \mathcal{M}$, and the backward direction of $(*)$ follows because for all the $u \in M$ we have $\mathcal{L}|_u \to \mathcal{M}$, also by definition. □

**Proposition 7.6.7**: Let $\mathcal{C}$ be a class of minions (i.e., a full subcategory of $\underline{\mathsf{Min}}$) such that
- The homomorphism order $\mathcal{C}/\leftrightarrow$ is a bounded lattice with $[\mathcal{M}] \vee [\mathcal{N}] = [\mathcal{M} + \mathcal{N}]$
- $\mathcal{C}$ is closed under forming subtractions.

Then $\mathcal{C}/\leftrightarrow$ is a co-Heyting algebra.

*Proof*: Analogous to the proof of Proposition 7.6.4. □

**Definition 7.6.4**: A lattice is called a *bi-Heyting algebra* if it is both a Heyting algebra and a co-Heyting algebra.

**Proposition 7.6.8**: The following bounded lattices are bi-Heyting algebras:
1. $\mathfrak{M}$, the homomorphism order of all (abstract) minions.
2. $\mathfrak{F}$, the homomorphism order of all finitely representable minions.
3. $\mathfrak{FG}$, the homomorphism order of all finitely generated minions.

*Proof*: The Heyting algebra structure has been established in Corollary 7.6.5. For the co-Heyting algebra structure we shall verify the premises of Proposition 7.6.7: That the classes are closed under taking sums is established in Proposition 2.3.1, Proposition 5.2.1, and Proposition 2.6.14. Since $\mathcal{M} \setminus \mathcal{N} \subseteq \mathcal{M}$ and all the classes are closed under taking subminions, they are closed under forming subtractions. □

### 7.6.1 Gaps, density, and duality pairs

As known for the homomorphism order of finite structures, one can utilize an exponential in order to obtain information about gaps by relating them to duality pairs. We will develop the same ideas as in [62], albeit with a little more abstract setting.

**Definition 7.6.5**: Let $P$ be a poset. A pair $(p, q) \in P^2$ is called a *duality pair* if for every $b \in P$ we have either $p \leq b$ or $b \leq q$.

Equivalently, $(p, q)$ is a duality pair if $P = p^{\uparrow} \cup q^{\downarrow}$, where $p^{\uparrow} := \{x \mid x \geq p\}$ is the *upward closure* of $p$ (dually, $p^{\downarrow}$ would be the so-called *downward closure* of $p$).

Recall that

- the *closed interval* $[a, c]$ is defined to be $[a, c] := \{b \mid a \leq b \leq c\}$
- an interval $[a, c]$ is called *proper* if $a < c$
- a *gap* in a poset is an interval $[a, c]$ with $a < c$ such that there exists no $b$ with $a < b < c$.

**Lemma 7.6.9**: Let $L$ be a distributive bounded lattice and $a \leq c$. Then the map $\pi(b) := a \vee (b \wedge c)$ is an idempotent lattice endomorphism whose image is the interval $[a, c]$ which maps $0$ to $a$ and $1$ to $c$.

*Proof*: Recall that in a lattice, $x \leq y$ if and only if $x \wedge y = x$ if and only if $x \vee y = y$.

First, let us verify that $\pi(b) \in [a, c]$. That $a \leq \pi(b)$ is clear because we join with $a$. For $\pi(b) \leq c$ note that

$$(a \vee (b \wedge c)) \wedge c = (a \wedge c) \vee (b \wedge c \wedge c) = a \vee (b \wedge c) = \pi(b).$$

To show idempotence, let $b \in [a, c]$. Then $b \leq c$ implies $b \wedge c = c$, and $a \leq b$ implies $b = a \vee b = a \vee (b \wedge c)$.

That this is a lattice homomorphism follows from distributivity: meeting with $c$ and joining with $a$ are both lattice homomorphisms, hence is their composite.

Finally, note that $\pi(0) = a \vee (0 \wedge c) = a \vee 0 = a$ and $\pi(1) = a \vee (1 \wedge c) = a \vee c = c$, proving the last assertion. □

**Proposition 7.6.10**: Let $(L, \wedge, \vee, \rightarrow, \setminus, 0, 1)$ be a bi-Heyting algebra, and $a < c$. Then $[a, c]$ is a gap if and only if $(c \setminus a, c \rightarrow a)$ is a duality pair.

*Proof*: Let $b \in [a, c]$, i.e., $b = \pi(b)$ with $\pi$ the projection from Lemma 7.6.9. Now note first that $b \leq a$ – since $b = b \wedge c$ – is equivalent to $b \leq (c \rightarrow a)$. Dually, we can observe that $c \leq b$ – since $b = a \vee b$ – is equivalent to $c \setminus a \leq b$.

First assume that $[a, c]$ is a gap. Then for all $b \in L$, we have either $\pi(b) = a$ or $\pi(b) = c$. The first possibility is equivalent to $\pi(b) = a \vee (b \wedge c) \leq a$, which is equivalent to $b \wedge c \leq a$, i.e., $b \leq (c \to a)$. The second possibility is equivalent to $c \leq \pi(b) = a \vee (b \wedge c)$, i.e., $c \setminus a \leq b \wedge c$. Since $c \leq a \vee c$, we always have $c \setminus a \leq c$, so the above is equivalent to $c \setminus a \leq b$.

Thus, if $[a, c]$ is a gap, then every element $b$ is either above $c \setminus a$ or below $c \to a$, proving that these form a duality pair.

For the other direction, assume that $(c \setminus a, c \to a)$ is a duality pair, and let $b \in [a, c]$. The previous lema implies that $b = \pi(b) = a \vee (b \wedge c)$. If $c \setminus a \leq b$, then $c \leq a \vee b$. Since $a \leq b$, this is equivalent to $c \leq b$, which implies that $a = c$. On the other hand, if $b \leq c \to a$, then $b \wedge c \leq a$. Since $b \leq c$, this is equivalent to $b \leq a$, so in that case $b = a$.

Therefore, $[a, c]$ must be a gap. □

Recall that in a poset, an element $y$ is a *cover* of $x$ if $[x, y]$ is a gap.

**Corollary 7.6.11**: In $\mathfrak{M}$, $\mathfrak{F}$, and $\mathfrak{F}\mathfrak{G}$, $[\mathcal{M}]$ is a cover of $[\mathcal{P}]$ if and only if $\big([\mathcal{M}] \setminus [\mathcal{P}], [\mathcal{P}^{\mathcal{M}}]\big)$ is a duality pair.

*Proof*: This follows from Proposition 7.6.10 because all of these posets are bi-Heyting algebras (Proposition 7.6.8). □

In light of the pp-formula characterization of the existence of a homomorphism (Proposition 2.5.4), we are looking for minions $\mathcal{M}$ which satisfy a minimal set of nontrivial minor conditions.

**Proposition 7.6.12**: Let $\mathcal{M}, \mathcal{L}$ be finitely generated minions. Then the following are equivalent:
1. The closed interval $[[\mathcal{M}], [\mathcal{L}]]$ is a gap in $\mathfrak{L}\mathfrak{F}$.
2. The closed interval $[[\mathcal{M}], [\mathcal{L}]]$ is a gap in $\mathfrak{F}\mathfrak{G}$.

*Proof*: Assume that $[[\mathcal{M}], [\mathcal{L}]]$ is a gap in $\mathfrak{F}\mathfrak{G}$, i.e., $(\mathcal{L} \setminus \mathcal{M}, \mathcal{M}^{\mathcal{L}})$ is a duality pair in $\mathfrak{F}\mathfrak{G}$. Now Let $\mathcal{K}$ be a locally finite minion which is not above $\mathcal{L} \setminus \mathcal{M}$. In particular, no finitely generated subminion $\mathcal{F} \leq \mathcal{K}$ can lie above $\mathcal{L} \setminus \mathcal{M}$. By duality in $\mathfrak{F}\mathfrak{G}$, we know that $\mathcal{F} \leq \mathcal{M}^{\mathcal{L}}$. By the compactness theorem for minion homomorphisms (Proposition 2.5.4), we must have a homomorphism $\mathcal{K} \to \mathcal{M}^{\mathcal{L}}$. Thus, $(\mathcal{L} \setminus \mathcal{M}, \mathcal{M}^{\mathcal{L}})$ is a duality pair in $\mathfrak{L}\mathfrak{F}$ as well, so $[\mathcal{M}, \mathcal{L}]$ must be a gap in $\mathfrak{L}\mathfrak{F}$.

The other direction is trivial since $\mathfrak{F}\mathfrak{G} \hookrightarrow \mathfrak{L}\mathfrak{F}$ (Proposition 5.1.1). □

In general, we might pose the bold question:

**Open Question 7.5**: Is the interval $[[\mathcal{P}], [\mathcal{I}_2]]$ *dense* in the order $\mathfrak{F}$ or $\mathfrak{L}\mathfrak{F}$?

If not, we should find gaps. One natural attempt would be to find a gap above $\mathcal{P}$:

**Question 7.6**: Does $\mathcal{P}$ have a cocover in $\mathfrak{F}$ or $\mathfrak{LF}$? Is there an infinite descending chain with infimum $\mathcal{P}$?

This question has a negative answer, as we will present briefly.

Another candidate for that might be the gaps in the Kazda-Moore embedding $\mathcal{P}(\{3,4,\ldots\}) \hookrightarrow \mathfrak{F}$; for instance, the smallest interesting gap would be $[\{3\},\{3,4\}]$. Since $\langle n_3, n_4\rangle \nrightarrow \langle n_3\rangle$ and both are only one unary, the subtraction $\langle n_3, n_4\rangle \setminus \langle n_3\rangle$ amounts to just $\langle n_3, n_4\rangle$. Hence, the question becomes:

**Open Question 7.7**: Which of the gaps of $\mathcal{P}(\{3,4,\ldots\})$ are preserved by Kazda and Moore's embedding

$$S \mapsto \langle n_s \mid s \in S\rangle \leq \mathcal{O}(2),$$

where $n_s(\underline{b}) = \left[\underline{b}^{-1}(0) \leq 1\right]$ is the $s$-ary near-unanimity term?

## 7.7 No gaps above $\mathcal{P}$

In a 2020 paper titled "ω-categorical structures avoiding height one identities" [63], the authors proved – as a corollary to some stronger claims – that there is no weakest nontrivial height-one identity.

More precisely, the central statement is as follows:

**Proposition 7.7.1** ([63, Thm. 1.3]): for every non-trivial height-1 condition $\Sigma$ there exists a structure $\mathbb{B}$ such that

1. $\mathbb{B}$ is a first-order reduct of a finitely bounded structure
2. $\mathrm{Pol}(\mathbb{B}) \not\models \Sigma$
3. $\mathrm{Pol}(\mathbb{B}) \models \Sigma'$ for some other non-trivial $\Sigma'$
4. $\mathrm{CSP}(\mathbb{B})$ is in $P$.

Paraphrasing in the language of this thesis and focusing on points 2 and 3, this means that

**Proposition**: For every finitely generated minion $\Sigma \nrightarrow \mathcal{P}$ there exists a structure $\mathbb{B}$ such that

- $\Sigma \nrightarrow \mathrm{Pol}(\mathbb{B})$
- there exists a finitely generated minion $\Sigma' \nrightarrow \mathcal{P}$ such that $\Sigma' \to \mathrm{Pol}(\mathbb{B})$.

To understand why this implies that $\mathcal{P}$ does not have a cocover, we must strengthen Proposition 2.5.3 to the following statement:

**Proposition 7.7.2**: Let $\mathcal{M}$ be an arbitrary minion and $\mathcal{N}$ locally finite. If all finitely generated subminions $\langle S\rangle \subseteq \mathcal{M}$ admit a homomorphism to $\mathcal{N}$, then $\mathcal{M}$ does as well.

*Proof*: We sketch a similar compactness argument as in the proof of Proposition 2.5.3. Let us view the set of all arity-preserving functions $\mathcal{M}$ to $\mathcal{N}$ as the topological product space

$$X := \prod_{\substack{i>0 \\ f \in \mathcal{M}_i}} \mathcal{N}_i,$$

where we recognize $\mathcal{N}_i$ as discrete finite spaces. By Tychonoff's theorem, $X$ is compact. For any finite set sequence $S \subseteq \mathcal{M}$, we define $C_S \subseteq X$ to be the set of all sequences $\left(h_{i,f}\right)_{i,f\in\mathcal{M}_i}$ such that the assignment $\langle S\rangle \ni f \mapsto h_{\mathrm{ar}\,f,f}$ gives rise to a minion homomorphism $\langle S\rangle \to \mathcal{N}$. Since every $C_S$ is defined by a condition referencing only finitely many factors, it is the preimage of a closed set under a projection, hence closed. The assumption that homomorphisms from every finitely generated subminion exist means precisely that every of the $C_S$ is nonempty. Since $C_{S_1} \cap C_{S_2} \supset C_{S_1 \cup S_2} \neq \emptyset$, the family $\{C_S \mid S \subset \mathcal{M} \text{ finite}\}$ satisfies the finite intersection property, so by compactness, the intersection $\bigcap_S C_S$ must be populated. Since sequences $\left(h_{i,f}\right)_{i>0, f\in\mathcal{M}_i}$ in this intersection correspond to homomorphisms $\mathcal{M} \to \mathcal{N}$, we are done. □

**Corollary 7.7.3**: Let $\mathcal{M}$ be any minion. If $\mathcal{M} > \mathcal{P}$, then there must be a nonempty finitely generated subminion $\mathcal{F} \leq \mathcal{M}$ such that $\mathcal{F} > \mathcal{P}$.

So far, the statement of Proposition 7.7.1 only guarantees that for each $\Sigma > \mathcal{P}$, there is another $\Sigma' > \mathcal{P}$ with $\Sigma' \not\to \Sigma$. However, to show that $\Sigma$ is no cocover of $\mathcal{P}$, we need to obtain a strictly weaker minion. The obvious candidate is the meet $\Sigma \times \Sigma'$, so it is left to show that this product admits no minion homomorphism to $\mathcal{P}$.

This step is not explicitly mentioned in [63], but deserves some attention. The argument, while known to the authors[58], is not quite trivial, and encodes a crucial property of $\mathcal{P}$.

**Lemma 7.7.4**: Let $\mathcal{M}$ and $\mathcal{N}$ be connected minions. If $\mathcal{M} \times \mathcal{N} \to \mathcal{P}$, then either $\mathcal{M} \to \mathcal{P}$ or $\mathcal{N} \to \mathcal{P}$.

To ease the proof, let us introduce some notation for certain elements of product minions.

**Definition 7.7.1** (Juxtaposed product): Let $\mathcal{M}$ and $\mathcal{N}$ be minions, $f \in \mathcal{M}_l$ and $g \in \mathcal{N}_k$. We define

$$f \times g := (f\iota_1, g\iota_2) \in (\mathcal{M} \times \mathcal{N})_{l+k},$$

where $\iota_1 \colon l \hookrightarrow l+k$ is the inclusion into the first summand and $\iota_2 \colon k \hookrightarrow l+k$ is the inclusion into the second summand, i.e., $\iota_1(i) = i$ and $\iota_2(i) = l + i$.

[58] The author is grateful to M. Bodirsky for providing a missing step in an attempted proof.

If $\alpha + \beta\colon l + k \to l' \to k'$ is the (sum-) juxtaposition of the minor operations $\alpha\colon l \to l'$ and $\beta\colon k \to k'$, it follows from the definition that $(f \times g)(\alpha + \beta) = f\alpha \times g\beta$. Viewing minor operations as tuples, this means

$$(f \times g)(\alpha_0, \ldots, \alpha_{l-1}, \alpha_l, \ldots, \alpha_{l+k-1}) = (f(\alpha_0, \ldots, \alpha_{l-1}), g(\alpha_l, \ldots, \alpha_{l+k-1})).$$

This is in contrast to the combination $(f, f')$ of elements with equal arity $l$, which satisfies

$$(f, f')(\alpha_0, \ldots, \alpha_{l-1}) = (f(\alpha_0, \ldots, \alpha_{l-1}), f'(\alpha_0, \ldots, \alpha_{l-1})).$$

*Remark*: The notation "$f \times g$" is motivated by viewing elements as homomorphisms: We can (product-) juxtapose $f\colon \mathcal{P}^l \to \mathcal{M}$ and $g\colon \mathcal{P}^k \to \mathcal{N}$ to a morphism $f \times g\colon \mathcal{P}^l \times \mathcal{P}^k \to \mathcal{M} \times \mathcal{N}$, and recognizing that $\mathcal{P}^l \times \mathcal{P}^k \cong \mathcal{P}^{l+k}$ gives us the corresponding element of arity $l + k$ in the product.

*Proof of Lemma 7.7.4*: Let $h\colon \mathcal{M} \times \mathcal{N} \to \mathcal{P}$ and let $u \in \mathcal{M}_1$ resp. $v \in \mathcal{N}_1$ be the unique unaries. A dichotomy will arise from the observation that $h(u \times v) \in \mathcal{P}_{1+1}$ is either $0$ or $1$. Indeed, for elements $f \in \mathcal{M}_l$ and $g \in \mathcal{N}_k$, it is equivalent that

1. $h(f \times g) < l$ and
2. $h(u \times v) = 0$:

Since $u \times v = (f \times g)(0, \ldots, 0, 1, \ldots, 1)$, we have $h(u \times v) = (h(f \times g)).(0, \ldots, 0, 1, \ldots, 1)$, so the the first clearly implies the latter; however if the first does not hold, i.e. $h(f \times g) \geq l$, the analogous observation shows that $h(u \times v) = 1$.

The case $h(u \times v) = 0$ however implies that $h(_ \times v)$ defines a minion homomorphism $\mathcal{M} \to \mathcal{P}$: indeed, if $f \in \mathcal{M}_l$, then $h(f \times v) \in \{0, \ldots, l-1\} = \mathcal{P}_l$, and for any $\alpha\colon l \to l'$ we have

$$h(f\alpha \times v) = \underbrace{h(f \times v)}_{\in \operatorname{cod} \alpha}(\alpha + \mathrm{id}_1) = h(f \times v)\alpha.$$

Analogously, if $h(u \times v) = 1$, the map $h(u \times _)$ induces a homomorphism $\mathcal{N} \to \mathcal{P}$, so in any case, one of the factors must map to $\mathcal{P}$. □

**Corollary 7.7.5**: In any of the minion homomorphism orders $\mathfrak{M}$, $\mathfrak{LF}$, $\mathfrak{F}$ and $\mathfrak{FG}$, the equivalence class of the minion $\mathcal{P}$ is *meet-prime*, i.e. $[\mathcal{M}] \wedge [\mathcal{N}] \leq [\mathcal{P}]$ implies that one of $[\mathcal{M}] \leq [\mathcal{P}]$ or $[\mathcal{N}] \leq [\mathcal{P}]$.

*Proof*: Let $\mathfrak{C}$ be any of the classes – all of which are closed under taking substructures and products. To argue contrapositively, let $\mathcal{M}, \mathcal{N}$ be minions from $\mathfrak{C}$ such that $\mathcal{M} \nrightarrow \mathcal{P}$ and $\mathcal{N} \nrightarrow \mathcal{P}$. Both premises must be witnessed on some connected component, i.e., there are unaries $u \in \mathcal{M}_1$ and $v \in \mathcal{N}_1$ such that $\mathcal{M}{\downarrow_u} \nrightarrow \mathcal{P}$ and $\mathcal{N}{\downarrow_v} \nrightarrow \mathcal{P}$. Applying Lemma 7.7.4, this implies that $\mathcal{M}{\downarrow_u} \times \mathcal{N}{\downarrow_v} \nrightarrow \mathcal{P}$, and being a subminion of $\mathcal{M} \times \mathcal{N}$, this obstructs a homomorphism $\mathcal{M} \times \mathcal{N} \to \mathcal{P}$. □

Not too surprisingly, meet-primeness and the exponential are closely related, as the following observation shows.

**Proposition 7.7.6**: A minion $\mathcal{M}$ is meet-prime in $\mathfrak{M}$ if and only if it does not contain a constant and for every $\mathcal{N}$, the exponential $\mathcal{M}^{\mathcal{N}}$ is either homomorphically equivalent to $\mathcal{M}$ or to $*$.

*Proof*: For one direction, assume $\mathcal{M}$ is meet-irreducible and let $\mathcal{N}$ be arbitrary. Then the evaluation homomorphism $\mathcal{N} \times \mathcal{M}^{\mathcal{N}} \to \mathcal{M}$ implies that either $\mathcal{N} \to \mathcal{M}$ or $\mathcal{M}^{\mathcal{N}} \to \mathcal{M}$. The first case tells us that $\mathcal{M}^{\mathcal{N}}$ has a constant, while the second case implies it is somomorphically equivalent to $\mathcal{M}$, since there is always a homomorphism $\mathcal{M} \to \mathcal{M}^{\mathcal{N}}$.

For the converse, let $\mathcal{N} \times \mathcal{K} \to \mathcal{M}$, or equivalently, $\mathcal{K} \to \mathcal{M}^{\mathcal{N}}$. If $\mathcal{N} \nrightarrow \mathcal{M}$, then the exponential has no constant, so by assumption it must be homomorphically equivalent to $\mathcal{M}$. We conclude that $\mathcal{K} \to \mathcal{M}$.

That $[\mathcal{M}] \neq [*]$ is precisely the demand that meet-prime elements shall not be the top element. □

As $\mathcal{P}$ is meet-prime, we conclude the following dichotomy for its powers:

**Corollary 7.7.7**: For any minion $\mathcal{M}$, we have

- $\mathcal{P}^{\mathcal{M}} \leftrightarrow \mathcal{P}$ if and only if $\mathcal{M} \nrightarrow \mathcal{P}$
- $\mathcal{P}^{\mathcal{M}} \leftrightarrow *$ if and only if $\mathcal{M} \nrightarrow \mathcal{P}$.

Finally, let us conclude that $\mathcal{P}$ has no cocover.

**Corollary 7.7.8**: In none of the orders $\mathfrak{M}, \mathfrak{LF}, \mathfrak{F}, \mathfrak{FG}$ does the minion $\mathcal{P}$ have a co-cover.

*Proof*: Let us first consider the class $\mathfrak{FG}$, so let $\mathcal{M}$ be finitely generated such that $\mathcal{M} \nrightarrow \mathcal{P}$. By Proposition 7.7.1, there exists a finitely generated minion $\mathcal{M}'$ such that $\mathcal{M}' \nrightarrow \mathcal{P}'$ and $\mathcal{M} \nrightarrow \mathcal{M}'$. by Lemma 7.7.4, $\mathcal{M} \times \mathcal{M}' \nrightarrow \mathcal{P}$, and since a homomorphism $\mathcal{M} \to \mathcal{M} \nrightarrow \mathcal{M}'$ would imply a homomorphism $\mathcal{M} \to \mathcal{M}'$, we have $[\mathcal{P}] < [\mathcal{M}] \wedge [\mathcal{M}'] < [\mathcal{M}]$, refuting the hypothesis.

For the other classes, we can apply the compactness argument (Proposition 7.7.2): If $\mathcal{M} \nrightarrow \mathcal{P}$, there must be a finitely generated minion $\mathcal{F}$ with $[\mathcal{P}] < [\mathcal{F}] \leq [\mathcal{M}]$, but from the argument above – and since all the other classes contain $\mathfrak{FG}$ – $[\mathcal{P}, \mathcal{F}]$ cannot be a gap. □

# 8 Epilogue

We have demonstrated how one can use abstract tools to analyze homomorphism orders of minions. Indeed, just by considering the canonical constructions on functors, and showing some preservation rules, we have shown that

- $\mathfrak{F}$ is a bi-Heyting algebra with coatom $[\mathcal{I}_2]$, for which we determined a core,
- $\mathfrak{LF}$ and $\mathfrak{M}$ are bi-Heyting algebras as well, and that
- gaps are characterizable by duality pairs in relation to the subtraction and the exponential minion.

We hope to have demonstrated the strength of the abstract approach, e.g. by being able to specify gadgets to primitive positive sentences like $\exists(f\colon 2).f(1,0) = f$ directly as a minion $\langle f\colon 2 \mid f(1,0) = f\rangle$.

As the scope of a master's thesis is rather limited, there is still much left to discover. The author in particular has selected some questions along the way, which he would like to see answered. In particular, it seems a good line of investigation to try and determine some exponential minions to determine whether they can be efficiently given a more explicit description, be it via generators and relations, or via invariant relation pairs.

Otherwise, the author is keen to see how current research regarding promise CSPs unfolds, and how this relates to the homomorphism order of minions. Perhaps it will turn out that some previously unconsidered class of minions precisely determines the boundary of a complexity dichotomy; or maybe even a non-dichotomy result will arise.

In the meantime, the author hopes that this thesis aids someone in advancing their understanding of minions.

## List of open questions

- **Question 5.1.** Are quotients of finitely representable minions finitely representable?
- **Question 5.2.** For finite sets $A, B$, is $\mathfrak{F}_{A,B}$ a lattice?
- **Question 5.3.** When is $\mathrm{res}_n$ surjective, i.e. when does an endomorphism of $\mathcal{T}_n \curvearrowright \mathcal{M}_n$ extend to an endomorphism of the full minion (assuming $n$-representability)?
- **Question 7.1.**
  1. Given a locally finite minion $\mathcal{M}$ and a finitely representable minion $\mathcal{N}$, under which conditions does $\mathcal{N}^{\mathcal{N}^{\mathcal{M}}} \to \mathcal{M}$?
  2. Given a locally finite minion $\mathcal{M}$, does there always exist a finitely representable minion $\mathcal{N}$ such that $\mathcal{N}^{\mathcal{N}^{\mathcal{M}}} \to \mathcal{M}$?
- **Question 7.2.** Is every locally finite minion homomorphically equivalent to a finitely representable minion?
- **Question 7.3.** Let $\underline{A}, \underline{B}$ be relational structures and let $n := |A|$. If $\mathcal{M}$ is a minion with $\mathcal{M}_n < \infty$, what does a relational basis of $\mathrm{Pol}(\underline{A}, \underline{B})^{\mathcal{M}}$ look like?
- **Question 7.4.** Does the exponential minion preserve being finitely related? In particular, given a (finite) relational basis for $\mathcal{M}$ and $\mathcal{N}$, is there an algorithm to produce a relational basis of $\mathcal{N}^{\mathcal{M}}$?
- **Question 7.5.** Is the interval $[[\mathcal{P}], [\mathcal{I}_2]]$ *dense* in the order $\mathfrak{F}$ or $\mathfrak{LF}$?
- **Question 7.7.** Which of the gaps of $\mathcal{P}(\{3, 4, \ldots\})$ are preserved by Kazda and Moore's embedding

$$S \mapsto \langle n_s \mid s \in S\rangle \leq \mathcal{O}(2),$$

where $n_s(\underline{b}) = \left[\underline{b}^{-1}(0) \leq 1\right]$ is the $s$-ary near-unanimity term?

# A Appendix

## A.1 Galois Connections

We follow [64, Chapter A5].

Recall that a partially ordered set, or *poset* for short, is a set $P$ equipped with a binary relation often denoted by $\leq$ which is reflexive ($\forall p\colon p \leq p$), transitive ($\forall p, q, r\colon p \leq q \wedge q \leq r \Rightarrow p \leq r$), and antisymmetric ($\forall p, q\colon p \leq q \wedge q \leq p \Rightarrow p = q$).

**Definition A.1.1**: A *antitone Galois connection* between two posets $(P, \leq)$ and $(Q, \leq)$ is a pair of maps $\varphi\colon P \to Q$ and $\psi\colon Q \to P$ which are antitone[59], such that for all $p \in P$ and $q \in Q$ we have

$$p \leq \varphi(q) \quad \text{if and only if} \quad q \leq \psi(p).$$

**Definition A.1.2**: A *closure operator* on a poset $(P, \leq)$ is an endomorphism $C$ which is extensive and idempotent, i.e., for all $p, q, r \in P$, we have
1. $p \leq C(p)$,
2. $p \leq q$ implies $C(p) \leq C(q)$, and
3. $C(C(p)) = C(p)$.

Note that under the assumption of the first two conditions, to verify the third it suffices for $C(C(p)) \leq C(p)$ to be upheld.

**Proposition A.1.1**: For an antitone Galois connection $(\varphi, \psi)$ between $(P, \leq)$ and $(Q, \leq)$, the following hold:
1. $\varphi \circ \psi$ and $\psi \circ \varphi$ are monotone and extensive.
2. $\varphi \circ \psi \circ \varphi = \varphi$ and $\psi \circ \varphi \circ \psi = \psi$.
3. $\varphi \circ \psi$ and $\psi \circ \varphi$ are closure operators.
4. $p$ is in the image of $\psi$ if and only if $p = \psi(\varphi(p))$, and $q$ is in the image of $\varphi$ if and only if $q = \varphi(\psi(q))$.

*Proof*: 1) Both $\varphi \circ \psi$ and $\psi \circ \varphi$ are compositions of two antitone maps, hence monotone. For extensivity, note that $q := \varphi(p) \leq \varphi(p)$ implies by the condition of a Galois connection that $p \leq \psi(q) = \psi(\varphi(p))$, and analogously for $\varphi$, proving extensivity. 2) By extensivity of $\psi \circ \varphi$, we have $p \leq \psi(\varphi(p))$, so applying the antitone map $\varphi$ yields $\varphi(\psi(\varphi(p))) \leq \varphi(p)$. Conversely, extensivity of $\varphi \circ \psi$ applied to $q := \varphi(p)$ yields $\varphi(p) \leq \varphi(\psi(\varphi(p)))$. The argument for $\psi = \psi \circ \varphi \circ \psi$ is analogous. 3) Idempotence follows since $(\varphi \circ \psi) \circ (\varphi \circ \psi) = (\varphi \circ \psi \circ \varphi) \circ \psi = \varphi \circ \psi$, and dually we get idempotence of $\psi \circ \varphi$. 4) If $p = \psi(q)$, then $\psi(\varphi(p)) = \psi(\varphi(\psi(q))) = \psi(q) = p$, and dually the argument for $p$ and $\operatorname{Im} \varphi$ follows. □

[59] Also known as "order-reversing". It means that they are homomorphisms from the domain poset to the oppositely-ordered codomain poset, e.g. $\varphi\colon (P, \leq) \to (Q, \leq^{\text{op}})$, where $q \leq^{\text{op}} q'$ if and only if $q' \leq q$.

**Proposition A.1.2**: Let $R \subseteq X \times Y$ be a binary relation. Then the maps

$$\varphi\colon \mathcal{P}(X) \to \mathcal{P}(Y), S \mapsto \{y \in Y \mid sRy \forall s \in S\}$$

and

$$\psi\colon \mathcal{P}(Y) \to \mathcal{P}(X), T \mapsto \{x \in X \mid xRt \forall t \in T\}$$

form an antitone Galois connection between $(\mathcal{P}(X), \subseteq)$ and $(\mathcal{P}(Y), \subseteq)$.

*Proof*: Assume $S \subseteq S'$ Then $t \in \varphi(S')$ if and only if $\forall s \in S'\colon tRs$, which by restriction implies $\forall s \in S\colon tRs$, i.e., $t \in \varphi(S)$. Hence $\varphi(S') \subseteq \varphi(S)$, i.e., $\varphi$ is antitone. $\psi$ is antitone by an analogous argument. The Galois property follows completely formally from

$$\begin{aligned} T \subseteq \varphi(S) &\Leftrightarrow \forall t \in T\colon t \in \varphi(S) && \text{Definition of } \subseteq \\ &\Leftrightarrow \forall t \in T\colon \forall s \in S\colon sRt && \text{Definition of } \varphi \\ &\Leftrightarrow \forall s \in S\colon \forall t \in T\colon sRt && \forall \text{ commutative} \\ &\Leftrightarrow \forall s \in S\colon s \in \psi(T) && \text{Definition of } \psi \\ &\Leftrightarrow S \subseteq \psi(T). && \text{Definition of } \subseteq \end{aligned}$$

□

## A.2 The named perspective

We shall prove here that every minion in the definition given here in this thesis can be extended to a minion with "named arities" in such a way that this extension is unique, and the restriction yields an isomorphic minion.

With some more concepts of category theory, this is easy to prove in the abstract[60], however we want to give an explicit construction using the tools of this thesis.

**Proposition A.2.1**: Let $\mathcal{M}$ be a minion. The assignment

$$X \mapsto \mathcal{M} \otimes X,$$

gives rise to a functor $\underline{\mathsf{FinSet}}_{>\emptyset} \to \underline{\mathsf{Set}}$ whose restriction to arities in $\mathbb{F}_{>0}$ yields an isomorphic minion.

Furthermore, first restricting a functor $F\colon \underline{\mathsf{FinSet}}_{>\emptyset} \to \underline{\mathsf{Set}}$ to $\mathbb{F}_{>0}$ and then applying this extension yields a naturally isomorphic functor.

*Proof*: Indeed, the functor structure is given by letting an operation $\alpha\colon X \to Y$ act pointwise on the right component, i.e., via

$$\mathcal{M} \otimes X \to \mathcal{M} \otimes Y, [f, \underline{x}] \to [f, \underline{x}].\alpha := [f, \underline{x};\alpha].$$

This assignment is well-defined because for every $\beta$ we have

[60]Namely, the essential statement here is that an equivalence of categories $\mathbb{F}_{>0} \simeq \underline{\mathsf{FinSet}}_{>\emptyset}$ extends to an equivalence of functor categories $[\mathbb{F}_{>0}, \underline{\mathsf{Set}}] \simeq [\underline{\mathsf{FinSet}}_{>\emptyset}, \underline{\mathsf{Set}}]$.

$$[f\beta, \underline{x}].\alpha = [f\beta, \underline{x};\alpha] = [f, \beta;\underline{x};\alpha] = [f, \beta;\underline{x}].\alpha,$$

i.e., because pre-composition and post-composition commute.

We will now sohow that the element-wise correspondence $\mathcal{M}_n \cong \mathcal{M} \otimes n$ which sends $f \mapsto [f, \mathrm{id}_n]$, is a minion homomorphism; indeed, for every $\alpha\colon n \to k$ we have

$$f\alpha \mapsto [f\alpha, \mathrm{id}_k] = [f, \alpha;\mathrm{id}_k] = [f, \alpha] = [f, \mathrm{id}_n].\alpha,$$

and bijectivity of this operation has been established in Lemma 3.3.3.

For the last point, let $F$ be a functor from $\underline{\mathrm{FinSet}}_{>\emptyset}$ to $\underline{\mathrm{Set}}$. Let us denote the action of a minor operation $\alpha\colon X \to Y$ of this functor on an element $f \in F_X$ as $f.\alpha$.

We have to show that we have a family of isomorphisms $F \otimes X \cong F_X$ which is natural in $X$. let $\underline{x} \in |X| \to X$ be a bijection, i.e., a tuple enumerating all the elements of the finite set $X$. We claim that the assignment $[f, \underline{x}] \mapsto f.\underline{x}$, where we interpret $\underline{x} \in X^{\mathrm{ar}\, f}$ as a minor operation $\mathrm{ar}\, f \to X$ acting on $f$, is natural in $X$ and a bijection. Naturality should be clear because $[f, \underline{x}].\alpha = [f, \underline{x};\alpha]$ maps to $(f.\underline{x}).\alpha$. Bijectivity follows from $F_n \cong F \otimes n$ and the naturality of this assignment, since an enumeration $\xi\colon n \to X$ induces isomorphisms $F_n \to F_X$ and $F \otimes n \to F \otimes X$. □

## A.3 A proof of uniqueness of the (right) adjoint

We start with a variant of the so-called yoneda lemma, a special case of which (focusing on posets) we have established in Lemma 7.6.1. It is an important principle which helps us prove many statements about universal constructions (such as adjunctions) with elegance.

Recall that for any category $\mathcal{D}$ and object $d \in \mathcal{D}$, we have an associated "hom-functor" $_ \to_{\mathcal{D}} d\colon \mathcal{D}^{\mathrm{op}} \to \underline{\mathrm{Set}}$, which assigns an arrow $f\colon e \to_{\mathcal{D}} e'$ the precomposition

$$f^* := _\circ f\colon (e' \to_{\mathcal{D}} d) \to (e \to_{\mathcal{D}} d).$$

**Proposition A.3.1** (Yoneda Lemma, cf. [6, Corollary 5.2.8]): If two hom-functors $_ \to_{\mathcal{D}} d$ and $_ \to_{\mathcal{D}} d'$ are naturally isomorphic, then $d \cong d'$.

*Proof*: Let $\eta\colon (_ \to_{\mathcal{D}} d) \to (_ \to_{\mathcal{D}} d')$ be a natural isomorphism, i.e.,

1. for every object $e \in \mathcal{D}_0$, we have a bijection $\eta_e\colon (e \to_{\mathcal{D}} d) \to (e \to_{\mathcal{D}} d')$
2. $\eta$ is a natural transformation: whenever $f\colon e \to e'$, we have $\eta_{e'}\circ f^* = f^*\circ\eta_e$. Using placeholder notation, this becomes $\eta_{e'}(_\circ f) = \eta_e(_)\circ f$.

Since $\eta$ bijective in every component, the component-wise inverse $\kappa_e := \eta_e^{-1}$ is easily seen to be natural as well, i.e., it satisfies $\kappa_{e'}(_\circ f) = \kappa_e(_)\circ f$.

We claim that $\eta_d(\mathrm{id}_d)\colon d \to_{\mathcal{D}} d'$ and $\kappa_{d'}(\mathrm{id}_{d'})\colon d' \to_{\mathcal{D}} d$ are inverses. Indeed, we have

$$\eta_d(\mathrm{id}_d)\circ\kappa_{d'}(\mathrm{id}_{d'}) = \eta_d(\mathrm{id}_d \circ\kappa_{d'}(\mathrm{id}_{d'})) = \eta_d(\kappa_{d'}(\mathrm{id}_{d'})) = \mathrm{id}_{d'}$$

by naturality and inverse-ness, and dually

$$\kappa_{d'}(\mathrm{id}_{d'})\circ\eta_d(\mathrm{id}_d) = \kappa_{d'}(\mathrm{id}_{d'} \circ\eta_d(\mathrm{id}_d)) = \kappa_{d'}(\eta_d(\mathrm{id}_d)) = \mathrm{id}_d\,.$$

□

**Proposition A.3.2**: If $R, R'\colon \mathcal{C} \to \mathcal{D}$ are both right adjoint to $L\colon \mathcal{D} \to \mathcal{C}$, then for every $c \in \mathcal{C}_0$, $Rc \cong R'c$. Furthermore, $R$ and $R'$ are even naturally isomorphic, i.e., isomorphic as functors $\mathcal{C} \to \mathcal{D}$.

*Proof*: For the first part, we shall note that by the definition of an adjunction, we have a chain of isomorphisms

$$d \to_{\mathcal{D}} Rc \cong Ld \to_{\mathcal{C}} c \cong d \to_{\mathcal{D}} R'c$$

which is natural in $d$ and in $c$, i.e., we have a natural isomorphism between the functors $_ \to_{\mathcal{D}} Rc$ and $_ \to_{\mathcal{D}} R'c$. By the yoneda lemma, these objects must be isomorphic.

For the second part, we use essentially the same argument, except we apply it instead to *functor category* $[\mathcal{D}, \mathcal{C}]$. Note that the hom-sets of the functor category are natural transformations. The correspondences then lift to assignments

$$F \to_{[\mathcal{D},\mathcal{C}]} R \quad \cong \quad L \circ F \to_{[\mathcal{D},\mathcal{D}]} \mathrm{id}_{\mathcal{D}} \quad \cong \quad F \to_{[\mathcal{D},\mathcal{C}]} R'.$$

These correspondences can be checked to be natural in $F$, allowing us once again to apply the yoneda lemma. □

# B Bibliography

# C Index

**R**

**S**

**T**

**U**

**W**

# D Symbols and notation

| | |
|---|---|
| $\mathcal{C}, \mathcal{D}$ | categories |
| $\mathcal{C}_0$ | objects of $\mathcal{C}$ |
| $\mathcal{C}_1$ | morphisms of $\mathcal{C}$ |
| $\mathcal{C}(X, Y)$, $X \to_{\mathcal{C}} Y$ | $\mathcal{C}$-morphisms from $X$ to $Y$ |
| $X \xrightarrow{f}_{\mathcal{C}} Y$ | a specific $\mathcal{C}$-morphism from $X$ to $Y$ |
| $\cong$ | isomorphism |
| $\mathrm{End}_{\mathcal{C}}(X)$ | Endomorphism monoid of $X$ |
| $\mathrm{Aut}_{\mathcal{C}}(X)$ | Automorphism group of $X$ |
| $\mathcal{C}^{\mathsf{op}}$ | opposite category |
| $FX$ | functor $F$ applied to object $X$ |
| $Ff$ | functor $F$ applied to morphism $f$ |
| $F \Rightarrow G$, $F \to G$ | Natural transformation between functors |
| $c \to_{\mathcal{C}} _$ | (covariant) hom-functor |
| $f_*$ | pushforward of $f$ |
| $c \to_{\mathcal{C}} _$ | (covariant) hom-functor |
| $f_*$ | pullback of $f$ |
| $c \to_{\mathcal{C}} _$ | (covariant) hom-functor |
| $[X]$ | equivalence class of $X$ |
| $X \to Y$ | there exists a morphism from $X$ to $Y$ |
| $\langle q_i \rangle_i$ | tupling |
| $\left(\Pi_i X_i, (\pi_i)_i\right)$ | product of objects with projections |
| $\pi_i^k$ | projection $X^k \twoheadrightarrow X$ onto the $i$th coordinate |
| $[s_i]_i$ | cotupling |
| $\left(\Sigma_i X_i, (\iota_i)_i\right)$ | sum of objects with injections |
| $F \dashv G$ | adjoint functors |
| $\underline{\mathsf{FinSet}}_{>\emptyset}$ | category of nonempty finite sets and functions |
| $\mathrm{Pol}_{\mathcal{C}}(A)$ | shorthand for $\mathrm{Pol}_{\mathcal{C}}(A, A)$ |
| $\mathcal{O}(A, B)$ | shorthand for $\mathrm{Pol}_{\underline{\mathsf{Set}}}(A, B)$ |
| $\mathrm{ar}\, h$ | arity of $h$ |
| $(s_i \colon k_i)_i$ | set sequence $S$ w/ $s_i \in S_{k_i}$ |
| $X \hookrightarrow Y$ | injection |
| $X \twoheadrightarrow Y$ | surjection |
| $\mathrm{Cong}(\mathcal{M})$ | congruences on $\mathcal{M}$ |
| $\Delta_{\mathcal{M}}$ | diagonal congruence |
| $\langle S \rangle$ | free minion generated by $S$ |
| $\langle f_i \colon n_i \rangle_i$ | free minion generated by $n_i$-ary elements $f_i$ |
| $[\![_]\!]$ | choice of witnesses |
| $(\alpha_0, \ldots, \alpha_{n-1})_k$ | A minor operation $\alpha \colon n \to k$ encoded as a tuple |
| $\iota$ | canonical inclusion |
| $\mathrm{Ess}\, f$ | essential coordinates of $f$ |
| $\mathrm{Iness}\, f$ | inessential coordinates of $f$ |
| $X^{\complement}$ | complement of the set $X$ |
| $*$ | terminal object (e.g. constant minion) |
| $\mathcal{B}_\infty$ | idempotent boolean functions preserving $\mathrm{NAZ}(k)$ |
| $\mathrm{NAZ}(k)$ | "not all zero"-relation |
| $\mathcal{C}$ [R] | see Lemma 3.2.2 |
| $\mathcal{M} \otimes X$ | height-one terms |
| $\mathrm{id}_{\mathcal{M}} \otimes m$ | induced map $\mathcal{M} \otimes X \to \mathcal{M} \otimes Y$ |
| $\mathfrak{M}$ | $\underline{\mathsf{Min}}\ /\!\leftrightarrow$ |
| $\underline{\mathsf{LFin}}$ | cat. of locally finite minions |
| $\mathfrak{LF}$ | $\underline{\mathsf{LFin}}\ /\!\leftrightarrow$ |
| $\underline{\mathsf{FinRep}}$ | cat. of finitely representable minions |
| $\mathfrak{F}$ | $\underline{\mathsf{FinRep}}\ /\!\leftrightarrow$ |
| $\underline{\mathsf{Rep}}_{A,B}$ | cat. of $(A, B)$-representable minions |
| $\mathfrak{F}_{A,B}$ | $\underline{\mathsf{Rep}}_{A,B}/\!\leftrightarrow$ |
| $\underline{\mathsf{Rep}}_A$ | cat. of $(A, A)$-representable minions |
| $\mathfrak{F}_A$ | $\underline{\mathsf{Rep}}_{A,A}/\!\leftrightarrow$ |
| $\underline{\mathsf{FinGen}}$ | cat. of finitely generated minions |
| $\mathfrak{FG}$ | $\underline{\mathsf{FinGen}}\ /\!\leftrightarrow$ |
| $\mathcal{M}{\downarrow_S}, \mathcal{M}{\downarrow_u}$ | restriction onto unary |
| $\mathcal{O}(n, k){\downarrow_{\mathrm{nc}}}$ | functions with non-constant unary |
| $\mathcal{I}_n$ | minion of idempotents over $n$ |
| $R_{l,r}$ | reflection along $(l, r)$ |
| $\mathrm{core}(X)$ | core |
| $\alpha_{\mathcal{M}}(n)$ | growth |
| $\gamma_{\mathcal{M}}(n)$ | essential growth |
| $\Omega$ | subobject classifier |
| $p^{\uparrow}$ | upward closure |
| $p^{\downarrow}$ | downward closure |
| $[a, c]$ | closed interval |
| $f \times g$ | Juxtaposed product of $f$ and $g$ |